\documentclass[11pt,oneside]{book}
\usepackage[a4paper,margin=2.6cm]{geometry}\usepackage{amsmath,amssymb,amsthm,booktabs,array,enumitem}\usepackage{xcolor}\definecolor{linkcol}{rgb}{0.10,0.18,0.42}
\usepackage[colorlinks=true,linkcolor=linkcol,citecolor=linkcol,urlcolor=linkcol]{hyperref}
\newenvironment{abstract}{\begin{quote}\small\textbf{Abstract.} }{\end{quote}}
\usepackage{mathtools,mathrsfs}
\usepackage{needspace}
\numberwithin{equation}{chapter}
\numberwithin{table}{chapter}
\numberwithin{figure}{chapter}
\makeatletter
\renewcommand\paragraph{\@startsection{paragraph}{4}{\z@}{1.3ex plus .5ex minus .2ex}{-.8em}{\normalfont\normalsize\bfseries}}
\makeatother

\newtheorem{theorem}{Theorem}[chapter]
\newtheorem{proposition}[theorem]{Proposition}
\newtheorem{lemma}[theorem]{Lemma}
\newtheorem{corollary}[theorem]{Corollary}
\theoremstyle{definition}
\newtheorem{definition}[theorem]{Definition}
\newtheorem{example}[theorem]{Example}
\theoremstyle{remark}
\newtheorem{remark}[theorem]{Remark}
\newtheorem{question}[theorem]{Question}
\newtheorem*{thmA}{Theorem A (the destination formula)}
\newtheorem*{thmB}{Theorem B (the stratified convergence principle)}
\newtheorem*{thmD}{Theorem D (the flag of limit sets)}
\newtheorem*{thmC}{Theorem C (limit sets and proper discontinuity)}
\newtheorem*{thmE}{Theorem E (principal images of $\PSL(2,\C)$)}
\newcommand{\restatedname}{}
\newtheorem*{restatedthm}{\restatedname}
\newenvironment{restated}[1]{\renewcommand{\restatedname}{#1}\begin{restatedthm}}{\end{restatedthm}}
\newcommand{\sectionpart}[1]{\par\bigskip\needspace{8\baselineskip}\begin{center}\normalfont\scshape #1\end{center}\phantomsection\addcontentsline{toc}{section}{\textsc{#1}}\nopagebreak\smallskip}

\newcommand{\C}{\mathbb C}\newcommand{\K}{\mathbb K}\newcommand{\R}{\mathbb R}\newcommand{\PP}{\mathbb P}
\newcommand{\g}{\mathfrak g}\newcommand{\ak}{\mathfrak a}\newcommand{\nk}{\mathfrak n}
\newcommand{\kk}{\mathfrak k}\newcommand{\pk}{\mathfrak p}

\DeclareMathOperator{\PSL}{PSL}\DeclareMathOperator{\PGL}{PGL}
\DeclareMathOperator{\Ad}{Ad}\DeclareMathOperator{\diag}{diag}

\newcommand{\D}{\mathscr D}
\DeclareMathOperator{\Eq}{Eq}
\newcommand{\Kul}{\mathrm{Kul}}\newcommand{\Myr}{\mathrm{Myr}}\newcommand{\CoG}{\mathrm{CoG}}
\title{Cartan Dynamics and Limit Sets on Flag Manifolds}
\author{\'Angel Cano, Marlon M.\ L\'opez-Flores\\Jos\'e Seade, Andrew Shaw}
\date{}
\providecommand{\doi}[1]{\href{https://doi.org/#1}{\nolinkurl{doi:#1}}}
\begin{document}

\maketitle

\begin{abstract}
The classical convergence principle says that a divergent sequence of M\"obius transformations, after passing to a subsequence, collapses the Riemann sphere minus one exceptional point onto a single point. We develop its counterpart for an arbitrary real semisimple Lie group $G$ acting on its flag manifolds, and use it to study the asymptotic dynamics of discrete subgroups of $G$: their limit sets and domains of discontinuity.

The key is to let the source move: the \emph{destination} of a point $x$ under a sequence $(g_m)$ is the set of all limits of $g_mx_m$ with $x_m\to x$. The \emph{destination formula} computes it for every divergent sequence, at every point of every flag manifold: it is the closure of a unipotent orbit, positioned by the Cartan data of the sequence and by a bounded Levi factor that has no counterpart in rank one. The \emph{stratified convergence principle} organises these destinations globally. After extraction, a divergent sequence converges to a boundary point of a compactification of $G$; this point cuts the flag manifold into cells and attaches to each cell a family of translated Schubert varieties, and every boundary point arises in this way. Over the complex numbers, the destination relation is the support of Brion's limiting graphs in the wonderful compactification, and this algebraic model gives scheme structures, multidegrees and diagonal intersections of the degenerating graphs.

For a discrete group, the union of these relations decides when a domain of proper discontinuity exists and when there is a largest one. It places the flag limit set of Quint, which through the Tits representation is a Conze--Guivarc'h set, the Kulkarni set and the Myrberg set in a single chain around the recurrent set, and balanced choices of cells produce domains of discontinuity for every discrete subgroup, with no Anosov hypothesis; the thickenings of Kapovich, Leeb and Porti fit into this framework.

These notions genuinely differ. Let a non-elementary Kleinian group whose limit set $\Lambda$ is not the whole sphere act on $\C\PP^N$ through the $N$-th symmetric power. Its largest domain consists of the divisors in which no point of $\Lambda$ has multiplicity at least $N/2$; for $N\ge3$ it strictly contains the Kulkarni domain, which here is the region of equicontinuity. Zariski-dense cusped Hitchin groups in $\mathrm{SL}_5(\R)$, which contain unipotent elements, have a nonempty largest domain in $\C\PP^4$, while even a strongly irreducible group may have no largest domain at all.

\end{abstract}

\clearpage

\chapter*{Preface}
\addcontentsline{toc}{chapter}{Preface}

This monograph springs from the work of the authors and others on complex Kleinian groups, a notion introduced in the 1990s by Seade and Verjovsky~\cite{SeadeVerjovsky2001,SeadeVerjovsky2002}, inspired by Sullivan's dictionary between Kleinian groups and iteration theory~\cite{Sullivan1985} and by the remarkable development of the latter in several complex variables.

The theory of complex Kleinian groups developed gradually, at first almost as an experimental science, through examples and results in particular cases; the bibliography gives a wide view of the subject, and~\cite{CNS} is a general reference. For a long time the theory seemed to face a complicated puzzle: several notions of limit set, and of related concepts, closely linked yet distinct. The article~\cite{BCNSFlag2026} finally showed us how to put the pieces of that puzzle together. Its trichotomy for divergent sequences in $\C\PP^2$, expressed through full and partial flags and the destinations of moving sources, is a direct precursor of our stratified convergence principle. That article suggested that the different limit sets are shadows of a common asymptotic geometry. We develop this perspective for real semisimple Lie groups acting on their flag manifolds and compare our results with the theories of Anosov and transverse groups and the work of Kapovich, Leeb and Porti.

\medskip
The text is addressed to researchers and graduate students interested in the geometry and dynamics of discrete subgroups of semisimple Lie groups acting on their flag manifolds. Our focus is on the asymptotic dynamics: destinations, limit sets and domains of discontinuity. The theory applies to Anosov representations and goes beyond them, for instance to groups with parabolic elements. The examples include classical Kleinian groups and their symmetric powers, complex hyperbolic groups, Hitchin representations (also cusped ones, which contain unipotent elements), Schottky groups and Nori groups.

We assume familiarity with discrete subgroups of $\PSL(2,\C)$ and their classical limit set, together with the basic structure of semisimple Lie groups; the required flag and representation theory is recalled where it enters. We work with arbitrary discrete subgroups: the general results do not require the domain of discontinuity to be nonempty.

The projective case, $\PSL(2,\C)$ acting on the Riemann sphere and, more generally, $\PSL(n+1,\C)$ acting on $\C\PP^n$, is treated in Chapter~\ref{ch:psl} and serves as a guide to the general statements.

\paragraph{How to read this text.} The Introduction states the main results, Theorems~A--E, and can be read on its own. Chapter~\ref{ch:cartan} is the core: it proves the destination formula (Theorem~A) and the stratified convergence principle (Theorem~B), which organises the destinations of a sequence over a stratification of the flag manifold. Every later chapter uses it. Chapters~\ref{ch:limit} and~\ref{ch:balanced} pass from one sequence to a discrete group and construct limit sets and domains of discontinuity. 
Chapters~\ref{ch:adjoint} and~\ref{ch:parabolic} are on the adjoint action and on general flag manifolds. Next, in  Chapter~\ref{ch:psl} we focus on the complex projective space. Chapter~\ref{r2gr:chapter} is addressed to readers interested in the complex algebraic geometry of degenerating graphs. It requires only Chapter~\ref{ch:cartan} and the definitions and properness criterion of Chapter~\ref{ch:limit}, and it also contains comparisons with Chapters~\ref{ch:parabolic} and~\ref{ch:psl}.

\clearpage
\tableofcontents
\clearpage

%==========================================================
\chapter*{Introduction}
\addcontentsline{toc}{chapter}{Introduction}

%==========================================================

Discrete groups of transformations have a tame part and a wild part. For a Kleinian group, a discrete subgroup of $\PSL(2,\C)$ acting on the Riemann sphere, the wild part is the limit set, and its complement is the region of discontinuity. In higher dimension this picture becomes more intricate. We work with a connected real semisimple Lie group $G$ with trivial centre and its flag manifolds $G/Q$, where $Q$ is parabolic. For a discrete subgroup of $G$, several notions of limit set describe different aspects of the action: Kulkarni's set, the Conze--Guivarc'h set, Quint's flag limit set, and the Myrberg-type~\cite{Myrberg} set, the closure of the points with two distinct destinations along a single divergent sequence. These sets need not agree. Whenever a largest domain of proper discontinuity exists, its complement is the recurrent set. On every non-trivial flag manifold of a simple group with trivial centre, the recurrent set is contained in the Kulkarni set, which in turn is contained in the Myrberg set. Quint's flag limit set, which is the Conze--Guivarc'h set of the Tits representation (Proposition~\ref{prop:Tits-comparison}), is contained in the recurrent set. These inclusions need not be equalities, and there may be no largest domain at all. Our aim is to understand their differences and, in many cases, to compute them.

The starting point is a question about a single divergent sequence: where can a point go if its source is allowed to move? Away from the exceptional strata, each source has a single destination. Near those strata, the answer can depend on the \emph{direction} of approach. Theorem~A describes all possible destinations by one formula: the closure of a unipotent orbit, positioned by the bounded Levi datum. The classical convergence principle is its first example.

For a discrete group, the union of the sequence relations determines recurrence and proper discontinuity. The data of individual sequences and the fixed-source set also give formulas for the limit sets (Chapter~\ref{ch:limit}). A finite combinatorial datum then yields domains of discontinuity, which we compare with the Anosov and transverse theories (Chapter~\ref{ch:balanced}). Chapters~\ref{ch:adjoint} and~\ref{ch:parabolic} treat the adjoint action and its relation to the action on $G/P$. Chapter~\ref{ch:psl} returns to $\C\PP^n$, where the critical wall and Kulkarni's set can be compared directly. Chapter~\ref{r2gr:chapter} returns to single sequences over $\C$, whose destinations are realised algebraically by Brion's graph family.

\medskip

We begin with two elementary examples. Let $(g_m)$ be a sequence of M\"obius transformations that leaves every compact subset of $\PSL(2,\C)$. The classical convergence principle says that, after passing to a subsequence, there are two points $x^+,x^-$ of the Riemann sphere $\C\PP^1$ such that
\[
g_m\to x^+\quad\text{uniformly on compact subsets of }\C\PP^1\setminus\{x^-\}.
\]
So the sequence collapses the complement of $x^-$ to the point $x^+$. The reason is visible in a single matrix. Write $g_m=k_m\,\diag(e^{t_m},e^{-t_m})\,k'_m$ with $k_m,k_m'$ unitary and $t_m\to\infty$, and let $k_m\to k$, $k'_m\to k'$. Then $x^+=k[e_1]$ and $x^-=k'^{-1}[e_2]$.

What happens \emph{at} $x^-$ itself? If we freeze the source at $x^-$ the answer depends on how fast $k_m'$ approaches its limit and cannot be read from the limit data alone. If the source is allowed to move, the answer is simple: for any point $y$ of the sphere there are sources $x_m\to x^-$ with $g_mx_m\to y$. Indeed, for $y=k[s:1]$ with $s\in\C$, take $x_m=(k_m')^{-1}[e^{-2t_m}s:1]$; then $g_mx_m=k_m[s:1]\to y$. For $y=k[e_1]$, take $x_m=(k_m')^{-1}[e^{-t_m}:1]$. So the set of limits of $g_mx_m$, over all sources $x_m\to x$, is $\{x^+\}$ if $x\ne x^-$ and the whole sphere if $x=x^-$.
We call this set the \emph{destination} of $x$ and write $\D_{(g_m)}(x)$. Taken together, these destinations form a closed relation on the sphere.

\medskip
The second example shows what changes in higher rank. Let $G=\PSL(3,\R)$ act on $\R\PP^2$ and take
\[
a_m=\diag\bigl(e^{2t_m},\,c\,e^{-t_m},\,c^{-1}e^{-t_m}\bigr),\qquad t_m\to\infty,\quad c\ge1 \text{ fixed}.
\]
If $x=[x_1:x_2:x_3]$ with $x_1\neq0$, then $a_mx_m\to[e_1]$ for every sequence $x_m\to x$. The exceptional set is now the whole line $\{x_1=0\}$, and on it the sources can approach from many directions. Take $x=[0:x_2:x_3]$ and sources $x_m=[s\,e^{-3t_m}:x_2:x_3]$. Then
\[
a_mx_m=[\,s:\,c\,x_2:\,c^{-1}x_3\,]\quad\text{for every }m .
\]
Varying $s$ gives every point of the line through $[e_1]$ and $[0:cx_2:c^{-1}x_3]$ except $[e_1]$, which is reached by taking $s=e^{t_m}$. For any moving source approaching $x$, its last two coordinates after applying $a_m$ remain in a ratio tending to $[cx_2:c^{-1}x_3]$ in $\R\PP^1$, so every destination lies on this line. Two features distinguish the example from the Riemann sphere.
\begin{itemize}[nosep]
\item the destination of an exceptional source is a \emph{line}, not the whole space, and its position depends on the source;
\item if $c>1$ and $x_2x_3\ne0$, the bounded part $\diag(1,c,c^{-1})$ of $a_m$ turns the line through $[e_1]$ and $x$ into the line through $[e_1]$ and $[0:cx_2:c^{-1}x_3]$. The diverging scale $e^{3t_m}$ alone cannot record this change. We call this bounded part the \emph{Levi datum}.
\end{itemize}
In $\PSL(2,\C)$ there is no Levi datum, and this is why the classical principle looks so simple.

\medskip We now state the general destination formula.
Let $G$ be a connected real semisimple Lie group with trivial centre, $G=K\exp(\ak^+)K$ a Cartan decomposition, $\Delta$ the simple restricted roots, and $X=G/Q$ a flag manifold. A sequence $g_m=k_m e^{\mu_m}k'_m$ leaving every compact set has, after extraction, a \emph{type} $\Theta\subset\Delta$ (the simple roots with $\alpha(\mu_m)\to\infty$), limit compact factors $k,k'\in K$ and a \emph{Levi class} $b$, the limit of the bounded part of $e^{\mu_m}$. To $\Theta$ belong a parabolic subgroup $P_\Theta=L_\Theta N_\Theta$ and a stratified retraction $\rho:X\to X^{A_\Theta}$, $\rho(x)=\lim_{z\to\infty}zx$, smooth on each source cell (Chapter~\ref{ch:cartan}). Then:

\medskip

\begin{thmA}
For every $x\in X$,
\[
\D_{(g_m)}(x)\;=\;k\;\overline{\,N_\Theta\cdot b\,\rho(k'x)\,}\,.
\]
That is, the destination of $x$ is the $k$-translate of the closure of the $N_\Theta$-orbit through $b\rho(k'x)$. Every point of the translated closure is reached along the \emph{whole extracted} sequence, and the real dimension of the orbit is $\sum_{\lambda}\dim\g_\lambda$, summed over the restricted roots $\lambda$ which expand $X$ at that cell. No discreteness, density, irreducibility or domain is assumed.
\end{thmA}

The formula contains the convergence principle and the second example above. For complex flag varieties, Brion's graph compactification gives the same destination relation: his degenerations of the diagonal describe the limiting graphs and their Schubert components \cite[Theorem~3 and Proposition~10]{r2grBrion2003}. Here we obtain it directly by constructing moving sources along the whole extracted sequence. The argument applies to every real flag manifold of $G$, including non-reduced restricted root systems and root spaces of any dimension.

\medskip
Theorem~A describes the destination of each point. After Cartan extraction, a divergent sequence determines a \emph{stratification} of $X$ with a family of Schubert varieties over each stratum. For $\PSL(d,\C)$ acting on $\C\PP^{d-1}$ it reads as follows (Theorem~\ref{thm:ext-boundary}). A divergent sequence has, after extraction, an ordered block decomposition $d=b_1+\dots+b_k$, a flag $F_1\supset\dots\supset F_k$ of kernels, a flag $G_1\subset\dots\subset G_k$ of images, and isomorphisms $A_j:F_j/F_{j+1}\to G_j/G_{j-1}$ up to scalar; this graded datum is a point of the exterior boundary (complete collineations), every datum occurs, and for $v\in F_j\setminus F_{j+1}$, $w\in G_j\setminus G_{j-1}$
\[
\D_{(g_m)}([v])=\PP\bigl(\pi_j^{-1}(A_j(\C\,q_jv))\bigr),\qquad\D_{(g_m^{-1})}([w])=\PP\bigl(q_j^{-1}(A_j^{-1}(\C\,\pi_jw))\bigr),
\]
projective subspaces of complex dimensions $b_1+\dots+b_{j-1}$ and $b_{j+1}+\dots+b_k$, respectively, forming holomorphic projective bundles over the corresponding strata. In general:

\begin{thmB}
Let $G$ be a connected real semisimple Lie group with trivial centre, $X=G/Q$ a flag manifold and $(g_m)$ a divergent sequence. Fix a projectively faithful representation $\varrho$ of regular highest restricted weight as in Theorem~\ref{thm:extG}. After extraction, $(g_m)$ has a limit $\xi$ in the exterior boundary associated with $\varrho$. It records the graded datum $(\Theta,\xi^+,\xi^-,[b])$, consisting of the type, attracting and repelling parabolics, and Levi class. It cuts $X$ into source cells $\mathcal O^-_\sigma(\xi)$ and target cells $\mathcal O^+_\tau(\xi)$, $\sigma,\tau\in W_\Theta\backslash W/W_\theta$. For the unique minimal-length representative $\dot\sigma$ of the double coset $\sigma$, put $X_\sigma=\overline{B\dot\sigma Q}$, where $B$ is the standard minimal parabolic and $\dot\sigma$ also denotes a lift to $N_K(\ak)$. For $x\in\mathcal O^-_\sigma(\xi)$ and $y\in\mathcal O^+_\tau(\xi)$,
\[
\D_{(g_m)}(x)=(k\,b\,l)\cdot X_\sigma,\qquad\D_{(g_m^{-1})}(y)=(k'^{-1}b^{-1}l')\cdot\overline{N^-_\Theta\tau Q},
\]
translates of a Schubert variety of dimension $\sum_{\Sigma_+(\sigma)}\dim\g_\lambda$ \textup{(}resp.\ of the closure of a cell of dimension $\sum_{\Sigma_-(\tau)}\dim\g_\lambda$\textup{)}, with $l,l'$ chosen from the retractions using the indicated representatives, forming a bundle over a flag manifold of the Levi group $L_\Theta$. The relation depends only on $\xi$, every datum occurs, and the case $\PSL(d,\C)$ is the display above.
\end{thmB}
The proof, with all the notation, is in Section~\ref{sec:strat}; it combines Theorem~A with the exterior boundary of Section~\ref{sec:extG}. Over $\C$, this is also the fibre description of Brion's limiting graphs. Chapter~\ref{r2gr:chapter} gives the comparison with Brion's graph family, including its scheme structure and Schubert components.

\medskip
A discrete subgroup $\Gamma$ of $G$ has many divergent sequences, and its various notions of limit sets are obtained by choosing which sources and which sequences to look at.
Chapter~\ref{ch:limit} organises these into one picture. The next two statements record the consequences used later. Their abstract properness criterion is classical; the calculation of its forbidden pairs is supplied by Theorem~A.

\begin{thmC}
Let $G$ be a simple real Lie group with trivial centre, $\Gamma<G$ discrete, $X=G/Q$ a non-trivial flag manifold. Then:

(a) The destination relation $\D_\Gamma$ is the union of the relations of Theorem~A over all divergent sequences with stabilised Cartan data, and a $\Gamma$-invariant open set is a domain of discontinuity iff it contains no pair of $\D_\Gamma$; there is a largest domain iff no two non-recurrent points are related. 

(b) Put $B_\Gamma=L_0\cup L_1$, the fixed-source set. Then Kulkarni's limit set is
\[
 \Lambda_{\Kul}=B_\Gamma\cup\overline{\D_\Gamma(X\setminus B_\Gamma)}:
\]
the fixed-source set together with the closure of the destinations of sources outside it. 

(c) Take the sources with at least two distinct destinations along some divergent sequence with stabilised Cartan data. Their closure $\Lambda_{\Myr}$ has complement a domain of discontinuity. 

(d) Quint's flag limit set $\Lambda^Q_\theta(\Gamma)$ is the set of attracting flags of sequences whose type contains $\theta$; if $\Gamma$ is $P_\theta$-divergent and $\theta$ is opposition invariant, let $\rho_\theta$ be the Tits representation with equivariant flag embedding $s_\theta:\mathcal F_\theta\hookrightarrow\PP(\mathcal V_\theta)$ as in Proposition~\ref{prop:Tits-comparison}. The trace of its equicontinuity region on the embedded flag manifold is given by
\[
 s_\theta^{-1}\!\left(\Eq(\rho_\theta(\Gamma))\right)
 =\{F\in\mathcal F_\theta:\ F\text{ is transverse to every }b\in\Lambda_\theta\}.
\] For every such flag manifold, $\Lambda_{\Kul}\subseteq\Lambda_{\Myr}$. For $\PSL(n+1,\C)$ on $\C\PP^n$, $\Lambda_{\Myr}$ is also the complement of the equicontinuity region.
\end{thmC}

\medskip

The limit sets also admit a natural filtration. For $\PSL(n+1,\C)$ these sets sit in one chain (Theorem~\ref{thm:flag-of-limit-sets}): the Conze--Guivarc'h set at the bottom, the Myrberg set at the top, the critical set and Kulkarni's set in between whenever the critical domain is the largest one. The corresponding chain on non-trivial flag manifolds of a centre-free simple real group is as follows.

\begin{thmD}
Let $G$ be a simple real Lie group with trivial centre, $\Gamma<G$ discrete, $P_\theta$ a proper parabolic subgroup, and $X=G/P_\theta$. The \emph{Cartan levels} $\Lambda_{\mathrm Z}^{(d)}$ ({\it i.e.}, the closure of the union of the target cells of dimension $\le d$ of divergent sequences in $\Gamma$)  form an increasing chain of closed invariant sets from Quint's flag limit set $\Lambda_\theta^Q$ $(d=0)$ to the Myrberg set $\Lambda_{\Myr}$ $(d=\dim X-1)$. For a dimension system as in Chapter~\ref{ch:balanced}, the wall $\Lambda_\delta$, the Kulkarni set $\Lambda_{\Kul}$ and $\Lambda_{\Myr}$ have complements which are domains, and contain the recurrent set $\mathscr R_\Gamma$. Either $\Lambda_\delta=\mathscr R_\Gamma$ \textup{(}then $\Omega_\delta$ is the largest domain\textup{)}, or $\mathscr R_\Gamma$ is itself a closed invariant cover and the largest domain is bigger, or two non-recurrent points are related and no largest domain exists. For every discrete $\Gamma<G$,
\[
 \Lambda^Q_\theta\subseteq\mathscr R_\Gamma
 \subseteq\Lambda_{\Kul}\subseteq\Lambda_{\Myr}.
\]
If $\Omega_\delta$ is the largest domain, then $\Lambda_\delta=\mathscr R_\Gamma$ and the dimension-system wall occupies the recurrent level of this chain. No density hypothesis is needed.
\end{thmD}
This is Theorem~\ref{thm:L}; the general comparison of Kulkarni and Myrberg is Proposition~\ref{prop:kulkarni-myrberg}. The kernel description of Myrberg and the identification of the levels with families of projective subspaces are specific to projective space.

The scalar Cartan levels do not in general recover Kulkarni's set: they record where the flag spaces lie, but not which sources are admissible. Keeping that condition gives exact descriptions of Kulkarni's set (Theorems~\ref{r2kul:ideal-level} and~\ref{r2kul:typed-levels}); the excluded sources must be removed before taking the closure. The principal images of Theorem~E below make this distinction explicit.

\medskip
The projective application that guides the second part of the monograph shows that, already for principal images, proper discontinuity and equicontinuity are governed by different multiplicity conditions. For an integer $N\ge1$, let $\rho_N:\PSL(2,\C)\to\PGL(N+1,\C)$ be the projective $N$-th symmetric power representation, acting on the space $X_N=\mathrm{Sym}^N\C\PP^1=\C\PP^N$ of effective divisors of degree $N$.

\medskip

\begin{thmE}
Let $\Gamma<\PSL(2,\C)$ be non-elementary and discrete, with proper classical limit set $\Lambda\subsetneq\C\PP^1$, and let $G=\rho_N(\Gamma)$. Then
$\Lambda_{\Kul}(G)=\Lambda_{\Myr}(G)=\{D\in X_N:\ \operatorname{supp}D\cap\Lambda\neq\varnothing\}$, the region of equicontinuity is $\mathrm{Sym}^N(\C\PP^1\setminus\Lambda)$, but the \emph{largest} open set on which $G$ acts properly discontinuously is
\[
\Omega_{\max}(G)=\{D:\ \operatorname{mult}_\xi D<\lceil N/2\rceil\ \text{for every }\xi\in\Lambda\},
\]
which is strictly larger than the Kulkarni domain for $N\ge3$. No convex cocompactness, geometric finiteness or absence of parabolics is needed.
\end{thmE}
Thus the largest domain need not be the Kulkarni domain. In Theorem~\ref{thm:C}, its computation uses divisors and the classical convergence principle directly. The general destination formula supplies the corresponding language on other flag manifolds. For $N\ge3$, Corollary~\ref{cor:principal-separation} also separates their complex geometry: the largest domain contains a projective line, whereas every component of the equicontinuity domain is complete hyperbolic and Stein.

Theorem~\ref{thm:crit-recurrence} gives a more general recurrence criterion in projective space, relating recurrence to the geometry of the critical flags; its precise hypotheses belong to Section~\ref{sec:crit-proj}. For Zariski-dense cusped Hitchin groups in $\mathrm{SL}_5(\R)$ it applies to their action on $\C\PP^4$, where it gives a nonempty largest domain (Theorem~\ref{r11:cusped-largest}). Such groups exist and contain principal unipotents (Proposition~\ref{r11:cusped-existence}); they are outside the traditional Borel--Anosov setting.

\medskip Four works serve as reference points for the relation with the literature. The table lists what each assumes and how it sits in our setting.

\begin{center}\small
\begin{tabular}{@{}p{.22\linewidth}p{.33\linewidth}p{.38\linewidth}@{}}
\toprule
Source & Hypothesis and output & Relation to this paper\\\midrule
Labourie~\cite{r2appLabourie2006} (Hitchin component) & Anosov, hyperconvex Frenet curve; discrete, faithful, loxodromic & Model family with $\Theta=\Delta$ for every divergent sequence; limit flags $\xi^\pm$ are the osculating flags.\\
Guichard--Wienhard~\cite{GuichardWienhard} & Anosov representations; domains of discontinuity in flag manifolds & Their domains are instances of ``complement of a union of fibres of $\D$''.\\
Kapovich--Leeb--Porti~\cite{KLP} & Arbitrary discrete groups; domains from thickenings over all face types; no antipodality for properness & Chapter~\ref{ch:balanced} expresses their construction through destination fibres; expansion, antipodality and balance enter their cocompactness theorem.\\
Canary--Zhang--Zimmer~\cite{mlCZZ}& transverse groups; convergence action on the limit set & Our setting has no antipodality; cusped examples are the test in Proposition~\ref{prop:unipotent}.\\\bottomrule
\end{tabular}
\end{center}
The questions developed here are: 

(a)~recurrence criteria for a \emph{largest} domain, and examples where it is larger than the Kulkarni domain;

(b)~the destinations of sequences of different Cartan types in the same group, including groups with parabolic elements;

(c)~a single explicit formula for what happens near the exceptional strata;

(d)~a fibre-by-fibre dictionary with the thickenings of Kapovich--Leeb--Porti, including their construction for arbitrary discrete groups (Section~\ref{sec:KLPfull}); powers of principal unipotents, whose attracting and repelling flags coincide, are treated by the same formula (Proposition~\ref{prop:unipotent}).

%==========================================================
%==========================================================
We assume familiarity with discrete subgroups of $\PSL(2,\C)$ and their classical limit sets, together with the basic structure of semisimple Lie groups. We refer to  \cite{CNS, r2appBCNSParabolic, BCNSFlag2026} for basics about limit sets for projective groups. The required flag and representation theory is recalled where it enters. Our concern is limit behaviour and proper discontinuity. Cocompactness, Patterson--Sullivan measures, Hausdorff dimension and the geometry of Hitchin representations are used only through the specific results cited below. The general results allow the domain of discontinuity to be empty.

\chapter{Cartan dynamics}
\label{ch:cartan}

In higher rank, different directions contract at different rates, and the exceptional locus has a geometry of its own. This chapter shows that the Cartan decomposition describes all the destinations that moving sources can reach. For a sequence $g_m=k_me^{\mu_m}k'_m$ of Cartan type $\Theta$ with limit flags $k,k'$ and Levi class $b$, and for every point $x$ of every flag manifold $X$ of $G$, Theorem~\ref{thm:dest} gives the set of limits of $g_mx_m$ over all sources $x_m\to x$:
\[
\D_{(g_m)}(x)=k\,\overline{N_\Theta\cdot b\,\rho(k'x)},
\]
the closure of an orbit of a nilpotent group. The datum $(\Theta,k,k',b)$ will become a point of a compact space (Sections~\ref{sec:ext}--\ref{sec:extG}); conversely, every such datum comes from a sequence.

The reason is visible in one linear-algebra fact: conjugation by $e^{\mu_m}$ multiplies the root space $\g_\lambda$ by the number $e^{\lambda(\mu_m)}$. Around each fixed point of the maximal torus, a flag manifold looks like a nilpotent group, whose root spaces are either expanded, contracted, or left alone by $e^{\mu_m}$; the destination is what the expanded part can reach. Sections~\ref{sec:classical} and~\ref{sec:sl3} give the precise statements and proofs for the two examples of the Introduction.

%--------------------------------------------------------------
\section{The classical principle, with a proof}\label{sec:classical}
We begin with the classical convergence principle. Its proof introduces the moving-source construction used below, so we give it in full.

\begin{proposition}[Convergence principle]\label{prop:classical}
Let $g_m\in\PSL(2,\C)$ leave every compact set. After extraction, $g_m=k_m\diag(e^{t_m},e^{-t_m})k'_m$ with $t_m\to\infty$, $k_m\to k$, $k'_m\to k'$ in $\mathrm{PSU}(2)$, and for every $x\in\C\PP^1$
\[
\D_{(g_m)}(x)=\begin{cases}\{k[e_1]\},&x\ne k'^{-1}[e_2],\\ \C\PP^1,&x=k'^{-1}[e_2].\end{cases}
\]
\end{proposition}
Here $\D_{(g_m)}(x)$ is the destination of $x$, the set of all limits $\lim g_mx_m$ over sequences of sources $x_m\to x$, as in the Introduction; the formal definition is in Section~\ref{sec:thmA}.
\begin{proof}
Let $a_m=\diag(e^{t_m},e^{-t_m})$ and write points of $\C\PP^1$ as $[u:1]$ or $[e_1]=[1:0]$. For $v_m\to[u:1]$ with $u\ne0$ the images are $a_mv_m=[e^{2t_m}u_m:1]\to[e_1]$, and the same holds for $v_m\to[e_1]$. If $v_m\to[e_2]=[0:1]$, write $v_m=[w_m:1]$ with $w_m\to0$, so $a_mv_m=[e^{2t_m}w_m:1]$. Given $s\in\C$, choose $w_m=e^{-2t_m}s$: then $a_mv_m=[s:1]$. To reach $[e_1]$, choose $w_m=e^{-t_m}$, so $w_m\to0$ and $e^{2t_m}w_m=e^{t_m}\to\infty$. Hence $\D_{(a_m)}([e_2])=\C\PP^1$. For general $g_m$, compose with $k_m\to k$ and $k'_m\to k'$, which converge uniformly on the compact space $\C\PP^1$.
\end{proof}
Three things are worth noticing. First, $a_m$ contracts vectors in $\C e_2$. The projective point $[e_2]$ is a repelling fixed point of $a_m$ and an attracting fixed point of $a_m^{-1}$. Second, near it, the image is controlled by the product $e^{2t_m}w_m$ of an enormous number and a tiny one, and that product can be anything. This is the ``indeterminacy'' that the rest of the monograph organises. Third, the answer ``anything'' is an orbit closure: the orbit of the upper unipotent group $\left\{\begin{pmatrix}1&s\\0&1\end{pmatrix}\right\}$ through $[0:1]$ is $\C\PP^1\setminus\{[e_1]\}$, whose closure is the whole sphere.

%--------------------------------------------------------------
\section{A second example: what the Levi datum does}\label{sec:sl3}
We record the full statement for the sequence $a_m=\diag(e^{2t_m},ce^{-t_m},c^{-1}e^{-t_m})$ on $\R\PP^2$ of the Introduction. It is the model for the general formula: the fibre is a line, not a point, but a line that depends on a finite matrix, namely $\diag(1,c,c^{-1})$.

\begin{proposition}\label{prop:sl3}
For this sequence, $\D_{(a_m)}([x_1:x_2:x_3])=\{[e_1]\}$ if $x_1\ne0$, and
\[
\D_{(a_m)}([0:x_2:x_3])=\PP\bigl(\R e_1+\R(0,cx_2,c^{-1}x_3)\bigr).
\]
\end{proposition}
\begin{proof}
If $x_1\neq0$ and $v_m\to x$, then $a_mv_m=[e^{3t_m}v_{m,1}:cv_{m,2}:c^{-1}v_{m,3}]\to[e_1]$. If $x_1=0$ and $v_m=[\varepsilon_m:x_2^{(m)}:x_3^{(m)}]\to x$, then $a_mv_m=[e^{3t_m}\varepsilon_m:cx_2^{(m)}:c^{-1}x_3^{(m)}]$. Every limit lies on the projective line $\PP(\R e_1+\R(0,cx_2,c^{-1}x_3))$. Its finite points $[s:cx_2:c^{-1}x_3]$, $s\in\R$, occur for $\varepsilon_m=se^{-3t_m}$ and $x_j^{(m)}=x_j$; its remaining point $[e_1]$ occurs for $\varepsilon_m=e^{-t_m}$, since $e^{3t_m}\varepsilon_m=e^{2t_m}\to\infty$.
\end{proof}
Both examples have the same shape: a unipotent group fixing a point ($[e_1]$), whose orbit closures are the destinations, and a finite diagonal matrix ($\diag(1,c,c^{-1})$) that acts on the source \emph{before} the unipotent group does. If $x_2x_3\ne0$, changing $c$ turns the line $\PP(\R e_1+\R(0,cx_2,c^{-1}x_3))$ around $[e_1]$ inside the plane. This is invisible in rank one, where the matrix $b$ below is trivial, and it is why the Levi class has to be part of the statement.

%--------------------------------------------------------------
\section{The setting}\label{sec:setting-sum}
The restricted roots replace the gaps between singular values, and the Levi factor carries the bounded motion. We now give names to these objects for a general semisimple group.

\paragraph{The group.} $G$ is the identity component of $\mathrm{Aut}(\g)$ for a real semisimple Lie algebra $\g$ (so $G$ has trivial centre). Let $\g=\kk\oplus\pk$ be a Cartan decomposition, $K$ the maximal compact subgroup, $\ak\subset\pk$ maximal abelian, $\Phi\subset\ak^*$ the \emph{restricted} roots with root spaces $\g_\lambda$ of any dimension (the system may be non-reduced, as for $\mathrm{SU}(p,q)$, $p\ne q$), $\Delta$ a set of simple roots, $\ak^+$ the closed chamber, and $G=K\exp(\ak^+)K$.

\paragraph{Types.} For $\Theta\subset\Delta$ put $\ak_\Theta=\bigcap_{\alpha\notin\Theta}\ker\alpha$, $A_\Theta=\exp\ak_\Theta$, $L_\Theta=Z_G(A_\Theta)$, $N_\Theta$ the group whose Lie algebra is spanned by the root spaces of the positive roots \emph{not} in the span of $\Delta\setminus\Theta$, $P_\Theta=L_\Theta N_\Theta$, and $P_\Theta^-=L_\Theta N_\Theta^-$ the opposite parabolic. Let $W=N_K(\ak)/Z_K(\ak)$ be the Weyl group and $W_\Theta\subset W$ the Weyl group of the Levi factor $L_\Theta$, generated by the reflections in the roots of $\Delta\setminus\Theta$. The flag manifold is $X=G/Q$ with $Q=P_\theta$ for some $\theta\subset\Delta$, and $W_\theta$ is defined likewise. \emph{We use $\theta$ for the flag manifold and $\Theta$ for the type of a sequence; they are unrelated.} (The extreme cases: $\Theta=\Delta$ gives $P_\Theta=B$, the minimal parabolic; $\theta=\Delta$ gives the full flag manifold $G/B$; $\theta=\{\alpha\}$ for a single $\alpha$ gives a ``Grassmannian-like'' $G/P_\alpha$.)

The Cartan formulas are stated for sequences in $G$. A finite-index extension can be used for the group action once an action on $X$ has been specified. For an extension by automorphisms of $G$, the induced flag action is defined when the automorphisms preserve the conjugacy class of $Q$: identify $X$ with the space of conjugates of $Q$. The relation for the extension is obtained by translating the relations from its finitely many cosets, as in Section~\ref{sec:Dgroup}.

\paragraph{The projective model.} $\PSL(n+1,\C)$, viewed as a real group, has $\ak$ the real diagonal matrices, restricted roots of type $A_n$ with multiplicity $2$, and $X_\theta$ the complex flag manifolds.
\begin{center}\small
\begin{tabular}{@{}p{.31\linewidth}p{.31\linewidth}p{.28\linewidth}@{}}
\toprule
$\PSL(n+1,\C)$ & General $G$ & $\PSL(2,\C)$\\\midrule
Gaps $\sigma_j/\sigma_{j+1}\to\infty$; block decomposition of $n+1$ & Type $\Theta=\{\alpha:\alpha(\mu(g_m))\to\infty\}$ & single gap, $\Theta=\Delta$\\
Source and target flags & Limit parabolics in the frames $k,k'\in K$ & points $x^\mp,x^+$\\
Graded isomorphisms of the blocks & Relative class in the Levi fibre $L_\Theta/Z_\Theta$ & none ($L$ is a torus)\\
$\C\PP^n$ & $X_\theta$ for all $\theta$ & $\C\PP^1$\\\bottomrule
\end{tabular}
\end{center}
This table is the dictionary of the paper. Whenever a general statement looks abstract, translate it into the first column.

%--------------------------------------------------------------
\section{The type, the Levi class and the compact factors}
We separate the divergent directions, which determine the unipotent group, from the bounded part that turns the destination.

Let $(g_m)$ leave every compact set and write $g_m=k_me^{\mu_m}k'_m$ with $\mu_m\in\ak^+$.
After extraction, $k_m\to k$, $k'_m\to k'$, and each $\alpha(\mu_m)$ converges in $[0,\infty]$. Put $\Theta=\{\alpha\in\Delta:\alpha(\mu_m)\to\infty\}$. Since $\Delta$ is a basis of $\ak^*$ and $(g_m)$ diverges, $\Theta\ne\varnothing$.

\begin{lemma}\label{lem:split}\
Write $\mu_m=\nu_m+\eta_m$ with $\nu_m\in\ak_\Theta$ and $\eta_m$ in the span of the coroots of $\Delta\setminus\Theta$. Then $\eta_m\to\eta$ and $\beta(\nu_m)\to\infty$ for $\beta\in\Theta$. Put $z_m=\exp\nu_m$ (central in $L_\Theta$), $b_m=\exp\eta_m\to b\in L_\Theta$ and $a_m=e^{\mu_m}=z_mb_m$. For a restricted root $\lambda$:
\[
\lambda(\nu_m)\to+\infty\ (\lambda\in\Phi(N_\Theta)),\qquad
\lambda(\nu_m)\to-\infty\ (\lambda\in\Phi(N_\Theta^-)),\qquad
\lambda(\nu_m)=0\ (\lambda\in\Phi(L_\Theta)).
\]
\end{lemma}
\begin{proof}
For $\alpha\notin\Theta$ we have $\alpha(\nu_m)=0$, so $\alpha(\eta_m)=\alpha(\mu_m)$ is bounded; the restrictions of these $\alpha$ to the span of their coroots are linearly independent (the Cartan matrix of the Levi factor is invertible), so $\eta_m$ converges because all these coordinates $\alpha(\mu_m)$ already converge. For $\beta\in\Theta$, $\beta(\nu_m)=\beta(\mu_m)-\beta(\eta_m)\to\infty$. A positive root $\lambda=\sum n_\alpha\alpha$ lies in $\Phi(N_\Theta)$ iff $n_\beta>0$ for some $\beta\in\Theta$, and then $\lambda(\nu_m)=\sum_{\beta\in\Theta}n_\beta\beta(\nu_m)\to\infty$. The other two cases follow in the same way.
\end{proof}
So $g_m=k_mz_mb_mk'_m$. The datum $(\Theta,k,k',b)$ remembers the divergent gaps, the limit frames and the bounded ratios. Theorem~\ref{thm:dest} shows that this is the only information that survives on a flag manifold: the speeds of divergence are forgotten.

%--------------------------------------------------------------
\section{Cells and the stratified retraction}
We cut the flag manifold into cells with the same contracting directions. Within each cell, contraction leaves a point on a flag manifold of the Levi factor, through which the expanded orbit will pass.

Let $X=G/Q$ and fix representatives $\sigma\in N_K(\ak)$ of the double cosets $W_\Theta\backslash W/W_\theta$, each of minimal length in $W_\Theta\sigma$. Since $A\subset\sigma Q\sigma^{-1}$, $A$ fixes each $\sigma Q$.

\begin{lemma}\label{lem:cells}\ (Bruhat; Bia\l ynicki-Birula)
\begin{enumerate}[label=(\roman*),nosep]
\item $X=\bigsqcup_\sigma O_\sigma$, with $O_\sigma=P_\Theta^-\,\sigma Q/Q$.
\item The fixed point set $X^{A_\Theta}=\bigsqcup_\sigma L_\Theta\sigma Q/Q$ is a disjoint union of compact flag manifolds of $L_\Theta$.
\item For $x=n\,l\,\sigma Q\in O_\sigma$ ($n\in N_\Theta^-$, $l\in L_\Theta$) the limit $\rho(x):=\lim zx=l\,\sigma Q$ exists as $z\to\infty$ in $A_\Theta$ (meaning $\beta(z)\to\infty$ for all $\beta\in\Theta$) and does not depend on how $z\to\infty$.
\item Let $\Sigma$ be the set of roots of the unipotent group $N_\sigma=\sigma N_\theta^-\sigma^{-1}$. The map $n\mapsto n\sigma Q$ is an $A$-equivariant chart of a neighbourhood of $\sigma Q$. Split $\Sigma=\Sigma_+\sqcup\Sigma_0\sqcup\Sigma_-$ according to $\Phi(N_\Theta)$, $\Phi(L_\Theta)$, $\Phi(N_\Theta^-)$, and put $M_+=\exp\bigoplus_{\Sigma_+}\g_\lambda$, $M_-=\exp\bigoplus_{\Sigma_0\cup\Sigma_-}\g_\lambda$. Both are subgroups, $N_\sigma=M_-M_+$ with unique decomposition, and $N_\Theta\,\sigma Q=M_+\,\sigma Q$.
\end{enumerate}
\end{lemma}
On each cell $O_\sigma$, the restriction of $\rho$ is a smooth retraction onto $L_\Theta\sigma Q/Q$, holomorphic when $G$ is complex. The map $\rho$ need not be continuous across cells; throughout, the term ``retraction'' refers to these restrictions.
\emph{Why the lemma holds.} Conjugation by $z$ contracts $N_\Theta^-$ and fixes $L_\Theta\sigma Q$, giving (iii). For (iv), the root sets $\Sigma_+$ and $\Sigma_0\cup\Sigma_-$ are closed, and the Baker--Campbell--Hausdorff formula, ordered by root height, gives the unique factorisation $N_\sigma=M_-M_+$. The full proof in Section~\ref{sec:proofcells} treats the whole restricted root spaces, including double roots; no split or reduced hypothesis is used.

\paragraph{A picture for $\R\PP^2$.} Take $\PSL(3,\R)$, $X=\R\PP^2$, and $\Theta=\{\alpha_1\}$ (the type of Proposition~\ref{prop:sl3}). There are two Weyl double cosets: $\sigma=1$ with $\sigma Q=[e_1]$, and $\sigma=s$ with $\sigma Q=$ a point of the line $\{x_1=0\}$. The cell $O_1=P_\Theta^-[e_1]$ is the affine plane $\{x_1\ne0\}$, which the map $\rho$ sends to $[e_1]$ (the sources whose limit is a single point). The cell $O_s=P^-_\Theta[e_2]$ is the line $\{x_1=0\}$, fixed pointwise by $\rho$. The set $\Sigma$ for $\sigma=s$ has two roots, one in $\Phi(N_\Theta)$ and one in $\Phi(L_\Theta)$: this is why the fibre over the line is one-dimensional.

%--------------------------------------------------------------
\section{The destination formula}\label{sec:thmA}
For a sequence $f_m$ of homeomorphisms of a compact space $X$, let
\[
\D_{(f_m)}(x)=\{y:\ \exists\,x_m\to x\text{ with }f_mx_m\to y\},
\]
where $m$ runs over the \emph{whole} given sequence. Extraction is made once and for all at the level of the Cartan data.

Allowing the sources to move brings us to the theory of dynamically related
points of Frances~\cite{Frances}. Cano--Liu--L\'opez explicitly define
$\D_{(g_m)}(x)$ and compute its projective fibres for complex hyperbolic groups
\cite[arXiv:1506.08113v5, Definition~3.2 and Proposition~3.6]{CLL}.
M\'endez studies the sequence and group relations, shows how they behave under
converging compact factors, and analyses them on a complex quadric
\cite[arXiv:1809.06953v1, \S4, Proposition~4.1 and Propositions~4.3--4.6]{Mendez2020}.
We now compute these fibres on every flag manifold, with no further extraction
depending on the source or its destination.

\begin{theorem}[Destination formula]\label{thm:dest}\
For a sequence $g_m=k_ma_mk'_m$ of type $\Theta$ with data $(k,k',b)$ as in Lemma~\ref{lem:split}, and every $x\in X$:
\begin{enumerate}[label=(\alph*),nosep]
\item $\D_{(a_m)}(x)=\overline{N_\Theta\cdot b\,\rho(x)}$; every point of the closure is reached along the whole sequence $(a_m)$.
\item $\D_{(g_m)}(x)=k\,\D_{(a_m)}(k'x)=k\,\overline{N_\Theta\,b\rho(k'x)}$.
\item If $k'x\in O_\sigma$ the fibre is the closure of a cell of real dimension $\sum_{\lambda\in\Sigma_+}\dim\g_\lambda$.
\item The graphs of $g_m$ converge in the Hausdorff topology to the closed relation $R=\{(x,y):y\in\D_{(g_m)}(x)\}$, and $\D_{(g_m^{-1})}$ is the transpose of $\D_{(g_m)}$.
\end{enumerate}
\end{theorem}

\begin{proof}
\emph{Step 1: the base point.} Take $x=\sigma Q$. A source $x_m\to\sigma Q$ is, by Lemma~\ref{lem:cells}(iv), $x_m=q_mp_m\sigma Q$ with $q_m\in M_-$, $p_m\in M_+$, and $q_m,p_m\to1$. On each root space of $\mathrm{Lie}(M_-)$, $\Ad(a_m)$ is multiplication by $e^{\lambda(\mu_m)}$, which is bounded (for $\Sigma_0$) or tends to $0$ (for $\Sigma_-$). There are only finitely many root spaces, so these restrictions have a common operator-norm bound; since $\log q_m\to0$, we have $\Ad(a_m)\log q_m\to0$. Hence $c_m:=a_mq_ma_m^{-1}\to1$, and since $a_m$ fixes $\sigma Q$,
\[
a_mx_m=c_m\,(a_mp_ma_m^{-1})\,\sigma Q.
\]
The second factor lies in $M_+\sigma Q=N_\Theta\sigma Q$, and $c_m\to1$ acts on the compact space $X$ uniformly close to the identity. So every limit point of $a_mx_m$ lies in $\overline{N_\Theta\sigma Q}$. This is true for any subsequence.

Conversely, let $y\in\overline{N_\Theta\sigma Q}$ and choose $p^{(j)}\in M_+$ with $p^{(j)}\sigma Q\to y$. For fixed $j$, $u^{(j)}_m:=a_m^{-1}p^{(j)}a_m\to1$ as $m\to\infty$, because $a_m^{-1}$ multiplies every root space of $\mathrm{Lie}(M_+)$ by $e^{-\lambda(\mu_m)}\to0$. Choose integers $M_1<M_2<\cdots$ with $d(u_m^{(j)},1)<1/j$ for $m\ge M_j$, and put $x_m=u^{(j)}_m\sigma Q$ for $M_j\le m<M_{j+1}$. Then $x_m\to\sigma Q$ and $a_mx_m=p^{(j)}\sigma Q\to y$. This is a construction along the whole sequence. Hence $\D_{(a_m)}(\sigma Q)=\overline{N_\Theta\sigma Q}$.

\emph{Step 2: an arbitrary point.} Let $x=h\sigma Q$ with $h=nl$, $n\in N_\Theta^-$, $l\in L_\Theta$. Since $z_m$ is central in $L_\Theta$ and contracts $N_\Theta^-$,
\[
a_mha_m^{-1}=(a_mna_m^{-1})(b_mlb_m^{-1})\longrightarrow c:=blb^{-1}.
\]
Sources $x_m\to x$ correspond to $h^{-1}x_m\to\sigma Q$, and $a_mx_m=(a_mha_m^{-1})\,a_m(h^{-1}x_m)$ with the first factor converging uniformly to $c$. So $\D_{(a_m)}(x)=c\,\D_{(a_m)}(\sigma Q)=c\,\overline{N_\Theta\sigma Q}$. Now $c\in L_\Theta$ normalises $N_\Theta$, and $b^{-1}$ fixes $\sigma Q$ (as $b\in A$), so
\[
c\,\overline{N_\Theta\sigma Q}=\overline{N_\Theta\,c\,\sigma Q}=\overline{N_\Theta\,b\,l\,\sigma Q}=\overline{N_\Theta\,b\rho(x)},
\]
which proves~(a).

\emph{Step 3: compact factors.} With $v_m=k'_mx_m$ we have $x_m\to x$ iff $v_m\to k'x$, and $g_mx_m=k_ma_mv_m$; since $k_m\to k$ uniformly, this gives (b).

\emph{Step 4: dimension and graphs.} The orbit $N_\Theta b\rho(x)$ is the translate by $bl\in L_\Theta$ of $N_\Theta\sigma Q=M_+\sigma Q$, which is a cell of dimension $\sum_{\Sigma_+}\dim\g_\lambda$ by the chart. It is semialgebraic by the polynomial nilpotent exponential and the algebraic flag action, so its closure has the same dimension. For (d), let $R^*$ be the set of cluster points of the graphs along arbitrary subsequences. In the compact metric space $X\times X$ it is closed, since $R^*=\bigcap_{M\ge1}\overline{\bigcup_{m\ge M}\operatorname{graph}(g_m)}$. The upper inclusion in Steps~1--3 applies to every subsequence, so $R^*\subset R$; the whole-sequence construction gives $R\subset R^*$. Thus $R$ is closed. If the graphs failed to lie in an $\varepsilon$-neighbourhood of $R$, compactness would give a cluster point outside $R$. Conversely, take a finite $\varepsilon$-net in $R$ and use the whole-sequence construction at its finitely many points: all are approached for every sufficiently large $m$. These two inclusions give Hausdorff convergence. Finally, $x_m\to x$ and $g_mx_m\to y$ iff $y_m=g_mx_m\to y$ and $g_m^{-1}y_m\to x$, which gives the transpose.
\end{proof}

\begin{remark}\label{rem:destA}
(1)~The relation depends on $(g_m)$ only through $\Theta$, $k$, $k'$ and $b$; the finer comparison of the rates $\alpha(\mu_m)$, $\alpha\in\Theta$, is invisible on a flag manifold. (2)~The convention matters: with $g_m=k_ma_mk'_m$ the compact factor sends the source to $k'x$, not $k'^{-1}x$. (3)~The theorem computes the closure of a strong unstable leaf; the closure need not be a manifold. (4)~Nothing is assumed about $\Gamma$: discreteness enters only in Chapter~\ref{ch:limit}. (5)~``Along the whole sequence'' is what lets us pass from sequences to \emph{groups} in Chapter~\ref{ch:limit}: if a statement holds for all pairwise distinct sequences, it holds for subsequences.
\end{remark}

\paragraph{Computing a destination.} Given $(g_m)$:
\begin{enumerate}[label=\arabic*.,nosep]
\item Write $g_m=k_ma_mk'_m$ and extract so that $k_m\to k$, $k'_m\to k'$ and the ratios converge. Read off $\Theta$ and $b$.
\item Move the source: $y=k'x$. Find the cell $O_\sigma$ containing $y$ and compute $\rho(y)=l\sigma Q$.
\item The destination is $k$ times the closure of the $N_\Theta$-orbit of $b\rho(y)$.
\end{enumerate}
In projective space step~3 is linear algebra, as Proposition~\ref{prop:proj} shows.

%--------------------------------------------------------------
\section{Two consequences: a conservation law and rank one}\label{sec:consequences}
The examples suggest that a smaller set of exceptional sources allows a larger destination. Part~(c) of the theorem makes this precise: the dimension lost by the source cell is exactly the dimension gained by the fibre.

\begin{corollary}[Conservation of dimension]\label{cor:conservation}
If $k'x\in O_\sigma$ then
\[
\dim_{\R}\D_{(g_m)}(x)+\dim_{\R}O_\sigma=\dim_{\R}X.
\]
In particular $\D_{(g_m)}(x)$ is a single point if and only if $k'x$ lies in the open cell, and the set $\mathcal E$ of \emph{exceptional} sources (those with non-trivial destination) is $k'^{-1}$ of the complement of the open cell, a proper algebraic subvariety.
\end{corollary}
\begin{proof}
The chart in Lemma~\ref{lem:cells}(iv) gives $\dim X=\sum_\Sigma\dim\g_\lambda$. The tangent-space calculation in Section~\ref{sec:proofcells} gives $\dim O_\sigma=\sum_{\Sigma_0\cup\Sigma_-}\dim\g_\lambda$, and Theorem~\ref{thm:dest}(c) gives the fibre dimension $\sum_{\Sigma_+}\dim\g_\lambda$. The complement of the open source cell is a finite union of opposite Schubert varieties, hence is algebraic.
\end{proof}
The law concerns \emph{dimension}; the orbit closure can still be singular. In $\PSL(2,\C)$ the exceptional point has fibre of dimension $2$, and all other fibres have dimension $0$.

\begin{corollary}[Rank one]\label{cor:rank-one}
If $G$ has real rank one and $P$ is a proper parabolic subgroup, then $X=G/P$ is the sphere at infinity of the symmetric space, $\Theta=\Delta$, there are two cells (a point $\{p_-\}$ and the open cell), and for every divergent sequence with Cartan data as in Lemma~\ref{lem:split}
\[
\D_{(g_m)}(x)=\begin{cases}\{k\,p_+\},&k'x\ne p_-,\\ X,&k'x=p_-.\end{cases}
\]
\end{corollary}
\begin{proof}
Here $Q=P_\Theta$ is the minimal parabolic, $p_+=eQ$ is fixed by $N$, and $p_-=wQ$ is the fixed point of $P^-$. The cells are $O_1=P^-eQ=X\setminus\{p_-\}$ (open) and $O_w=\{p_-\}$. The group $N$ acts simply transitively on $X\setminus\{p_+\}$, so $\rho(x)=p_+$ for $x\ne p_-$ and $\overline{N\,b\,p_-}=X$. In the normalisation of Lemma~\ref{lem:split}, $\Theta=\Delta$ gives $\eta_m=0$, so $b=1$.
\end{proof}
This is the convergence property of the boundary of a rank-one symmetric space. The same destination formula holds on the perfect boundary of a non-elementary hyperbolic group after passing to a convergence subsequence with attractor $x^+$ and repeller $x^-$. For $y\ne x^+$ take $x_m=g_m^{-1}y\to x^-$. To reach $x^+$, choose targets $y_j\ne x^+$ tending to $x^+$ and use $x_m=g_m^{-1}y_{j(m)}$ with $j(m)\to\infty$ sufficiently slowly. In higher rank, several divergent root directions and a residual Levi class replace this two-point picture.

%--------------------------------------------------------------
\section{Projective space}\label{sec:proj}
For $\PSL(n+1,\mathbb K)$ acting on $\mathbb K\PP^n$ ($\mathbb K=\R$ or $\C$), we can see the destination directly in the coordinates of a vector. We give this computation both as an independent proof of Theorem~\ref{thm:dest} in projective space and because these coordinates will return throughout the paper.

Let $V=\mathbb K^{n+1}=V_1\oplus\cdots\oplus V_r$ be a block decomposition into coordinate subspaces of dimensions $d_1,\dots,d_r$. Let $a_m$ be the diagonal matrix with $a_m|_{V_i}=z_{i,m}\,b_{i,m}$, where $z_{i,m}>0$ are scalars with $z_{i+1,m}/z_{i,m}\to0$ and $b_{i,m}$ are positive diagonal matrices converging to invertible positive diagonal matrices $b_i$. For a Cartan sequence one can take $z_{i,m}$ to be the largest singular value of its $i$th block. Bounded ratios within each block then keep the normalised entries uniformly away from zero, and after extraction the matrices $b_{i,m}$ converge. Thus this form covers every extracted Cartan sequence of type $\Theta=\{\alpha_{d_1},\alpha_{d_1+d_2},\dots\}$, with the simple roots at the block boundaries.

\begin{proposition}[Destinations on projective space]\label{prop:proj}
Let $x=[y]\in\mathbb K\PP^n$ with $y=(y_1,\dots,y_r)$, $y_i\in V_i$, and let $i_0=\min\{i:y_i\ne0\}$. Then
\[
\D_{(a_m)}(x)=\PP\bigl(V_1\oplus\cdots\oplus V_{i_0-1}\oplus\mathbb K\,b_{i_0}y_{i_0}\bigr).
\]
In particular the fibre is a projective subspace of $\mathbb K$-dimension $d_1+\cdots+d_{i_0-1}$; it is a point exactly when $i_0=1$.
\end{proposition}
\begin{proof}
Let $x_m\to x$ and choose representatives $y^{(m)}\to y$. Divide $a_my^{(m)}$ by $z_{i_0,m}$. The components of $y^{(m)}$ in blocks $j>i_0$ are multiplied by $z_{j,m}/z_{i_0,m}\to0$ (times bounded factors) and disappear. The component in the block $i_0$ converges to $b_{i_0}y_{i_0}\ne0$. The components in blocks $j<i_0$ are $(z_{j,m}/z_{i_0,m})b_{j,m}y^{(m)}_j$, a product of a number tending to infinity with a vector tending to $0$, with no constraint. Any limit point of $[a_my^{(m)}]$ is therefore of the form $[w+b_{i_0}y_{i_0}]$ with $w$ in the limit of $V_{<i_0}$, possibly at infinity, and so lies in the displayed projective subspace, which is closed.

Conversely, let $w\in V_{<i_0}$. Put $y^{(m)}_j=(z_{i_0,m}/z_{j,m})\,b_{j,m}^{-1}w_j\to0$ for $j<i_0$ and $y^{(m)}_j=y_j$ for $j\ge i_0$. Then $a_my^{(m)}/z_{i_0,m}\to w+b_{i_0}y_{i_0}$ along the whole sequence. Points at infinity of this affine chart are obtained by letting $|w|\to\infty$ slowly (diagonal choice), exactly as in Theorem~\ref{thm:dest}, Step~1.
\end{proof}
\begin{remark}
This is Theorem~\ref{thm:dest} in the following dictionary: $\rho(x)=[0:\dots:0:y_{i_0}:0\dots]$, the group $N_\Theta$ consists of the block-upper-unipotent matrices, and the orbit of $b_{i_0}y_{i_0}$ is $V_{<i_0}+b_{i_0}y_{i_0}$, an affine subspace whose projective closure adds $\PP(V_{<i_0})$. The conservation law of Corollary~\ref{cor:conservation} reads: the set of sources with a given $i_0$ is a locally closed subset of $\mathbb K$-dimension $n-(d_1+\cdots+d_{i_0-1})$, and the fibre has dimension $d_1+\cdots+d_{i_0-1}$. The case $\PSL(3,\R)$, with blocks $(1,2)$, is Proposition~\ref{prop:sl3}.
\end{remark}

%--------------------------------------------------------------
\section{Examples in higher rank}\label{sec:examples}
In projective space every destination is a projective subspace. On other flag manifolds the same root calculation produces different shapes, sometimes singular ones. We begin with the two small examples and then watch what changes when the source is a plane or a complete flag.

\begin{example}[$\PSL(2,\mathbb K)$ on $\mathbb K\PP^1$, $\mathbb K=\R,\C$]
One positive restricted root, of multiplicity $1$ or $2$; $N_\Theta$ is the upper unipotent group, which fixes $[e_1]$ and is transitive on the complement. $\rho(x)=[e_1]$ for $x\ne[e_2]$ and $\rho([e_2])=[e_2]$. The theorem gives Proposition~\ref{prop:classical}, with fibres of real dimension $0$ and $1$ or $2$.
\end{example}
\begin{example}[$\PSL(3,\R)$ on $\R\PP^2$]
Here $\Theta=\{\alpha_1\}$, $b=\diag(1,c,c^{-1})$, $N_\Theta=\{I+ue_{12}+ve_{13}\}$. For $x=[0:x_2:x_3]$ we have $\rho(x)=x$, and the orbit of $b\rho(x)=[0:cx_2:c^{-1}x_3]$ is $\{[s:cx_2:c^{-1}x_3]\}$, whose closure adds $[e_1]$. This is Proposition~\ref{prop:sl3}. Over $\C$ the same calculation holds with root multiplicity $2$ and twice the real dimensions.
\end{example}
\begin{example}[$\PSL(3,\R)$ on the manifold of complete flags]\label{ex:flags}
Let $X=\{(p\in\ell)\}$ be the manifold of flags in $\R^3$, of dimension $3$, and take $a_m=\diag(e^{2t_m},1,e^{-2t_m})$ with $t_m\to+\infty$, so $\Theta=\Delta$. There are six Bruhat cells and the fibre dimensions are
\[
\begin{array}{c|cccccc}
\text{cell dimension}&0&1&1&2&2&3\\\hline
\text{fibre dimension}&3&2&2&1&1&0
\end{array}
\]
The point of the smallest cell (the flag $(e_3\in\langle e_2,e_3\rangle)$, once $k'=1$) has as destination the whole flag manifold. A generic flag is sent to the single flag $(e_1\in\langle e_1,e_2\rangle)$. The four intermediate cells interpolate, in the manner of Corollary~\ref{cor:conservation}.
\end{example}
\begin{example}[$\mathrm{SU}(2,1)$ on its null-line sphere]\label{ex:su21}
Realise $\mathrm{SU}(2,1)$ with the Hermitian form $J$ having ones on the anti-diagonal, and $a_m=\diag(e^{t_m},1,e^{-t_m})$ with $t_m\to+\infty$. The restricted roots are $\alpha,2\alpha$ with multiplicities $2,1$, and $N_\Theta$ is the Heisenberg group $\{n(z,u)\}$ with $a_mn(z,u)a_m^{-1}=n(e^{t_m}z,e^{2t_m}u)$. It fixes $[e_1]$ and is transitive on the rest of the null-line sphere $S^3$. So the destination of the repelling point $[e_3]$ is the whole sphere (dimension $2+1=3$), and every other source has destination $\{[e_1]\}$. This is the picture with a double root: the real dimension is the sum of the root multiplicities, not the number of roots. It is also Corollary~\ref{cor:rank-one}.
\end{example}
\begin{example}[$\mathrm{Sp}(4,\R)$]\label{ex:sp4}
The root system is $C_2$, with simple roots $\alpha_1=e_1-e_2$ (short) and $\alpha_2=2e_2$ (long), all multiplicities $1$. There are three nontrivial flag manifolds: $\R\PP^3$ ($\theta=\{\alpha_1\}$), the Lagrangian Grassmannian $\Lambda(2,4)$ of dimension $3$ ($\theta=\{\alpha_2\}$), and the manifold of complete isotropic flags, of dimension $4$. Table~\ref{tab:sp4} lists the distinct cell and fibre dimensions for each type $\Theta$ of sequence, ordered by increasing cell dimension; multiplicities are suppressed.
\begin{table}[ht]
\centering\small
\begin{tabular}{@{}lllll@{}}
\toprule
$X$ & $\Theta$ & cell dimensions & fibre dimensions\\\midrule
$\R\PP^3$ & $\{\alpha_1,\alpha_2\}$ & $0,1,2,3$ & $3,2,1,0$\\
 & $\{\alpha_1\}$ & $0,2,3$ & $3,1,0$\\
 & $\{\alpha_2\}$ & $1,3$ & $2,0$\\\midrule
$\Lambda(2,4)$ & $\{\alpha_1,\alpha_2\}$ & $0,1,2,3$ & $3,2,1,0$\\
 & $\{\alpha_1\}$ & $1,3$ & $2,0$\\
 & $\{\alpha_2\}$ & $0,2,3$ & $3,1,0$\\\midrule
complete flags & $\{\alpha_1,\alpha_2\}$ & $0,1,2,3,4$ & $4,3,2,1,0$\\
 & $\{\alpha_1\}$, $\{\alpha_2\}$ & $1,2,3,4$ & $3,2,1,0$\\\bottomrule
\end{tabular}
\caption{Distinct cell and fibre dimensions for $\mathrm{Sp}(4,\R)$, with multiplicities suppressed. In every row, cell dimension $+$ fibre dimension $=\dim X$.}\label{tab:sp4}
\end{table}
The rows for $\R\PP^3$ are exactly what Proposition~\ref{prop:proj} predicts: the type $\{\alpha_1\}$ corresponds to blocks $(1,2,1)$ (the first and last coordinates are separated from the middle plane), the type $\{\alpha_2\}$ to blocks $(2,2)$. The Lagrangian rows agree with those for $\R\PP^3$ after exchanging the two simple reflections. This is the symmetry of the Coxeter system $I_2(4)$ with equal restricted root multiplicities. The two fundamental representations have dimensions $4$ and $5$.
\end{example}
\begin{example}[$\PSL(4,\mathbb K)$ on the Grassmannian $\mathrm{Gr}(2,\mathbb K^4)$]\label{ex:gr24}
Let $X=\mathrm{Gr}(2,\mathbb K^4)$, of dimension $4$ over $\mathbb K$ (real dimension $4$ for $\R$ and $8$ for $\C$), and take $a_m=\mathrm{diag}(e^{t_m}b_1,e^{-t_m}b_2)$ with $t_m\to\infty$ and $b_i\in\mathrm{GL}(E_i)$ positive diagonal, $E_1=\langle e_1,e_2\rangle$, $E_2=\langle e_3,e_4\rangle$: type $\Theta=\{\alpha_2\}$, blocks $(2,2)$. For a plane $W$ let $s=\dim(W\cap E_2)$. A direct computation gives three kinds of fibres.
\begin{itemize}[nosep]
\item $s=0$: $W$ is the graph of a map $\varphi:E_1\to E_2$. Every nearby source is a graph $W_m=\operatorname{graph}(\varphi_m)$ with $\varphi_m\to\varphi$, and $a_mW_m$ is the graph of $e^{-2t_m}b_2\varphi_m b_1^{-1}\to0$. Thus the destination is the single point $E_1$. These planes form the open cell (dimension $4$ over $\mathbb K$).
\item $s=1$: $W=\langle u,w\rangle$ with $0\ne u\in E_2$ and $w=w_1+w_2$, $w_1\in E_1\setminus0$. Put $\ell=\mathbb Kb_1w_1$ and $T=E_1\oplus\mathbb Kb_2u$. For arbitrary $W_m\to W$, choose frames tending to $(w,u)$. The three components of their Pl\"ucker vectors in $\bigwedge^2E_1$, $E_1\wedge E_2$ and $\bigwedge^2E_2$ are multiplied by $e^{2t_m}$, $1$ and $e^{-2t_m}$, respectively, with the fixed maps induced by $b_1,b_2$. The middle component tends to $b_1w_1\wedge b_2u\ne0$, and the last tends to zero. Every projective limit therefore lies in the projective line spanned by $\bigwedge^2E_1$ and $b_1w_1\wedge b_2u$. Its planes are exactly
\[
\{P:\ \mathbb K\,b_1w_1\subset P\subset E_1\oplus\mathbb K\,b_2u\}.
\]
We can reach every one of them. Write $P=\langle b_1w_1,v+c\,b_2u\rangle$, where $v\in E_1$, $c\in\mathbb K$ and $v+c\,b_2u\notin\ell$; here $c$ is the coefficient of $b_2u$ in the chosen second generator. If $c\ne0$, take $\delta_m=e^{-2t_m}b_1^{-1}(v/c)$ and $W_m=\langle w,u+\delta_m\rangle$. Then
\[
e^{-t_m}a_mw\longrightarrow b_1w_1,\qquad
c\,e^{t_m}a_m(u+\delta_m)=v+c\,b_2u.
\]
If $c=0$, the target is $P=E_1$; take $v\notin\ell$ and $\delta_m=e^{-t_m}b_1^{-1}v$, so $a_m(u+\delta_m)=v+e^{-t_m}b_2u\to v$. Both constructions have $\delta_m\to0$ and realise $P$ along the whole sequence. The destination is a projective line, the planes in the $3$-space $T$ through $\ell$. These sources form a codimension-one cell, so the conservation law $3+1=4$ holds over $\mathbb K$ (and $6+2=8$ in real dimensions over $\C$). The maps $b_1,b_2$ \emph{turn the fibre}, exactly as in Proposition~\ref{prop:sl3}.
\item $s=2$: $W=E_2$ is a point. For any map $F:E_2\to E_1$, take $W_m$ to be the graph over $E_2$ of $e^{-2t_m}b_1^{-1}Fb_2$. Then $W_m\to E_2$ and $a_mW_m=\operatorname{graph}(F)$. These target graphs form the open dense orbit of $E_2$ under $N_\Theta=\mathrm{Hom}(E_2,E_1)$. For a target on its boundary choose graphs $\operatorname{graph}(F_j)$ tending to it and let $j=j(m)\to\infty$ slowly enough that $e^{-2t_m}\|b_1^{-1}F_{j(m)}b_2\|\to0$, as in Step~1 of Theorem~\ref{thm:dest}. Thus every plane in $X$ is reached along the whole sequence.
\end{itemize}
Table~\ref{tab:gr24} gives the distinct cell and fibre dimensions for all seven types of sequence, with multiplicities suppressed.
\begin{table}[ht]
\centering\small
\begin{tabular}{@{}lll@{}}
\toprule
$\Theta$ & cell dimensions & fibre dimensions\\\midrule
$\{\alpha_1\}$, $\{\alpha_3\}$ & $2,4$ & $2,0$\\
$\{\alpha_2\}$ & $0,3,4$ & $4,1,0$\\
$\{\alpha_1,\alpha_3\}$ & $1,2,4$ & $3,2,0$\\
$\{\alpha_1,\alpha_2\}$, $\{\alpha_2,\alpha_3\}$ & $0,2,3,4$ & $4,2,1,0$\\
$\{\alpha_1,\alpha_2,\alpha_3\}$ & $0,1,2,3,4$ & $4,3,2,1,0$\\\bottomrule
\end{tabular}
\caption{Distinct cell and fibre dimensions for $\PSL(4,\R)$ on $\mathrm{Gr}(2,\R^4)$, with multiplicities suppressed. The row $\{\alpha_2\}$ is the case described above.}\label{tab:gr24}
\end{table}
\end{example}

\begin{proposition}[Fibres are Schubert varieties when $\Theta=\Delta$]\label{prop:schubert}
Suppose $\Theta=\Delta$ (all simple roots diverge: ``fully regular'' sequences) and let $X=G/Q$ be any flag manifold. If $k'x\in O_\sigma$ then
\[
\D_{(g_m)}(x)=k\,\overline{B\sigma Q/Q}=:k\,X_\sigma,
\]
the translate by $k$ of the Schubert variety $X_\sigma$, which has dimension $\sum\dim\g_\lambda$ over the positive roots $\lambda$ with $\sigma^{-1}\lambda\in\Phi(N_\theta^-)$. In particular the fibres of the destination relation are exactly the $k$-translates of Schubert varieties, one for each Weyl double coset, and they do not depend on the Levi class $b$.
\end{proposition}
\begin{proof}
For $\Theta=\Delta$, $P_\Theta=B=MAN$, $N_\Theta=N$ and $L_\Theta=MA$. By Lemma~\ref{lem:cells}(iii), $\rho(k'x)=l\sigma Q$ with $l\in MA$. The element $b\in A$ and the $A$-part of $l$ fix $\sigma Q$, and $m\in M$ normalises $N$ with $m\sigma Q=\sigma(\sigma^{-1}m\sigma)Q=\sigma Q$ since $\sigma^{-1}M\sigma=M\subset Q$. Hence $N\,b\,\rho(k'x)=N\sigma Q=B\sigma Q/Q$, and Theorem~\ref{thm:dest}(b) gives the claim. The dimension is that of $M_+\sigma Q$, the number of roots in $\Sigma_+=\Sigma\cap\Phi^+$ counted with multiplicity, where $\Sigma=\sigma\Phi(N_\theta^-)$.
\end{proof}
Proposition~\ref{prop:schubert} illustrates the observation that the closure of an orbit ``need not be a manifold'' (Remark~\ref{rem:destA}(3)): Schubert varieties are singular in general, so destinations of single points can be singular varieties, and they are in rank two.

\begin{example}[The Maslov cycle]\label{ex:maslov}
Take $\mathrm{Sp}(4,\R)$ acting on the Lagrangian Grassmannian $\Lambda=\Lambda(2,4)$ (a $3$-manifold) and a regular sequence, $a_m=\diag(e^{3t_m},e^{t_m},e^{-t_m},e^{-3t_m})$ with $t_m\to+\infty$, $k=k'=1$ and symplectic form pairing $e_1\leftrightarrow e_4$, $e_2\leftrightarrow e_3$. The four cells have dimensions $3,2,1,0$ and the fibres are the Schubert varieties of $\Lambda$ relative to the flag $\langle e_1\rangle\subset E=\langle e_1,e_2\rangle\subset\langle e_1\rangle^\perp$:
\begin{center}\small
\begin{tabular}{@{}lll@{}}
\toprule
source $L$ & dimension of fibre & destination\\\midrule
generic ($L\cap\langle e_3,e_4\rangle=0$) & $0$ & the point $E$\\
$L=\langle e_3,e_1\rangle$ & $1$ & the circle $\{L':\ \langle e_1\rangle\subset L'\subset\langle e_1\rangle^\perp\}$\\
$L=\langle e_2,e_4\rangle$ & $2$ & the \emph{Maslov cycle} $\{L':\ L'\cap E\ne0\}$\\
$L=\langle e_3,e_4\rangle$ & $3$ & all of $\Lambda$\\\bottomrule
\end{tabular}
\end{center}
The Maslov cycle is singular exactly at $L'=E$: in the chart of Lagrangian planes $L'$ that are graphs over $E$ of symmetric maps $S=\begin{pmatrix}x&y\\y&z\end{pmatrix}:E\to E'$ ($E'=\langle e_3,e_4\rangle$), the condition $L'\cap E\ne0$ means $\det S=0$, so the cycle is the cone $\{xz=y^2\}\subset\R^3$, with vertex $S=0$ ($L'=E$). Every $L'$ lies in such a chart for a Lagrangian complement $E'$ transverse to both $E$ and $L'$. These charts exist because the transverse loci are open dense. The determinant has nonzero differential at every symmetric matrix of rank one. So already for $\mathrm{Sp}(4,\R)$ the destination of a single point can be a singular surface, and the Levi class plays no role. In $\PSL(4,\R)$ acting on $\mathrm{Gr}(2,\R^4)$ the same happens for the Schubert divisor $\{W:W\cap E\ne0\}$, a cone over the quadric $\{x_1x_4-x_2x_3=0\}\subset\PP^3$; and in full flags of $\PSL(4,\R)$ the fibres run through the $24$ Schubert varieties $X_\sigma$, $\sigma\in S_4$, of dimension $\ell(\sigma)$ (there are $1,3,5,6,5,3,1$ of dimensions $0,\dots,6$), of which exactly $X_{3412}$ and $X_{4231}$ (dimensions $4$ and $5$) are singular (Lakshmibai--Sandhya).
\end{example}

\begin{remark}
The stabiliser of a coordinate source $\sigma Q$ in $N_\Theta$ has
Lie algebra $\nk_\Theta\cap\Ad_\sigma\mathfrak q$. Hence its orbit has
dimension
\[
\dim\nk_\Theta-\dim(\nk_\Theta\cap\Ad_\sigma\mathfrak q).
\]
Together with Theorem~\ref{thm:dest}, this gives the dimensions in
Example~\ref{ex:flags} and Tables~\ref{tab:sp4} and~\ref{tab:gr24}.
\end{remark}

\begin{remark}[The Levi class need not lie in $A$]\label{rem:generalb}
The formula of Theorem~\ref{thm:dest} also applies to $g_m=k_mz_mb\,k'_m$ with an arbitrary fixed $b\in L_\Theta$. Since $b$ commutes with $z_m$ and normalises $N_\Theta$, we have $\D_{(z_mb)}(x)=b\D_{(z_m)}(x)=\overline{N_\Theta\,b\rho(x)}$. The converging compact factors give the formula for $g_m$.
\end{remark}

%--------------------------------------------------------------
\section{The exterior boundary and complete collineations}\label{sec:ext}
The moving sources have given us their destinations. We now turn the question
around: how much of a divergent sequence must we remember to recover them all?
For induced actions of discrete subgroups of $\mathrm{PU}(1,n)$ on
Grassmannians and their Pl\"ucker projective spaces, Z\'u\~niga
\cite{Zuniga,Zuniga2} studies limit sets and equicontinuity regions under
the complex-hyperbolic hypotheses. Here we keep all Cartan types before
imposing any such hypothesis on a group.
For $\PSL(d,\mathbb K)$ we shall recover the type, the two flags and the graded isomorphisms from a compact space: exterior powers restore the slower scales lost by the ordinary projective limit. The statements of Chapters~\ref{ch:limit} and~\ref{ch:balanced} can be read without the details of this construction.

Throughout, $\mathbb K\in\{\R,\C\}$, $V=\mathbb K^d$ with its standard Euclidean or Hermitian structure, $d\ge2$. We write $g_m=k_ma_mk'_m$ for Cartan decompositions, so the compact factor $k'_m$ acts directly on sources, in agreement with Theorem~\ref{thm:dest}.

%..............................................................
\subsection{The pseudo-projective compactification}\label{sec:SP}
We first allow a matrix to become singular while retaining its class up to scalar. This is the meaning of ``limits of matrices'' in Chapters~\ref{ch:limit}--\ref{ch:adjoint}. The pseudo-projective compactification used by Cano--Seade~\cite{CanoSeade2014} and developed in~\cite{CNS} describes boundary points as \emph{pseudo-projective transformations} $\mathrm{QP}(n+1,\C)$: the maps $[v]\mapsto[Mv]$ off $\PP(\ker M)$, with convergence uniform on compact sets there. We use the description $\PP(\mathrm{End})$ over both $\R$ and $\C$. To recover what happens at the kernel we pass to complete collineations (Vainsencher~\cite{Vainsencher1984}, Thaddeus~\cite{rwThaddeus1999}, Neretin~\cite{Neretin}), the wonderful compactification of De~Concini--Procesi~\cite{DCP1983}.

\begin{definition}[Pseudo-projective compactification]\label{def:SP}
Let $V$ be a vector space of dimension $d\ge2$ over $\mathbb K\in\{\R,\C\}$. Put
\[
\mathrm{SP}(V)=\PP(\mathrm{End}_{\mathbb K}V)=(\mathrm{End}\,V\setminus\{0\})/\mathbb K^*\ \cong\ \mathbb K\PP^{d^2-1},
\]
with the quotient topology; it is compact and metrisable. The group $\PGL(V)=\{[A]:\det A\ne0\}$ is an open dense subset (the complement is the hypersurface $\{\det=0\}$). For $\tau=[A]\in\mathrm{SP}(V)$ the \emph{rank} $\operatorname{rk}\tau=\operatorname{rk}A$, the \emph{kernel} $\operatorname{Ker}\tau=\PP(\ker A)\subset\PP(V)$ and the \emph{image} $\operatorname{Im}\tau=\PP(\operatorname{Im}A)$ are well defined, and $\tau$ acts as a projective map on $\PP(V)\setminus\operatorname{Ker}\tau$ by $[v]\mapsto[Av]$.
\end{definition}

\begin{proposition}[Structure of $\mathrm{SP}(V)$]\label{prop:SP}
\begin{enumerate}[label=\textup{(\roman*)},nosep]
\item $\PGL(V)\times\PGL(V)$ acts continuously on $\mathrm{SP}(V)$ by $(g,h)\cdot[A]=[gAh^{-1}]$, extending the left--right action on $\PGL(V)$.
\item The orbits are the rank strata $S_r=\{\operatorname{rk}\tau=r\}$, $1\le r\le d$, with $S_d=\PGL(V)$ and $\overline{S_r}=\bigcup_{s\le r}S_s$.
\item $S_r$ is a bundle over $\mathrm{Gr}_{d-r}(V)\times\mathrm{Gr}_r(V)$ (kernel, image) whose fibre is the space of isomorphisms $V/\ker\to\mathrm{Im}$ modulo scalars, a copy of $\PGL_r(\mathbb K)$; so $\dim_{\mathbb K}S_r=2dr-r^2-1$. For $d=2$, $S_1=\PP^1\times\PP^1$.
\item \textbf{Continuity off the kernel.} If $\tau_m\to\tau$ in $\mathrm{SP}(V)$ and $x_m\to x$ in $\PP(V)\setminus\operatorname{Ker}\tau$, then $\tau_mx_m\to\tau x$; the convergence is uniform on compact subsets of $\PP(V)\setminus\operatorname{Ker}\tau$.
\item \textbf{Cartan reading.} Let $g_m=k_ma_mk'_m$ with $k_m\to k$, $k'_m\to k'$, $a_m=\diag(e^{t_{m,1}},\dots,e^{t_{m,d}})$, $t_{m,1}\ge\dots\ge t_{m,d}$, and suppose the first block of $(g_m)$ is $E_1=\langle e_1,\dots,e_{b_1}\rangle$ (that is, $t_{m,1}-t_{m,j}$ stays bounded for $j\le b_1$ and tends to $\infty$ for $j>b_1$). After extraction, $g_m\to\tau^+=[kDk']$ where $D=\diag(c_1,\dots,c_{b_1},0,\dots,0)$ with all $c_j>0$. Hence
\[
\operatorname{rk}\tau^+=b_1,\qquad\operatorname{Im}\tau^+=\PP(kE_1),\qquad\operatorname{Ker}\tau^+=\PP\bigl(k'^{-1}(E_2\oplus\cdots\oplus E_k)\bigr).
\]
\end{enumerate}
\end{proposition}
\begin{proof}
(i) The action is the projectivisation of a linear automorphism of $\mathrm{End}\,V$. (ii) Row and column operations carry every rank-$r$ matrix to $\diag(I_r,0)$; rank is invariant. Letting some of its nonzero diagonal entries tend to zero gives the lower ranks, and the vanishing of the $(r+1)$-minors excludes higher ranks from the closure. (iii) The equivariant map to (kernel, image) is a homogeneous bundle with fibre $\PP(\mathrm{Iso}(V/K,I))\cong\PGL_r$. Its dimension is $2r(d-r)+r^2-1$. (iv) Choose representatives $A_m\to A\ne0$ and unit representatives $v_m\to v$ for $x_m\to x$. Since $Av\ne0$, $A_mv_m\to Av\ne0$ and hence $[A_mv_m]\to[Av]$. A failure of uniformity on a compact set would give, after extraction, a convergent sequence of sources contradicting this assertion. (v) Normalising $a_m$ by $e^{-t_{m,1}}$ gives entries tending to $c_j\in(0,1]$ for $j\le b_1$ and to $0$ for $j>b_1$. Thus $e^{-t_{m,1}}g_m\to kDk'$.
\end{proof}

\begin{remark}\label{rem:SPwarn}
\begin{enumerate}[label=\textup{(\arabic*)},nosep]
\item \emph{Inversion separates dimension two from higher dimensions.} For $d\ge3$, $[A]\mapsto[A^{-1}]$ on $\PGL(V)$ has no continuous extension to $\mathrm{SP}(V)$, so $\tau^-=\lim g_m^{-1}$ is additional information. For $d=3$ take $g_m=\diag(m,1,m^{-1})$ and $g'_m=\diag(m,m^{-1/2},m^{-1/2})$. Both converge to $\tau^+=[e_1e_1^*]$, but $g_m^{-1}\to[e_3e_3^*]$ has rank $1$ while $g_m'^{-1}\to[\mathrm{pr}_{\langle e_2,e_3\rangle}]$ has rank $2$; adjoining identity blocks gives the same obstruction for every $d>3$. For $d=2$, in contrast,
\[
 A=\begin{pmatrix}a&b\\c&d\end{pmatrix}\longmapsto
 \operatorname{adj}A=\begin{pmatrix}d&-b\\-c&a\end{pmatrix}
\]
is a linear involution of $\mathrm{End}(V)$, with $A\operatorname{adj}A=(\det A)I$. Its projectivisation is therefore a continuous involution of $\mathrm{SP}(V)$ agreeing with inversion on $\PGL(V)$.
\item \emph{A limit sees only the fastest scale} (Proposition~\ref{prop:SP}(v)): only the leading image $G_1$, the leading kernel $F_2$ and the projective map $[A_1:V/F_2\to G_1]$ survive. For $d\ge3$ this does not determine all destinations. The pair $(\tau^+,\tau^-)$ sees the two end blocks and suffices in dimensions $2$ and $3$, but can miss the middle for $d\ge4$ (Corollary~\ref{cor:noquot}). The exterior boundary records every scale; for $d=2$ it is already $\mathrm{SP}(V)$.
\end{enumerate}
\end{remark}

\begin{remark}[Beyond $\PSL(d)$]\label{rem:SPgeneral}\ A representation $\rho$ of a covering group, whose central covering kernel acts by scalars, gives a projective action $\bar\rho:G\to\PGL(V_\rho)$. Its closure in $\PP(\mathrm{End}\,V_\rho)$ compactifies the effective group $G_{\rm eff}=G/\ker\bar\rho$; it compactifies $G$ itself when $\bar\rho$ is faithful. For a complex adjoint group, a regular highest weight gives the wonderful compactification of De~Concini--Procesi. The comparison for $\PGL(V)$ in detail and for general $G$ in outline is made in Section~\ref{sec:SPvsWond}.\end{remark}

\paragraph{A single projective limit.}
The first collapse has a projective limit in each direction.

\begin{lemma}\label{lem:qp}
Every sequence $(g_m)$ in $\PSL(d,\mathbb K)$ has a subsequence with $g_m\to\tau^+$ and $g_m^{-1}\to\tau^-$ in $\mathrm{SP}(V)$, uniformly on compact sets away from $\operatorname{Ker}\tau^\pm$. If $(g_m)$ leaves every compact set, both limits are singular and $\operatorname{Im}\tau^+\subset\operatorname{Ker}\tau^-$, $\operatorname{Im}\tau^-\subset\operatorname{Ker}\tau^+$.
\end{lemma}
\begin{proof}
Compactness and Proposition~\ref{prop:SP}(iv) give both limits and uniform convergence off their kernels. For determinant-one lifts $A_m$, divergence implies $\|A_m\|\|A_m^{-1}\|\to\infty$: both norms are $\ge1$, and a bounded product would bound the lifts and their inverses. Extract so that the unit normalisations tend to $A\ne0$, $B\ne0$. Then $AB=\lim A_mA_m^{-1}/(\|A_m\|\|A_m^{-1}\|)=0$, and likewise $BA=0$, giving the inclusions and singularity.
\end{proof}
For $d=2$ the singular boundary is the quadric of pairs (image point, kernel point), so this is the convergence principle of Proposition~\ref{prop:classical}. In higher dimension the kernel can itself carry slower dynamics, invisible in the first limit; Example~\ref{ex:d3fibre} will display the missing destinations.

\begin{remark}[Twistor lines; known]
The twistor fibration $p:\PP^3\to S^4$ realises $\PP^3=\PP(\C^4)$ as the twistor space of the conformal $4$-sphere, with fibres $\PP^1$. An orientation-preserving conformal group of $S^4$ acts holomorphically on $\PP^3$; if it is discrete, the action is properly discontinuous on $p^{-1}(\Omega)$ (Seade--Verjovsky~\cite{SV2,SeadeVerjovsky2002}) and the complement $p^{-1}(\Lambda)$ is a union of complete projective lines. This is the geometric form of the obstruction above: even in the most classical case, point-like attracting data are replaced by positive-dimensional ones.
\end{remark}

%..............................................................
\paragraph{Successive kernels.}
Restricting to the kernel of the first limit and normalising again gives a flag with convergence on each stratum. At every step the kernel is fixed: the source has not yet been allowed to move.

\begin{lemma}\label{lem:Furst}
Let $(g_m)\subset\PSL(d,\mathbb K)$ leave every compact set. After extraction there are a flag $0=V_0\subsetneq V_1\subsetneq\cdots\subsetneq V_k=V$ with $k\ge2$ and linear maps $A_j:V_j\to V$ with $\ker A_j=V_{j-1}$ such that $g_m[v]\to[A_jv]$ for $[v]\in\PP(V_j)\setminus\PP(V_{j-1})$, uniformly on compact subsets of each stratum.
\end{lemma}
\begin{proof}
Choose determinant-one lifts $M_m$. Put $W_0=V$ and extract so that $M_m|_{W_0}/\|M_m|_{W_0}\|\to B_0\ne0$. The limit is singular: otherwise $M_m$ and $M_m^{-1}$ would be bounded. Put $W_1=\ker B_0\subsetneq W_0$. Inductively, if $W_i\ne0$ then $M_m|_{W_i}\ne0$ (as $M_m$ is invertible), and after extraction $M_m|_{W_i}/\|M_m|_{W_i}\|\to B_i\ne0$ and we put $W_{i+1}=\ker B_i$. The $W_i$ strictly decrease, so the procedure ends at some $W_k=0$ after finitely many extractions. Set $V_j=W_{k-j}$, $A_j=B_{k-j}$. On a compact set of $\PP(W_i)\setminus\PP(W_{i+1})$ unit representatives satisfy $\|B_iv\|\ge c>0$, so convergence of the normalised linear maps gives uniform convergence of their projectivisations.
\end{proof}
These maps record fixed-source limits only. The singular subspaces of a Cartan decomposition move with the sequence, and faster directions can magnify a small motion of the compact factors (Example~\ref{ex:drift}). To recover every destination we must keep those moving subspaces until the limit is taken.

%..............................................................
\paragraph{Cartan blocks and types.}
We keep every singular block, including those killed by the first projective limit: ratios stay bounded within a block and diverge between consecutive blocks.

Let $g_m=k_ma_mk'_m$ with $a_m=\diag(e^{\lambda_{1,m}},\dots,e^{\lambda_{d,m}})$, $\lambda_{1,m}\ge\cdots\ge\lambda_{d,m}$, $\sum_i\lambda_{i,m}=0$; extract so that $k_m\to k$, $k'_m\to k'$, and each gap $\lambda_{i,m}-\lambda_{i+1,m}$ converges in $[0,\infty]$. Let $r_1<\cdots<r_{k-1}$ be the indices of the divergent gaps, $r_0=0$, $r_k=d$, and put
\[
E_j=\mathrm{Span}(e_{r_{j-1}+1},\dots,e_{r_j}),\qquad b_j=\dim E_j,\qquad j=1,\dots,k,
\]
the blocks, ordered from fastest to slowest. In the language of Section~\ref{sec:setting-sum}, the type is $\Theta=\{\alpha_{r_1},\dots,\alpha_{r_{k-1}}\}$.

\begin{lemma}[Bounded multiplication does not change the type]\label{lem:bdd}
Fix $g\in\mathrm{GL}(d,\mathbb K)$. There is $C_g$ with $|\log\sigma_i(gh)-\log\sigma_i(h)|\le C_g$ and $|\log\sigma_i(hg)-\log\sigma_i(h)|\le C_g$ for all $h$ and all $i$. The same holds if the Hermitian norm is replaced by another one. Hence left or right multiplication by a sequence contained in a compact subset of $\mathrm{GL}(d,\mathbb K)$ preserves the divergent gaps and the block dimensions. This is the meaning of \emph{bounded} here: both the matrices and their inverses are uniformly bounded.
\end{lemma}
\begin{proof}
By the min--max characterisation of singular values, $\|g^{-1}\|^{-1}\sigma_i(h)\le\sigma_i(gh)\le\|g\|\sigma_i(h)$; the right-hand version follows by passing to adjoints. If $c\|v\|\le\|v\|'\le C\|v\|$ then $(c/C)\sigma_i\le\sigma_i'\le(C/c)\sigma_i$. Adjacent differences of logarithms change by at most a constant, so an infinite gap stays infinite and a bounded gap stays bounded.
\end{proof}

\begin{lemma}[Types]\label{lem:types}
The block types of divergent sequences in $\PSL(d,\mathbb K)$ are exactly the ordered compositions $d=b_1+\cdots+b_k$ with $k\ge2$; there are $2^{d-1}-1$ of them. If $F_j=k'^{-1}(E_j\oplus\cdots\oplus E_k)$ and $G_j=k(E_1\oplus\cdots\oplus E_j)$, then $\dim F_j=\sum_{i\ge j}b_i$ and $\dim G_j=\sum_{i\le j}b_i$. The first forward and inverse pseudo-projective limits have ranks $b_1$ and $b_k$, and, ordered from the generic stratum to the deepest, the projective dimensions of the forward destinations are
\[
0,\ b_1,\ b_1+b_2,\ \dots,\ b_1+\cdots+b_{k-1},
\]
and those of the inverse sequence are obtained by reversing the block order.
\end{lemma}
\begin{proof}
The divergent gaps are the proper partial sums of the $b_i$, giving the bijection and count. The flag dimensions follow from their definitions, the end ranks from Proposition~\ref{prop:SP}(v) applied in both directions, and the destination dimensions from Proposition~\ref{prop:proj}. Conversely, choose $\eta_1>\cdots>\eta_k$, put $c_i=\eta_i-\frac1d\sum_ab_a\eta_a$, and take $a_m=\diag(e^{mc_1}I_{b_1},\dots,e^{mc_k}I_{b_k})$: its divergent gaps are exactly the prescribed ones.
\end{proof}
\begin{table}[ht]
\centering\small
\begin{tabular}{@{}ccccc@{}}
\toprule
type & $\operatorname{Im}\tau^+$ & $\operatorname{Im}\tau^-$ & forward dimensions & inverse dimensions\\\midrule
$(1,3)$ & point & plane & $0,1$ & $0,3$\\
$(3,1)$ & plane & point & $0,3$ & $0,1$\\
$(2,2)$ & line & line & $0,2$ & $0,2$\\
$(1,1,2)$ & point & line & $0,1,2$ & $0,2,3$\\
$(1,2,1)$ & point & point & $0,1,3$ & $0,1,3$\\
$(2,1,1)$ & line & point & $0,2,3$ & $0,1,2$\\
$(1,1,1,1)$ & point & point & $0,1,2,3$ & $0,1,2,3$\\\bottomrule
\end{tabular}
\caption{The seven Cartan types in $\PSL(4,\mathbb K)$ acting on $\mathbb K\PP^3$. Dimensions are projective dimensions, from the generic stratum to the deepest.}\label{tab:cp3-types}
\end{table}
The two types with first images (point, point) have different destination dimensions: the end images do not determine the type.

\begin{example}\label{ex:two-four}
For $a_m=\diag(e^m,e^m,e^{-m},e^{-m})$ there are two blocks $E_1=\langle e_1,e_2\rangle$, $E_2=\langle e_3,e_4\rangle$. Off $\PP(E_2)$ the destination is the point given by the $E_1$-component; for $x=[\xi]\in\PP(E_2)$ it is $\PP(E_1\oplus\mathbb K\xi)$, a projective plane. For $a_m=\diag(e^{3m},e^m,e^{-m},e^{-3m})$ every gap diverges and the destination dimensions run through $0,1,2,3$; in particular $\D_{(a_m)}([e_4])=\PP^3$.
\end{example}

%..............................................................
\paragraph{Flags and graded maps.}
To locate the destinations we add the source and target flags and the finite maps between corresponding blocks: the Levi datum that turned the line in the second example.

Put $s_{j,m}=e^{\lambda_{r_{j-1}+1,m}}$, the largest scale in $E_j$, and $D_{j,m}=s_{j,m}^{-1}a_m|_{E_j}$. All gaps inside $E_j$ stay bounded, so after one more extraction $D_{j,m}\to D_j\in\mathrm{GL}(E_j)$ (a positive diagonal matrix). Let $\pi_j:G_j\to G_j/G_{j-1}$ and $q_j:F_j\to F_j/F_{j+1}$ be the quotient maps and set
\[
A_j:F_j/F_{j+1}\longrightarrow G_j/G_{j-1},\qquad A_j(v+F_{j+1})=k\,D_j\,\mathrm{pr}_{E_j}(k'v)+G_{j-1}.
\]
The maps induced by $k'$, $D_j$ and $k$ are isomorphisms, so $A_j$ is well defined and invertible. The \emph{full graded asymptotic datum} is $\mathcal L=(F_\bullet,G_\bullet,([A_j])_{j=1}^k)$, with each $A_j$ retained up to scalar. This is $(\Theta,k,k',b)$ in the language of flags: the flags record the type and the compact factors modulo their action inside blocks, while $[A_j]$ records $b$ together with those inner factors.

\begin{proposition}[Destination fibres in terms of the graded datum]\label{prop:graded}
For $x=[\xi]$ with $\xi\in F_j\setminus F_{j+1}$,
\[
\D_{(g_m)}(x)=\PP\bigl(\pi_j^{-1}(A_j(\mathbb K\,q_j\xi))\bigr)=\PP\bigl(G_{j-1}\oplus\mathbb K\,k\,D_j\,\mathrm{pr}_{E_j}(k'\xi)\bigr),
\]
of projective dimension $\dim G_{j-1}=b_1+\cdots+b_{j-1}$, and for $w\in G_j\setminus G_{j-1}$ the inverse sequence has $\D_{(g_m^{-1})}([w])=\PP\bigl(q_j^{-1}(A_j^{-1}(\mathbb K\,\pi_jw))\bigr)$.
\end{proposition}
\begin{proof}
For $\xi\in F_j\setminus F_{j+1}$ the first nonzero block of $k'\xi$ is $j$. Proposition~\ref{prop:proj}, transported by the compact factors as in Theorem~\ref{thm:dest}(b), gives the forward formula. Applying the same argument to $g_m^{-1}=k_m'^{-1}a_m^{-1}k_m^{-1}$ reverses the blocks and inverts the graded maps, giving the inverse formula.
\end{proof}

%..............................................................
The projective block argument already appears in Toledo-Acosta
\cite[\S\S3.1, 3.3, Proposition~3.3.3]{ToledoAcosta2019}.
In adapted source and target coordinates, his block limits $D_i$
are the maps $D_i$ above, and his $V_i$ corresponds to $\PP(G_i)$.
The maps $[A_i]$ retain these limits on the quotient flags.
Theorem~\ref{thm:ext-boundary} records the datum in all exterior
degrees; Theorem~\ref{thm:dest} treats arbitrary real flag manifolds
and all restricted-root multiplicities.

\paragraph{The exterior boundary.}
The exterior powers recover the blocks hidden in the kernel of the first limit. Define
\[
\jmath:\PGL(d,\mathbb K)\longrightarrow\prod_{q=1}^{d-1}\PP\bigl(\mathrm{End}\,\wedge^qV\bigr),\qquad[g]\longmapsto\bigl([\wedge^qg]\bigr)_{q=1}^{d-1},
\]
and let $X_{\rm ext}=\overline{\jmath(\PGL(d,\mathbb K))}$, the \emph{exterior compactification}. Here we compactify the full projective group; the closure of its subgroup $\PSL(d,\R)$ is a separate real subset when $d$ is even, as discussed below. Over $\C$ this is the classical space of complete collineations in its simultaneous exterior-power realisation (Vainsencher~\cite{Vainsencher1984}; Thaddeus~\cite{rwThaddeus1999}), and for $d=2$ it is $\mathrm{SP}(V)$. The projectivisations must be taken \emph{independently} for each $q$: for $g_m=\diag(e^m,1,e^{-m})$ the normalised first power has a rank-one limit and zero second exterior power, while the normalised second powers have a nonzero rank-one limit, which is the datum recorded by the second factor.

\begin{proposition}[Boundary points and graded data]\label{prop:ext-data}
Let $g_m=k_ma_mk'_m$ be as above. For $1\le q\le d-1$ write $q=r_{j-1}+t$, $1\le t\le b_j$.
\begin{enumerate}[label=(\alph*),nosep]
\item Put $c_{q,m}=\exp(\lambda_{1,m}+\cdots+\lambda_{q,m})$. Then $c_{q,m}^{-1}\wedge^qg_m\to T_q\ne0$ and $\operatorname{rank}T_q=\binom{b_j}{t}$. In particular the set of indices $q$ with $\operatorname{rank}T_q=1$ is $\{r_1,\dots,r_{k-1}\}$.
\item The point $([T_q])_q$ of $X_{\rm ext}$ determines, and is determined by, the graded datum $\mathcal L$.
\end{enumerate}
\end{proposition}
\begin{proof}
(a) In the wedge basis $e_J$ ($J\subset\{1,\dots,d\}$, $|J|=q$) the diagonal coefficient of $c_{q,m}^{-1}\wedge^qa_m$ is $\exp\bigl(\sum_{a\in J}\lambda_{a,m}-\sum_{a\le q}\lambda_{a,m}\bigr)$. By Abel summation the negative of its exponent is
\[
\sum_{i=1}^{d-1}\bigl(\min(i,q)-n_i(J)\bigr)(\lambda_{i,m}-\lambda_{i+1,m}),\qquad n_i(J)=|J\cap\{1,\dots,i\}|,
\]
and every term is $\ge0$. The coefficient has a nonzero limit exactly when every divergent gap $i=r_s$ has $n_{r_s}(J)=\min(r_s,q)$. For $q=r_{j-1}+t$ this means that $J$ contains $\{1,\dots,r_{j-1}\}$, contains exactly $t$ indices of $E_j$, and none of the later ones; there are $\binom{b_j}{t}$ such $J$, each with a nonzero (finite) limit because the gaps inside $E_j$ are bounded. Multiplying by $\wedge^qk_m\to\wedge^qk$ and $\wedge^qk'_m\to\wedge^qk'$ gives $T_q$ of that rank. The rank equals $1$ iff $t=b_j$, i.e. $q=r_j$.

(b) \emph{From the point to the datum.} At $q=r_j$ the limit is rank one with image the Pl\"ucker line of $G_j$ (the line spanned by $\wedge^{r_j}k(e_1\wedge\cdots\wedge e_{r_j})$), and its adjoint has image the Pl\"ucker line of $\wedge^{r_j}(F_{j+1}^\perp)$ since $F_{j+1}^\perp=k'^{-1}(E_{\le j})$; the adjoint is determined by $[T_q]$ and the Hermitian structure. A decomposable vector determines its subspace, so all $G_j$ and $F_{j+1}$ are recovered. For a one-dimensional block the graded map is unique up to scalar. For $b_j>1$, put $U_j=F_j\cap F_{j+1}^\perp\cong F_j/F_{j+1}$, $Z_j=G_j\cap G_{j-1}^\perp\cong G_j/G_{j-1}$ and $r=r_{j-1}$. Choose nonzero $\omega^-\in\wedge^rF_j^\perp$ and $\omega^+\in\wedge^rG_{j-1}$ (decomposable lines, as above, with $\wedge^0=\mathbb K$) and define injections $\iota^-:U_j\to\wedge^{r+1}V$, $u\mapsto\omega^-\wedge u$ and $\iota^+:Z_j\to\wedge^{r+1}V$, $z\mapsto\omega^+\wedge z$. In the diagonal coordinates of (a) with $q=r+1$, i.e.\ $t=1$, the only nonzero coefficients are $e_{\{1..r\}}\wedge e_a$ with $a\in E_j$, equal to $e^{\lambda_a-\lambda_{r+1}}$, so $T_{r+1}\circ\iota^-$ takes values in $\operatorname{Im}\iota^+$ and the induced map $B_j=(\iota^+)^{-1}T_{r+1}\iota^-:U_j\to Z_j$ equals $kD_j(k')|_{U_j}$ up to a scalar, that is, $[B_j]=[A_j]$.

\emph{From the datum to the point.} By the computation in (a), the nonzero part of $T_q$ is, in the coordinates $E_i$ transported by $k,k'$, the map $\mathrm{Id}_{\wedge^{r_{j-1}}(E_{<j})}\otimes\wedge^tD_j$ on $\wedge^{r_{j-1}}(E_{<j})\otimes\wedge^tE_j$, and zero on the other weight summands, up to a scalar. The domain and target summands are determined by the flags and the map by $[A_j]$. Hence $\mathcal L$ determines each $[T_q]$.
\end{proof}

\begin{theorem}[Stratified convergence principle for $\PSL(d,\mathbb K)$]\label{thm:ext-boundary}
Let $(g_m)$ be a sequence in $\PSL(d,\mathbb K)$ leaving every compact set.
\begin{enumerate}[label=(\roman*),nosep]
\item After extraction $\jmath(g_m)$ converges in $X_{\rm ext}$ to a boundary point $\xi$, which corresponds to a unique graded datum $\mathcal L(\xi)=(F_\bullet,G_\bullet,([A_j]))$ with $d=b_1+\cdots+b_k$, $k\ge2$.
\item For this subsequence, and for every other sequence $\jmath(h_m)\to\xi$, the whole-sequence destination relations are
\begin{align*}
&\D_{(h_m)}([v])=\PP\bigl(\pi_j^{-1}(A_j(\mathbb K q_jv))\bigr)\quad(v\in F_j\setminus F_{j+1}),\\
&\D_{(h_m^{-1})}([w])=\PP\bigl(q_j^{-1}(A_j^{-1}(\mathbb K\pi_jw))\bigr)\quad(w\in G_j\setminus G_{j-1}).
\end{align*}
\item \textup{(Strata and bundles.)} $\PP(V)$ is the disjoint union of the source strata $\PP(F_j\setminus F_{j+1})$ and of the target strata $\PP(G_j\setminus G_{j-1})$. The destination of $[v]\in\PP(F_j\setminus F_{j+1})$ depends only on $[q_jv]\in\PP(F_j/F_{j+1})$, and is a projective subspace of dimension $\dim G_{j-1}$: over $\PP(F_j/F_{j+1})$ these subspaces form the projective bundle $\PP(G_{j-1}\oplus\mathcal O(-1))$ pulled back by $[A_j]$ from $\PP(G_j/G_{j-1})$ \textup{(}holomorphic if $\mathbb K=\C$\textup{)}.
\item Every composition of $d$ with $k\ge2$ occurs as a type, and every datum $\mathcal L$ is realised by a divergent sequence in $\PGL(d,\mathbb K)$.
\end{enumerate}
For an arbitrary group this is Theorem~\ref{thm:S}.
\end{theorem}
\begin{proof}
(i) $X_{\rm ext}$ is compact, and Proposition~\ref{prop:ext-data} identifies the limit.

(ii) Every subsequence of $\jmath(h_m)\to\xi$ has a further Cartan extraction. Proposition~\ref{prop:ext-data}(b) gives the same datum $\mathcal L(\xi)$ on each extraction, and Proposition~\ref{prop:graded} and Theorem~\ref{thm:dest}(d) give Hausdorff convergence of its graphs to $R_{\mathcal L}$. Thus the graphs of the whole sequence converge to $R_{\mathcal L}$. Each $(x,y)\in R_{\mathcal L}$ can therefore be approached by points $(x_m,h_mx_m)$ of these graphs, with $x_m\to x$ and $h_mx_m\to y$; conversely every such limit lies in $R_{\mathcal L}$. This proves the whole-sequence assertion. The inverse relation is the transpose.

(iii) By (ii) the destination depends on $v$ only through $q_jv$, and $\PP(\pi_j^{-1}(\ell))$ for a line $\ell\subset G_j/G_{j-1}$ is the projective span of $G_{j-1}$ and any $w$ with $\pi_jw\in\ell$, which is the fibre of $\PP(G_{j-1}\oplus\mathcal O(-1))$ over $\ell$. (iv) Types: Lemma~\ref{lem:types}. Realisation: given $\mathcal L$, write each $A_j$, after identifying $F_j/F_{j+1}$ and $G_j/G_{j-1}$ with $U_j$ and $Z_j$, in singular-value form $A_j=V_j\Sigma_jW_j$ with $W_j:U_j\to E_j$ and $V_j:E_j\to Z_j$ isometries and $\Sigma_j$ positive diagonal. Choose unitary (orthogonal) $k'$ with $k'|_{U_j}=W_j$ and $k$ with $k|_{E_j}=V_j$, and put $a_m=\diag(e^{mc_1}\Sigma_1,\dots,e^{mc_k}\Sigma_k)$ with $c_1>\cdots>c_k$ and $\sum b_jc_j=0$ (after normalising the determinant). The sequence $k a_m k'$ has the datum $\mathcal L$ by construction.
\end{proof}
\begin{remark}
Over $\R$ and $d$ even, $\PSL(d,\R)$ is the positive-determinant component of $\PGL(d,\R)$. The sign condition for realising an exterior datum in this component is given after Table~\ref{tab:SPvsX}.
\end{remark}

Moving compact factors can change the limit of a fixed source even when the Cartan and fixed-kernel flags agree; they can also erase an intermediate step of the latter flag. We begin with the first possibility.

\begin{example}[Fixed kernels see only one point of a fibre]\label{ex:drift}
Let $g_m=\diag(e^{2m},1,e^{-2m})\,R_{\theta_m}$, where $R_\theta$ is the rotation by $\theta$ in the plane $\langle e_1,e_2\rangle$ and $\theta_m=e^{-3m/2}$. Then $k_m=1$, $k'_m=R_{\theta_m}\to1$, and all gaps diverge, so the type is $(1,1,1)$ with $D_j=1$. For the fixed source $x=[e_2]$,
\[
g_me_2=\bigl(-e^{2m}\sin\theta_m,\ \cos\theta_m,\ 0\bigr)\sim\bigl(-e^{m/2},1,0\bigr),
\]
so $g_mx\to[e_1]$. The fixed-kernel construction of Lemma~\ref{lem:Furst} has the same flags $\langle e_3\rangle\subset\langle e_2,e_3\rangle\subset\C^3$ and records exactly this: its map on $\langle e_2,e_3\rangle$ sends $e_2$ to a nonzero multiple of $e_1$. The Cartan graded map is obtained using the moving compact frames. It sends $e_2+F_3$ to $e_2+G_1$ with $G_1=\langle e_1\rangle$, and Proposition~\ref{prop:graded} gives $\D_{(g_m)}([e_2])=\PP\langle e_1,e_2\rangle$. The limit $[e_1]$ is one point of this line, the one reached by the fixed source because the rotation leaks the fast direction $e_1$ into the image faster than the slow scale of $e_2$ can compete. Moving sources reach the whole line.
\end{example}

There is a simple sufficient condition for the two flags to agree.

\begin{proposition}[Fixed kernels and the Cartan source flag]\label{r3fix:comparison}
For the Cartan extraction above, put $F_{k+1}=0$, and write
$V=W_0\supsetneq W_1\supsetneq\cdots\supsetneq W_\ell=0$
for the fixed-kernel tower of Lemma~\ref{lem:Furst}. Then $W_1=F_2$.
If both compact factors $k_m,k'_m$ are constant, then $\ell=k$ and
$W_i=F_{i+1}$ for $0\le i\le k$. Thus the whole tower agrees with
the source flag recovered from the independently normalised exterior
operators: $\operatorname{Im}T_{r_j}^*=\bigwedge^{r_j}F_{j+1}^{\perp}$ for $1\le j<k$
(Proposition~\ref{prop:ext-data}). In general its
intermediate kernels need not agree with that flag.
\end{proposition}
\begin{proof}
The leading normalisation is invertible on $E_1$ and vanishes on the
later blocks, so its kernel is $k'^{-1}(E_2\oplus\cdots\oplus E_k)=F_2$.
For constant compact factors, restricting to this fixed kernel leaves
exactly the later singular blocks; successive normalisations remove
them one at a time. This gives $W_i=F_{i+1}$ by induction.
\end{proof}
With moving compact factors a later restriction can have an injective limit and stop the construction too soon.

\begin{example}[A strictly coarser fixed-kernel flag]\label{r3fix:coarser}
In $\mathrm{SL}(3,\C)$ let $g_m=\diag(m,1,m^{-1})R_m$, where $R_m$
fixes $e_2$ and is the rotation by $1/m$ in $\langle e_1,e_3\rangle$,
with $R_me_3=-\sin(1/m)e_1+\cos(1/m)e_3$.
The singular values are exactly $m,1,m^{-1}$ and $R_m\to I$, so the
Cartan source flag is $0\subset\C e_3\subset\langle e_2,e_3\rangle
\subset\C^3$. But $m^{-1}g_m\to E_{11}$ gives
$W_1=\langle e_2,e_3\rangle$, and on this fixed plane
\[
 g_me_2=e_2,\qquad
 g_me_3=-m\sin(1/m)e_1+m^{-1}\cos(1/m)e_3\longrightarrow-e_1.
\]
The restriction has an injective limit, so the fixed-kernel flag is
only $0\subset W_1\subset\C^3$. Unlike Example~\ref{ex:drift}, here a
Cartan flag step itself is lost. This is a sequence-level example;
no common discrete subgroup is asserted.
\end{example}

%..............................................................
\paragraph{Fibres as projective bundles.}
Within a source stratum the destinations fit together holomorphically; their dimensions jump between strata.

Work over $\C$ with a realised datum $\mathcal L$. Put $S_j=\PP(F_j)\setminus\PP(F_{j+1})$, $r_{j-1}=\dim G_{j-1}$ and
\[
m_j:S_j\to\PP(G_j/G_{j-1}),\ m_j([v])=[A_jq_jv],\qquad\mathcal E_j=\{(x,w)\in S_j\times G_j:\pi_j(w)\in m_j(x)\}.
\]
Since $q_jv\ne0$ on $S_j$, $m_j$ is holomorphic and $\mathcal E_j$ is a holomorphic bundle with exact sequence $0\to S_j\times G_{j-1}\to\mathcal E_j\to m_j^*\mathcal O(-1)\to0$. Proposition~\ref{prop:graded} identifies the destination relation over $S_j$ with the locally trivial holomorphic $\PP^{r_{j-1}}$-bundle $\PP(\mathcal E_j)$. Across strata the projective dimensions, hence the top integral homology of the fibres, change; the global relation is not a bundle.

\begin{example}[Exponents are not enough]\label{ex:exponents}
Let $t_m=m$ and
\[
a_m=\diag(e^{\sqrt m},e^{\sqrt m},e^{-2\sqrt m}),\quad
b_m=\diag(e^{2\sqrt m},1,e^{-2\sqrt m}),\quad
c_m=\diag(2e^{\sqrt m},e^{\sqrt m},\tfrac12e^{-2\sqrt m}).
\]
All three have normalised exponents $\lambda_{i,m}/t_m\to0$, yet their block types are $(2,1)$, $(1,1,1)$, $(2,1)$. Moreover $\D_{(a_m)}([e_2])=\{[e_2]\}$ whereas $\D_{(b_m)}([e_2])=\PP\langle e_1,e_2\rangle$, because near $[e_2]$ the sequence $b_m$ has the faster line $\C e_1$: $[ue^{-2\sqrt m}e_1+e_2]\mapsto[ue_1+e_2]$. The sequences $a_m$ and $c_m$ have the same normalised limiting exponents, flags and cuts, but the active map on $\langle e_1,e_2\rangle$ is $\mathrm{Id}$ for $a_m$ and $\diag(2,1)$ for $c_m$, so the destinations over $[e_1+e_2]$ are $[e_1+e_2]$ and $[2e_1+e_2]$. Thus Lyapunov-type exponents do not determine moving-source destinations: the flags, the graded maps and the sub-exponential cuts are additional data. Also, $a_m$ illustrates the failure of local triviality: its fibres are points off $\PP\langle e_3\rangle$ and all of $\PP^2$ over $[e_3]$.
\end{example}

\paragraph{Relation with the wonderful compactification.}\ $X_{\rm ext}$ is the wonderful compactification of $\PGL(V)$; the precise comparison with $\mathrm{SP}(V)$, and the fact that $X_{\rm ext}$ is the closure of $\PGL(V)$ in the space of correspondences of $\PP(V)$, are the subject of Section~\ref{sec:SPvsWond}.

%--------------------------------------------------------------
\section{The pseudo-projective and the wonderful compactifications}\label{sec:SPvsWond}
Forgetting the slower scales gives a map from the exterior compactification $X_{\rm ext}$ to $\mathrm{SP}(V)$. Its fibres describe what the first collapse leaves undetermined. We shall also recover the boundary datum from the destinations, identifying $X_{\rm ext}$ with a space of correspondences of $\PP(V)$.

Throughout, $\mathbb K\in\{\R,\C\}$ and $\dim V=d\ge2$; we write $X=X_{\rm ext}(V)$ and $\pi:X\to\mathrm{SP}(V)$ for the projection to the first factor $q=1$ of $\prod_q\PP(\mathrm{End}\wedge^qV)$.

\paragraph{Orbits and the forgetful map.}
The boundary orbits remember every cut between blocks; the projection remembers only the first block.

\begin{proposition}[Orbits of $X_{\rm ext}$]\label{prop:orbitsX}
$\PGL(V)\times\PGL(V)$ acts on $X$ with exactly $2^{d-1}$ orbits $O_{\mathbf b}$, one for each ordered composition $\mathbf b=(b_1,\dots,b_k)$ of $d$ ($k=1$ is the open orbit $\PGL(V)$). The orbit $O_{\mathbf b}$ is the set of points whose graded datum has block sizes $\mathbf b$; it fibres over the product of two partial flag varieties of type $\mathbf b$ (the flags $F_\bullet,G_\bullet$) with fibre $\prod_j\PGL_{b_j}$ (the classes $[A_j]$), so
\[
\dim_{\mathbb K}O_{\mathbf b}=d^2-k,
\]
the closed orbit being $k=d$, i.e.\ $\mathrm{Fl}(V)\times\mathrm{Fl}(V)$. Moreover $O_{\mathbf b}\subset\overline{O_{\mathbf c}}$ iff $\mathbf c$ is a coarsening of $\mathbf b$.
\end{proposition}
\begin{proof}
By Theorem~\ref{thm:ext-boundary}(i),(iv), points are graded data. The group is transitive on pairs of flags of a fixed type, and their stabiliser induces arbitrary automorphisms of each graded piece, hence acts transitively on the $[A_j]$. The dimension is $2\dim\mathrm{Fl}_{\mathbf b}+\sum(b_j^2-1)=d^2-k$. To split a block of $\mathbf c$ of size $p+q$, take its graded map to be $B_t=\diag(I_p,tI_q)$, $t>0$, keeping the other maps equal to the identity in standard flags. In an exterior factor taking $u$ vectors from this block, divide independently by $t^{\max(0,u-p)}$. The surviving wedges take $\min(u,p)$ vectors from the first subblock and the rest from the second, giving exactly the split datum. Repeating splits and using equivariance gives the asserted inclusion when $\mathbf b$ refines $\mathbf c$.

Conversely, at every cut $r$ of $\mathbf c$ the exterior operator of degree $r$ has rank one. Under a projective limit choose nonzero convergent representatives; their $2\times2$ minors remain zero, so the limit has rank at most one and, being nonzero, rank exactly one. Proposition~\ref{prop:ext-data}(a) says that $r$ is still a cut of the limiting type $\mathbf b$. Thus every cut of $\mathbf c$ is a cut of $\mathbf b$, which is precisely the asserted coarsening condition.
\end{proof}
The orbit calculation identifies the boundary components: $\partial X=D_1\cup\dots\cup D_{d-1}$, where $D_r$ is the closure of the orbit of type $(r,d-r)$ and a point lies in $D_r$ iff $r$ is one of the partial sums $r_j=b_1+\dots+b_j$; the orbit closures are the intersections of the $D_r$'s. Smoothness of $X$ and normal crossings of its boundary follow from the complete-collineation model~\cite{Vainsencher1984,rwThaddeus1999}, recalled over $\C$ in Theorem~\ref{thm:blowup}. The determinantal blow-up construction is defined over $\R$ as well; its real graded data are all realised by Theorem~\ref{thm:ext-boundary}(iv), so its real points give $X_{\rm ext}(\R^d)$ with the same smoothness and crossing properties. These are the boundary and orbit properties of the \emph{wonderful} compactification. The boundary of $\mathrm{SP}(V)$, by contrast, is the single irreducible hypersurface $\{\det=0\}$, whose singular locus is the rank-$\le d-2$ locus. Among its $d$ orbit closures, $\overline{S_r}$ is singular exactly for $2\le r\le d-1$; $\overline{S_1}$ is the smooth Segre variety and $\overline{S_d}$ is projective space. For $d=2$ both orbit closures are smooth and the singular locus of the boundary is empty.

\begin{proposition}[The forgetful map]\label{prop:forget}
The map $\pi:X_{\rm ext}(V)\to\mathrm{SP}(V)$ is a $\PGL\times\PGL$-equivariant continuous surjection with the following properties.
\begin{enumerate}[label=\textup{(\alph*)},nosep]
\item $\pi$ sends $\xi$ with graded datum $(F_\bullet,G_\bullet,[A_j])$ to $\tau^+=[A_1\oplus0]$: the class of the endomorphism $V\to V/F_2\xrightarrow{A_1}G_1\hookrightarrow V$. Thus $\pi(O_{\mathbf b})=S_{b_1}$.
\item $\pi$ is a homeomorphism over $\PGL(V)\cup S_{d-1}$, that is, over the matrices of rank $\ge d-1$.
\item The fibre over $\tau\in S_r$, $r<d$, is in canonical bijection with the exterior compactification $X_{\rm ext}(\ker\tau,\ \mathrm{coker}\,\tau)$ of the projective space of isomorphisms from $\ker\tau$ to $\mathrm{coker}\,\tau=V/\mathrm{Im}\,\tau$. It has dimension $(d-r)^2-1$, with open part $\PGL(\ker\tau\to\mathrm{coker}\,\tau)$ (the datum with only two blocks) and boundary strata the finer types.
\end{enumerate}
\end{proposition}
\begin{proof}
(a) is Proposition~\ref{prop:SP}(v) and Proposition~\ref{prop:ext-data}. (c) Once $(G_1,F_2,[A_1])$ is fixed, the remaining flags in $F_2$ and $V/G_1$, with the maps $[A_j]$, $j\ge2$, are exactly the graded data of a complete collineation $F_2\to V/G_1$. All are realised by Theorem~\ref{thm:ext-boundary}(iv); the single-block data form the open part, of dimension $(d-r)^2-1$. This also proves surjectivity. (b) The fibre is a point for rank $d-1$, as it plainly is for rank $d$. The proper map $\pi$ is therefore a homeomorphism over the open rank-$\ge d-1$ locus.
\end{proof}
\begin{theorem}[Vainsencher; Thaddeus; Laksov]\label{thm:blowup}\cite{Vainsencher1984,rwThaddeus1999,Laksov1988}
Over $\C$, $\pi$ is the composition of $d-2$ blow-ups $X_{\rm ext}=Y_{d-2}\to\cdots\to Y_1\to Y_0=\PP(\mathrm{End}\,V)$, where $Y_{i}\to Y_{i-1}$ is the blow-up along the proper transform of the locus $\overline{S_i}$ of matrices of rank $\le i$, taken in increasing order $i=1,2,\dots,d-2$ (so first the Segre variety of rank-one matrices, then the proper transform of rank $\le2$, and so on). Each centre is smooth, the exceptional divisor of the $i$-th blow-up is $D_i$, and the proper transform of $\{\det=0\}$ is $D_{d-1}$.
\end{theorem}

For $d=2$ there are no blow-ups, so $X_{\rm ext}=\mathrm{SP}(V)=\PP^3$; for $d=3$ there is one, along the Segre variety $\PP^2\times\PP^{2*}\subset\PP^8$. In every dimension exterior duality identifies $[\wedge^qg^{-1}]$ with the transpose of $[\wedge^{d-q}g]$, extending inversion to $X_{\rm ext}$, exchanging the flags and sending $D_r$ to $D_{d-r}$. Its different behaviour on $\mathrm{SP}(V)$ is described in Remark~\ref{rem:SPwarn}(1).

\paragraph{The case $d=3$ in detail.}
Table~\ref{tab:d3} lists the four orbits. The preimage of the rank-one locus $S_1$ (dimension $4$, codimension $4$) is $D_1$, a $\PP^3$-bundle over $S_1$, and it is the union of the two orbits of types $(1,2)$ and $(1,1,1)$.
\begin{table}[ht]
\centering\small
\begin{tabular}{@{}lcccl@{}}
\toprule
Type $\mathbf b$ & $\dim X_{\rm ext}$-orbit & Image in $\mathrm{SP}$ & Fibre of $\pi$ & Description of $\xi$\\
\midrule
$(3)$ & $8$ & $\PGL_3$ & point & a collineation\\
$(2,1)$ & $7$ & $S_2$ ($\dim 7$) & point & rank $2$; kernel line, image plane\\
$(1,2)$ & $7$ & $S_1$ ($\dim4$) & $\PGL_2$ ($\dim3$) & rank $1$; plus a collineation $\ker\to\mathrm{coker}$\\
$(1,1,1)$ & $6$ & $S_1$ ($\dim4$) & $\PP^1\times\PP^1$ & rank $1$; plus a rank-one map $\ker\to\mathrm{coker}$\\
\bottomrule
\end{tabular}
\caption{Orbits of $X_{\rm ext}(\C^3)$. The fibre over each point of $S_1$ is $\PP(\mathrm{End}\,\C^2)=\PP^3\supset\PGL_2\sqcup(\PP^1\times\PP^1)$.}\label{tab:d3}
\end{table}

\begin{example}[The fibre $\PP^3$ over a rank-one point]\label{ex:d3fibre}
Let $\tau=[e_1e_1^*]\in S_1$ (image $\langle e_1\rangle$, kernel $\langle e_2,e_3\rangle$). The fibre $\pi^{-1}(\tau)$ is $\PP(\mathrm{End}\langle e_2,e_3\rangle)\cong\PP^3$, and by Theorem~\ref{thm:ext-boundary}(ii) every point $[B]$ of it prescribes a destination relation $R_B$:
\begin{itemize}[nosep]
\item $[B]\in\PGL_2$ (type $(1,2)$). $R_B([x])=\{e_1\}$ for $x\notin\PP\langle e_2,e_3\rangle$, and for $y\in\PP\langle e_2,e_3\rangle$, $R_B(y)=\PP\langle e_1,By\rangle$, the \emph{line} through $e_1$ and the point $[By]$. All the lines through $e_1$ occur, and the attachment is specified by the projective isomorphism $[B]$ from the source kernel to the target cokernel.
\item $[B]=[e_2e_2^*]$ (type $(1,1,1)$ with flags $\langle e_1\rangle\subset\langle e_1,e_2\rangle$ and $\langle e_3\rangle\subset\langle e_2,e_3\rangle$). $R_B(y)=\PP\langle e_1,e_2\rangle$ for $y\in\PP\langle e_2,e_3\rangle\setminus\{e_3\}$ and $R_B(e_3)=\PP^2$.
\end{itemize}
The second relation is the limit of the first as $B_t=\diag(1,t)\to[e_2e_2^*]$: for $y\ne e_3$ the lines $\PP\langle e_1,(y_2:ty_3)\rangle$ tend to $\PP\langle e_1,e_2\rangle$, while $R_{B_t}(e_3)=\PP\langle e_1,e_3\rangle$ for all $t\ne0$ jumps in the limit to the whole plane. The relations vary with $B$ in the Hausdorff topology; the values at a fixed source need not.
\end{example}
The two sequences in Remark~\ref{rem:SPwarn}(1) realise these two parts of the same fibre: $g'_m$ gives $[B]=[\mathrm{Id}]$, while $g_m$ gives $[B]=[e_2e_2^*]$.

\paragraph{$X_{\rm ext}$ is the graph compactification.}
The dimensions of the destination fibres recover the strata, and the fibres recover the graded maps. Thus the destinations determine the point of $X_{\rm ext}$.

\begin{theorem}[Completeness for projective space]\label{thm:graphcomp}
Let $\mathcal K(\PP(V)^2)$ be the space of closed subsets of $\PP(V)\times\PP(V)$ with the Hausdorff metric, and $\Gamma_g=\{(x,gx)\}$ the graph of $g\in\PGL(V)$.
\begin{enumerate}[label=\textup{(\roman*)},nosep]
\item For $\xi\in\partial X_{\rm ext}$ with datum $\mathcal L$ let $R_\xi=\{(x,y):y\in\D_{\mathcal L}(x)\}$ be the destination relation of Theorem~\ref{thm:ext-boundary}(ii) (and $R_g=\Gamma_g$ for $g\in\PGL$). Then $\xi\mapsto R_\xi$ is continuous, and \emph{injective}: $R_\xi$ determines $\mathcal L(\xi)$.
\item Consequently $\xi\mapsto R_\xi$ is a homeomorphism of $X_{\rm ext}$ onto the closure of $\{\Gamma_g:g\in\PGL(V)\}$ in $\mathcal K(\PP(V)^2)$.
\item The single-valued part of $R_\xi$ is $\pi(\xi)$: $R_\xi(x)=\{\tau^+x\}$ for $x\notin\PP(F_2)$.
\end{enumerate}
\end{theorem}
\begin{proof}
\emph{Continuity.} For $g_m\in\PGL$ with $\jmath(g_m)\to\xi$ the graphs $\Gamma_{g_m}$ converge to $R_\xi$ in the Hausdorff topology, by Theorem~\ref{thm:ext-boundary}(ii) (and by continuity of the action if $\xi\in\PGL$). For $\xi_m\to\xi$ in $X_{\rm ext}$ choose $g_m\in\PGL$ with $\jmath(g_m)$ close to $\xi_m$ and $\Gamma_{g_m}$ Hausdorff-close to $R_{\xi_m}$ (possible since $R_{\xi_m}$ is a Hausdorff limit of graphs) within $1/m$; then $g_m\to\xi$, so $R_{\xi_m}$ is within $1/m$ of $\Gamma_{g_m}$, which tends to $R_\xi$.

\emph{Injectivity.} On $\PGL$ it is faithfulness; a boundary relation has positive-dimensional fibres, so cannot be the graph of a bijection. By Proposition~\ref{prop:graded}, the distinct fibre dimensions $r_0<r_1<\cdots$ recover the strata $S_j=\PP(F_j)\setminus\PP(F_{j+1})$ and hence $\PP(F_j)=\overline{S_j}$. Next, $G_j=\mathrm{span}\bigcup_{x\in S_j}R(x)$, because $A_j$ is surjective on the quotient. Finally, quotienting the fibres by $G_{j-1}$ recovers the projective isomorphism $[v]\mapsto[A_jq_jv]$ from $\PP(F_j/F_{j+1})$ to $\PP(G_j/G_{j-1})$. This determines $A_j$ up to scalar, and hence the entire datum.

(ii) A continuous injection from a compact space to a Hausdorff space is a homeomorphism onto its image, and the image is the closure of the image of the dense subset $\PGL$. (iii) is Proposition~\ref{prop:graded} with $j=1$ ($G_0=0$).
\end{proof}
\begin{corollary}[No coarser compactification works]\label{cor:noquot}
Let $q:X_{\rm ext}\to Y$ be a continuous equivariant quotient with $\PGL$ embedded in $Y$ and suppose that the destination relation $\D_{(g_m)}$ on $\PP(V)$ is a function of the limit of $(g_m)$ in $Y$. Then $q$ is injective. In particular $\mathrm{SP}(V)$ cannot compute the destination relation for $d\ge3$, and a compactification in which only $(\tau^+,\tau^-)$ survive cannot do so for $d\ge4$. For $d=2$ or $3$, the pair of first forward and inverse pseudo-projective limits does determine the destination relation.
\end{corollary}
\begin{proof}
If $q(\xi)=q(\xi')$ then by hypothesis $R_\xi=R_{\xi'}$, hence $\xi=\xi'$ by Theorem~\ref{thm:graphcomp}. That $\mathrm{SP}$ is not injective for $d\ge3$ is Example~\ref{ex:d3fibre} (and, for $d\ge4$, the same pair embedded in a block).

To see what the forward--inverse pair forgets, let $d\ge4$ and, for $a=1,2$, take the determinant-one matrices
\[
g_m^{(a)}=\diag\bigl(e^m,a,a^{-1},\underbrace{1,\dots,1}_{d-4},e^{-m}\bigr).
\]
Both sequences have type $(1,d-2,1)$, with the bounded middle block $\diag(a,a^{-1},1,\dots,1)$, and
\[
\tau^+=\lim_m[g_m^{(a)}]=[e_1e_1^*],\qquad
\tau^-=\lim_m[(g_m^{(a)})^{-1}]=[e_de_d^*]
\]
for either value of $a$. But Proposition~\ref{prop:proj} gives
\[
\D_{(g_m^{(a)})}([e_2+e_3])=
\PP\bigl(\mathbb K e_1+\mathbb K(ae_2+a^{-1}e_3)\bigr).
\]
These are different lines for $a=1$ and $a=2$. Thus the pair does not determine the destination relation for any $d\ge4$. These examples concern sequences; no common discrete subgroup is asserted.

For $d=3$, the canonical identification
$\bigwedge^2V\cong V^*\otimes\det V$ identifies
$\bigwedge^2g$ with $(\det g)g^{-T}$.
The two exterior factors of degrees one and two are therefore precisely
the forward limit and the dual of the inverse limit. The pair determines
the exterior datum, hence the destination relation by
Theorem~\ref{thm:ext-boundary}. In dimension two only the first exterior
factor is needed, as in the classical convergence principle.
\end{proof}
\begin{remark}
Thus $X_{\rm ext}(V)$ is the space of correspondences on $\PP(V)$ arising as Hausdorff limits of collineations; inversion makes it closed under transposition. The space $\mathrm{SP}(V)$ instead retains the partially defined collineations, with convergence off the kernel (Proposition~\ref{prop:SP}(iv)).
\end{remark}

\medskip
Table~\ref{tab:SPvsX} summarises the comparison.
\begin{table}[ht]
\centering\small
\renewcommand{\arraystretch}{1.15}
\begin{tabular}{@{}p{3.3cm}p{5.6cm}p{5.6cm}@{}}
\toprule
 & $\mathrm{SP}(V)=\PP(\mathrm{End}\,V)$ & $X_{\rm ext}(V)$ (wonderful)\\
\midrule
Definition & one projective space; matrices up to scalar & closure of $\PGL$ in $\prod_{q<d}\PP(\mathrm{End}\wedge^qV)$; iterated blow-up of $\mathrm{SP}$\\
Dimension and smoothness & $d^2-1$, smooth & $d^2-1$, smooth\\
Boundary & irreducible hypersurface $\det=0$; singular for $d\ge3$, smooth quadric for $d=2$ & normal-crossings divisor $D_1\cup\dots\cup D_{d-1}$\\
Orbits of $\PGL\times\PGL$ & $d$ (rank strata); closures singular only for $2\le r\le d-1$ & $2^{d-1}$ (compositions of $d$); closures smooth\\
Closed orbit & rank one: $\PP^{d-1}\times\PP^{d-1*}$ & $\mathrm{Fl}(V)\times\mathrm{Fl}(V)$\\
Inversion $g\mapsto g^{-1}$ & extends for $d=2$; not for $d\ge3$ & extends\\
Boundary point records & $(G_1,F_2,[A_1])$: fastest scale only & whole graded datum $(F_\bullet,G_\bullet,[A_j])$\\
Determines $\D_{(g_m)}$ on $\PP(V)$ & yes for $d=2$; no for $d\ge3$ (Ex.~\ref{ex:d3fibre}) & yes, and conversely (Thm.~\ref{thm:graphcomp})\\
Determines kernels, images and the exceptional set for stabilised Cartan data & yes (Prop.~\ref{prop:kernel}) & yes (via $\pi$)\\
Real case $d$ even, $\PSL(d,\R)$ & closure contains the whole singular boundary & contains the whole boundary for $d=2$; omits some points for even $d\ge4$\\
\bottomrule
\end{tabular}
\caption{The two compactifications of $\PGL(V)$.}\label{tab:SPvsX}
\end{table}

\medskip For discrete groups this comparison leads to the kernel and image descriptions of limit sets in Chapter~\ref{ch:limit} (Propositions~\ref{prop:kernel} and~\ref{prop:Tits-comparison}), and to the finer flag $\Lambda_{\CoG}\subseteq\cdots\subseteq\Lambda_{\Myr}$ of Chapter~\ref{ch:psl}. The rank-one first limits give the first level; the full datum keeps the subsequent quotients and their destinations.

\medskip Before passing to a general group, we must keep track of one real distinction. Let $\mathbb K=\R$ and $d$ even. The group $\PSL(d,\R)$ is the index-two subgroup of $\PGL(d,\R)$ of classes of positive determinant. Every boundary point of $\mathrm{SP}(\R^d)$ is a limit of matrices of either sign of determinant (bring a rank-$r$ representative to $\diag(I_r,0)$ by fixed changes of basis and perturb the zero block by $t\diag(\pm1,1,\dots,1)$), so $\overline{\PSL(d,\R)}=\PSL(d,\R)\cup\{[A]:\det A=0\}$ inside $\mathrm{SP}(\R^d)$; the negative-determinant invertible component is omitted. In $X_{\rm ext}(\R^d)$ the answer depends on the blocks. The determinant lines of the two flags identify the determinant of $\bigoplus A_j$ with a nonzero real number up to the factors $c_j^{b_j}$ coming from $A_j\mapsto c_jA_j$. If every $b_j$ is even, its sign is therefore well defined and locally constant on the orbit $O_{\mathbf b}$; nearby invertible matrices have that sign. If some $b_j$ is odd, changing its representative changes the sign without changing the datum, and the realisation in Theorem~\ref{thm:ext-boundary}(iv) gives approaching matrices of either sign. Thus for even $d\ge4$ the closure of $\PSL(d,\R)$ omits exactly the boundary data whose blocks are all even and whose orientation is negative. For $d=2$ the only boundary type is $(1,1)$, so the whole boundary is present, in agreement with $X_{\rm ext}=\mathrm{SP}$. Theorem~A applies to sequences in $\PSL(d,\R)$ and the destination relations occurring are those of the points of this closure only.

\medskip
The distinction between the first limit and the full datum also appears in representations. Let $G$ be a complex semisimple adjoint group and let $\rho_\lambda:\widetilde G\to\mathrm{GL}(V_\lambda)$ be an irreducible representation of its simply connected cover. The centre acts by scalars, so $[\rho_\lambda]$ descends to $\bar\rho_\lambda:G\to\PGL(V_\lambda)$. Put $G_{\rm eff}=G/\ker\bar\rho_\lambda$. The closure $\mathrm{SP}_\lambda(G)$ of the projective image is the \emph{simple} group compactification~\cite{BraviGandiniMaffeiRuzzi2011} of $G_{\rm eff}$. It compactifies $G$ itself precisely when the projective action is faithful, or equivalently when no simple factor of $G$ is omitted by $\lambda$. For $G=\PGL(V)$ and $\lambda=\omega_1$, the standard representation is linear on the cover $\mathrm{SL}(V)$ and projective on $G$; its closure is $\mathrm{SP}(V)$.

When $\lambda$ is regular, every simple factor is seen and $\mathrm{SP}_\lambda(G)$ is the wonderful compactification $W(G)$ of De~Concini--Procesi, smooth, with boundary a normal-crossings divisor with $|\Delta|$ components and $2^{|\Delta|}$ orbits. For $\PGL(V)$ this is $X_{\rm ext}(V)$: the product of independently projectivised exterior factors encodes the fundamental representations, and the representation of the cover of regular weight $\sum_{q=1}^{d-1}\omega_q$ gives the same compactification. Neither the individual fundamental representations nor their sum of weights need descend linearly to the adjoint group. For a weight with support $\theta\subset\Delta$, $\mathrm{SP}_\lambda(G)$ is a $G\times G$-equivariant quotient of $W(G)$, with the omitted factors absorbed in $G_{\rm eff}$. The weight-space normalisation gives a limit $[\rho_\lambda(k)D\rho_\lambda(k')]$; the surviving weights are those for which $\lambda-\nu$ has no contribution from a divergent simple root, and $D$ also retains the finite Levi factors. The orbits of $W(G)$ correspond to subsets of $\Delta$, i.e.\ to the types $\Theta$ of Theorem~A (with open orbit $\Theta=\varnothing$), and the orbit of type $\Theta$ fibres over $G/P_\Theta\times G/P^-_\Theta$ with fibre the adjoint group of the Levi factor $L_\Theta$. In the compact frames $k,k'$ the point of this fibre is represented by $[b]$. Changing those frames to $kq_1,q_2k'$ changes the representative to $[q_1^{-1}bq_2^{-1}]$, with $q_1,q_2\in K\cap L_\Theta$.

So the pattern of the present section is general: Theorem~A records the destination relation as a function on $W(G)$ (a point of the orbit of type $\Theta$ records $kP_\Theta$, $k'^{-1}P^-_\Theta$ and the relative Levi class over this pair), while $\mathrm{SP}_\rho$, for a Tits representation supported on $\alpha$, retains the top operator, with its source and target data, and may forget slower scales. Its rank-one limits record the two flags in $\mathcal F_\alpha$ and the opposite flag manifold. The Tits representations of Section~\ref{sec:Quint} have highest restricted weights that are positive multiples of the fundamental weights; their associated pseudo-projective spaces play this role, and Quint's set is the set of images of the rank-one points.

For a connected centre-free real semisimple group $G$, let $\mathbf G$ be
its real adjoint algebraic group, so $G=\mathbf G(\R)^\circ$. Chanfi's construction
gives the algebraic wonderful compactification $W(\mathbf G)$ over $\R$,
including for non-split forms~\cite[\S1.3]{Chanfi2024}. Take the closure
of $G$ in $W(\mathbf G)(\R)$ for the real topology. To identify this closure
with a model of the next section, choose a representation of the simply
connected cover of $\mathbf G_\C$ with regular complex highest weight,
for instance the sum of the complex fundamental weights. Its projective
closure is the complexification
of $W(\mathbf G)$, and its restricted highest weight is regular. Under this
identification the real closure of $G$ is $X_\varrho$ for the restricted
representation. Theorem~\ref{thm:extG} and Theorem~A therefore recover the
destinations on every real flag manifold from each attained boundary point.
For a discrete $\Gamma<G$, the union of the relations over the boundary
points of its closure is $\D_\Gamma$ by Proposition~\ref{prop:Dunion}.
The supports of limiting real graphs can forget part of this boundary
datum.

%--------------------------------------------------------------
%--------------------------------------------------------------
\section{The exterior boundary of a general group}\label{sec:extG}
For $\PSL(d)$ the flags and graded maps associated with $(\Theta,k,k',b)$ give a point of the space of complete collineations (Theorem~\ref{thm:ext-boundary}). To retain the same information for an arbitrary $G$ of Section~\ref{sec:setting-sum}, we choose a representation that is faithful projectively and whose highest restricted weight sees every simple root. Its weight decomposition will recover the type, the two parabolics and the finite Levi class, and hence every destination. The regular complex highest-weight choice above gives the real wonderful model.

\paragraph{The representation.} Let $\varrho:\tilde G\to\mathrm{GL}(V)$ be a finite-dimensional complex irreducible representation of the universal cover whose highest \emph{restricted} weight $\lambda$ is regular, $\langle\lambda,\alpha^\vee\rangle>0$ for all $\alpha\in\Delta$. We can obtain such a $\varrho$ by taking the complex highest weight to be the sum of all fundamental weights, with a Borel subalgebra of $\g_\C$ compatible with $\ak$ and the chosen positive system.

The centre of $\tilde G$ acts by scalars, so $\varrho$ induces $\bar\varrho:G\to\PGL(V)$. We require it to act non-trivially on every simple factor, including the compact ones. Its differential then has zero projective kernel: on a simple ideal a scalar action would be zero, while a nonzero representation is faithful at the Lie-algebra level. The projective kernel in the connected centre-free group $G$ is consequently discrete and central, hence trivial. The construction using all complex fundamental weights has this property, so no factor is hidden.

Choose a Hermitian form compatible with a Cartan involution $\vartheta$, so $d\varrho(\vartheta X)=-d\varrho(X)^*$; it is obtained by making the compact form $\kk\oplus i\pk$ unitary. Thus $K$ acts unitarily and the restricted weight spaces are orthogonal. The projective image is closed in $\PGL(V)$. Indeed, $s(H)=\max_\nu\nu(H)-\min_\nu\nu(H)$ is positive for every $0\ne H\in\ak$: otherwise $d\varrho(H)$ would be scalar, contrary to the zero projective differential kernel. By homogeneity and compactness of the unit sphere, $s(H)\ge c\|H\|$ for some $c>0$. If $\|\mu_m\|\to\infty$, the condition number $e^{s(\mu_m)}$ tends to infinity, so the projective operators cannot converge to an invertible one. Hence $\bar\varrho$ is an embedding, and its closure in $\PP(\mathrm{End}\,V)$ compactifies $G$. In all formulas $\varrho(g)$ for $g\in G$ means any lift to $\tilde G$; the resulting projective operator is independent of that choice.

Write $V=\bigoplus V_\nu$ for the weight decomposition under $A$ (restricted weights); every weight is $\nu=\lambda-\sum_{\alpha\in\Delta}n_\alpha\alpha$ with $n_\alpha\in\mathbb Z_{\ge0}$. For $\Theta\subseteq\Delta$ put
\[
V^\Theta=\bigoplus\{V_\nu:\ n_\alpha(\nu)=0\ \ \forall\alpha\in\Theta\},\qquad V=V^\Theta\oplus V^{\Theta,\perp},
\]
$V^{\Theta,\perp}$ the sum of the other weight spaces, and $\Pi_\Theta$ the projection onto $V^\Theta$ along $V^{\Theta,\perp}$. $V^\Theta$ is the eigenspace of $\varrho(A_\Theta)$ of largest weight, hence $L_\Theta$-stable, and $\Pi_\Theta$ commutes with $\varrho(L_\Theta)$.

\begin{lemma}\label{lem:Pi}\
(i) Use the canonical exponential lifts of $N_\Theta$ and $N_\Theta^-$. Then $\varrho(n)\Pi_\Theta=\Pi_\Theta$ for $n\in N_\Theta$ (a raising operator applied to $V^\Theta$ leaves the set of weights) and $\Pi_\Theta\varrho(n^-)=\Pi_\Theta$ for $n^-\in N_\Theta^-$ (lowering operators have no component in $V^\Theta$). (ii) The stabiliser of the subspace $V^\Theta$ in $G$ is $P_\Theta$, and that of $V^{\Theta,\perp}$ is $P_\Theta^-$.
\end{lemma}
\begin{proof}
(i) If $\gamma$ is a root of $N_\Theta$, it has a positive coefficient at some $\alpha\in\Theta$. For a weight in $V^\Theta$, raising by $\gamma$ would make that coefficient of $\lambda-\nu$ negative, so the raising operator is zero. Lowering cannot send any weight outside $V^\Theta$ into $V^\Theta$. Exponentiating gives the two identities.

(ii) The stabiliser $H$ of $V^\Theta$ contains $P_\Theta$, since $L_\Theta$ preserves it and $N_\Theta$ fixes it pointwise. In particular $H$ contains the minimal real parabolic $B$. Closed overgroups of $B$ are the standard real parabolics, obtained by adjoining the negative simple restricted root spaces \cite[Propositions~7.76, 7.82(d) and~7.83(b)]{KnappLG}. To prove $H=P_\Theta$ it is therefore enough to exclude $\g_{-\alpha}$ for each simple $\alpha\in\Theta$; no assertion about an arbitrary negative root is needed.

Here is the exclusion also for non-split real groups. Choose $0\ne X\in\g_\alpha$ and put $Y=-\vartheta X\in\g_{-\alpha}$. Then $[X,Y]=cH_\alpha$ with $c>0$ and $H_\alpha$ the positive restricted coroot direction: the bracket belongs to $\pk\cap\g_0=\ak$, and its positivity follows from the positive form $-B(\cdot,\vartheta\cdot)$. For $0\ne v\in V_\lambda$, $d\varrho(X)v=0$ and $d\varrho(Y)=d\varrho(X)^*$. Hence
\[
 \|d\varrho(Y)v\|^2
 =\langle d\varrho([X,Y])v,v\rangle
 =c\lambda(H_\alpha)\|v\|^2>0.
\]
The resulting vector has restricted weight $\lambda-\alpha$, outside $V^\Theta$. Thus $\g_{-\alpha}$ does not preserve $V^\Theta$. This uses neither one-dimensional restricted root spaces nor a reduced restricted root system. The classification now gives $H=P_\Theta$. Finally $V^{\Theta,\perp}$ is the orthogonal complement for the chosen form, and $\varrho(g)^*=\varrho(\vartheta(g)^{-1})$ up to scalar, so its stabiliser is $\vartheta(H)=P_\Theta^-$.
\end{proof}

\begin{theorem}[Limits of $\varrho(g_m)$ and the exterior boundary]\label{thm:extG}\
Let $g_m=k_me^{\mu_m}k'_m$ have Cartan data $(\Theta,k,k',b)$ in the sense of Lemma~\ref{lem:split}. Then in $\PP(\mathrm{End}\,V)$
\[
\lim_m[\varrho(g_m)]=[\varrho(k)\,\Pi_\Theta\varrho(b)\,\varrho(k')].
\]
Let $X_\varrho$ be the closure of $\bar\varrho(G)$ in $\PP(\mathrm{End}\,V)$ and $\partial X_\varrho=X_\varrho\setminus\bar\varrho(G)$.
\begin{enumerate}[label=\textup{(\alph*)},nosep]
\item $\partial X_\varrho=\bigsqcup_{\varnothing\ne\Theta\subseteq\Delta}\mathcal S_\Theta$, with $\mathcal S_\Theta=\{[\varrho(k)\Pi_\Theta\varrho(\ell)\varrho(k')]:k,k'\in K,\ \ell\in L_\Theta\}$, and $\mathcal S_\Theta$ is the $G\times G$-orbit of $[\Pi_\Theta]$ for the action $(g,h)\cdot\Psi=\varrho(g)\Psi\varrho(h)^{-1}$.
\item The map $\mathcal S_\Theta\to G/P_\Theta\times G/P^-_\Theta$, $[\Psi]\mapsto(\text{stabiliser of }\mathrm{Im}\Psi,\ \text{stabiliser of }\ker\Psi)$, is a $G\times G$-equivariant fibre bundle with fibre $L_\Theta/Z_\Theta$, where $Z_\Theta$ is the subgroup of $L_\Theta$ acting by scalars on $V^\Theta$ (it contains $A_\Theta$). For the limit of $(g_m)$, the stabiliser of the image is the attracting parabolic $kP_\Theta k^{-1}$ and the stabiliser of the kernel is the repelling parabolic $k'^{-1}P_\Theta^-k'$ of Section~\ref{sec:balanced}.
\item \textup{(Realisation)} Every point of $\mathcal S_\Theta$ is the limit of $g_m=k\exp(m\nu_\Theta)\,\ell\,k'$, where $\nu_\Theta\in\ak_\Theta$ has $\alpha(\nu_\Theta)=1$ for $\alpha\in\Theta$ and $0$ for $\alpha\notin\Theta$; this sequence has type exactly $\Theta$. Hence $\partial X_\varrho$ is exactly the set of limits of divergent sequences, and \emph{every datum $(\Theta,\xi^+,\xi^-,[\ell])$ occurs}.
\item The destination relation $\D_{(g_m)}$ depends on $(g_m)$ only through the limit point in $\partial X_\varrho$.
\end{enumerate}
\end{theorem}
\begin{proof}
\emph{The limit.} On $V_\nu$, $\varrho(e^{\mu_m})$ is the scalar $e^{\nu(\mu_m)}=e^{\lambda(\mu_m)}e^{-\sum n_\alpha\alpha(\mu_m)}$. If $n_\alpha>0$ for some $\alpha\in\Theta$ this is $o(e^{\lambda(\mu_m)})$, because $\alpha(\mu_m)\to\infty$ and all $\alpha'(\mu_m)\ge0$. If $n_\alpha=0$ for $\alpha\in\Theta$, then only $\alpha\notin\Theta$ occur, $\alpha(\mu_m)=\alpha(\eta_m)\to\alpha(\eta)$, and the scalar divided by $e^{\lambda(\mu_m)}$ tends to $e^{-\sum n_\alpha\alpha(\eta)}$, which is $\varrho(b)|_{V_\nu}$ divided by $e^{\lambda(\eta)}$ ($b=\exp\eta$). Since $k_m\to k$, $k'_m\to k'$, the claim follows.
\emph{(a)} By compactness every sequence in $G$ leaving compacts has a subsequence with a Cartan datum, so limits lie in the right-hand side; conversely (c) shows that every point is attained. The orbit statement: by Lemma~\ref{lem:Pi}(i), if $g=kp$ with $p=ln\in P_\Theta$ and $h^{-1}=p'^-k'$ with $p'^-=n^-l'\in P_\Theta^-$ (Iwasawa-type decompositions $G=KP_\Theta=P_\Theta^-K$), then $[\varrho(g)\Pi_\Theta\varrho(h)^{-1}]=[\varrho(k)\varrho(l)\Pi_\Theta\varrho(l')\varrho(k')]$.
\emph{(b)} By Lemma~\ref{lem:Pi}(ii), $\mathrm{Im}=\varrho(k)V^\Theta$ has stabiliser $kP_\Theta k^{-1}$, and $\ker=\varrho(k')^{-1}V^{\Theta,\perp}$ has stabiliser $k'^{-1}P^-_\Theta k'$. The $G\times G$-stabiliser of $[\Pi_\Theta]$ is $\{(ln,\,l'n^-)\}$ with $\varrho(l)\Pi_\Theta=c\,\Pi_\Theta\varrho(l')$, i.e.\ $l'\equiv l$ modulo $Z_\Theta$; this gives the fibre $L_\Theta/Z_\Theta$.
\emph{(c)} $\varrho(\exp(m\nu_\Theta))$ acts on $V^\Theta$ (top eigenspace of $\varrho(A_\Theta)$) by the scalar $e^{m\lambda(\nu_\Theta)}$ and on the other eigenspaces by strictly smaller powers; $\ell\in L_\Theta$ commutes with $A_\Theta$ and preserves the eigenspaces. Hence $e^{-m\lambda(\nu_\Theta)}\varrho(\exp(m\nu_\Theta)\ell)\to\Pi_\Theta\varrho(\ell)$. For the type, write $\ell=u\exp(\eta)v$ with $u,v\in K\cap L_\Theta$ and $\eta$ dominant for the Levi roots. This polar decomposition stays in the centraliser of $A_\Theta$: the positive polar factor and its logarithm commute with $A_\Theta$, and the compact factor does as well. The logarithm can be conjugated into $\ak$ by $K\cap L_\Theta$ and then into the closed Levi chamber. Since $u,v$ commute with $\exp(m\nu_\Theta)$ and $m\nu_\Theta+\eta$ belongs to $\ak^+$ for all sufficiently large $m$, this is exactly the Cartan projection of $k\exp(m\nu_\Theta)\ell k'$. Thus the complementary root values are constant and the root values on $\Theta$ are $m+\alpha(\eta)$.
\emph{(d)} Replace $(k,k',b)$ by $(kq_1,q_2k',\,q_1^{-1}bq_2^{-1})$ with $q_1,q_2\in K\cap L_\Theta$; the limit point is unchanged since $\Pi_\Theta$ commutes with $\varrho(L_\Theta)$. By Theorem~\ref{thm:dest} and the $L_\Theta$-equivariance of $\rho$ and the normalisation of $N_\Theta$ by $L_\Theta$, $k q_1\,\overline{N_\Theta\,q_1^{-1}b q_2^{-1}\rho(q_2k'x)}=k\,\overline{N_\Theta\,b\rho(k'x)}$ (Remark~\ref{rem:generalb} for $b\notin A$). Two data with the same limit point have the same $\Theta$, parabolics and (by (b)) Levi classes modulo $Z_\Theta$, hence differ by such a change and by an element $z\in Z_\Theta$; the scalar ambiguity acts trivially on every $L_\Theta$-flag manifold. Indeed, for $J\subseteq\Delta\setminus\Theta$, the subspace $V^{\Theta\cup J}\subset V^\Theta$ has stabiliser $L_\Theta\cap P_{\Theta\cup J}$ in $L_\Theta$, by Lemma~\ref{lem:Pi}(ii). An element $z\in Z_\Theta$ fixes every subspace $lV^{\Theta\cup J}$, $l\in L_\Theta$, because it is scalar on $V^\Theta$. It therefore acts trivially on $L_\Theta/(L_\Theta\cap P_{\Theta\cup J})$. The classification of parabolics of the reductive Levi gives all its flag manifolds this way, up to conjugacy; compact simple factors are already contained in their parabolics. Thus $z$ fixes $\rho(y)$ and normalises $N_\Theta$, so it does not change $\overline{N_\Theta\,b\rho(y)}$.
\end{proof}
\begin{remark}
(1) For $\PSL(d)$ the corresponding construction in the \emph{product} $\prod_{q=1}^{d-1}\PP(\mathrm{End}\,\wedge^qV)$, with an independent normalisation in each factor, is Theorem~\ref{thm:ext-boundary}. On the linear cover $\mathrm{SL}(d)$, a direct sum $\bigoplus_q\wedge^qV$ inside one $\PP(\mathrm{End})$ has only one scalar normalisation and can lose entire summands. The centre acts by different characters on the summands in general, whereas independent projectivisation gives a well-defined action of $\PSL(d)$. For $V=\g_\C$ and $\mathrm{Ad}$, which need not be irreducible or regular, the weight-space normalisation gives the first pseudo-projective limit of Chapter~\ref{ch:adjoint}, after a further extraction if the relative gaps between simple factors require it. The recovery assertions (a)--(d) are not claimed without the hypotheses above. (2)~When $G$ is a complex adjoint group and $\varrho$ is a holomorphic algebraic representation of its simply connected cover with regular complex highest weight, $X_\varrho$ is the wonderful compactification of De Concini--Procesi, and for two regular $\lambda$ the identity of $G$ extends to an isomorphism $X_\varrho\cong X_{\varrho'}$. This is not used: the statements above are proved for each $\varrho$ separately. The closure in the product $\prod_i\PP(\mathrm{End}\,V_{\omega_i})$ over the fundamental representations has the same strata, but the identification with $X_\varrho$ is not proved here. (3)~For real $G$ the weight spaces $V_\nu$ may have dimension $>1$, and the argument above allows these multiplicities. Proximality requires a one-dimensional top space, as in the Tits representations used in Proposition~\ref{prop:highest-weight}.
\end{remark}

%--------------------------------------------------------------
\section{The stratified convergence principle for a general group}\label{sec:strat}
A boundary point now remembers all the destinations. We return to the flag manifold and let the source vary. The classical convergence principle has an open set of sources with one common destination and a single exceptional source whose destination is the whole sphere. In higher rank the space $X$ has finitely many \emph{strata}. On each stratum the destinations are translates of one fixed Schubert variety, and the retraction tells us how the translate varies: the family is parametrised by a flag manifold of the Levi factor. The inverse sequence gives the opposite picture. For $\PSL(d,\mathbb K)$ this is Theorem~\ref{thm:ext-boundary}; we may now state it for an arbitrary $G$ and an arbitrary flag manifold.

\paragraph{Notation.} For $\Theta\subseteq\Delta$ non-empty and $\sigma\in\mathcal C_\Theta=W_\Theta\backslash W/W_\theta$ let $O^-_\sigma=P^-_\Theta\sigma Q$ and $O^+_\sigma=P_\Theta\sigma Q$ (source and target cells), and with $\Sigma_+,\Sigma_0,\Sigma_-$ as in Lemma~\ref{lem:cells}(iv) put
\[
d_\pm(\sigma)=\sum_{\lambda\in\Sigma_\pm(\sigma)}\dim\g_\lambda,\qquad d_0(\sigma)=\sum_{\lambda\in\Sigma_0(\sigma)}\dim\g_\lambda,\qquad d_++d_0+d_-=\dim X.
\]
Then $\dim O^+_\sigma=d_++d_0$, $\dim O^-_\sigma=d_0+d_-$, and $d_0(\sigma)$ is the dimension of the flag manifold $L_\Theta\sigma Q$ of $L_\Theta$. If $\sigma$ is the representative of minimal length in its double coset, put $X_\sigma=\overline{B\sigma Q}$, the \emph{Schubert variety} of $\sigma$. Finally $\rho^+:X\to X^{A_\Theta}$ is the stratified retraction of Lemma~\ref{lem:cells}(iii) for the opposite direction, $\rho^+(x)=\lim_{z\to\infty}z^{-1}x$, with smooth restrictions to the cells $O^+_\tau$.

\Needspace{8\baselineskip}
\begin{theorem}[Stratified convergence principle]\label{thm:S}
Let $G$ be as in Section~\ref{sec:setting-sum}, $X=G/Q$ a flag manifold, $\varrho$ the representation of Theorem~\ref{thm:extG}, and $(g_m)$ a sequence leaving every compact set of $G$.
\begin{enumerate}[label=\textup{(\roman*)},nosep]
\item \textup{(The asymptotic datum.)} After extraction $[\varrho(g_m)]$ converges to a point $\xi$ of the exterior boundary $\partial X_\varrho$, which lies in a stratum $\mathcal S_\Theta$, $\Theta\ne\varnothing$. The point $\xi$ is the same thing as the \emph{graded datum} $\mathcal L(\xi)=(\Theta,\ \xi^+,\ \xi^-,\ [b])$: the type $\Theta$, the attracting parabolic $\xi^+=kP_\Theta k^{-1}$, the repelling parabolic $\xi^-=k'^{-1}P^-_\Theta k'$, and a point $[b]$ of the fibre $L_\Theta/Z_\Theta$ over $(\xi^+,\xi^-)$ \textup{(}the Levi class\textup{)}.
\item \textup{(Strata.)} The sets
\[
\mathcal O^-_\sigma(\xi)=k'^{-1}O^-_\sigma\quad(\text{sources}),\qquad\mathcal O^+_\tau(\xi)=k\,O^+_\tau\quad(\text{targets}),\qquad\sigma,\tau\in\mathcal C_\Theta,
\]
decompose $X$ into finitely many strata, which depend only on $\xi^-$ and $\xi^+$ respectively.
\item \textup{(Destination formula, source side.)} Let $x\in\mathcal O^-_\sigma(\xi)$, with $\sigma$ of minimal length, and $\rho(k'x)=l\,\sigma Q$, $l\in L_\Theta$. Then
\[
\D_{(g_m)}(x)=(k\,b\,l)\cdot X_\sigma,
\]
a translate of one Schubert variety of dimension $d_+(\sigma)$, which depends on $x$ only through the point $\rho(k'x)$ of the flag manifold $L_\Theta\sigma Q$. It lies in $\bigcup_{\tau\le\sigma}\mathcal O^+_\tau(\xi)$, and its intersection with $\mathcal O^+_\sigma(\xi)$ is open and dense in the destination fibre. The fibre is a single point iff $\Sigma_+(\sigma)=\varnothing$.
\item \textup{(Inverse formula, target side.)} Let $y\in\mathcal O^+_\tau(\xi)$ and $\rho^+(k^{-1}y)=l'\tau Q$. Then
\[
\D_{(g_m^{-1})}(y)=k'^{-1}\,\overline{N_\Theta^-\,b^{-1}l'\,\tau Q}=(k'^{-1}b^{-1}l')\cdot\overline{N^-_\Theta\tau Q},
\]
the closure of a cell of dimension $d_-(\tau)$; and $\D_{(g_m^{-1})}$ is the transpose of $\D_{(g_m)}$.
\item \textup{(Bundles.)} The map $x\mapsto\rho(k'x)$ makes $\mathcal O^-_\sigma(\xi)$ a fibre bundle over $L_\Theta\sigma Q$ whose fibres are cells of dimension $d_-(\sigma)$, and $x\mapsto\D_{(g_m)}(x)$ is constant on its fibres. The translates of $X_\sigma$ form a locally trivial family over $L_\Theta\sigma Q$, and the destination incidence over $\mathcal O^-_\sigma(\xi)$ is its pullback. For complex $G$ the bundles and local trivialisations are holomorphic, with complex dimensions equal to half the displayed real dimensions.
\item \textup{(Uniqueness and realisation.)} The relation $\D_{(g_m)}$ depends on $(g_m)$ only through $\xi$: every sequence with the same limit point has the same relation. Every $\Theta\ne\varnothing$ and every datum $(\Theta,\xi^+,\xi^-,[b])$ occurs, for a sequence of type exactly $\Theta$.
\end{enumerate}
\end{theorem}
\begin{proof}
(i) and the realisation in (vi) follow from Theorem~\ref{thm:extG}(b),(c). For uniqueness in (vi), consider any sequence whose projective operators converge to $\xi$. Every subsequence has a further Cartan extraction with the same exterior limit. By Theorem~\ref{thm:extG}(d) all these extracted relations agree, and Theorem~\ref{thm:dest}(d) gives Hausdorff convergence of their graphs to this common relation. The graphs of the full sequence therefore converge to the same relation. Each pair in it is realised along the full sequence by choosing points in the graphs approaching the pair, and every full-sequence limit belongs to it. (ii) follows from Lemma~\ref{lem:cells}(i) and its mirror image for $P_\Theta$, translated by $k'^{-1}$ and $k$.
(iii) By Theorem~\ref{thm:dest}(b), $\D_{(g_m)}(x)=k\overline{N_\Theta\,b\rho(k'x)}=k\overline{N_\Theta\,bl\sigma Q}$. As $bl\in L_\Theta$ normalises $N_\Theta$, this is $kbl\,\overline{N_\Theta\sigma Q}$. For $\sigma$ minimal, the proof of Proposition~\ref{prop:Eorder} shows $N_\Theta\sigma Q=N\sigma Q=B\sigma Q$, so $\overline{N_\Theta\sigma Q}=X_\sigma$, of dimension $d_+(\sigma)$ by Theorem~\ref{thm:dest}(c), and $B\sigma Q$ is open and dense in its closure; it lies in $O^+_\sigma$, while the rest of $X_\sigma$ is a union of $B\tau'Q$ with $\tau'<\sigma$, which lie in the strata $O^+_{[\tau']}$, $[\tau']<\sigma$ by the minimality of $\sigma$ in its double coset and Lemma~\ref{lem:doublecosets}. The elements $b,l\in L_\Theta\subset P_\Theta$ preserve every $O^+_\tau$. The point criterion follows from the cell dimension $d_+(\sigma)$.
(iv) Write $g_m^{-1}=k_m'^{-1}a_m^{-1}k_m^{-1}$. Now $a_m^{-1}\in\exp(-\ak^+)$ is a Cartan sequence for the opposite positive system: the roles of $N_\Theta,N_\Theta^-$ and of $\rho,\rho^+$ are exchanged, the Levi class is $b^{-1}$, and the proof of Theorem~\ref{thm:dest}(a),(b) is symmetric. This gives the formula; the last sentence is Theorem~\ref{thm:dest}(d), and the dimension is as in (iii) with $N_\Theta^-$ in place of $N_\Theta$.
(v) Put $S_0=L_\Theta\cap\sigma Q\sigma^{-1}$ and $U_0=N_\Theta^-\cap\sigma Q\sigma^{-1}$. The stabiliser of $\sigma Q$ in $P_\Theta^-$ is $U_0S_0$. Indeed, if $n l$ fixes $\sigma Q$ with $n\in N_\Theta^-$ and $l\in L_\Theta$, conjugation by $\exp(t\nu_\Theta)$ contracts $n$ and fixes $l$. It preserves the stabiliser because $\sigma Q$ is fixed by $A_\Theta$. Taking $t\to+\infty$ shows $l\in S_0$ and then $n\in U_0$. Consequently the source stratum, in the source frame, is the associated bundle $L_\Theta\times^{S_0}(N_\Theta^-/U_0)$ over $L_\Theta/S_0$. Its fibre is a cell of dimension $d_-(\sigma)$ by Lemma~\ref{lem:cells}(iv) for the opposite system. Since $S_0$ fixes $\sigma Q$ and normalises $N_\Theta$, it preserves $X_\sigma=\overline{N_\Theta\sigma Q}$. The bundle $L_\Theta\times^{S_0}X_\sigma$ therefore parametrises the destination fibres by $[l,y]\mapsto (l\sigma Q,kbl\,y)$, by (iii). Pulling it back along $x\mapsto\rho(k'x)$ gives the claimed incidence. Local smooth sections of $L_\Theta\to L_\Theta/S_0$, holomorphic for complex $G$, give the local trivialisations.
\end{proof}

\begin{remark}[Rank one and the classical principle]
For the non-trivial rank-one flag manifold and $\Theta=\Delta=\{\alpha\}$, $\mathcal C_\Theta=\{1,w\}$: $\sigma=1$ has $\Sigma_+=\varnothing$ and $X_1=\{\mathrm{pt}\}$, $\sigma=w$ has $X_w=X$. So there are two source strata, the open one (destination the point $x^+$) and the point $x^-$ (destination everything): the convergence principle. In $\PSL(2,\C)$ there is no Levi datum, $L_\Theta/Z_\Theta$ is a point and $\xi$ is just the pair $(x^+,x^-)$.
\end{remark}

\begin{example}\label{ex:SPSL} Consider  the model case $\PSL(d,\mathbb K)$. 
Let $X=\PP(V)$ and let $(g_m)$ have Cartan blocks $d=b_1+\dots+b_k$ and flags $F_\bullet$, $G_\bullet$ as in Section~\ref{sec:ext}. The type $\Theta$ is the set of the $k-1$ divergent gaps; $\mathcal C_\Theta=\{1,\dots,k\}$; the source stratum of index $j$ is $\PP(F_j\setminus F_{j+1})$ and the target stratum is $\PP(G_j\setminus G_{j-1})$; if $\sigma Q=[e_\sigma]$ in the $j$th standard block, then $X_\sigma$ is the projectivisation of the earlier standard blocks together with $\mathbb K e_\sigma$, so $kX_\sigma=\PP(G_{j-1}\oplus\mathbb K ke_\sigma)$ has projective dimension over $\mathbb K$ equal to $b_1+\dots+b_{j-1}=\dim_{\mathbb K}G_{j-1}$; the source frame identifies the base $L_\Theta\sigma Q$ with $\PP(F_j/F_{j+1})$, and the datum $[b]$ is the family $([A_j])$. Parts (iii)--(vi) of Theorem~\ref{thm:S} become the two displayed formulas of Theorem~\ref{thm:ext-boundary}(ii).
\end{example}

\begin{remark}
Everything in Theorem~\ref{thm:S} holds for every discrete or non-discrete sequence, every semisimple $G$ with trivial centre and every flag manifold. What the classical statement for $\PSL(d,\C)$ gives in addition is explicit \emph{linear} formulas, because the Schubert varieties are projective subspaces and the graded maps $A_j$ are linear; for a general $X$ the translates of $X_\sigma$ are Schubert varieties of a flag manifold and the formulas are the ones above.
\end{remark}

\section{The proof of the cell lemma}\label{sec:proofcells}
The proof of Lemma~\ref{lem:cells} completes the root calculation behind the destination formula. We keep the same indexing throughout: the double cosets $W_\Theta\backslash W/W_\theta$, the retraction $\rho$ and the decomposition $N_\sigma=M_-M_+$. Although the lemma is classical, giving the proof in this form makes clear which roots contract the source and which expand its destination. We use the real Bruhat and parabolic structure described in
\cite[4.8, 4.10(d), 5.15--5.18 and~5.20]{BorelTits} and
\cite[Chapter~VII]{KnappLG}.
The root-coordinate argument in (B2) is included below and allows
non-reduced restricted root systems.

\begin{description}[leftmargin=2em,nosep]
\item[(B1) BN-pair.] Let $B=MAN$ be the minimal parabolic subgroup ($N$ the unipotent group of the positive restricted roots) and $W=N_K(\ak)/Z_K(\ak)$. Then $G=\bigsqcup_{w\in W}BwB$, and for every simple reflection $s$ and every $w\in W$,
\[
BsB\cdot BwB\subset BswB\ \cup\ BwB,\qquad BwB\cdot BsB\subset BwsB\ \cup\ BwB.
\]
\item[(B2) Root groups.] Call $\Psi\subset\Phi$ \emph{closed} if $\lambda,\mu\in\Psi$ and $\lambda+\mu\in\Phi$ imply $\lambda+\mu\in\Psi$ (take $\lambda=\mu$ to see that $\lambda\in\Psi,2\lambda\in\Phi$ implies $2\lambda\in\Psi$). If moreover $\Psi$ is contained in a positive system, then $\bigoplus_{\Psi}\g_\lambda$ is a nilpotent subalgebra and $U_\Psi=\exp\bigoplus_\Psi\g_\lambda$ is a closed unipotent subgroup. For \emph{any} ordering $\lambda_1,\dots,\lambda_r$ of $\Psi$, the product of exponential slices $\exp\g_{\lambda_1}\times\cdots\times\exp\g_{\lambda_r}\to U_\Psi$ is a diffeomorphism. To see the last statement, choose a positive root system containing $\Psi$ and filter by root height. The BCH formula is triangular: the coordinate in a given height is the corresponding input plus a polynomial in smaller heights. Solve successively by height; nilpotence makes the inverse polynomial and global. This works for any ordering of the slices. In a non-reduced system $\exp\g_\alpha$ need not itself be a subgroup when $2\alpha$ is a root; its brackets are included at the $2\alpha$ level of the filtration.
\item[(B3) Big cell.] For a standard parabolic $Q=P_\theta$, the product $N_\theta^-\times Q\to G$ is a diffeomorphism onto an open subset. Equivalently $n\mapsto nQ$ is a diffeomorphism of $N_\theta^-$ onto an open neighbourhood of the base point of $G/Q$.
\end{description}

Write $W_\Theta\subset W$ for the Weyl group of the Levi factor $L_\Theta$, generated by the reflections in the roots of $\Delta\setminus\Theta$; similarly $W_\theta$. Representatives of $W$ in $N_K(\ak)$ will be denoted by the same letters, and we choose each double-coset representative $\sigma$ of minimal length in $W_\Theta\sigma$. This means $\sigma^{-1}\lambda>0$ for every positive root $\lambda$ in the span of $\Delta\setminus\Theta$.

\paragraph{Step 0: Bruhat decomposition for parabolics.}
Let $I,J\subset\Delta$ and let $P_I,P_J\supset B$ be standard parabolics, with Weyl groups $W_I,W_J$ of their Levi factors. We claim
\[
G=\bigsqcup_{w\in W_I\backslash W/W_J}P_IwP_J.
\]
Since $P_I=\bigsqcup_{u\in W_I}BuB$ (apply (B1) in $L_I$ and use $N_I\subset B$), we have
\[
 P_IwP_J=\bigcup_{u\in W_I,v\in W_J}BuB\,wBvB.
\] Choose a reduced expression $u=s_1\cdots s_a$ with all $s_i$ simple reflections of $L_I$. The product $Bs_1B\cdots Bs_aB$ is bi-$B$-invariant and contains a representative of $u$, so it contains $BuB$. Applying (B1) $a$ times to $Bs_1B\cdots Bs_aBwB$ shows that $BuBwB$ is contained in $\bigcup_{u'\in W_I}Bu'wB$, and the same on the right gives $P_IwP_J\subset\bigcup_{u'\in W_I,v'\in W_J}Bu'wv'B$. The reverse inclusion is obvious. Distinct double cosets of $W_I\backslash W/W_J$ give disjoint unions of $B$-double cosets by (B1), which proves the claim.

\paragraph{Step 1: part (i).}
Let $w_0$ be the longest element of $W$ and $\Theta^*=-w_0\Theta$. Then $w_0$ maps $\Phi^+$ to $\Phi^-$ and $\mathrm{span}(\Delta\setminus\Theta^*)$ to $\mathrm{span}(\Delta\setminus\Theta)$, so $P_\Theta^-=w_0P_{\Theta^*}w_0^{-1}$ and $w_0^{-1}W_\Theta w_0=W_{\Theta^*}$. By Step~0,
\begin{align*}
G&=w_0\Bigl(\bigsqcup_{\tau\in W_{\Theta^*}\backslash W/W_\theta}P_{\Theta^*}\,\tau\,P_\theta\Bigr)\\
&=\bigsqcup_{\sigma\in W_\Theta\backslash W/W_\theta}P_\Theta^-\,\sigma\,P_\theta,\qquad\sigma=w_0\tau.
\end{align*}
Dividing by $Q=P_\theta$ gives $X=\bigsqcup_\sigma O_\sigma$. The set $P_\Theta^-\sigma Q$ depends only on the double coset because $W_\Theta\subset L_\Theta\subset P_\Theta^-$ and $W_\theta\subset Q$.

\paragraph{Step 2: part (ii).}
Since $\sigma$ normalises $\ak$, $A=\sigma A\sigma^{-1}\subset\sigma Q\sigma^{-1}$, so $\sigma Q$ is fixed by $A_\Theta\subset A$; and $L_\Theta$ centralises $A_\Theta$, so every point of $L_\Theta\sigma Q$ is fixed. Conversely, let $x\in X^{A_\Theta}$. By Step~0 applied to the \emph{positive} parabolic $P_\Theta$, $x=nl\sigma Q$ with $n\in N_\Theta$, $l\in L_\Theta$. Choose $z_t\in A_\Theta$ with $\beta(z_t)\to\infty$ for $\beta\in\Theta$ (as $t\to\infty$). Then $z_t^{-1}x=x$, and also $z_t^{-1}x=(z_t^{-1}nz_t)\,l\sigma Q\to l\sigma Q$, because $\mathrm{Ad}(z_t^{-1})$ multiplies $\g_\lambda$ by $e^{-\lambda(\log z_t)}\to0$ for $\lambda\in\Phi(N_\Theta)$. Hence $x=l\sigma Q$, and the fixed set is $\bigcup_\sigma L_\Theta\sigma Q$; the pieces are disjoint because they lie in distinct $P_\Theta$-double cosets.

Compactness of $L_\Theta\sigma Q/Q$: by the Iwasawa decomposition of the reductive group $L_\Theta$, $L_\Theta=K_\Theta\,A\,N_L$ with $K_\Theta=K\cap L_\Theta$ and $N_L=N\cap L_\Theta$. We have $A\subset\sigma Q\sigma^{-1}$, and $N_L\subset\sigma Q\sigma^{-1}$ because for a positive root $\lambda$ of $L_\Theta$ our choice of $\sigma$ gives $\sigma^{-1}\lambda>0$, a root of $Q$. Hence $L_\Theta\sigma Q=K_\Theta\sigma Q$ is the image of the compact group $K_\Theta$.

\paragraph{Step 3: part (iii).}
By Step~1 and $P_\Theta^-=N_\Theta^-L_\Theta$, a point of $O_\sigma$ is $x=nl\sigma Q$ with $n\in N_\Theta^-$, $l\in L_\Theta$. For $z\in A_\Theta$ with $\beta(z)\to\infty$ ($\beta\in\Theta$), every root $\lambda$ of $N_\Theta^-$ has a non-positive expansion $\lambda=-\sum n_\alpha\alpha$ with $n_\beta>0$ for some $\beta\in\Theta$, so $\lambda(\log z)=-\sum_{\beta\in\Theta}n_\beta\beta(\log z)\to-\infty$, and $znz^{-1}\to1$ (uniformly for $n$ in compact sets). Since $z$ commutes with $l$ and fixes $\sigma Q$,
\[
zx=(znz^{-1})\,l\,\sigma Q\longrightarrow l\sigma Q.
\]
The limit is intrinsic to $x$, hence $\rho(x)=l\sigma Q$ is well defined and independent of how $z\to\infty$.

\paragraph{Step 4: part (iv).}
Put $n_\sigma=\mathrm{Ad}(\sigma)\,\mathfrak n_\theta^-=\bigoplus_{\lambda\in\Sigma}\g_\lambda$ with $\Sigma=\sigma\Phi(N_\theta^-)$. The set $\Phi(N_\theta^-)$ is closed (coefficients add), hence so is $\Sigma$. By (B3) translated by $\sigma$, the map $n\mapsto n\sigma Q$ is a diffeomorphism of $N_\sigma=\sigma N_\theta^-\sigma^{-1}$ onto an open neighbourhood of $\sigma Q$, and it is $A$-equivariant for the conjugation action on $N_\sigma$: for $a\in A$, $a\,n\sigma Q=(ana^{-1})\,a\sigma Q=(ana^{-1})\sigma Q$.

Define the $\Theta$-height of a root $\lambda=\sum_\alpha n_\alpha\alpha$ by $h(\lambda)=\sum_{\beta\in\Theta}n_\beta$. It is additive, $h>0$ exactly on $\Phi(N_\Theta)$, $h=0$ on $\Phi(L_\Theta)$, $h<0$ on $\Phi(N_\Theta^-)$. Then $\Sigma_+=\{\lambda\in\Sigma:h>0\}$ and $\Sigma_{\le0}=\{h\le0\}\cap\Sigma$ are closed (the height of a sum is the sum of heights), so (B2) applies to both, and
\[
N_\sigma=\prod_{\Sigma_{\le0}}\exp\g_\lambda\cdot\prod_{\Sigma_+}\exp\g_\lambda=M_-M_+
\]
with unique decomposition, taking in (B2) the ordering that lists $\Sigma_{\le0}$ first.

It remains to show $N_\Theta\sigma Q=M_+\sigma Q$. The set $\Phi(N_\Theta)$ is closed with $\Phi(N_\Theta)\cap-\Phi(N_\Theta)=\varnothing$. Split it as $\Sigma_+\sqcup\Psi$ with $\Psi=\{\lambda\in\Phi(N_\Theta):\sigma^{-1}\lambda\in\Phi(Q)\}$. Indeed, if $\lambda\in\Phi(N_\Theta)$ and $\sigma^{-1}\lambda\notin\Phi(Q)$ then $\sigma^{-1}\lambda\in\Phi(N_\theta^-)$, so $\lambda\in\Sigma\cap\Phi(N_\Theta)=\Sigma_+$. The exponential slices of the root spaces in $\Psi$ lie in $\sigma Q\sigma^{-1}$. By (B2), listing $\Sigma_+$ first, $N_\Theta=M_+\cdot\prod_\Psi\exp\g_\mu$, and the second factor fixes $\sigma Q$. Therefore $N_\Theta\sigma Q=M_+\sigma Q$.

\paragraph{Dimensions.}
We also record the two dimension counts used in Corollary~\ref{cor:conservation} and in Theorem~\ref{thm:dest}(c). The orbit map of $M_+$ at $\sigma Q$ is injective (since $M_+\subset N_\sigma$ and $N_\sigma\cap\sigma Q\sigma^{-1}=1$) and an immersion, so $\dim M_+\sigma Q=\sum_{\Sigma_+}\dim\g_\lambda$. For the cell, $O_\sigma$ is a homogeneous space of $P_\Theta^-$ and its tangent space at $\sigma Q$ is $\mathfrak p_\Theta^-/(\mathfrak p_\Theta^-\cap\mathrm{Ad}(\sigma)\mathfrak q)$, spanned by the root spaces $\g_\lambda$ with $h(\lambda)\le0$ and $\sigma^{-1}\lambda\notin\Phi(Q)$, that is, with $\lambda\in\Sigma_0\cup\Sigma_-$. So
\[
\dim O_\sigma=\sum_{\Sigma_0\cup\Sigma_-}\dim\g_\lambda,\qquad \dim X=\sum_\Sigma\dim\g_\lambda=\dim O_\sigma+\dim M_+\sigma Q.
\]
This is the identity behind the conservation law of Corollary~\ref{cor:conservation}.
\qed

%--------------------------------------------------------------
\medskip
\begin{remark}\label{sec:notsay}
We have followed one divergent sequence from its moving sources to its boundary datum. Before passing to all sequences in a group, four distinctions will be useful.
\begin{enumerate}[label=(\roman*),nosep]
\item \emph{It is a statement about one sequence.} The sets $\D_{(g_m)}(x)$ for different sequences in the same group fit together only after the work of Chapter~\ref{ch:limit}.
\item \emph{It forgets rates.} Two sequences with the same $(\Theta,k,k',b)$ but different speeds have the same relation on every flag manifold, yet they can differ on other homogeneous spaces (for instance on $G/H$ for non-parabolic $H$) and in the way they act on measures.
\item \emph{It is not a statement about a neighbourhood of $\mathcal E$.} The theorem gives the destination of each point; it does not say how the images of a ball around an exceptional point look for finite $m$. The exterior boundary determines the limiting graph. Finite-time images require the actual Cartan scales and frames.
\item \emph{It uses the Cartan decomposition of $G$, not of a subgroup.} If $\Gamma<G$ is discrete, the sequences $(\gamma_m)$ in $\Gamma$ form a subset of the sequences of $G$, and the data $(\Theta,k,k',b)$ of $\gamma_m$ encode how $\Gamma$ sits in $G$; this is the subject of Chapter~\ref{ch:limit}.
\end{enumerate}
\end{remark}

%==========================================================
%==========================================================
\chapter{Limit sets and domains}
\label{ch:limit}

The destination formula describes individual divergent sequences. Taking the union of their relations for a discrete group $\Gamma$ gives the group relation $\D_\Gamma$, which determines recurrence and the domains of proper discontinuity. For $\PSL(2,\C)$ the classical limit set already describes both. In higher rank the individual sequences carry further information. Myrberg's construction tests the size of each sequence's destination fibres, and the definitions of Quint and Conze--Guivarc'h retain limiting flags or projective images. The Kulkarni formula uses $\D_\Gamma$ together with the stabiliser and fixed-source information in $B_\Gamma$.

\medskip
The first result says exactly when a nonempty domain of discontinuity exists and when there is a largest domain. It holds for any action of a countable group by homeomorphisms on a compact metrizable space, with no Lie theory.

\begin{theorem}[Relational criterion]\label{thm:relational}
Let $\Gamma$ be a countable group acting by homeomorphisms on a compact metrizable space $X$, let $\D_\Gamma\subset X\times X$ be its destination relation (Definition~\ref{def:Dgamma}), let $\mathscr R_\Gamma=\{x:(x,x)\in\D_\Gamma\}$ be its recurrent set, and $\mathscr W_\Gamma=X\setminus\mathscr R_\Gamma$.
\begin{enumerate}[label=\textup{(\roman*)},nosep]
\item A $\Gamma$-invariant open set $U$ is a domain of proper discontinuity iff $(U\times U)\cap\D_\Gamma=\varnothing$.
\item $\mathscr W_\Gamma$ is the union of all domains of proper discontinuity, and every domain is contained in a maximal one.
\item A \emph{largest} domain exists iff $(\mathscr W_\Gamma\times\mathscr W_\Gamma)\cap\D_\Gamma=\varnothing$, that is, iff no two non-recurrent points are related by $\D_\Gamma$; it is then $\mathscr W_\Gamma$.
\item If no largest domain exists, there are at least two incomparable maximal domains.
\end{enumerate}
\end{theorem}

\medskip

The proof is a few lines (a diagonal argument and Zorn's lemma). For discrete subgroups acting on flag manifolds, Theorem~A turns the relation $\D_\Gamma$ into a union of explicit nilpotent orbit closures. We may then ask which sources each limit set remembers. Kulkarni starts with stabilisers and fixed sources, then allows sources to move in compact sets outside the resulting set; Myrberg records sources with non-trivial destination; Quint records attracting flags. We work out the rank-one family, where everything reduces to the classical picture, then the whole theory for a discrete group in an arbitrary simple $G$, and one example in $\PSL(3,\C)$ worked out completely, in which the Kulkarni domain is strictly smaller than both of two incomparable maximal domains.

%--------------------------------------------------------------
\section{The destination relation of a group}\label{sec:Dgroup}
\begin{definition}\label{def:Dgamma}
Let $\Gamma$ be a countable group acting on a compact metrizable space $X$ by homeomorphisms. Its \emph{destination relation} is
\[
\D_\Gamma=\{(x,y)\in X\times X:\ \exists\,\gamma_m\in\Gamma\text{ pairwise distinct},\ x_m\to x,\ \gamma_mx_m\to y\}.
\]
The \emph{recurrent set} is $\mathscr R_\Gamma=\{x:(x,x)\in\D_\Gamma\}$. A \emph{domain of proper discontinuity} is a $\Gamma$-invariant open set $U$ such that $\{\gamma:\gamma K\cap K\ne\varnothing\}$ is finite for every compact $K\subset U$.
\end{definition}
Note the difference with the fixed-source relation: the sources $x_m$ move, and they may move while $\gamma_mx_m\to x$ as well. A point is recurrent as soon as it is related to itself along moving sources. For a single sequence the same definition is $\D_{(\gamma_m)}$ of Section~\ref{sec:thmA}, and $\D_\Gamma$ is the union of these over all sequences of distinct elements.

\begin{lemma}\label{lem:Dbasic}
$\D_\Gamma$ is closed, symmetric, and invariant under $\Gamma\times\Gamma$: if $(x,y)\in\D_\Gamma$ and $g,h\in\Gamma$ then $(gx,hy)\in\D_\Gamma$. In particular $\mathscr R_\Gamma$ is closed and $\Gamma$-invariant.
\end{lemma}
\begin{proof}
\emph{Closed.} Let $(x_r,y_r)\in\D_\Gamma$ converge to $(x,y)$, with witnesses $\gamma_{r,m}$, $x_{r,m}\to x_r$, $\gamma_{r,m}x_{r,m}\to y_r$. For each $r$ choose $m(r)$ so large that $\gamma_{r,m(r)}\notin\{\gamma_{s,m(s)}:s<r\}$ and $d(x_{r,m(r)},x_r)<1/r$, $d(\gamma_{r,m(r)}x_{r,m(r)},y_r)<1/r$. Then $x_{r,m(r)}\to x$, $\gamma_{r,m(r)}x_{r,m(r)}\to y$, with pairwise distinct elements. \emph{Symmetric.} Put $y_m=\gamma_mx_m$; then $y_m\to y$ and $\gamma_m^{-1}y_m=x_m\to x$. \emph{Invariant.} Use the witnesses $h\gamma_mg^{-1}$, $gx_m$.
\end{proof}
Two consequences are worth stating at once. First, $\mathscr R_\Gamma$ contains every point with infinite stabiliser (take $x_m=x$ and $\gamma_m\in\mathrm{Stab}(x)$). Second, if $\Gamma'<\Gamma$ then $\D_{\Gamma'}\subset\D_\Gamma$, so $\mathscr R_{\Gamma'}\subset\mathscr R_\Gamma$.

If $H\le\Gamma$ has finite index and $R$ is a set of representatives with $\Gamma=\bigsqcup_{r\in R}rH$, then $\D_\Gamma=\bigcup_{r\in R}(1\times r)\D_H$. Extract a fixed left coset in any sequence witnessing a pair and write its elements as $r h_m$ with $h_m\in H$. Conversely, multiplying witnesses for $\D_H$ on the left by $r$ gives the displayed inclusion into $\D_\Gamma$.

\medskip

Suppose now that $\Gamma$ is a discrete subgroup of $G$ and $X=G/Q$.

\begin{proposition}\label{prop:Dunion}
$\D_\Gamma$ is the union of the destination relations of all
sequences of pairwise distinct elements of $\Gamma$ obtained after
an extraction satisfying the Cartan hypotheses of
Theorem~\ref{thm:dest}. Each of these extracted relations is given
by the formula of that theorem.
\end{proposition}
\begin{proof}
A sequence of pairwise distinct elements of the discrete group
$\Gamma$ leaves every compact subset of $G$. Given a pair in
$\D_\Gamma$, extract Cartan data from its witnessing sequence;
both witnessing limits persist, so the pair belongs to the
destination relation of this extracted sequence. Conversely,
every pair in the relation of such a sequence is witnessed along
that sequence by Theorem~\ref{thm:dest}(a)--(b), and hence belongs
to $\D_\Gamma$.
\end{proof}
So computing $\D_\Gamma$ means collecting the Cartan data $(\Theta,k,k',b)$ of all divergent sequences in $\Gamma$ and taking the union of the corresponding orbit closures. This relation contains \emph{all} the information needed to decide properness. Before reading off the different limit sets, we ask what happens when we forget some of the geometry. The answer only uses compactness, so it applies well beyond flag manifolds.

%--------------------------------------------------------------
\begin{theorem}[Compact factors]\label{r3fac:compact-factor}
Let $f:X\to Y$ be a continuous equivariant surjection between compact
metrizable spaces on which the same countable group $\Gamma$ acts. Then
\[
 \D^Y_\Gamma=(f\times f)(\D^X_\Gamma),\qquad
 \D^Y_\Gamma(y)=f\bigl(\D^X_\Gamma(f^{-1}(y))\bigr),
\]
where a relation on a subset means the union of its fibres. In particular,
\[
 y\in\mathscr R^Y_\Gamma
 \quad\Longleftrightarrow\quad
 \D^X_\Gamma\cap\bigl(f^{-1}(y)\times f^{-1}(y)\bigr)\ne\varnothing.
\]
For every invariant open $U\subset Y$, the action on $U$ is proper if
and only if the action on $f^{-1}(U)$ is proper.
\end{theorem}
\begin{proof}
Projecting witnessing sequences gives one inclusion. Conversely, lift
$y_j\to y$, $\gamma_jy_j\to z$ to arbitrary $x_j\in f^{-1}(y_j)$.
Compactness gives a subsequence with $x_j\to x$ and
$\gamma_jx_j\to w$, where $f(x)=y$ and $f(w)=z$. This proves the
relation formula, hence the fibre and recurrence tests. A pair upstairs
has both endpoints in $f^{-1}(U)$ exactly when its image pair lies in
$U\times U$; Theorem~\ref{thm:relational}(i) proves the last assertion.
\end{proof}
Thus a return downstairs may come from two different points upstairs:
the endpoints over a recurrent point of $Y$ need not coincide in
$X$. Forgetting can create recurrence by bringing the two ends of a
destination pair into the same fibre. Properness, however, passes
exactly between an open set and its full inverse image.

For locally trivial flag projections we can say more. Instead of
choosing both endpoints by compactness, we may prescribe the initial
lift; Theorem~\ref{thm:forget} gives that stronger fibre formula. We
return to it when comparing full and partial flags.

\section{First tests: \texorpdfstring{$\PSL(2,\C)$}{PSL(2,C)} and rank one}\label{sec:rankone}
Let $\Gamma<\PSL(2,\C)$ be non-elementary and discrete, $X=\C\PP^1$, and $\Lambda$ its limit set.

\begin{proposition}\label{prop:D-psl2}
$\D_\Gamma=(X\times\Lambda)\cup(\Lambda\times X)$. Hence $\mathscr R_\Gamma=\Lambda$, and the largest domain of proper discontinuity is $X\setminus\Lambda$.
\end{proposition}
\begin{proof}
By Proposition~\ref{prop:classical}, after extracting Cartan data from a divergent sequence we have $\D(x)=\{x^+\}$ if $x\neq x^-$ and $\D(x^-)=X$, with $x^\pm\in\Lambda$. Any witnessing pair persists on the extracted sequence. This gives $\subseteq$. Conversely let $y\in\Lambda$ and $x\in X$. There is a sequence $\gamma_m$ with attractor $y$ and some repeller $r$. Composing on the right with a fixed $h\in\Gamma$ gives a sequence with attractor $y$ and repeller $h^{-1}r$; since $\Gamma$ is non-elementary, $r$ is not fixed by $\Gamma$, so we can choose $h$ with $h^{-1}r\neq x$. Then $(x,y)\in\D_\Gamma$. This shows $X\times\Lambda\subseteq\D_\Gamma$, and $\Lambda\times X$ follows by symmetry. Finally $(x,x)\in\D_\Gamma$ iff $x\in\Lambda$, and $(X\setminus\Lambda)^2\cap\D_\Gamma=\varnothing$, so Theorem~\ref{thm:relational}(iii) gives the largest domain.
\end{proof}
The proof uses only Corollary~\ref{cor:rank-one}, so it extends verbatim to every real rank-one group.

\begin{proposition}[Rank one]\label{prop:D-rank1}
Let $G$ have real rank one, $X=G/P$ its sphere at infinity, and $\Gamma<G$ an infinite discrete subgroup with no global fixed point in $X$. Let $\Lambda$ be the set of attracting points of sequences of distinct elements. Then $\Lambda$ is closed, $\Gamma$-invariant, equal to the set of repelling points, and
\[
\D_\Gamma=(X\times\Lambda)\cup(\Lambda\times X),\qquad\mathscr R_\Gamma=\Lambda,
\]
so $X\setminus\Lambda$ is the largest domain of proper discontinuity.
\end{proposition}
\begin{proof}
Corollary~\ref{cor:rank-one} gives the same two fibres as Proposition~\ref{prop:D-psl2}. Attracting points form a closed invariant set by diagonal extraction, and inversion identifies them with repelling points. The proof of Proposition~\ref{prop:D-psl2} now applies: given attractor $y$ and repeller $r$, right multiplication by $h\in\Gamma$ replaces $r$ by $h^{-1}r$, which can be chosen different from any prescribed source $x$ because $\Gamma$ has no global fixed point. Symmetry gives the other half of the relation, and Theorem~\ref{thm:relational}(iii) gives the largest domain.
\end{proof}
This is the classical theory in one line: $\Lambda$ is at the same time the recurrent set, the minimal closed invariant set (for non-elementary groups), and the unique minimal ``vertex cover'' of $\D_\Gamma$ (see below). In higher rank these roles separate.

%--------------------------------------------------------------
\section{Recurrence and domains: proof of Theorem~\ref{thm:relational}}\label{sec:criterion}
\begin{proposition}\label{prop:recur-equiv}
For $x\in X$ the following are equivalent: (i) $x\notin\mathscr R_\Gamma$; (ii) $x$ has a neighbourhood $N$ with $\{\gamma:\gamma N\cap N\ne\varnothing\}$ finite; (iii) $x$ lies in a domain of proper discontinuity.
\end{proposition}
\begin{proof}
(i)$\Rightarrow$(ii). If (ii) fails, then for each $r$ the ball $N_r=B(x,1/r)$ satisfies $\gamma N_r\cap N_r\neq\varnothing$ for infinitely many $\gamma$; choose pairwise distinct $\gamma_r$ and $x_r\in N_r$ with $\gamma_rx_r\in N_r$. Then $x_r\to x$ and $\gamma_rx_r\to x$, so $x\in\mathscr R_\Gamma$.

(ii)$\Rightarrow$(iii). Let $U=\Gamma N$, open and invariant. Given a compact $K\subset U$, cover it by finitely many $g_iN$. If $\gamma K\cap K\neq\varnothing$ then $(g_j^{-1}\gamma g_i)N\cap N\neq\varnothing$ for some $i,j$, so $\gamma$ lies in a finite set.

(iii)$\Rightarrow$(i). If $x\in U$ and $x$ were recurrent, a compact neighbourhood $K\subset U$ of $x$ would contain $x_m$ and $\gamma_mx_m$ for all large $m$, giving infinitely many $\gamma$ with $\gamma K\cap K\neq\varnothing$.
\end{proof}

\begin{proof}[Proof of Theorem~\ref{thm:relational}]
(i) If $(x,y)\in(U\times U)\cap\D_\Gamma$, take compact neighbourhoods $K_x,K_y\subset U$; a witnessing sequence gives $\gamma_mK_x\cap K_y\ne\varnothing$ for infinitely many $m$, so the action on $K_x\cup K_y$ is not proper. Conversely, if the action on $U$ is not proper, there are a compact $K\subset U$, pairwise distinct $\gamma_m$ and $x_m\in K$ with $\gamma_mx_m\in K$; extracting, $(x_m,\gamma_mx_m)\to(x,y)\in(U\times U)\cap\D_\Gamma$.

(ii) The first claim is Proposition~\ref{prop:recur-equiv}. For the second, let $U$ be a domain and put $S_U=\mathscr R_\Gamma\cup\overline{\D_\Gamma(U)}$, which is disjoint from $U$ by (i). If $x\notin S_U$, by closedness there is an open $W\ni x$ with $(W\times W)\cap\D_\Gamma=\varnothing$ and $W\cap\overline{\D_\Gamma(U)}=\varnothing$. By invariance $\Gamma W$ has the same properties, so $U\cup\Gamma W$ is again a domain by (i). Hence a domain $U$ is maximal iff $X\setminus U=S_U$. A union of a chain of domains is a domain, since two points lie in one member; Zorn gives the maximal extension.

(iii) A largest domain contains every domain, hence equals $\mathscr W_\Gamma$, which is invariant and open; by (i) it is a domain iff $(\mathscr W_\Gamma^2)\cap\D_\Gamma=\varnothing$.

(iv) Call a closed invariant set $L$ a \emph{vertex cover} if every pair of $\D_\Gamma$ has an endpoint in $L$. By (i), $U\mapsto X\setminus U$ is an inclusion-reversing bijection between domains and vertex covers. By (ii), every domain extends to a maximal one; if this maximal domain were unique, it would contain every domain and hence be the largest. Thus without a largest domain there are at least two maximal domains, necessarily incomparable.
\end{proof}
The proof shows a little more, and we record it because it explains the shape of everything that follows.

\begin{corollary}[Vertex covers]\label{cor:covers}
The recurrent set is the intersection of all closed invariant vertex covers of $\D_\Gamma$, and $\mathscr R_\Gamma$ is itself a vertex cover iff a largest domain exists. If a largest domain exists, its complement is the smallest vertex cover and equals $\mathscr R_\Gamma$.
\end{corollary}
\begin{proof}
Every cover contains $\mathscr R_\Gamma$, since $(x,x)\in\D_\Gamma$ must be covered. If $x\notin\mathscr R_\Gamma$, Proposition~\ref{prop:recur-equiv} puts $x$ in a domain whose complement is a cover missing $x$. So $\mathscr R_\Gamma$ is the intersection. It is a cover iff $\mathscr W_\Gamma$ is a domain, which is part~(iii).
\end{proof}
\begin{remark}
Criterion~(i) and the wandering-point criterion behind~(ii) are the specialisation to groups of known criteria for properness of actions (Kapovich--Leeb--Porti; Kapovich~\cite{KapovichProper}, \emph{Proper actions}; Frances~\cite{Frances} for the relation defined by moving points); the nonuniqueness of maximal domains is already seen in the cyclic diagonal action of Kapovich--Leeb--Porti. The new content of this monograph is that, for flag manifolds, $\D_\Gamma$ is \emph{computed} by Theorem~A. We use the elementary proofs above because nothing else is needed.
\end{remark}

%--------------------------------------------------------------
\section{Kulkarni's set and the exceptional set}\label{sec:KulMyr}
Kulkarni's definitions apply to any action of a countable group $\Gamma$ by homeomorphisms on a compact metrizable space $X$. They use stabilisers, orbits and compact sets.

\begin{definition}
\begin{align*}
L_0&=\overline{\{x\in X:\ |\mathrm{Stab}_\Gamma(x)|=\infty\}},\\
L_1&=\overline{\{y:\ \exists\,x\notin L_0,\ \exists\,\gamma_m\text{ distinct},\ \gamma_mx\to y\}}.
\end{align*}
Put $B_\Gamma=L_0\cup L_1$ (the \emph{fixed-source set}), let $L_2$ be the closure of the cluster points of $\gamma_mK$, $K\subset X\setminus B_\Gamma$ compact and $(\gamma_m)$ distinct, and set
\[
\Lambda_{\Kul}=B_\Gamma\cup L_2,\qquad\Omega_{\Kul}=X\setminus\Lambda_{\Kul}.
\]
\end{definition}

\begin{proposition}\label{prop:L2}\ (a)~$L_2=\overline{\D_\Gamma(X\setminus B_\Gamma)}$, hence
\[
\Lambda_{\Kul}=B_\Gamma\cup\overline{\D_\Gamma(X\setminus B_\Gamma)}:
\]
Kulkarni's set is the fixed-source set together with the closure of the destinations of all sources outside it. (b)~$\Gamma$ acts properly discontinuously on $\Omega_{\Kul}$. (c)~$\mathscr R_\Gamma\subseteq\Lambda_{\Kul}$.
\end{proposition}
\begin{proof}
(a) A generating cluster point is $y=\lim\gamma_mx_m$ with $x_m\in K\subset X\setminus B_\Gamma$; extracting, $x_m\to x\in K$ and $(x,y)\in\D_\Gamma$. Conversely if $(x,y)\in\D_\Gamma$ with $x\notin B_\Gamma$ (a closed set), take a compact neighbourhood $K\subset X\setminus B_\Gamma$ of $x$; the moving sources lie in $K$ for large $m$, so $y$ is a cluster point of $\gamma_mK$. Take closures.
(b) By (a), every destination of a source in $\Omega_{\Kul}\subset X\setminus B_\Gamma$ lies in $\Lambda_{\Kul}$, so Theorem~\ref{thm:relational}(i) gives properness.
(c) If $x\in\mathscr R_\Gamma$ then either $x\in B_\Gamma$, or $x\notin B_\Gamma$ and $x\in\D_\Gamma(x)\subset\overline{\D_\Gamma(X\setminus B_\Gamma)}$.
\end{proof}
For discrete subgroups acting on flag manifolds, Proposition~\ref{prop:Dunion} writes the set inside the closure in~(a) as an explicit union of nilpotent orbit closures from Theorem~A. Note that $\Omega_{\Kul}$ is a domain, but nothing says that it is maximal: in Section~\ref{sec:toral} it is strictly smaller than each of two incomparable maximal domains.

\paragraph{The exceptional set.} The classical convergence principle has a single exceptional source. On a compact metrizable $\Gamma$-space the following definition makes sense for any action by homeomorphisms.
\begin{definition}
For a sequence $(\gamma_m)$ of pairwise distinct elements of $\Gamma$, set
\[
E_{(\gamma_m)}=\{x:\ \D_{(\gamma_m)}(x)\text{ contains at least two distinct points}\}.
\]
Put $\Lambda_{\mathrm{exc}}(\Gamma)=\bigcup E_{(\gamma_m)}$, over all such sequences, and put $\Lambda_{\Myr}=\overline{\Lambda_{\mathrm{exc}}}$ and $\Omega_{\Myr}=X\setminus\Lambda_{\Myr}$.
\end{definition}
The set $\Lambda_{\mathrm{exc}}(\Gamma)$ is $\Gamma$-invariant. For $\delta\in\Gamma$, two destinations of $x$ for $(\gamma_m)$ give two destinations of $\delta x$ for $(\gamma_m\delta^{-1})$, by applying $\delta$ to both witnessing source sequences. Invariance of $\Lambda_{\Myr}$ and $\Omega_{\Myr}$ follows.

For a discrete subgroup of a group $G$ as in Section~\ref{sec:setting-sum} acting on a flag manifold, the union can be taken over sequences with extracted Cartan data because the two witnesses persist under extraction. Theorem~A then gives
\[
 E_{(g_m)}=\{x:\ N_\Theta\text{ does not fix }b\rho(k'x)\}.
\]
For the exceptional-domain results below, let $G$ be simple with trivial centre, $\Gamma<G$ discrete, and $X=G/Q$ a non-trivial flag manifold, so the action of $G$ on $X$ is effective.

\begin{lemma}\label{lem:cellfix}\
Let $\sigma$ be a cell and $\Sigma_\pm(\sigma)$ as in Lemma~\ref{lem:cells}. $N_\Theta$ fixes $\sigma Q$ iff $\Sigma_+(\sigma)=\varnothing$; $N_\Theta^-$ fixes $\sigma Q$ iff $\Sigma_-(\sigma)=\varnothing$; and, as $G$ acts effectively, they are not both empty.
\end{lemma}
\begin{proof}
The orbits $N_\Theta\sigma Q$ and $N_\Theta^-\sigma Q$ are cells of respective dimensions
\[
 \sum_{\Sigma_+}\dim\g_\lambda,
 \qquad \sum_{\Sigma_-}\dim\g_\lambda
\]
by Lemma~\ref{lem:cells}(iv) and its mirror image for $N_\Theta^-$. If both sets were empty, the chart at $\sigma Q$ would lie in\linebreak $L_\Theta\sigma Q\subset X^{A_\Theta}$, so $A_\Theta\ne1$ would fix an open set of the connected analytic manifold $X$ and hence act trivially, contradicting effectiveness.
\end{proof}

\begin{theorem}[The exceptional domain]\label{thm:exc}\
$\Omega_{\Myr}$ is a $\Gamma$-invariant open set on which $\Gamma$ acts properly discontinuously.
\end{theorem}
\begin{proof}
By Theorem~\ref{thm:relational}(i) it suffices to show $(\Omega_{\Myr}\times\Omega_{\Myr})\cap\D_\Gamma=\varnothing$. Suppose $x,y\in\Omega_{\Myr}$ and $(x,y)\in\D_\Gamma$, witnessed by distinct $\gamma_m$; after extraction they have Cartan data $(\Theta,k,k',b)$. Since $x\notin\Lambda_{\mathrm{exc}}$, $\D_{(\gamma_m)}(x)$ is a point, so $N_\Theta$ fixes $p:=b\rho(k'x)$ and $y=kp$. If $\rho(k'x)=l\sigma Q$, then $p=bl\sigma Q$ and, since $bl$ normalises $N_\Theta$, $N_\Theta$ fixes $\sigma Q$: $\Sigma_+(\sigma)=\varnothing$. By Lemma~\ref{lem:cellfix}, $\Sigma_-(\sigma)\neq\varnothing$. Now apply Theorem~\ref{thm:dest}(d) to $(\gamma_m^{-1})$:
\[
\D_{(\gamma_m^{-1})}(y)=\{x':\ y\in\D_{(\gamma_m)}(x')\}\ \supseteq\ k'^{-1}\bigl(N_\Theta^-\,l\sigma Q\bigr),
\]
since every $x'$ with $\rho(k'x')=l\sigma Q$ satisfies $b\rho(k'x')=p$, whose $N_\Theta$-orbit is $\{p\}$. This set has dimension $\sum_{\Sigma_-(\sigma)}\dim\g_\lambda>0$. Take two distinct destinations in this set and extract Cartan data for the inverse sequence, preserving both witnesses. This puts $y$ in $\Lambda_{\mathrm{exc}}$, contradicting $y\in\Omega_{\Myr}$.
\end{proof}
The proof uses the positive-dimensional family of sources with the same singleton destination. That family gives at least two destinations for the inverse sequence. The set $\Lambda_{\Myr}$ is explicit through Theorem~A: it is the closure of the union, over Cartan data of divergent sequences of $\Gamma$, of $k'^{-1}\{x:\ N_\Theta\text{ does not fix }b\rho(x)\}$, a union of Bruhat-type cells. Its relation to the weak-convergence and ``uncertainty'' sets of Kapovich--Leeb--Porti is taken up in Chapter~\ref{ch:balanced}.

\begin{remark}\label{rem:PSL2exc}
For $\PSL(2,\C)$, an extracted sequence with repelling point $x^-$ has $E_{(g_m)}=\{x^-\}$. Thus $\Lambda_{\mathrm{exc}}=\Lambda_{\Myr}=\Lambda$ and Theorem~\ref{thm:exc} is the classical statement.
\end{remark}

\begin{corollary}\label{cor:R-in-both}
$\mathscr R_\Gamma\subseteq\Lambda_{\Kul}\cap\Lambda_{\Myr}$.
\end{corollary}
\begin{proof}
$\Omega_{\Myr}$ is a domain (Theorem~\ref{thm:exc}), hence disjoint from $\mathscr R_\Gamma$ by Proposition~\ref{prop:recur-equiv}; and $\mathscr R_\Gamma\subseteq\Lambda_{\Kul}$ is Proposition~\ref{prop:L2}(c).
\end{proof}
So the recurrent set sits inside both limit sets, and the interesting question in every example is how much bigger each of them is. In Chapter~\ref{ch:balanced} we will see that Theorem~\ref{thm:relational} turns this into a trichotomy.

\subsection*{Projective space over \texorpdfstring{$\mathbb R$ and $\mathbb C$}{R and C}}
The first pseudo-projective limit determines the exceptional set once the full Cartan data have been extracted.

\begin{proposition}[Exceptional set and equicontinuity in projective space]\label{prop:kernel}
Let $V$ have dimension $d\ge2$ over $\mathbb K\in\{\R,\C\}$. For a divergent sequence $(g_m)$ in $\PGL(V)$ with extracted Cartan data and first pseudo-projective limit $\tau^+$, one has $E_{(g_m)}=\operatorname{Ker}\tau^+=\PP(F_2)$. Consequently, for any discrete $\Gamma<\PGL(V)$,
\[
\Lambda_{\mathrm{exc}}(\Gamma)=\Lambda_{\Myr}(\Gamma)
 =\bigcup\{\operatorname{Ker}\tau:\ \tau=\lim_m\gamma_m,\ \gamma_m\in\Gamma\text{ pairwise distinct}\},
\qquad
\Eq(\Gamma)=\Omega_{\Myr}(\Gamma).
\]
Here $\Eq(\Gamma)$ denotes the set of points at which the family $\Gamma$ is equicontinuous.
\end{proposition}
\begin{proof}
The proof of Proposition~\ref{prop:graded} uses singular-value decompositions, so it applies in $\PGL(V)$, allowing the compact factors to be orthogonal over $\R$ and unitary over $\C$. By that proposition, the destination of $[\xi]$ for $\xi\in F_j\setminus F_{j+1}$ has projective dimension $\dim G_{j-1}=b_1+\cdots+b_{j-1}$. This is zero precisely for $j=1$. Hence $E_{(g_m)}=\PP(F_2)$, and Proposition~\ref{prop:SP}(v) identifies $\PP(F_2)$ with $\operatorname{Ker}\tau^+$. A sequence converging to any first limit $\tau$ admits a further extraction with full Cartan data and retains the same $\tau$. This proves the kernel union formula.

The set of limits $\mathcal L_\Gamma\subset\PP(\operatorname{End}V)$ of sequences of distinct elements of $\Gamma$ is closed by diagonal extraction, choosing a new element at each step. It is therefore compact, and discreteness excludes invertible limits. The incidence set
\[
 \{([T],[v])\in\mathcal L_\Gamma\times\PP(V):Tv=0\}
\]
is compact. Its projection to $\PP(V)$ is exactly the kernel union, which is consequently closed.

A point in this union has two different destinations along one extracted sequence. Equicontinuity at the point would force the images of the two witnessing source sequences to approach one another, which is impossible. Conversely, failure of equicontinuity at $x$ gives sources $x_m\to x$ and elements $\gamma_m\in\Gamma$ with the distance between $\gamma_mx_m$ and $\gamma_mx$ bounded below by a positive constant. Continuity of each fixed element lets us extract pairwise distinct $\gamma_m$. By compactness of $\mathrm{SP}(V)$ we may assume $\gamma_m\to\tau\in\mathcal L_\Gamma$. If $x$ is outside the kernel union, both image sequences converge to $\tau x$ by Proposition~\ref{prop:SP}(iv), a contradiction.
\end{proof}
Over $\C$ the same region occurs in complex Kleinian theory~\cite{CNS}. For discrete subgroups of $\PSL(3,\C)$ its pseudo-projective description and comparison with the Kulkarni region are established by Barrera--Cano--Navarrete~\cite{BCNLimit2011}. Cano--Seade~\cite{CanoSeade2010} treat the complex-hyperbolic case.

The inverse sequence gives the general comparison. Outside the exceptional set its destination is a point; a compact-source limit there would therefore create a forbidden fixed-source accumulation point.

\begin{proposition}[Kulkarni and Myrberg limit sets]
\label{prop:kulkarni-myrberg}\label{q:KulMyr}
Let $G$ be a simple real Lie group with trivial centre, let
$X=G/Q$ be a non-trivial flag manifold, and let $\Gamma<G$ be discrete.
Then
\[
 B_\Gamma\subseteq\mathscr R_\Gamma
 \subseteq\Lambda_{\Kul}(\Gamma)
 \subseteq\Lambda_{\Myr}(\Gamma).
\]
In particular $\Omega_{\Myr}(\Gamma)\subseteq\Omega_{\Kul}(\Gamma)$.
\end{proposition}
\begin{proof}
We first show that $B_\Gamma$ is disjoint from every domain of proper
discontinuity $U$. Such a domain contains no point with infinite
stabiliser, and $X\setminus U$ is closed, so $L_0\cap U=\varnothing$.
Suppose a generating point $y=\lim_m\gamma_mx$ of $L_1$ belongs to $U$.
Since $U$ is invariant, $\gamma_mx\in U$ for large $m$ forces $x\in U$.
Choose a compact neighbourhood $C\subset U$ of $y$. For all large $m$,
the compact set $C\cup\{x\}\subset U$ meets its translate by $\gamma_m$,
contradicting properness. Thus the generating points of $L_1$ lie in
the closed set $X\setminus U$, and so does $L_1$ itself.
Proposition~\ref{prop:recur-equiv} now gives
$B_\Gamma\subseteq\mathscr R_\Gamma$.

Put $\Omega=\Omega_{\Myr}(\Gamma)$, which is a domain by
Theorem~\ref{thm:exc}. In particular $B_\Gamma\cap\Omega=\varnothing$.
Suppose a generating cluster point $y$ of $L_2$ belongs to $\Omega$.
There are a compact set $K\subset X\setminus B_\Gamma$, pairwise
distinct $\gamma_m\in\Gamma$, and points $x_m\in K$ such that
$\gamma_mx_m\to y$. First extract so that
$x_m\to x\in K$ and $\gamma_m^{-1}y\to z$. Then extract Cartan data
for $(\gamma_m^{-1})$; both limits persist.
Writing $y_m=\gamma_mx_m$, we have
\[
 y_m\longrightarrow y,\qquad
 \gamma_m^{-1}y_m=x_m\longrightarrow x,
 \qquad
 \gamma_m^{-1}y\longrightarrow z.
\]
Thus both $x$ and $z$ belong to
$\D_{(\gamma_m^{-1})}(y)$, along this same extracted sequence:
the first uses the moving sources $y_m$, the second the fixed source
$y$. Since $y\notin\Lambda_{\Myr}(\Gamma)$, this fibre is a single
point, so $z=x$. As $y\notin L_0$ and the inverses $\gamma_m^{-1}$
are pairwise distinct, the fixed-source limit
$\gamma_m^{-1}y\to x$ puts $x$ in $L_1\subseteq B_\Gamma$,
contradicting $x\in K\subset X\setminus B_\Gamma$.

Every generating cluster point of $L_2$ therefore lies in the closed
set $X\setminus\Omega$, hence so does $L_2$. Together with
$B_\Gamma\cap\Omega=\varnothing$, this proves
$\Lambda_{\Kul}(\Gamma)\subseteq\Lambda_{\Myr}(\Gamma)$.
The inclusion $\mathscr R_\Gamma\subseteq\Lambda_{\Kul}(\Gamma)$ is
Proposition~\ref{prop:L2}(c).
\end{proof}
For complex projective groups, the inclusion
$\Lambda_{\Kul}\subseteq\Lambda_{\Myr}$ is established in
Cano--Navarrete--Seade~\cite{CNS}. The argument above applies to every
flag manifold in the present effective setting.

%--------------------------------------------------------------
\section{Examples in higher rank}\label{sec:hr-examples}
Rank one is too small to see the distinctions of the previous section. We now compute $\D_\Gamma$, $\mathscr R_\Gamma$, $\Lambda_{\Kul}$, $\Lambda_{\Myr}$ and the maximal domains in four families. In the first three, which are cyclic or torus-like, Theorem~A reduces the calculation to the cells and the pairs between them. The fourth is the principal image of a Kleinian group, which is the model for Chapter~\ref{ch:psl}.

\paragraph{A loxodromic cyclic group in $\PSL(3,\C)$.}
Let $g=\diag(a,b,c)$ with $|a|>|b|>|c|$ and $\Gamma=\langle g\rangle$ acting on $\C\PP^2$, with fixed points $e_1,e_2,e_3$ and lines $L_{12}=\PP\langle e_1,e_2\rangle$, $L_{23}$, $L_{13}$. Divergent sequences are, up to a subsequence, of the form $g^m$ or $g^{-m}$ with $m\to+\infty$. For $g^m$ the type is $(1,1,1)$ and the source and target flags are the coordinate flags; the diagonal unitary phases do not change the projective maps on the one-dimensional graded quotients. Proposition~\ref{prop:proj} (or the direct computation $[a^m\varepsilon_1:b^m\varepsilon_2:c^m x_3]$) gives, with sources $x_m\to x$:
\[
\D_{(g^m)}(x)=\begin{cases}\{e_1\}& x\notin L_{23},\\ L_{12}& x\in L_{23}\setminus\{e_3\},\\ \C\PP^2& x=e_3,\end{cases}
\]
and $\D_{(g^{-m})}$ is the transpose relation, with the roles of $e_1,e_3$ and $L_{12},L_{23}$ exchanged.

\begin{proposition}\label{prop:cyc3}
For this $\Gamma$: (a)~$\mathscr R_\Gamma=B_\Gamma=\{e_1,e_2,e_3\}$. (b)~$\Lambda_{\Kul}=\Lambda_{\Myr}=L_{12}\cup L_{23}$, and $\Omega_{\Kul}=\Omega_{\Myr}=\Eq=\{[x_1:x_2:x_3]:x_1x_3\ne0\}$. (c)~The only pairs of $\D_\Gamma$ between non-recurrent points are the pairs $(x,y)$ with $x\in L_{23}\setminus\{e_2,e_3\}$, $y\in L_{12}\setminus\{e_1,e_2\}$, and their transposes. Hence there is no largest domain, and there are exactly two maximal domains,
\[
U_1=\C\PP^2\setminus(L_{12}\cup\{e_3\}),\qquad U_2=\C\PP^2\setminus(L_{23}\cup\{e_1\}),
\]
each strictly larger than $\Omega_{\Kul}=U_1\cap U_2$.
\end{proposition}
\begin{proof}
(a) The displayed fibres give $x\in\D_\Gamma(x)$ exactly at $e_1,e_2,e_3$. These are also the fixed points and all fixed-source limits, so $B_\Gamma=\mathscr R_\Gamma$.

(b) Outside $B_\Gamma$ the destinations are $\{e_1\},\{e_3\}$, $L_{12}$ and $L_{23}$; Proposition~\ref{prop:L2} gives $\Lambda_{\Kul}=L_{12}\cup L_{23}$. The exceptional sets are the kernels $L_{23},L_{12}$ of the first limits $[e_1e_1^*]$, $[e_3e_3^*]$ (Proposition~\ref{prop:kernel}).

(c) After removing recurrent endpoints, the displayed relation and its transpose leave the complete bipartite relation between $E_{23}=L_{23}\setminus\{e_2,e_3\}$ and $E_{12}=L_{12}\setminus\{e_1,e_2\}$. A vertex cover must contain one whole side, together with $\mathscr R_\Gamma$; otherwise two points outside it are related. Closedness therefore gives exactly the minimal covers $L_{12}\cup\{e_3\}$ and $L_{23}\cup\{e_1\}$, whose complements are the two maximal domains by Theorem~\ref{thm:relational}.
\end{proof}
Here the pattern of the toral example is seen in its simplest form: $\mathscr R_\Gamma$ consists of three points, while the limit set fills two lines. Most points of those lines do not return to themselves, but the pairs between them force us to remove one line or the other. Each maximal domain is the complement of a line and one extra point.

\paragraph{The cyclic group in $\PSL(n+1,\C)$.}
Let $g=\diag(a_1,\dots,a_{n+1})$ with $|a_1|>\cdots>|a_{n+1}|$. Torus cells in $\C\PP^n$ are indexed by the non-empty supports $S\subset\{1,\dots,n+1\}$ of coordinates. A source in the cell $S$ has, by Proposition~\ref{prop:proj} and the computation above, $\D_{(g^m)}(x)=\PP\langle e_1,\dots,e_{\min S}\rangle$, and $\D_{(g^{-m})}(x)=\PP\langle e_{\max S},\dots,e_{n+1}\rangle$.
\begin{proposition}\label{prop:cycn}
(a)~$\mathscr R_\Gamma=\{e_1,\dots,e_{n+1}\}$. (b)~$\Lambda_{\Kul}=\Lambda_{\Myr}=\{x_1=0\}\cup\{x_{n+1}=0\}$ (two hyperplanes), with $\Omega_{\Kul}=\{x_1x_{n+1}\ne0\}$. (c)~Relations between non-recurrent cells $S,T$ ($|S|,|T|\ge2$) occur iff $\max T\le\min S$ or $\max S\le\min T$. For $n=2$ there are exactly $2$ maximal domains, and for $n=3$ exactly $3$.
\end{proposition}
\begin{proof}
(a),(b) follow from the support formulas exactly as in Proposition~\ref{prop:cyc3}: recurrence means $|S|=1$, and the exceptional kernels are $\{x_1=0\}$ and $\{x_{n+1}=0\}$. Destinations of non-fixed sources lie in these hyperplanes and attain them when $\max S=2$ and $\min S=n$, respectively; Proposition~\ref{prop:L2} gives the Kulkarni set. (c) The relation is complete bipartite between whole cells. In a cover, every neighbour of a cell not contained in full must be contained in full. Retaining its full cells and the fixed points, then taking closures, therefore still gives a cover. Hence every minimal closed invariant cover is a union of cells closed under taking sub-supports, and minimal covers are minimal solutions of a finite combinatorial problem on the $2^{n+1}-n-2$ cells with $|S|\ge2$. For $n=2$ the cells are $\{1,2\},\{2,3\},\{1,3\},\{1,2,3\}$ and only $\{1,2\}\sim\{2,3\}$ are related. For $n=3$, abbreviate a support by its indices: every edge has an endpoint $12$ or $34$, and $12\sim34$.
If a cover omits $34$, it contains all its neighbours $12,13,23,123$; if it omits $12$, it contains $23,24,34,234$.
Each of these families, together with the fixed points, is closed under taking sub-supports and covers every edge, while a cover containing both $12$ and $34$ contains the cover $\{12,34\}$.
These three covers are minimal under inclusion, so they are exactly the closed unions of cells with supports
\[
\{\{1,2\},\{1,3\},\{2,3\},\{1,2,3\}\},\qquad\{\{1,2\},\{3,4\}\},\qquad\{\{2,3\},\{2,4\},\{3,4\},\{2,3,4\}\}.
\]
\end{proof}
The middle cover corresponds to the maximal domain obtained by discarding the two lines $\PP\langle e_1,e_2\rangle$ and $\PP\langle e_3,e_4\rangle$. None of the three maximal domains contains the complement of a projective hyperplane. Each domain omits all four coordinate points, and no hyperplane contains all four.

\paragraph{The Grassmannian $\mathrm{Gr}_2(\C^4)$.}
Let $g=\diag(a_1,a_2,a_3,a_4)$ with $|a_1|>\cdots>|a_4|$ act on $X=\mathrm{Gr}_2(\C^4)\subset\PP^5$. This is $\Theta=\Delta$ and the stabiliser is not a Borel, so the cells are indexed by $W_\Theta\backslash W/W_\theta$ with six elements and the complex fibre dimensions are $4,3,2,2,1,0$ (Proposition~\ref{prop:schubert}). Put
\[
D^+=\{W:\ W\cap\langle e_3,e_4\rangle\ne0\},\qquad D^-=\{W:\ W\cap\langle e_1,e_2\rangle\ne0\},
\]
the two Schubert divisors (cones over the quadric $x_1x_4=x_2x_3$, singular at $\langle e_3,e_4\rangle$ and $\langle e_1,e_2\rangle$).
\begin{proposition}\label{prop:gr24}
For this cyclic group: (a)~$\mathscr R_\Gamma$ is the set of six coordinate planes. (b)~$\Lambda_{\Myr}=D^+\cup D^-$ and $\Lambda_{\Kul}=D^+\cup D^-$, so $\Omega_{\Kul}=\Omega_{\Myr}=\{W:\ W\text{ transverse to }\langle e_1,e_2\rangle\text{ and }\langle e_3,e_4\rangle\}\cong\mathrm{GL}_2(\C)$, the space of graphs of invertible maps $\langle e_1,e_2\rangle\to\langle e_3,e_4\rangle$. (c)~There is no largest domain.
\end{proposition}
\begin{proof}
(a) Index the coordinate planes by $I=(i_1,i_2)$, $1\le i_1<i_2\le4$, with componentwise Bruhat order. For positive powers, write $C_I^-$ for the opposite source cell through $e_I$ and $X_I^+$ for its Schubert destination. On every Cartan extraction of positive powers, Proposition~\ref{prop:schubert} gives the destination $X_I^+$ for $W\in C_I^-$. The limiting diagonal compact factors preserve these cells and varieties. Hence every subsequence of positive powers has a further extraction whose graphs converge in Hausdorff distance to the same relation, by Theorem~\ref{thm:dest}(d). The whole sequence of graphs converges to that relation, so $\D_{(g^m)}(W)=X_I^+$ for $W\in C_I^-$. The same argument applies to any sequence of exponents tending to $+\infty$. In Pl\"ucker coordinates, $C_I^-$ has $p_I\ne0$ and $p_J=0$ unless $J\ge I$, whereas $X_I^+$ has $p_J=0$ unless $J\le I$. Hence
\[
 C_I^-\cap X_I^+=\{e_I\}.
\]
Thus only coordinate planes can be recurrent along positive powers. Reversing the two Borels proves the same assertion for negative powers. Every sequence of distinct elements of this cyclic group has a subsequence of exponents tending to $+\infty$ or to $-\infty$; a moving-source return persists on that subsequence. Consequently every recurrent plane is coordinate. Conversely, all six coordinate planes are fixed by $g$, whose powers are distinct, so all six are recurrent.

(b) The unique source cell with $0$-dimensional destination consists of the planes transverse to $\langle e_3,e_4\rangle$. Therefore the exceptional set is $D^+$ on a Cartan extraction of positive powers, and $D^-$ on one of negative powers. Every infinite sequence of distinct powers has an extraction in one of these directions, and two distinct destinations persist on extraction. Thus $\Lambda_{\Myr}=D^+\cup D^-$. All coordinate planes have infinite stabiliser. By (a) and Proposition~\ref{prop:kulkarni-myrberg}, $B_\Gamma$ is exactly this set of six planes.

For a non-fixed source, the largest positive-power destination is the unique $3$-dimensional Schubert variety $D^-$. It is the closure of the target cell through $e_{24}$; the corresponding source cell has dimension $1$ and consists of
\[
 W_t=\langle e_2+t e_3,e_4\rangle,\qquad t\in\C.
\]
Its points with $t\ne0$ are non-fixed. All other destinations of non-fixed sources are contained in $D^-$: the Bruhat poset of $\mathrm{Gr}_2(\C^4)$ has one element each of dimensions $0,1,3,4$ and two of dimension $2$, and its unique $4$-dimensional destination comes from the fixed source $e_{34}$. The inverse statement has $D^+$ in place of $D^-$. Proposition~\ref{prop:L2}(a) now gives $\Lambda_{\Kul}\subseteq D^+\cup D^-$. Conversely, $W_1\notin B_\Gamma$ has destination $D^-$, so that same proposition puts $D^-$ in $L_2$; the inverse argument puts $D^+$ in $L_2$. Hence $\Lambda_{\Kul}=D^+\cup D^-$. Their common complement is the set of graphs of invertible maps $\langle e_1,e_2\rangle\to\langle e_3,e_4\rangle$, and properness follows from Proposition~\ref{prop:L2}(b).

(c) Choose explicitly
\[
 x=\langle e_2+e_3,e_4\rangle,\qquad
 y=\langle e_1+e_2,\ e_2+e_3+e_4\rangle.
\]
Neither plane is coordinate, so (a) shows that both are non-recurrent. The point $x=W_1$ lies in the $1$-dimensional source cell above and has $\D_{(g^m)}(x)=D^-$. Since $y$ contains $e_1+e_2$, it meets $\langle e_1,e_2\rangle$ nontrivially, and therefore $y\in D^-$. We have obtained a pair of related non-recurrent points. Theorem~\ref{thm:relational}(iii) excludes a largest domain.
\end{proof}
The domain $\Omega_{\Kul}=\mathrm{GL}_2(\C)$ is the set of planes transverse to both the attracting and the repelling plane, as in the transverse-group literature; the complement is the union of two singular cones, and this domain is not the largest.

\subsection{Principal images of Kleinian groups}\label{sec:principal-pre}
Let $\Gamma_0<\PSL(2,\C)$ be non-elementary discrete with limit set $\Lambda_0$, and $\Gamma=\mathrm{Sym}^2\Gamma_0<\PSL(3,\C)$ its image under the irreducible representation on $\mathrm{Sym}^2\C^2=\C^3$. The Veronese map $\nu:\C\PP^1\to C\subset\C\PP^2$ embeds $\C\PP^1$ as the conic $C$, and $\C\PP^2\setminus C$ is the set of unordered pairs $\{p,q\}\subset\C\PP^1$, $p\ne q$, so that $\C\PP^2=\mathrm{Sym}^2(\C\PP^1)$: a point $z\notin C$ lies on the tangent line $T_{\nu(p)}C$ iff $p$ belongs to the pair of $z$.
\begin{proposition}\label{prop:sym2}
\begin{enumerate}[label=\textup{(\alph*)},nosep]
\item For $g=kak'$ with $a=\diag(t,t^{-1})$ and $t\ge1$, one has the Cartan decomposition $\mathrm{Sym}^2g=\mathrm{Sym}^2k\cdot\diag(t^2,1,t^{-2})\cdot\mathrm{Sym}^2k'$. If $(g_m)$ is a divergent sequence with extracted Cartan data and attracting and repelling limits $p,q$, then $(\mathrm{Sym}^2g_m)$ has block type $(1,1,1)$ with attracting and repelling limits $\nu(p),\nu(q)$.
\item $E_{(\mathrm{Sym}^2g_m)}=T_{\nu(q)}C$, and $\D_{(\mathrm{Sym}^2g_m)}(x)=T_{\nu(p)}C$ for $x\in T_{\nu(q)}C\setminus\nu(q)$.
\item $\Lambda_{\Kul}(\Gamma)=\Lambda_{\Myr}(\Gamma)=\bigcup_{p\in\Lambda_0}T_{\nu(p)}C$, and $\Omega_{\Kul}(\Gamma)=\Omega_{\Myr}(\Gamma)=\mathrm{Sym}^2(\C\PP^1\setminus\Lambda_0)$.
\end{enumerate}
\end{proposition}
\begin{proof}
(a) With the induced Hermitian structure and weighted monomial basis, $\mathrm{Sym}^2$ preserves unitarity and sends $\diag(t,t^{-1})$ to $\diag(t^2,1,t^{-2})$. The highest and lowest weight lines are $\nu(p)$ and $\nu(q)$.
(b) By Proposition~\ref{prop:graded} and (a) the kernel of the first limit is $F_2=\langle e_2,e_3\rangle$ in the adapted basis, which is the tangent line to $C$ at the lowest weight point $e_3=\nu(q)$; and the destination of $x\in F_2\setminus e_3$ is $G_2=\langle e_1,e_2\rangle=T_{\nu(p)}C$.
(c) By (a), the image of a source sequence escaping compact sets has singular-value ratio $\sigma_1/\sigma_3=t^4\to\infty$. The projective symmetric-square map is therefore proper, so $\Gamma$ is discrete. Put $\Omega_0=\C\PP^1\setminus\Lambda_0$. We first claim
$B_\Gamma=L_0=\mathrm{Sym}^2\Lambda_0$.
An infinite stabiliser has a finite-index subgroup fixing every
root of the corresponding pair; classical rank-one convergence
then puts every root in $\Lambda_0$. Conversely, loxodromic
fixed-point pairs are dense in $\Lambda_0^2$, and each such pair
has an infinite stabiliser. To see the density, in any disjoint
neighbourhoods $U,V$ meeting $\Lambda_0$, choose loxodromics $s,h$
with both endpoints in $U,V$, respectively, using
Lemma~\ref{lem:principal-coalescing-fixed-pairs}. For large $n$,
$s^nh^{-n}$ is loxodromic with attracting point tending to $s^+$
and repelling point tending to $h^+$; this gives density off the
diagonal, and the same lemma supplies the diagonal.
Finally, for any fixed-source orbit limit every root other than
the repeller tends to the attractor, while a root equal to the
repeller remains in the closed invariant set $\Lambda_0$.
Thus $L_1\subseteq\mathrm{Sym}^2\Lambda_0$, proving the claim.

The set
$T=\bigcup_{p\in\Lambda_0}T_{\nu(p)}C$
is closed: in pair coordinates it is the image of
$\Lambda_0\times\C\PP^1$ under $(p,r)\mapsto\{p,r\}$.
By (b) and Proposition~\ref{prop:kernel}, $\Lambda_{\Myr}=T$.
If $\Omega_0\ne\varnothing$, for each $p\in\Lambda_0$ choose
an extracted sequence with attractor $p$ and repeller $q$
(Proposition~\ref{prop:D-psl2}), and choose $r\in\Omega_0$.
The mixed source $x=\{q,r\}$ lies outside $B_\Gamma$, and (b)
gives $\D_{(\mathrm{Sym}^2g_m)}(x)=T_{\nu(p)}C$.
Every witnessing sequence of moving sources is eventually in
a compact neighbourhood of $x$ disjoint from $B_\Gamma$.
Proposition~\ref{prop:L2} therefore gives $T\subseteq\Lambda_{\Kul}$.
Conversely, along any extracted sequence, a source outside $B_\Gamma$
has at least one root different from the repeller; that root tends to the attractor
even under moving-source perturbations. Hence all its
destinations lie in $T$. Since $B_\Gamma\subseteq T$ and $T$
is closed, Proposition~\ref{prop:L2} gives
$\Lambda_{\Kul}\subseteq T$.
If $\Lambda_0=\C\PP^1$, the claim already gives $L_0=\C\PP^2$,
so $\Lambda_{\Kul}=T=\C\PP^2$.
In both cases the complement of $T$ consists exactly of pairs
whose roots both lie in $\Omega_0$, namely
$\mathrm{Sym}^2(\C\PP^1\setminus\Lambda_0)$.
\end{proof}
In words: the Kulkarni and Myrberg limit sets of $\mathrm{Sym}^2\Gamma_0$ are the cone of tangent lines to the conic at the limit set, and the domain is the symmetric square of the Kleinian domain. This is the model for Theorem~E of Chapter~\ref{ch:psl}, where osculating flags of the rational normal curve replace the tangent lines for the principal embedding $\PSL(2,\C)\to\PSL(n+1,\C)$. The higher-degree equality of Kulkarni and Myrberg there assumes that the classical limit set is proper.

%--------------------------------------------------------------
\section{A worked example: the toral group}\label{sec:toral}
We compute the objects introduced above for a discrete subgroup of $\PSL(3,\C)$ acting on $\C\PP^2$. Hyperbolic toral groups occur in the four-line analysis of Barrera--Cano--Navarrete~\cite{BCNFour2011,BCNLines2016}; the presentation used here is also given in \cite[Example~3.5]{CCL2026}. For cyclic projective dynamics, see \cite{Navarrete2008,CLU}. Let $M=\begin{pmatrix}2&1\\1&1\end{pmatrix}$ and let $\Gamma_M<\PSL(3,\C)$ be the discrete group of matrices $\begin{pmatrix}M^k&a\\0&1\end{pmatrix}$, $k\in\mathbb Z$, $a\in\mathbb Z^2$. It is discrete because it lies in $\mathrm{SL}(3,\mathbb Z)$.

\paragraph{Coordinates.} The eigenvalues of $M$ are $\lambda>1$ and $\lambda^{-1}$, with eigenvectors of irrational slope. Let $(z,w)$ be affine coordinates in the eigenbasis, so that
\[
h_{k,a}(z,w)=(\lambda^kz+a_+,\ \lambda^{-k}w+a_-),\qquad k\in\mathbb Z,\ a=(a_+,a_-)\in L,
\]
where $L\subset\R^2$ is the image of $\mathbb Z^2$, a lattice preserved by $D=\diag(\lambda,\lambda^{-1})$. Both coordinate projections of $L$ are dense in $\R$ (a rank-two group $\mathbb Z\alpha+\mathbb Z\beta$ with $\alpha/\beta$ irrational). Write $\ell_\infty$ for the line at infinity, $\R^2=\{\operatorname{Im}z=\operatorname{Im}w=0\}$ for the real points, and
\[
U_+=\{\operatorname{Im}z\ne0\},\qquad U_-=\{\operatorname{Im}w\ne0\},\qquad B=\ell_\infty\cup(\R\times\C)\cup(\C\times\R).
\]

\begin{proposition}\label{prop:toral}
For $\Gamma=\Gamma_M$:
\begin{enumerate}[label=\textup{(\alph*)},nosep]
\item $\mathscr R_\Gamma=B_\Gamma=L_0=\R^2\cup\ell_\infty$.
\item $\Lambda_{\Kul}(\Gamma)=\Lambda_{\Myr}(\Gamma)=B$, hence $\Omega_{\Kul}=\Eq=\Omega_{\Myr}=U_+\cap U_-$, whose four components are products of half-planes and so complete Kobayashi hyperbolic.
\item $U_+$ and $U_-$ are two distinct, incomparable maximal domains of proper discontinuity, and there is no largest one. The obstruction is a heteroclinic edge: $p=(0,i)\in U_-$ and $q=(i,0)\in U_+$ are non-recurrent and $(p,q)\in\D_\Gamma$.
\item Each of $U_\pm$ contains affine complex lines, so neither of the two maximal domains is hyperbolic.
\end{enumerate}
\end{proposition}
Part~(b) says that the Kulkarni domain, the equicontinuity region and the exceptional domain all coincide, and are strictly smaller than both maximal domains. Part~(a) says $\mathscr R_\Gamma\subsetneq\Lambda_{\Kul}$.

\begin{proof}
\emph{$U_\pm$ are domains.} Let $K\subset U_+$ be compact; its points have $|\operatorname{Im}z|\in[c,C]$ with $c>0$ and bounded coordinates. If $h_{k,a}K\cap K\ne\varnothing$, then $|\operatorname{Im}z'|=\lambda^k|\operatorname{Im}z|$ forces $\lambda^k\in[c/C,C/c]$, so $k$ ranges in a finite set. The real parts then bound $a_+$ and $a_-$, which leaves finitely many lattice vectors. The argument for $U_-$ is identical.

\emph{(a)} Points of $\ell_\infty$ have the whole translation group in their stabiliser. A real point $r$ with rational coordinates in the original $\mathbb Z^2$-basis also has an infinite stabiliser: for a common denominator $q$ some power $M^j\equiv I\pmod q$, and then $(I-M^j)r\in\mathbb Z^2$, so $x\mapsto M^jx+(I-M^j)r$ belongs to the group and fixes $r$. Such points are dense in $\R^2$, so $\R^2\cup\ell_\infty\subseteq L_0\subseteq\mathscr R_\Gamma$. Conversely every point of $\C^2\setminus\R^2=U_+\cup U_-$ lies in a domain, so it is non-recurrent and has finite stabiliser. Hence $L_0=\mathscr R_\Gamma=\R^2\cup\ell_\infty$. Finally $L_1\subseteq\R^2\cup\ell_\infty$: take $x\in U_+$ and $\gamma_mx$ with $\gamma_m=h_{k_m,a_m}$. If $k_m\to+\infty$, $|\operatorname{Im}z_m'|=\lambda^{k_m}|\operatorname{Im}z|\to\infty$ and the limit is at infinity. If $k_m\to-\infty$, a finite limit has real first coordinate (since $\operatorname{Im}z'\to0$), and the same for the second coordinate: $|\operatorname{Im}w'|=\lambda^{-k_m}|\operatorname{Im}w|\to\infty$ unless $\operatorname{Im}w=0$, in which case it is $0$. If $k_m$ is bounded, $a_m\to\infty$ and the limit is at infinity. The case $x\in U_-$ is symmetric. So $B_\Gamma=\R^2\cup\ell_\infty$.

\emph{(b), $\Lambda_{\Kul}\subseteq B$.} Take a compact $K\subset\C^2\setminus\R^2$ and $x_m\in K$ with $h_{k_m,a_m}x_m\to(u,v)\in\C^2$ finite. If $k_m\to+\infty$, then $\operatorname{Im}z'_m=\lambda^{k_m}\operatorname{Im}z_m$ is bounded, so $\operatorname{Im}z_m\to0$; the limit source lies in $U_-$ (not in $\R^2$), so $\operatorname{Im}w_m$ stays away from $0$ and $\operatorname{Im}w'_m=\lambda^{-k_m}\operatorname{Im}w_m\to0$: $v\in\R$. Thus $(u,v)\in\C\times\R$. If $k_m\to-\infty$, symmetrically $(u,v)\in\R\times\C$. If $k_m$ is bounded the limits lie at infinity. With $B_\Gamma\subset B$ this gives $\Lambda_{\Kul}\subseteq B$.

\emph{(b), $B\subseteq\Lambda_{\Kul}$.} $\ell_\infty\cup\R^2\subseteq B_\Gamma$. Fix $(u,v)\in\C\times\R$. Choose $a_j\in L$ with $a_{-,j}\to v$ (the projection of $L$ is dense) and integers $k_j\to+\infty$ so fast that $|a_{+,j}|\lambda^{-k_j}\to0$. Put $x_j=((u-a_{+,j})\lambda^{-k_j},\,i)$. Then $x_j\to(0,i)\in U_-$, so the $x_j$ eventually lie in a compact subset of $\C^2\setminus\R^2$, and
\[
h_{k_j,a_j}x_j=\bigl(u,\ a_{-,j}+i\lambda^{-k_j}\bigr)\longrightarrow(u,v).
\]
So $\C\times\R\subseteq L_2$, and $\R\times\C\subseteq L_2$ by taking $k_j\to-\infty$.

\emph{(b), $\Lambda_{\Myr}=B$.} By Proposition~\ref{prop:kernel} it suffices to compute the kernels of the pseudo-projective limits. Normalise the matrix of $h_{k_m,a_m}$ by its largest absolute entry. If $k_m$ is bounded, $a_m\to\infty$ and the limit is $\begin{pmatrix}0&0&\alpha\\0&0&\beta\\0&0&0\end{pmatrix}$ with $(\alpha,\beta)\in\R^2\setminus0$, whose kernel is $\ell_\infty$. If $k_m\to+\infty$ the limit has the form $\begin{pmatrix}t&0&b\\0&0&c\\0&0&0\end{pmatrix}$ with $t\ge0$, $b,c\in\R$ and $(t,b,c)\ne(0,0,0)$; its kernel is the line $\{z=-b/t\}$ if $t>0=c$, the point $[0:1:0]\in\ell_\infty$ if $t>0$ and $c\ne0$, and $\ell_\infty$ if $t=0$. So all kernels lie in $B$ for $k_m\to+\infty$, and symmetrically for $k_m\to-\infty$ (lines $w=y$, $y\in\R$). Conversely, given $x\in\R$ choose $a_m\in L$ within bounded distance of $(-x\lambda^m,0)$ (a lattice is cocompact); then
\[
\lambda^{-m}\begin{pmatrix}\lambda^m&0&a_{+,m}\\0&\lambda^{-m}&a_{-,m}\\0&0&1\end{pmatrix}\longrightarrow\begin{pmatrix}1&0&-x\\0&0&0\\0&0&0\end{pmatrix},
\]
whose kernel is the line $z=x$. Taking $m\to-\infty$ gives all lines $w=y$, $y\in\R$, and unbounded translations give $\ell_\infty$. The union of these kernels is $B$.

\emph{(c) Maximality.} Suppose $\Omega\supsetneq U_+$ were an invariant domain. It cannot meet $\R^2\cup\ell_\infty=\mathscr R_\Gamma$, so it contains a point $(x,w)$ with $x\in\R$ and $\beta=\operatorname{Im}w\ne0$. Since $\Omega$ is open, there is a real point $r=(r_+,r_-)$ with rational original coordinates such that $(r_+,r_-+i\beta)\in\Omega$. By (a) some $h\in\Gamma$ fixes $r$ and has linear part $D^k$ with $k>0$. For any $\alpha\in\R\setminus0$ put $z_m=(r_++i\alpha\lambda^{-km},\ r_-+i\beta)$. Then $z_m\to(r_+,r_-+i\beta)\in\Omega$ and
\[
h^mz_m=\bigl(r_++i\alpha,\ r_-+i\beta\lambda^{-km}\bigr)\longrightarrow(r_++i\alpha,r_-)\in U_+\subseteq\Omega.
\]
This is a pair of $\D_\Gamma$ inside $\Omega$, contradicting Theorem~\ref{thm:relational}(i). So $U_+$ is maximal, and by symmetry so is $U_-$. They are incomparable. For the heteroclinic edge take $p_m=(i\lambda^{-m},i)\to p=(0,i)$ and $h_{m,0}$: $h_{m,0}p_m=(i,i\lambda^{-m})\to q=(i,0)$. Both are non-recurrent, being in $U_-$ and $U_+$. By Theorem~\ref{thm:relational}(iii) there is no largest domain.

\emph{(d)} $U_+$ contains the line $\{z=i\}$ and $U_-$ the line $\{w=i\}$, which are complex affine lines.
\end{proof}

\paragraph{Consequences of the example.}
The heteroclinic pair $(p,q)$ leaves $\mathscr R_\Gamma=\R^2\cup\ell_\infty$ uncovered, so no largest domain exists (Corollary~\ref{cor:covers}). The canonical region $\Omega_{\Kul}=\Omega_{\Myr}$ is smaller than both maximal domains, yet has the better geometry: its components are hyperbolic and Stein. The extra wall
\[
\Lambda_{\Kul}\setminus\mathscr R_\Gamma=(\R\times(\C\setminus\R))\cup((\C\setminus\R)\times\R)
\]
consists of non-recurrent endpoints of heteroclinic pairs. The rank-one picture of attractors and repellers has acquired a second layer: such an endpoint can be a source of a moving-source pair for one sequence and a destination for the inverse sequence.

\begin{remark}
Inoue's $S^0$ projective groups give a complementary warning. They act freely and properly on $\mathbb H^+\times\C$ with compact quotient, yet their projective equicontinuity region is empty (Cano--Seade). Thus equicontinuity alone cannot replace the global relation $\D_\Gamma$. 
\end{remark}

%--------------------------------------------------------------
\section{Quint's flag limit set}\label{sec:Quint}
Quint's~\cite{qfQuint2002} limit set is defined by concentration of measure rather than by sequences of points, and it is the limit set used for transverse and Anosov groups. In our language it is the set of attracting flags of sequences whose type contains $\theta$.

Let $\varnothing\ne\theta\subseteq\Delta$, $\mathcal F_\theta=G/P_\theta$, and $\nu_\theta$ the $K$-invariant probability measure on $\mathcal F_\theta$. For each $\alpha\in\theta$ let $(\rho_\alpha,V_\alpha)$ be a Tits representation~\cite{Tits1971} (one-dimensional highest restricted weight space, highest weight a positive multiple of the fundamental weight), with an adapted $K$-invariant norm. Then
\[
\frac{\sigma_2(\rho_\alpha(g))}{\sigma_1(\rho_\alpha(g))}=e^{-\alpha(\kappa(g))},
\]
where $\kappa(g)\in\ak^+$ is the Cartan projection, and when these gaps are positive the top singular lines of the $\rho_\alpha(g)$ define a flag $U_\theta(g)\in\mathcal F_\theta$: if $g=k\,a\,k'$ then $U_\theta(g)=kP_\theta$, the attracting flag of $g$.

\begin{definition}
For a discrete $\Gamma<G$,
\[\Lambda^Q_\theta(\Gamma)=\Bigl\{a\in\mathcal F_\theta:\ (\gamma_n)_*\nu_\theta\to\delta_a\ \text{weakly, some distinct }\gamma_n\in\Gamma\Bigr\}.\]
\end{definition}

\phantomsection\label{par:measure-return}
\paragraph{Concentration gives recurrence.}
Let $\Gamma$ act by homeomorphisms on a compact metrizable space $X$,
and let $\nu$ be a probability measure with full support. If
$(\gamma_m)_*\nu\to\delta_a$ weakly for pairwise distinct
$\gamma_m\in\Gamma$, then $a\in\mathscr R_\Gamma$.
Indeed, for every open neighbourhood $U$ of $a$, full support and
weak convergence give
\[
 \nu(U)>0,\qquad \nu(\gamma_m^{-1}U)\longrightarrow1,
 \qquad
 \nu(U\cap\gamma_m^{-1}U)
 \ge\nu(U)+\nu(\gamma_m^{-1}U)-1>0
\]
for all sufficiently large $m$. Choose a decreasing neighbourhood
basis $U_j$ at $a$, increasing indices $m_j$, and points
$x_j\in U_j\cap\gamma_{m_j}^{-1}U_j$. Then $x_j\to a$ and
$\gamma_{m_j}x_j\to a$, proving recurrence. The $K$-invariant
probability $\nu_\theta$ has full support on $\mathcal F_\theta$;
thus the definition alone gives, for every discrete $\Gamma<G$,
\begin{equation}\label{eq:quint-recurrent}
 \Lambda^Q_\theta(\Gamma)
 \subseteq\mathscr R_\Gamma^{\mathcal F_\theta}.
\end{equation}

\begin{theorem}[Quint's limit set from the destination data]\label{thm:Quint}\
\textup{(a)} $(g_n)_*\nu_\theta\to\delta_a$ iff $\min_{\alpha\in\theta}\alpha(\kappa(g_n))\to\infty$ and $U_\theta(g_n)\to a$. Hence, for every discrete $\Gamma$,
\[
\Lambda^Q_\theta(\Gamma)=\{\lim_n k_nP_\theta:\ \gamma_n=k_na_nk'_n\text{ distinct},\ \theta\subseteq\Theta(\gamma_n)\},
\]
a compact set. \textup{(b)} If every sequence of distinct elements of $\Gamma$ has $\theta\subseteq\Theta$ ($\Gamma$ is $P_\theta$-divergent) and $\theta$ is invariant under opposition, then $\Lambda^Q_\theta(\Gamma)$ is the Cartan flag limit set of Canary--Zhang--Zimmer~\cite{mlCZZ}; no Zariski density, transversality, finite generation or Anosov hypothesis is needed. For $\Gamma$ also transverse and non-elementary it is a perfect compact set on which $\Gamma$ acts minimally and as a convergence group~\cite[author version, Proposition~2.8]{mlCZZ}. \textup{(c)} If $\Gamma$ is Zariski dense over $\R$ in the real adjoint algebraic group, the minimality theorem of Benoist and Guivarc'h (see~\cite[Theorem~7.2]{qfQuint2002} and~\cite{Benoist97,BenoistQuint}) identifies $\Lambda^Q_\theta(\Gamma)$ as the unique minimal closed invariant subset of $\mathcal F_\theta$, also when $\Gamma$ is not $P_\theta$-divergent.
\end{theorem}
\begin{proof}[Proof of (a)]
All normalised subsequential limits $T$ of $\rho_\alpha(g_n)$ have rank one iff the gap $\alpha(\kappa(g_n))\to\infty$. The $K$-orbit of the highest weight line is a connected real-analytic homogeneous manifold whose cone spans $V_\alpha$, so $\PP(\ker T)$ meets it in a proper analytic subset and hence in a set of $\nu$-measure zero. Along the subsequence defining $T$, the images $[\rho_\alpha(g_{n_j})v]$ converge to $[Tv]$ off this set, and dominated convergence gives convergence of the corresponding pushforward measures to $T_*\nu$, a Dirac mass when $T$ has rank one, concentrated at $\mathrm{Im}\,T$. Conversely, if $T_*\nu$ is concentrated on a line, continuity and full support force the whole orbit cone into its inverse image, so the spanning property gives $\operatorname{rank}T=1$. Applying this to every normalised subsequential limit, simultaneously for $\alpha\in\theta$, gives the equivalence, and the displayed description follows from $U_\theta(g_n)=k_nP_\theta$ and the definition of the type. Closedness follows by diagonal extraction in the weak topology, choosing distinct group elements at each step; the flag manifold is compact.
\end{proof}
The identification of flag limit sets in (b) follows from (a); the convergence, minimality and perfectness assertions are supplied by the cited result of Canary--Zhang--Zimmer. In terms of Theorem~A, after extraction of Cartan data, a sequence with $\Theta\supseteq\theta$ has one common destination $kP_\theta$ for $\nu_\theta$-almost every source. If $\Theta\not\supseteq\theta$, at least one normalised Tits limit has rank greater than one, and the spanning argument in the proof of (a) shows that the generic limit map is nonconstant. Such a sequence does not concentrate $\nu_\theta$ at a point. Thus Quint's set records sequences with one common destination for almost every source: the singleton destination must be independent of the generic source.

\paragraph{Comparison in the Tits representation.}
Assume $\varnothing\ne\theta\subseteq\Delta$ is opposition invariant and $\Gamma$ is $P_\theta$-divergent; write $\Lambda_\theta=\Lambda^Q_\theta(\Gamma)$. Let $\mathcal V_\theta$ be the highest-weight component of $\bigotimes_{\alpha\in\theta}V_\alpha$ (complexified), $\rho_\theta$ the resulting representation of $G$, $s_\theta:\mathcal F_\theta\to\PP(\mathcal V_\theta)$ the equivariant Segre embedding. Let $\phi_\theta$ be the coefficient of a fixed highest-weight vector, vanishing on the other weight spaces. For $b=hP^-_\theta$, identifying the opposite flag space with $\mathcal F_\theta$, put $H_b=\PP(\rho_\theta(h)\ker\phi_\theta)$. The opposite parabolic preserves this kernel, so the definition is independent of $h$; its flag trace is $s_\theta(F)\in H_b\iff F$ is not transverse to $b$. Write $\Lambda_{\CoG}$ for the set of images $\operatorname{Im}\tau$ of the rank-one pseudo-projective limits (the Conze--Guivarc'h set).
\begin{proposition}\label{prop:Tits-comparison}\
$\rho_\theta(\Gamma)$ is discrete, and
\begin{align*}
s_\theta\Lambda_\theta&=\Lambda_{\CoG}(\rho_\theta\Gamma),\qquad
\Lambda_{\Myr}(\rho_\theta\Gamma)=\bigcup_{b\in\Lambda_\theta(\Gamma)}H_b,\\
s_\theta^{-1}\Eq(\rho_\theta\Gamma)&=\{F:\ F\text{ transverse to every }b\in\Lambda_\theta(\Gamma)\}.
\end{align*}
\end{proposition}
\begin{proof}
The projective image is discrete and the fibres on $\Gamma$ are finite: otherwise a sequence of distinct elements would have images in a compact subset of the projective linear group, where $\sigma_2/\sigma_1$ is bounded away from zero. This contradicts $P_\theta$-divergence, since the singular values of $\rho_\theta(g)$ satisfy $\sigma_2/\sigma_1=e^{-\min_{\alpha\in\theta}\alpha(\kappa(g))}$, so by $P_\theta$-divergence every normalised limit $T$ of $\rho_\theta(\gamma_n)$ has rank one. Write $\gamma_n=k_na_nk'_n$ with limits $k,k'$. Then $T=\rho_\theta(k)\,\Pi\,\rho_\theta(k')$, where $\Pi$ is the projection onto the highest-weight line $\mathbb Cv_\theta$ along the other weight spaces. Hence $\PP(\operatorname{Im}T)=s_\theta(kP_\theta)$, which gives the first identity. The kernel of $T$ is the hyperplane $\{v:\langle\rho_\theta(k')v,v_\theta\rangle=0\}$. For a flag $F=hP_\theta$ one has $\langle\rho_\theta(k'h)v_\theta,v_\theta\rangle\ne0$ iff $k'hP_\theta$ lies in the open cell $N_\theta^-P_\theta/P_\theta$ of flags transverse to the opposite base flag $P^-_\theta$ (the matrix coefficient of a highest weight vector with itself vanishes off the big cell), i.e.\ iff $F$ is transverse to $b=k'^{-1}P_\theta^{\mathrm{opp}}$, the repelling flag of $(\gamma_n)$. Thus $\PP(\ker T)\cap s_\theta(\mathcal F_\theta)=\{s_\theta F:\ F\text{ not transverse to }b\}$ and directly $\PP(\ker T)=\PP(\rho_\theta(k')^{-1}\ker\phi_\theta)=H_b$. Since $\theta$ is invariant under opposition, repelling flags of sequences in $\Gamma$ are attracting flags of the inverse sequences, so they run over exactly $\Lambda_\theta$. Proposition~\ref{prop:kernel} identifies the kernel union with $\Lambda_{\Myr}$ and its complement with the equicontinuity region. This proves the second and third identities.
\end{proof}
Thus the trace of the equicontinuity region is the locus of flags transverse to every limit flag. This trace is obtained without an antipodality assumption, and it may be empty. The domains constructed by Guichard--Wienhard and Kapovich--Leeb--Porti depend on the chosen thickening and can strictly contain this transverse locus. The Kulkarni set of $\rho_\theta\Gamma$ need not equal the union of hyperplanes. The comparison reads the Quint set, the Conze--Guivarc'h image and the Myrberg wall from data of individual sequences. Taking the union of the destination relations can erase the information distinguishing those sequences. The Kulkarni formula also uses the fixed-source set $B_\Gamma$ through Proposition~\ref{prop:L2}(a).
\begin{remark}
The separate marginal limit sets must not be multiplied: in $\mathrm{SL}(4,\R)$ the discrete diagonal group with logarithmic entries $(m+n,m-n,-m-n,-m+n)$ has limit point, plane and hyperplane families containing $\langle e_1\rangle$, $\langle e_1,e_2\rangle$, $\langle e_1,e_2,e_3\rangle$ separately, but no sequence realises all three at once ($e_1$ needs $m,n\to+\infty$, while $e_4$ as the slowest direction needs $n\to-\infty$). The joint type $\Theta\supseteq\theta$ in Theorem~\ref{thm:Quint}(a) excludes this flag.
\end{remark}

\section{What the geometric supports omit}\label{r2bd:supports}
The four readings of the limit set above describe where sources and
destinations can lie. How much do they tell us about the speed of the
motion? The destination relation determines the open sets on which
the action is proper, but it does not retain rates of divergence.
Even keeping every boundary graph and every fixed-source limit does
not recover those rates. We begin with two diagonal matrices.

\begin{example}[The same graphs with different Jordan directions]
\label{r2bd:rates}
In $\PSL(3,\R)$ put
\[
 a=\diag(e^3,e,e^{-4}),\qquad
 b=\diag(e^4,e,e^{-5}).
\]
The cyclic groups $\langle a\rangle$ and $\langle b\rangle$ have the
same boundary graph families on $\R\PP^2$ and the same fixed-source
limit set, but different pairs of asymptotic Jordan rays.
\end{example}
\begin{proof}
For positive powers the two matrices have the same one-dimensional Cartan blocks, coordinate flags and graded projective maps; Theorem~\ref{thm:ext-boundary} therefore gives the same graph limit. Negative powers give its transpose, and every sequence of distinct powers has a subsequence at one of these two ends. As in Proposition~\ref{prop:cyc3}, their fixed-source limit set is $\{[e_1],[e_2],[e_3]\}$. Their Jordan rays, however, are generated respectively by
\[
 (3,1,-4),\ (4,-1,-3)
 \quad\hbox{and}\quad
 (4,1,-5),\ (5,-1,-4),
\]
two different unordered pairs.
\end{proof}
Both sequences have the same geometric shape, but they travel in
different Jordan directions. The example makes no Zariski-density
assertion: it shows what the unrestricted geometric relation forgets,
and leaves limit-cone theorems under their own hypotheses untouched.

Keeping some spectral information does not remove every ambiguity.
Complementary spectral widths, for instance, cannot distinguish a
representation from its contragredient. Use the standard Euclidean norm on $\R^d$
and their induced exterior norms.
For $g\in\mathrm{SL}_d(\R)$ write $g^\vee=(g^{-1})^t$ and set
\[
 H_k(g)=\log\|\textstyle\bigwedge^kg\|
          +\log\|\textstyle\bigwedge^{d-k}g\|,\qquad
 \ell_k(g)=\sum_{i=1}^k\lambda_i(g)
                 -\sum_{i=d-k+1}^d\lambda_i(g),
\]
where $\lambda_1(g)\ge\cdots\ge\lambda_d(g)$ are the logarithms
of the eigenvalue moduli.
\begin{proposition}\label{r2bd:widths}
For $1\le k<d$, $H_k(g^\vee)=H_k(g)$ and
$\ell_k(g^\vee)=\ell_k(g)$. Nevertheless a Zariski-dense
representation $\rho:F_2\to\mathrm{SL}_5(\R)$ is not conjugate
to $\rho^\vee$ in $\mathrm{GL}_5(\R)$ or $\mathrm{PGL}_5(\R)$.
\end{proposition}
\begin{proof}
If the logarithmic singular values are
$\kappa_1\ge\cdots\ge\kappa_d$, those of $g^\vee$ are
$-\kappa_d\ge\cdots\ge-\kappa_1$. Their sum is zero, hence
$\log\|\bigwedge^kg^\vee\|=\log\|\bigwedge^{d-k}g\|$.
The same reversal of the $\lambda_i$ proves the Jordan identity.
If $C\rho C^{-1}=\rho^\vee$, then
$\rho(\gamma)^tC\rho(\gamma)=C$ for every $\gamma$.
Zariski density extends this polynomial identity to all
$g\in\mathrm{SL}_5(\R)$. A diagonal determinant-one matrix scales
$C_{ij}$ by $t_it_j$. For $i=j$, vary $t_i$ and a second coordinate;
for $i\ne j$, vary $t_i$ and a third coordinate. Thus every entry
vanishes, contradicting invertibility of $C$. A projective conjugacy
lifts to $C\rho(\gamma)C^{-1}=c_\gamma\rho^\vee(\gamma)$;
determinants give $c_\gamma^5=1$, hence $c_\gamma=1$ over $\R$.
\end{proof}
For example $g=\diag(e^4,e^3,e,e^{-2},e^{-6})$ has dominant
eigenvalue $e^4$, whereas $g^\vee$ has dominant eigenvalue $e^6$;
both have $\ell_1=10$ and $\ell_2=15$. Thus the widths remember the
separation of the two ends of the spectrum while losing its signed
information. This is a statement about complementary widths; the
marked geometric destination relations have not been identified.

There is a further distinction when measures enter the picture.
Quint's definition tests concentration by translating a fixed probability;
a stationary measure also depends on a probability law on the group.
For a probability $\mu$ on a countable $\Gamma$, stationarity means
$\nu=\sum_\gamma\mu(\gamma)\gamma_*\nu$. The identity tells us
something about the support. If $\operatorname{supp}\mu$ generates
$\Gamma$ as a semigroup, this identity implies that
$\operatorname{supp}\nu$ is $\Gamma$-invariant: an open
neighbourhood $O$ of $\gamma x$, with
$x\in\operatorname{supp}\nu$ and $\mu(\gamma)>0$, satisfies
$\nu(O)\ge\mu(\gamma)\nu(\gamma^{-1}O)>0$. Identifying this
support with a particular minimal set needs further hypotheses.

The support alone cannot give back the measure. Distinct
probabilities with positive weights on a two-point space are
stationary for the trivial action, although their supports agree.
Thus the geometric identifications in this section leave further
questions about probability and growth: stationary measure uniqueness,
critical exponents and dimension--entropy identities require their
own arguments.

\chapter{Balanced systems}
\label{ch:balanced}

The examples of Chapter~\ref{ch:limit} required us to find every pair in $\D_\Gamma$. To construct a domain, we can often do less: it is enough to remove a set meeting every pair, and the cells tell us how to choose it. Theorem~A tells us that the destination of a source is the closure of an orbit of $N_\Theta$ through a point of a Bruhat cell, and an orbit closure meets only cells below the cell where it started. So one may declare, for each Cartan type $\Theta$, some cells to be \emph{bad as sources} and some to be \emph{bad as targets}, in such a way that every cell is either a bad source or sends only into bad targets. We call such a choice a balanced system.

\medskip
We make this choice separately for each type and then take the closure of the union over the sequences of a group. Theorem~\ref{thm:balanced} says that, for every discrete subgroup $\Gamma$ of every simple $G$ acting on every flag manifold, the resulting set $\Lambda_{\mathfrak B}(\Gamma)$ of bad sources and targets of the sequences of $\Gamma$ is closed and invariant, and its complement is a domain of proper discontinuity. Corollary~\ref{cor:trichotomy} tests whether too much has been removed: the domain is the largest one iff $\Lambda_{\mathfrak B}$ is contained in the recurrent set.

The proof is simple: it is Theorem~A read cell by cell. The interest is in the choices and in what makes their removed points recurrent. We begin with the exceptional set of Chapter~\ref{ch:limit}, which gives a minimal vertex cover of each cell graph. The balanced thickenings of Kapovich--Leeb--Porti and the domains of Guichard--Wienhard give familiar examples supported on one type; in Sections~\ref{sec:compare} and~\ref{sec:KLPfull} we express the thickening construction through the destination fibres for all types and show that the cell graph is the Bruhat order. The same calculation treats sequences with coincident repelling and attracting flags, such as powers of a principal unipotent (Proposition~\ref{prop:unipotent}). In complex hyperbolic geometry, these coincident flags explain why, for every non-elementary discrete $\Gamma<\mathrm{PU}(2,1)$, the complement of the tangent lines at the limit set is the largest domain on $\C\PP^2$; for a lattice it is the ball (Proposition~\ref{prop:PU21}). The same return mechanism underlies the critical level of $\PSL(n+1,\C)$ (Chapter~\ref{ch:psl}).

%--------------------------------------------------------------
\section{Cells and the cell graph}\label{sec:cellgraph}
Fix $\Theta\subset\Delta$ nonempty, a flag manifold $X=G/Q$, $Q=P_\theta$, and the set of cells $\mathcal C_\Theta=W_\Theta\backslash W/W_\theta$ of Lemma~\ref{lem:cells}. For $\sigma\in\mathcal C_\Theta$ (a representative of minimal length) put
\[
O^-_\sigma=P_\Theta^-\,\sigma Q/Q,\qquad O^+_\sigma=P_\Theta\,\sigma Q/Q .
\]
Each family partitions $X$ into locally closed cells (Lemma~\ref{lem:cells}(i) and its mirror image). The closure of $O^+_\sigma$ is a union of cells $O^+_\tau$; write $\tau\le\sigma$ if $O^+_\tau\subset\overline{O^+_\sigma}$. This is the \emph{closure} (Bruhat) \emph{order} on $\mathcal C_\Theta$. Recall that $O^+_\sigma=N_\Theta L_\Theta\sigma Q$ and that Theorem~A says
\[
\D_{(g_m)}(x)=k\,\overline{N_\Theta\,b\,\rho(k'x)},\qquad k'x\in O^-_\sigma,\ \rho(k'x)=l\sigma Q.
\]

\begin{definition}[The cell graph]\label{def:cellgraph}
For $\sigma\in\mathcal C_\Theta$ let
\[
E(\sigma)=\{\tau\in\mathcal C_\Theta:\ O^+_\tau\cap\overline{N_\Theta\,\sigma Q}\ne\varnothing\},
\]
the set of cells in which the destination of a source in the cell $\sigma$ can lie. The \emph{cell graph} $\mathcal G_\Theta$ is the bipartite graph with two copies of $\mathcal C_\Theta$ (``source cells'' and ``target cells'') and an edge $\sigma\to\tau$ iff $\tau\in E(\sigma)$.
\end{definition}

\begin{lemma}\label{lem:cellgraph}
\textup{(i)} $\sigma\in E(\sigma)$ for every $\sigma$ (a loop at every vertex). \textup{(ii)} If $\tau\in E(\sigma)$ then $\tau\le\sigma$. \textup{(iii)} $E(\sigma)=\{\sigma\}$ iff $\Sigma_+(\sigma)=\varnothing$ iff $O^+_\sigma=L_\Theta\sigma Q$ is compact (a \emph{closed cell}). \textup{(iv)} $E(\sigma)$ does not depend on the representative of the double coset.
\end{lemma}
\begin{proof}
(i) $\sigma Q\in N_\Theta\sigma Q\cap O^+_\sigma$. (ii) $N_\Theta\sigma Q\subset P_\Theta\sigma Q=O^+_\sigma$, and the closure of a $P_\Theta$-orbit is $P_\Theta$-invariant, hence a union of cells. (iii) If $\Sigma_+(\sigma)=\varnothing$ then $N_\Theta$ fixes $\sigma Q$ (Lemma~\ref{lem:cells}(iv)), so $N_\Theta\sigma Q=\{\sigma Q\}$ and $E(\sigma)=\{\sigma\}$; moreover $O^+_\sigma=N_\Theta L_\Theta\sigma Q=L_\Theta\sigma Q$, since $N_\Theta$ fixes every point of $L_\Theta\sigma Q$ ($nl\sigma Q=l\cdot(l^{-1}nl)\sigma Q=l\sigma Q$), and this is compact by Step~2 of the proof of Lemma~\ref{lem:cells}. Conversely, if $\Sigma_+(\sigma)\ne\varnothing$, the orbit $N_\Theta\sigma Q=M_+\sigma Q$ is diffeomorphic to $\R^{\dim}$ with positive dimension (Lemma~\ref{lem:cells}(iv)). It is closed in $O^+_\sigma$: it is the fibre over $\sigma Q$ of the retraction $O^+_\sigma\to L_\Theta\sigma Q$, $nl\sigma Q\mapsto l\sigma Q$ (the mirror image of Step~3: $z_t^{-1}nl\sigma Q\to l\sigma Q$ for $z_t\in A_\Theta$ with $\beta(z_t)\to\infty$, $\beta\in\Theta$). Being noncompact it is not closed in the compact space $X$, so its closure contains points outside $O^+_\sigma$, i.e.\ some $\tau\ne\sigma$ lies in $E(\sigma)$. (iv) is clear since $N_\Theta$ and the closure order are intrinsic to the double coset.
\end{proof}
\begin{remark}
The graph is computable in each case from the root data of Lemma~\ref{lem:cells}: $\overline{N_\Theta\sigma Q}$ is the closure of the cell $M_+\sigma Q\cong\R^{\sum_{\Sigma_+}\dim\g_\lambda}$ and its intersections with the cells $O^+_\tau$ are Bruhat intersections. For projective space and for $\PSL(3)$ it is the complete order graph $\{\sigma\to\tau:\tau\le\sigma\}$ (Examples~\ref{ex:PSL3cells} and~\ref{ex:rankone-cells} below); it holds in general (Proposition~\ref{prop:Eorder}).
\end{remark}

%--------------------------------------------------------------
\section{Balanced systems and the theorem}\label{sec:balanced}
Given a divergent sequence $g_m=k_ma_mk'_m$ of type $\Theta$ with Cartan data $(\Theta,k,k',b)$ and subsets $\mathfrak S_\Theta,\mathfrak T_\Theta\subset\mathcal C_\Theta$, put
\[
S_{\mathfrak S}(g_m)=k'^{-1}\!\bigcup_{\sigma\in\mathfrak S_\Theta}O^-_\sigma,\qquad
T_{\mathfrak T}(g_m)=k\bigcup_{\tau\in\mathfrak T_\Theta}O^+_\tau .
\]
These sets depend only on the repelling and attracting parabolic subgroups of the sequence, $k'^{-1}P_\Theta^-k'$ and $kP_\Theta k^{-1}$.

\begin{definition}[Balanced systems]\label{def:balanced}
A \emph{balanced system} is a family $\mathfrak B=(\mathfrak S_\Theta,\mathfrak T_\Theta)_{\varnothing\ne\Theta\subset\Delta}$ such that for every $\Theta$ the pair $(\mathfrak S_\Theta,\mathfrak T_\Theta)$ is a \emph{vertex cover} of $\mathcal G_\Theta$:
\begin{description}[nosep]
\item[(B)] if $\sigma\notin\mathfrak S_\Theta$ and $\tau\in E(\sigma)$, then $\tau\in\mathfrak T_\Theta$.
\end{description}
In words: a cell which is not a \emph{bad source} can be sent only into \emph{bad targets}. For a discrete $\Gamma<G$ put
\[
\Lambda_{\mathfrak B}(\Gamma)=\overline{\bigcup_{(\gamma_m)}\bigl(S(\gamma_m)\cup T(\gamma_m)\bigr)},\qquad\Omega_{\mathfrak B}=X\setminus\Lambda_{\mathfrak B},
\]
the union over all sequences of pairwise distinct elements of $\Gamma$, each taken after extraction of Cartan data.
\end{definition}
Because of the loops (Lemma~\ref{lem:cellgraph}(i)), every balanced system has $\mathfrak S_\Theta\cup\mathfrak T_\Theta=\mathcal C_\Theta$.

\begin{theorem}[Balanced systems give domains]\label{thm:balanced}\
Let $\mathfrak B$ be a balanced system. Then $\Lambda_{\mathfrak B}(\Gamma)$ is closed and $\Gamma$-invariant, every pair of $\D_\Gamma$ has an endpoint in it, and $\Gamma$ acts properly discontinuously on $\Omega_{\mathfrak B}$.
\end{theorem}
\begin{proof}
\emph{Invariance.} Write $P^+=kP_\Theta k^{-1}$ and $P^-=k'^{-1}P_\Theta^-k'$. In the projective model of Theorem~\ref{thm:extG}, let $[\Psi]$ be the boundary limit of $(\gamma_m)$. For fixed $\gamma,\delta\in\Gamma$, the sequence $\gamma\gamma_m\delta^{-1}$ has limit $[\varrho(\gamma)\Psi\varrho(\delta)^{-1}]$, in the same stratum $\mathcal S_\Theta$. Its image and kernel stabilisers are $\gamma P^+\gamma^{-1}$ and $\delta P^-\delta^{-1}$, respectively, by Theorem~\ref{thm:extG}(a),(b). After extracting Cartan data for this sequence, its bad target and source sets are therefore $\gamma T(\gamma_m)$ and $\delta S(\gamma_m)$: these sets depend only on the corresponding parabolic and the fixed cell labels. Their union over all sequences, and then its closure, is $\Gamma$-invariant. \emph{Cover.} Let $(x,y)\in\D_\Gamma$, witnessed by $(\gamma_m)$ with Cartan data $(\Theta,k,k',b)$, and let $k'x\in O^-_\sigma$ with $\rho(k'x)=l\sigma Q$. If $\sigma\in\mathfrak S_\Theta$ then $x\in S(\gamma_m)$. Otherwise, by Theorem~\ref{thm:dest} and because $bl\in L_\Theta$ normalises $N_\Theta$,
\[
y\in k\,\overline{N_\Theta\,bl\sigma Q}=k\,bl\,\overline{N_\Theta\,\sigma Q}.
\]
Each cell $O^+_\tau$ is $P_\Theta$-invariant, hence $bl$-invariant, so every point of $\overline{N_\Theta\sigma Q}$ lies in a cell $O^+_\tau$ with $\tau\in E(\sigma)\subset\mathfrak T_\Theta$ by (B). Hence $y\in T(\gamma_m)$. \emph{Properness.} Theorem~\ref{thm:relational}(i): the complement of a closed invariant vertex cover is a domain.
\end{proof}

\begin{corollary}[Trichotomy]\label{cor:trichotomy}
Let $\mathfrak B$ be balanced, $\Lambda=\Lambda_{\mathfrak B}(\Gamma)$ and $\mathscr R=\mathscr R_\Gamma$. Then $\mathscr R\subseteq\Lambda$ and exactly one of the following holds.
\begin{enumerate}[label=\textup{(\roman*)},nosep]
\item $\Lambda=\mathscr R$: $\Omega_{\mathfrak B}$ is the largest domain, and $\Lambda$ is the smallest vertex cover.
\item $\mathscr R\subsetneq\Lambda$ and $\mathscr R$ is a vertex cover: the largest domain exists and equals $X\setminus\mathscr R\supsetneq\Omega_{\mathfrak B}$.
\item $\mathscr R$ is not a vertex cover: two non-recurrent points are related by $\D_\Gamma$, and there is no largest domain.
\end{enumerate}
In particular $\Omega_{\mathfrak B}$ is the largest domain iff $\Lambda_{\mathfrak B}\subset\mathscr R_\Gamma$.
\end{corollary}
\begin{proof}
Corollary~\ref{cor:covers} and Theorem~\ref{thm:relational}(ii)--(iii), since $\Lambda$ is a closed invariant vertex cover.
\end{proof}
So a balanced system produces a domain, and whether it is the largest one is decided by the equality $\Lambda_{\mathfrak B}=\mathscr R_\Gamma$, a recurrence statement. We already know that every recurrent point must be removed. What remains is to make each removed point return. The smaller the balanced system, the less there is to prove.

\begin{remark}[Vertex covers of the cell graph]\label{rem:covers}
Balanced systems for a given $\Theta$ are exactly the vertex covers of the finite bipartite graph $\mathcal G_\Theta$, partially ordered by inclusion. Their minimal elements are the minimal vertex covers. For the order graph of Proposition~\ref{prop:Eorder}, every minimal cover is a partition: if $v\in\mathfrak S\cap\mathfrak T$, minimality gives $w\notin\mathfrak T$ with $w\le v$ and $u\notin\mathfrak S$ with $v\le u$, and transitivity leaves $u\to w$ uncovered. Thus every minimal cover has $|\mathcal C_\Theta|$ vertices, the minimum certified by the diagonal matching and K\"onig's theorem; equivalently, they are the sublevel systems of Proposition~\ref{prop:sublevel}. Different minimal covers can give the same wall after taking the union over sequences (Example~\ref{ex:PSL3cells}); in general one still has to choose; in the sequel the preferred choice is the \emph{exceptional} system, which is the cover generated by the compact cells (Section~\ref{sec:sublevel}).
\end{remark}

%--------------------------------------------------------------
\section{Sublevel systems and the exceptional system}\label{sec:sublevel}
A subset $\mathcal J\subset\mathcal C_\Theta$ is a \emph{sublevel set} (or down-set) if $\sigma\in\mathcal J$ and $\tau\le\sigma$ imply $\tau\in\mathcal J$; these are the ``unions of sublevels of the Bruhat order'' of the literature.

\begin{proposition}[Sublevel systems]\label{prop:sublevel}
For every family $\mathcal J=(\mathcal J_\Theta)$ of sublevel sets, $\mathfrak B(\mathcal J)=(\mathfrak S_\Theta,\mathfrak T_\Theta)=(\mathcal C_\Theta\setminus\mathcal J_\Theta,\ \mathcal J_\Theta)$ is a balanced system. For it $\mathfrak S_\Theta\sqcup\mathfrak T_\Theta=\mathcal C_\Theta$ (a partition).
\end{proposition}
\begin{proof}
If $\sigma\notin\mathfrak S_\Theta$ then $\sigma\in\mathcal J_\Theta$, and $\tau\in E(\sigma)$ gives $\tau\le\sigma$ by Lemma~\ref{lem:cellgraph}(ii), so $\tau\in\mathcal J_\Theta$.
\end{proof}
\begin{example}[The exceptional system]\label{ex:excsys}
Let $\mathcal J^{\mathrm{exc}}_\Theta$ be the set of closed cells ($\Sigma_+(\sigma)=\varnothing$). It is a sublevel set: a closed cell is minimal for the closure order, since its closure is itself. The corresponding system is $\mathfrak S^{\mathrm{exc}}_\Theta=\{\sigma:\Sigma_+(\sigma)\ne\varnothing\}$ (sources $x$ for which $N_\Theta$ moves $b\rho(k'x)$: the exceptional sources of Chapter~\ref{ch:limit}), $\mathfrak T^{\mathrm{exc}}_\Theta=\{\sigma:\Sigma_+(\sigma)=\varnothing\}$ (targets fixed by $N_\Theta$). When $G$ acts effectively on $X$, each bad target $kO^+_\tau$ with $\tau$ closed is exceptional for the inverse sequence. Indeed, $O^+_\tau$ is fixed pointwise by $N_\Theta$, and the argument of Theorem~\ref{thm:exc} gives more than one inverse destination. The bad-source union has closure $\Lambda_{\Myr}$, so $\Lambda_{\mathfrak B}=\Lambda_{\Myr}$ and Theorem~\ref{thm:exc} is the case $\mathfrak B^{\mathrm{exc}}$ of Theorem~\ref{thm:balanced}.
\end{example}
\begin{proposition}[Partition systems are minimal covers]\label{prop:excmin}
If $\mathfrak S_\Theta\sqcup\mathfrak T_\Theta=\mathcal C_\Theta$ is a partition which is a vertex cover of $\mathcal G_\Theta$ (for instance any sublevel system), then it is a \emph{minimal} vertex cover: no proper subpair is a cover. In particular $\mathfrak B^{\mathrm{exc}}$ is minimal.
\end{proposition}
\begin{proof}
Each loop $\sigma\to\sigma$ must be covered, so $\sigma\in\mathfrak S$ or $\sigma\in\mathfrak T$. In a partition $\sigma$ lies in exactly one of them, so removing it uncovers its loop.
\end{proof}

\begin{example}[Rank one and $\PSL(2,\C)$]\label{ex:rankone-cells}
$\Theta=\Delta=\{\alpha\}$, $X=\C\PP^1$, $\mathcal C_\Theta=\{1,2\}$ with standard representatives $p_0=[e_1]$ and $q_0=[e_2]$. The group $N_\Theta$ fixes $p_0$ and its orbit through $q_0$ is $X\setminus\{p_0\}$, so $E(1)=\{1\}$ and $E(2)=\{1,2\}$. The minimal covers are $(\{2\},\{1\})=\mathfrak B^{\mathrm{exc}}$, $(\varnothing,\{1,2\})$ and $(\{1,2\},\varnothing)$. For a sequence with limiting factors $k,k'$, the exceptional cover removes the repelling source $q=k'^{-1}q_0$ and the attracting target $p=kp_0$, which may coincide. For $\mathfrak B^{\mathrm{exc}}$, $\Lambda_{\mathfrak B}=\Lambda$ is the classical limit set, and by Proposition~\ref{prop:D-rank1} it equals $\mathscr R_\Gamma$ in the rank-one setting of that proposition.
\end{example}

\begin{example}[$\PSL(3,\C)$ on $\C\PP^2$: all minimal covers]\label{ex:PSL3cells}
There are three types.
\begin{itemize}[nosep]
\item $\Theta=\Delta$ (blocks $(1,1,1)$, attracting flag $p\in L$, repelling flag $q\in H$). The cells are $\mathcal C_\Delta=\{1,2,3\}$ with $kO^+_1=\{p\}$, $kO^+_2=L\setminus\{p\}$, $kO^+_3=\C\PP^2\setminus L$, and $k'^{-1}O^-_1=\C\PP^2\setminus H$, $k'^{-1}O^-_2=H\setminus\{q\}$, $k'^{-1}O^-_3=\{q\}$. In standard coordinates the orbit closures are $\overline{N[e_1]}=\{[e_1]\}$, $\overline{N[e_2]}=\PP\langle e_1,e_2\rangle$ and $\overline{N[e_3]}=\C\PP^2$. The target frame $k$ takes $[e_1]$ to $p$ and $\PP\langle e_1,e_2\rangle$ to $L$, so $E(1)=\{1\}$, $E(2)=\{1,2\}$, $E(3)=\{1,2,3\}$: the graph is the complete order graph of the chain $1<2<3$. Its minimal covers are exactly the four chain systems $T=\{1,\dots,t\}$, $S=\{t+1,\dots,3\}$, $t=0,1,2,3$. For $t=1$, $\mathfrak B^{\mathrm{exc}}$: bad sources $H=\C\PP^2\setminus(k'^{-1}O^-_1)$ (the kernel line), bad target $\{p\}$. For $t=2$, the \emph{mirror} system: bad source $\{q\}$, bad target $L$. The two are exchanged by inversion; since $\Gamma=\Gamma^{-1}$ they give the same $\Lambda(\Gamma)$.
\item $\Theta=\{\alpha_1\}$ (blocks $(1,2)$: attractor point $p$, kernel line $H$). $\mathcal C=\{1,2\}$, $E(1)=\{1\}$, $E(2)=\{1,2\}$. Minimal covers: $(\{2\},\{1\})$ with bad source $H$ and bad target $\{p\}$ (the exceptional system), and the two trivial ones.
\item $\Theta=\{\alpha_2\}$ (blocks $(2,1)$: attractor line $L$, kernel point $q$). Likewise: bad source $\{q\}$, bad target $L$ (exceptional), and the trivial ones.
\end{itemize}
In all three, the nontrivial minimal covers are the exceptional one and (for $\Theta=\Delta$) its mirror.
\end{example}
\begin{remark}[A larger cover]
For $\Theta=\Delta$ one can also choose bad source $H$ and bad target $L$. This is a balanced system. Once $H$ is a bad source, one can discard the target cell $L\setminus\{p\}$. The target $p$ is still needed to cover the edge from the open source cell.
\end{remark}

\begin{remark}[The graph is always the order graph]\label{rem:Eorder}
The small cases suggest that $E(\sigma)=\{\tau:\tau\le\sigma\}$ in general: every cell below the source cell is reached. Proposition~\ref{prop:Eorder} below proves this for every type.
\end{remark}

%--------------------------------------------------------------
\section{The trichotomy at work}\label{sec:trichotomy-examples}
\begin{enumerate}[label=\textup{(\arabic*)},nosep]
\item \emph{Rank one} (Proposition~\ref{prop:D-rank1}): case (i), with $\Lambda_{\mathfrak B^{\mathrm{exc}}}=\mathscr R_\Gamma=\Lambda$.
\item \emph{A loxodromic cyclic group in $\PSL(3,\C)$} (Proposition~\ref{prop:cyc3}): $\mathscr R=\{e_1,e_2,e_3\}\subsetneq\Lambda=L_{12}\cup L_{23}$, and two non-recurrent points are related: case (iii), two maximal domains.
\item \emph{The toral group} (Proposition~\ref{prop:toral}): case (iii).
\item \emph{Discrete subgroups of $\mathrm{PU}(2,1)$}: case (i) for every non-elementary $\Gamma$, case (iii) for a loxodromic cyclic group (Proposition~\ref{prop:PU21}).
\end{enumerate}

\subsection{Discrete subgroups of \texorpdfstring{$\mathrm{PU}(2,1)$}{PU(2,1)}}\label{sec:PU21}
Let $h$ be a Hermitian form on $\C^3$ with two positive directions and one negative direction, $B=\{[z]:h(z,z)<0\}$ the complex hyperbolic ball, $S=\partial B$, and $G=\mathrm{PU}(2,1)<\PSL(3,\C)$. For $p\in S$, the \emph{complex tangent line} is $T_p=\PP(p^\perp)$ ($p^\perp$ the $h$-orthogonal complement); $T_p\cap\overline B=\{p\}$. Let $\Gamma<G$ be discrete with Chen--Greenberg limit set $\Lambda_0\subset S$.

\begin{proposition}[$\mathrm{PU}(2,1)$]\label{prop:PU21}
\begin{enumerate}[label=\textup{(\alph*)},nosep]
\item The set $\mathcal P$ of pairs (attractor, repeller) obtained by Cartan extraction from sequences of distinct elements of $\Gamma$ is contained in $\Lambda_0\times\Lambda_0$, and equals it if $\Gamma$ is non-elementary.
\item $\Lambda_{\mathrm{exc}}(\Gamma)=\Lambda_{\Myr}(\Gamma)=\bigcup_{q\in\Lambda_0}T_q$ if $\Gamma$ is non-elementary, so $\Omega_{\Myr}=\C\PP^2\setminus\bigcup_{q\in\Lambda_0}T_q$; when $\Lambda_0=S$ (e.g.\ for a lattice) this is the ball $B$.
\item If $\Gamma$ is non-elementary, then $\mathscr R_\Gamma=\Lambda_{\mathrm{exc}}$: every point of every tangent line $T_q$, $q\in\Lambda_0$, is recurrent. Hence Corollary~\ref{cor:trichotomy}\textup{(i)} holds and $\Omega_{\Myr}=\mathrm{Eq}(\Gamma)$ is the largest domain of proper discontinuity of $\Gamma$ on $\C\PP^2$.
\item If $\Gamma=\langle g\rangle$ with $g$ loxodromic with fixed points $p,q$, then $\mathscr R_\Gamma=\{p,q,e_2\}$ ($e_2=T_p\cap T_q$), $\Lambda_{\mathrm{exc}}=T_p\cup T_q$, and case \textup{(iii)} holds (Proposition~\ref{prop:cyc3}).
\end{enumerate}
\end{proposition}
\begin{proof}
Choose coordinates with $h(z)=2\operatorname{Re}(z_1\bar z_3)+|z_2|^2$ and the Cartan involution of $G$ standard. Then every element has a Cartan decomposition $g=kak'$ with
\[
 K=\{[u]:u\in\mathrm{SU}(h),\ u^*u=I\}\subset\mathrm{PU}(3),
 \qquad k,k'\in K,\qquad a=\diag(e^t,1,e^{-t}).
\]
The subgroup before projectivisation is unitarily conjugate to $\mathrm{S}(\mathrm U(2)\times\mathrm U(1))$. This is also a Cartan decomposition in $\PSL(3,\C)$. For a divergent sequence, $t_m\to+\infty$, so the block type is $(1,1,1)$. After extraction let $k,k'$ be the limiting compact factors and put $p=ke_1$, $q=k'^{-1}e_3$, both in $S$. Since $k,k'$ preserve $h$ and the standard norm, $G_2=k\langle e_1,e_2\rangle=p^\perp$ and $F_2=k'^{-1}\langle e_2,e_3\rangle=q^\perp$, so $\PP(G_2)=T_p$, $\PP(F_2)=T_q$. By Proposition~\ref{prop:graded},
\[
\D_{(g_m)}(x)=\begin{cases}\{p\}&x\notin T_q,\\ T_p&x\in T_q\setminus\{q\},\\ \C\PP^2&x=q,\end{cases}\tag{$*$}
\]
and this holds also when $p=q$.

(a) Attractors and repellers are accumulation points of orbits, so $\mathcal P\subset\Lambda_0\times\Lambda_0$. The set $\mathcal P$ is closed (diagonal extraction and compactness of $K$) and contains the pair of fixed points $(g^+,g^-)$ of every loxodromic $g\in\Gamma$. For non-elementary $\Gamma$ these pairs are dense in $\Lambda_0\times\Lambda_0$ (standard for non-elementary discrete groups of isometries of $\C\mathbb H^2$; Chen--Greenberg, Cano--Seade \S4), including near the diagonal; hence $\mathcal P=\Lambda_0\times\Lambda_0$, the diagonal pairs $(p,p)$ being realised by sequences with attractor and repeller both equal to $p$.

(b) By $(*)$ the exceptional set of a sequence is $T_q$; by (a) $q$ ranges over $\Lambda_0$. The union of the $T_q$ is closed (a continuous family over the compact set $\Lambda_0$). A point lies in $B$ iff it lies on none of the tangent lines $T_q$ with $q\in S$. This proves the assertion when $\Lambda_0=S$. For $n=2$ this is the equicontinuity theorem of Cano--Seade~\cite{CanoSeade2010} and the description of the Kulkarni region by Navarrete~\cite{Navarrete2006}. For all $n$, Cano--Liu--L\'opez prove the tangent-hyperplane formula for $\Lambda_{\Kul}$ and, for non-elementary $\Gamma<\mathrm{PU}(1,n)$, that $\Omega_{\Kul}=\Eq(\Gamma)$ is the largest domain \cite[arXiv:1506.08113v5, Theorems~0.1 and~5.1]{CLL}. Their proof uses moving-source destinations and control groups; the calculation here places the $n=2$ case in the general flag-manifold relation.

(c) Let $x\in T_q$, $q\in\Lambda_0$. If $x=q$, use a sequence with attractor $q$ and repeller $q'\ne q$ (it exists by (a), $\Lambda_0$ being infinite): $q\notin T_{q'}$, so $\D(q)=\{q\}$ and $(q,q)\in\D_\Gamma$. If $x\in T_q\setminus\{q\}$, use a sequence with attractor and repeller both equal to $q$: $(*)$ gives $\D(x)=T_q\ni x$. Hence $\Lambda_{\mathrm{exc}}\subset\mathscr R_\Gamma$, and the reverse inclusion is Corollary~\ref{cor:R-in-both}. Corollary~\ref{cor:trichotomy}(i) gives the rest.

(d) Here $\mathcal P=\{(p,q),(q,p)\}$, so only $T_q$ and $T_p$ are exceptional, and the second case of $(*)$ never gives $x\in\D(x)$ for $x\ne e_2$ off $\{p,q\}$; the point $e_2$ lies in $T_p\cap T_q$ and $\D(e_2)=T_p\ni e_2$.
\end{proof}
The cyclic loxodromic example therefore differs from the non-elementary case: as soon as the group is non-elementary, the \emph{parabolic-like} sequences with attractor equal to repeller (which exist by closure) make every exceptional point recurrent. The same mechanism reappears in Chapter~\ref{ch:psl}, for $\mathrm{PU}(1,s)$ and for suitable subgroups of $\PSL(n+1,\C)$, where the middle block may carry a non-trivial residual action.

%--------------------------------------------------------------
%--------------------------------------------------------------
\section{The flag of limit sets of a general group}\label{sec:flagLSgen}
The exceptional system took the closed cells as bad targets. We can also gather cells in increasing order of dimension and ask how the familiar limit sets sit among them. For $\PSL(n+1,\C)$ this gives the chain ordered by the dimensions of the target spaces of the Cartan flags (Theorem~\ref{thm:flag-of-limit-sets}). The same question makes sense on any flag manifold: we use the cells of Theorem~\ref{thm:S} in place of projective subspaces. Throughout this subsection, $G$ is simple with trivial centre, $X=G/P_\theta$ is nontrivial, and $\Gamma<G$ is discrete.

\paragraph{Cartan levels.} In this subsection $\dim$ denotes real dimension. For a divergent sequence with Cartan data $(\Theta,k,k',b)$ and an integer $d\ge0$ let
\[
T_{\le d}(g_m)=k\bigcup\{O^+_\tau:\ \tau\in\mathcal C_\Theta,\ \dim O^+_\tau\le d\},\qquad
\Lambda_{\mathrm Z}^{(d)}(\Gamma)=\overline{\bigcup_{(\gamma_m)}T_{\le d}(\gamma_m)},
\]
the union over sequences of pairwise distinct elements of $\Gamma$, each after extraction. The \emph{Cartan level} $\Lambda_{\mathrm Z}^{(d)}$ is thus obtained by taking the closure of the points in target cells of dimension $\le d$ of divergent sequences of $\Gamma$. It is closed and increasing in $d$. For $\gamma\in\Gamma$, the equivariance of attracting parabolics established in the proof of Theorem~\ref{thm:balanced} carries each target cell of $(\gamma_m)$ to the corresponding target cell of $(\gamma\gamma_m)$, after extraction. Cell dimensions are preserved, so $T_{\le d}(\gamma\gamma_m)=\gamma T_{\le d}(\gamma_m)$ and $\Lambda_{\mathrm Z}^{(d)}$ is $\Gamma$-invariant.

\paragraph{Dimension systems.} The levels order target cells, but a domain also requires a choice of bad sources. For a function $\delta:\Theta\mapsto\delta_\Theta\in\mathbb Z_{\ge0}$ put $\mathcal J_\Theta=\{\tau\in\mathcal C_\Theta:\dim O^+_\tau\le\delta_\Theta\}$. It is a down-set of the Bruhat order (the closure of $O^+_\sigma$ contains $O^+_\tau$ for $\tau\le\sigma$, and a smaller orbit has smaller dimension), so by Proposition~\ref{prop:sublevel} the pair $\mathfrak B^\delta=(\mathcal C_\Theta\setminus\mathcal J_\Theta,\ \mathcal J_\Theta)_\Theta$ is a balanced system, the \emph{dimension system} of $\delta$; write $\Lambda_\delta=\Lambda_{\mathfrak B^\delta}$, $\Omega_\delta=X\setminus\Lambda_\delta$.

One natural cut lies near the middle. The \emph{middle-dimension} system takes $\delta_\Theta$ to be the smallest real dimension of a cell $O^+_\tau$, $\tau\in\mathcal C_\Theta$, at least $\lfloor\frac12\dim_{\mathbb R}X\rfloor$ for a real flag manifold. For a complex flag manifold we use the complex midpoint: $\delta_\Theta$ is twice the smallest complex cell dimension at least $\lfloor\frac12\dim_{\mathbb C}X\rfloor$. The bad targets run up to this first available level in $X$, and the bad sources are the remaining cells. Write $\Lambda_{\mathrm{mid}}$ for this dimension-system wall. This construction uses every Cartan type and a dimension cut in each cell space. The projective critical wall of Section~\ref{sec:critsys} makes a different choice: it uses admissible orientations and keeps their complementary levels paired. We keep the two walls distinct; no equality between them is assumed.

\begin{proposition}[The two ends of the chain]\label{prop:ends}\
\textup{(a)} $\Lambda_{\mathrm Z}^{(0)}(\Gamma)=\Lambda^Q_\theta(\Gamma)$, Quint's flag limit set.
\textup{(b)} $\Lambda_{\mathrm Z}^{(\dim X-1)}(\Gamma)=\Lambda_{\Myr}(\Gamma)$.
\end{proposition}
\begin{proof}
(a) A zero-dimensional cell $O^+_\tau$ is a point: it is a union of $B$-orbits and the only zero-dimensional $B$-orbit in $G/Q$ is $\{Q\}$, by the Bruhat dimension formula. A cell $O^+_\tau$ is a point iff $P_\Theta$ fixes $\tau Q$, i.e. iff $P_\Theta\subseteq\tau P_\theta\tau^{-1}$. The right-hand side then contains $B$, so it is a standard parabolic subgroup, conjugate to $P_\theta$, hence equal to $P_\theta$ (two standard parabolic subgroups which are conjugate are equal), and $\tau Q=Q$. And $P_\Theta\subseteq P_\theta$ iff $\theta\subseteq\Theta$. So a sequence has a point cell iff $\theta\subseteq\Theta$, and then the point is $kQ=kP_\theta$, the attracting flag; Theorem~\ref{thm:Quint}(a) gives the claim.
(b) The cells $O^+_\tau$ of dimension $<\dim X$ are all those other than the open $P_\Theta$-orbit, which is unique. So $T_{\le\dim X-1}(g_m)=k\,(X\setminus O^+_{\mathrm{open}})$. On the other hand, by Lemma~\ref{lem:cellfix} and the dimension count $d_++d_0+d_-=\dim X$, the exceptional sources of $(g_m)$ are $k'^{-1}\bigl(X\setminus O^-_{\mathrm{open}}\bigr)$, since $\Sigma_+(\sigma)=\varnothing$ iff $O^-_\sigma$ has full dimension. Now $g_m^{-1}=(k_m'^{-1}\dot w_0)\,e^{\iota\mu_m}\,(\dot w_0^{-1}k_m^{-1})$ with $\iota\mu=-w_0\mu\in\ak^+$ is a Cartan decomposition of type $\iota\Theta$, and $\dot w_0P_{\iota\Theta}\dot w_0^{-1}=P^-_\Theta$; so the target set of dimension $<\dim X$ of $(g_m^{-1})$ is $k'^{-1}\dot w_0\bigl(X\setminus O^+_{\mathrm{open}}(\iota\Theta)\bigr)=k'^{-1}(X\setminus O^-_{\mathrm{open}}(\Theta))$. Hence $E_{(g_m)}=T_{\le\dim X-1}(g_m^{-1})$, and since $\Gamma=\Gamma^{-1}$, $\Lambda_{\mathrm{exc}}=\bigcup T_{\le\dim X-1}$. Take closures.
\end{proof}

\begin{theorem}[Flag of limit sets]\label{thm:L}
Let $G$ be a simple real Lie group with trivial centre, $X=G/P_\theta$ a non-trivial flag manifold, $\Gamma<G$ discrete, and $\delta$ a function as above \textup{(}for instance the middle-dimension choice\textup{)}.
\begin{enumerate}[label=\textup{(\alph*)},nosep]
\item \textup{(The chain.)} $\Lambda^Q_\theta=\Lambda_{\mathrm Z}^{(0)}\subseteq\Lambda_{\mathrm Z}^{(1)}\subseteq\cdots\subseteq\Lambda_{\mathrm Z}^{(\dim X-1)}=\Lambda_{\Myr}$, all closed and $\Gamma$-invariant; $\Lambda_{\mathrm Z}^{(d)}\subseteq\Lambda_\delta$ whenever $d\le\min_\Theta\delta_\Theta$.
\item \textup{(Domains.)} $\Omega_\delta$, $\Omega_{\Kul}$, $\Omega_{\Myr}$ are domains of proper discontinuity, and $\mathscr R_\Gamma\subseteq\Lambda_\delta\cap\Lambda_{\Kul}\cap\Lambda_{\Myr}$.
\item \textup{(Trichotomy.)} Exactly one holds: $\Lambda_\delta=\mathscr R_\Gamma$ and $\Omega_\delta$ is the largest domain; or $\mathscr R_\Gamma\subsetneq\Lambda_\delta$ and the largest domain $X\setminus\mathscr R_\Gamma$ exists and strictly contains $\Omega_\delta$; or two non-recurrent points are related by $\D_\Gamma$ and there is no largest domain.
\item \textup{(Position.)} $B_\Gamma\subseteq\mathscr R_\Gamma\subseteq\Lambda_{\Kul}\subseteq\Lambda_{\Myr}$. If $\Omega_\delta$ is the largest domain, then $\Lambda_\delta=\mathscr R_\Gamma$ and $\Lambda_\delta\subseteq\Lambda_{\Kul}\subseteq\Lambda_{\Myr}$.
\item \textup{(Bottom.)} For every discrete $\Gamma<G$,
\[
 \Lambda^Q_\theta\subseteq\mathscr R_\Gamma
 \subseteq\Lambda_{\Kul}\subseteq\Lambda_{\Myr},
 \qquad \mathscr R_\Gamma\subseteq\Lambda_\delta.
\]
If $\Omega_\delta$ is the largest domain, then
\[
 \Lambda^Q_\theta\subseteq\mathscr R_\Gamma=\Lambda_\delta
 \subseteq\Lambda_{\Kul}\subseteq\Lambda_{\Myr}.
\]
\end{enumerate}
\end{theorem}
\begin{proof}
(a) Monotonicity and invariance are clear from the definition, the two ends are Proposition~\ref{prop:ends}, and the last assertion holds because every cell of dimension $\le d\le\delta_\Theta$ is a target of $\mathfrak B^\delta$. (b) Theorems~\ref{thm:balanced}, \ref{thm:exc}, Proposition~\ref{prop:L2} and Corollary~\ref{cor:R-in-both}. (c) Corollary~\ref{cor:trichotomy}. (d) Proposition~\ref{prop:kulkarni-myrberg} gives the unconditional chain, and (c) gives $\Lambda_\delta=\mathscr R_\Gamma$ under the largest-domain hypothesis. (e) The \hyperref[par:measure-return]{concentration argument} of Section~\ref{sec:Quint} gives $\Lambda^Q_\theta\subseteq\mathscr R_\Gamma$ directly from the definition. Combine this with Proposition~\ref{prop:kulkarni-myrberg}, with (b), and with $\Lambda_\delta=\mathscr R_\Gamma$ when $\Omega_\delta$ is the largest domain, by (c).
\end{proof}
\begin{remark}[What is proved in general and what is not]
Theorem~\ref{thm:L} uses neither Zariski density nor minimality of Quint's set: its bottom inclusion follows directly from concentration of a full-support measure. Proposition~\ref{prop:kulkarni-myrberg} compares Kulkarni and Myrberg on every flag manifold in the present setting, without a largest-domain hypothesis. When a dimension-system wall gives a largest domain, it equals the recurrent set and therefore lies below Kulkarni. The geometric descriptions remain different: the kernel formula for Myrberg and the identification of the levels $\Lambda_{\mathrm Z}^{(k)}$ with families of projective subspaces are specific to projective space. For $\PSL(n+1,\C)$, $T_{\le d}$ is the union of the spaces $\PP(G_j)$ with $2\dim_{\mathbb C}\PP(G_j)\le d$. Chapter~\ref{ch:psl} indexes this same chain by complex dimension.
\end{remark}

\section{Comparison with the Anosov and transverse theories}\label{sec:compare}
The cells have given us domains for arbitrary discrete groups. We now
compare this construction with the thickenings used in the Anosov and
transverse theories, keeping track of the question each theorem answers.
For a $\Delta$-Anosov representation of a non-elementary word-hyperbolic group, Stecker proves that every maximal
domain is obtained from a minimal fat ideal
\cite[arXiv version, Proposition~3.15 and Corollary~3.16]{SteckerBalanced};
no cocompactness is required. In $X=\PP(\mathbb K^d)$, $\mathbb K=\R$ or
$\C$, the relative positions form the chain $1<\cdots<d$, with opposition
$j\mapsto d+1-j$. Its unique minimal fat ideal is
$\{1,\ldots,\lceil d/2\rceil\}$. Writing $\xi$ for the complete
boundary flag map, the complement of
$\bigcup_{z\in\partial\Gamma}\PP(\xi^{\lceil d/2\rceil}(z))$
is the largest projective domain. Over $\C$, Proposition~\ref{r3an:borel-return}
gives a direct moving-source proof and identifies every point of the
removed set as recurrent.
Let $\tau\subset\Delta$ be opposition-invariant. A $\tau$-convergence group admits the attracting and repelling flags used in the corresponding flag convergence principle. KLP prove properness for fat thickenings without antipodality or expansion, and cocompactness under the additional convergence, expansion, antipodality and balance hypotheses \cite[Theorems~6.8, 6.13 and~7.12]{KLP}. For arbitrary discrete groups they use thickenings over all face types \cite[\S6.4, Propositions~6.15--6.16]{KLP}. The properness construction also restricts to partial flag orbits. Right Weyl-group invariance allows the thickening itself to be expressed as a union of double cosets. Fat, slim and balanced are defined using the complementary thickening
\[
 \mathrm{Th}^{c}=w_0(W\setminus\mathrm{Th}),
 \qquad W=\mathrm{Th}\sqcup w_0\mathrm{Th}^{c},
\]
by $\mathrm{Th}\supseteq,\subseteq,=\mathrm{Th}^{c}$, respectively \cite[\S3.4.1, Definition~3.18]{KLP}. These are statements about Bruhat down-sets.

\paragraph{The dictionary.} Our sublevel sets are the Bruhat down-sets used for thickenings. The loops force every balanced system to satisfy $\mathfrak S\cup\mathfrak T=\mathcal C_\Theta$, the analogue of fatness; complementary pairs give the balanced case (Proposition~\ref{prop:sublevel}). The proposition below makes the correspondence precise for regular sequences: bad sources are thickened repelling flags, bad targets are thickened attracting flags, and opposition invariance puts both flag families in $\Lambda_\tau$. Fatness supplies properness, while slimness enters cocompactness under convergence, expansion and antipodality; balanced thickenings satisfy both conditions.

\begin{proposition}[The dictionary in the regular case]\label{prop:KLPdict}
Assume $\Delta\ne\varnothing$ and let $\Gamma<G$ be discrete. Suppose every divergent sequence of $\Gamma$ has a subsequence of type $\Delta$ (``Borel-divergent'': Borel--Anosov groups, and the regular-unipotent example of Proposition~\ref{prop:unipotent}), let $X=G/Q$, $Q=P_\theta$, $W^\theta=W/W_\theta$ with the Bruhat order, and let $\Lambda\subset G/B$ be the closed set of attracting flags of sequences of $\Gamma$ (it contains the repelling flags, which are the attracting flags of the inverse sequences). For a Bruhat down-set $\mathrm{Th}\subset W^\theta$ and $\xi=gB$ put $\mathrm{Th}(\xi)=\{x\in X:\ x\in gB\tau Q\text{ for some }\tau\in\mathrm{Th}\}$.
\begin{enumerate}[label=\textup{(\roman*)},nosep]
\item The cell graph $\mathcal G_\Delta$ is the Bruhat order: $E(\sigma)=\{\tau:\tau\le\sigma\}$.
\item The pair $(\mathfrak S,\mathfrak T)=(w_0\mathrm{Th},\mathrm{Th})$ is a vertex cover (so $\mathfrak B$, completed by taking all cells bad on other types, is a balanced system) if and only if $\mathrm{Th}$ contains the Bruhat down-set of $W^\theta\setminus w_0\mathrm{Th}$. This is the ``fat'' condition on thickenings; it holds with equality exactly when $\mathrm{Th}$ is a down-set and $w_0\mathrm{Th}=W^\theta\setminus\mathrm{Th}$ (the ``balanced'' case).
\item In that case
\[
\Lambda_{\mathfrak B}(\Gamma)=\bigcup_{\xi\in\Lambda}\mathrm{Th}(\xi),\qquad \Omega_{\mathfrak B}=X\setminus\bigcup_{\xi\in\Lambda}\mathrm{Th}(\xi),
\]
and $\Gamma$ acts properly discontinuously on $\Omega_{\mathfrak B}$ (Theorem~\ref{thm:balanced}). \textup{(}No antipodality of $\Lambda$ is used.\textup{)}
\end{enumerate}
\end{proposition}
\begin{proof}
(i) is Proposition~\ref{prop:schubert}: $\overline{N\sigma Q}=\overline{B\sigma Q}=\bigcup_{\tau\le\sigma}B\tau Q$. Thus (B) is exactly the down-set inclusion in (ii). For (iii), the bad targets are $T=\mathrm{Th}(kB)$. Conjugating $B^-$ by $w_0$ gives
\[
 k'^{-1}B^-\sigma Q=(k'^{-1}w_0)B(w_0^{-1}\sigma)Q,
\]
so the bad sources are $S=\mathrm{Th}(k'^{-1}w_0B)$, the thickening of the inverse sequence's attracting flag. Both flags lie in $\Lambda$, and every flag of $\Lambda$ occurs. The union is closed by compactness of $\Lambda$ and closedness of the Schubert thickenings, as detailed for all types in Theorem~\ref{thm:KLPfull}(ii); properness is Theorem~\ref{thm:balanced}.
\end{proof}
\noindent For $A_2$ and $\mathrm{Th}=\{e,s_1,s_2\}$ one finds $w_0\mathrm{Th}=\{w_0,s_1s_2,s_2s_1\}=W\setminus\mathrm{Th}$, so this is balanced; the domain is that of the Hitchin example of Section~\ref{sec:lost}. The case of sequences of type $\Theta\subsetneq\Delta$ is Theorem~\ref{thm:KLPfull}.
\paragraph{The role of the destination formula.}
Theorem~\ref{thm:dest} keeps the exact destinations and the residual Levi action, even when a thickening records only which cells they meet. We can therefore use the relational criterion both to test properness and to decide whether a largest domain exists. Section~\ref{r2bd:absorption} spells out the quantifier in this test: every destination of a source outside the removed set must lie inside it. The vertex-cover condition alone gives neither cocompactness nor nonemptiness. Those questions have their own geometric hypotheses in the Anosov and transverse theories of Guichard--Wienhard and Canary--Zhang--Zimmer \cite{GuichardWienhard,mlCZZ}.

There is no need for the attracting and repelling flags of a sequence to be distinct in this calculation. The principal unipotent gives a useful test: both flags coincide, while the destination fibres still run through the Schubert varieties.

\begin{proposition}[Unipotent sequences: the cusp test]\label{prop:unipotent}
Let $\iota:\PSL(2,\K)\to G$ be a principal homomorphism ($\K=\R$ or $\C$). We use compatible Cartan data in which its principal element
\[
 H=d\iota\bigl(\diag(1,-1)\bigr)\in\ak^+
 \quad\text{satisfies}\quad \alpha(H)=2\quad(\alpha\in\Delta).
\]
This specifies the convention also for non-reduced restricted root systems. Put $u=\iota\begin{pmatrix}1&1\\0&1\end{pmatrix}$, and let $F=\iota_*(\infty)\in G/B$ be the flag fixed by the image under $\iota$ of the upper triangular group. For $g_n=u^n$:
\begin{enumerate}[label=\textup{(\alph*)},nosep]
\item $g_n$ has type $\Theta=\Delta$ (all $\alpha(\mu(u^n))\to\infty$, like $2\log n$), the Levi class is trivial, and the attracting flag and the repelling flag of $(g_n)$ coincide: both are $F$.
\item Let $X=G/Q$, $Q=P_\theta$, and $F_\theta$ the image of $F$ in $X$. Write $F=hB$ and let $\tau\in W/W_\theta$ be the relative position of $x\in X$ to $F$, so $x\in hB\tau Q/Q$. Then
\[
\D_{(u^n)}(x)=X_{w_0\tau}(F):=h\,\overline{Bw_0\tau Q/Q}.
\]
\end{enumerate}
In particular $\D_{(u^n)}(x)=\{F_\theta\}$ if and only if $x$ is transverse to $F$ (the open cell), and $\D_{(u^n)}(F_\theta)=X$.
\end{proposition}
\begin{proof}
(a) In $\PSL(2,\K)$, $\begin{pmatrix}1&n\\0&1\end{pmatrix}=\kappa_n\,\mathrm{diag}(s_n,s_n^{-1})\,\kappa'_n$ with $s_n\sim n$, $\kappa_n\to\kappa$, $\kappa_n'\to\kappa'$ in $\mathrm{SO}(2)$ or $\mathrm{SU}(2)$. With the compatible Cartan decompositions ($\iota(K_2)\subset K$, $\iota(A_2)\subset A$),
\[
 u^n=\iota(\kappa_n)\,\exp((\log s_n)H)\,\iota(\kappa'_n)
\]
is a Cartan decomposition, and $\alpha(H)=2$ for every simple restricted root $\alpha$. So $\Theta=\Delta$, $b=1$, and $k=\iota(\kappa)$, $k'=\iota(\kappa')$. The attracting and repelling points of the matrices $\left(\begin{smallmatrix}1&\pm n\\0&1\end{smallmatrix}\right)$ on $\PP^1$ both equal $\infty$, and $\iota$ is equivariant for the boundary map $\xi:\PP^1\to G/B$, $\xi(\infty)=F$; the attracting flag of $(g_n)$ is $kB=\xi(\infty)$ and the repelling one is $k'^{-1}w_0B=\xi(\infty)$. Hence $k'^{-1}w_0B=kB$.
(b) Write $v=k^{-1}x\in B\tau Q$. From $k'^{-1}w_0B=kB$ we get $k'^{-1}w_0=kq_0$ with $q_0\in B$, so $k'x=w_0q_0^{-1}v\in w_0B\tau Q=B^-\,w_0\tau Q$ (as $w_0Bw_0^{-1}=B^-$), i.e.\ $\sigma=w_0\tau$. By Proposition~\ref{prop:schubert} (the case $\Theta=\Delta$), $\D(x)=k\,\overline{Bw_0\tau Q/Q}=X_{w_0\tau}(F)$, since $kB=F=hB$.
\end{proof}
\begin{example}[The cusp in $\PSL(3,\R)$ on $\R\PP^2$]\label{ex:cusp}
Take $u=I+E_{12}+E_{23}$ (the principal unipotent), so $u^n(x_1,x_2,x_3)=(x_1+nx_2+\binom n2x_3,\;x_2+nx_3,\;x_3)$ and $F_\theta=[e_1]$, $F=([e_1],\langle e_1,e_2\rangle)$. If $x_3\ne0$, $u^nx_m\to[e_1]$ for every $x_m\to x$ (the term $\binom n2x_3$ dominates): the open cell has a single destination. If $x_3=0$, $x_2\ne0$, put $x_m=(x_1,x_2,-2x_2/n+\delta_m)$ with $n^2\delta_m/2\to c$: the image tends to $[x_1+x_2+c:-x_2:0]$, so the destination is the line $\{y_3=0\}=\PP(\langle e_1,e_2\rangle)$. At $x=[e_1]$ a source $(1,x_2,\varepsilon)$ with $\varepsilon\sim2/n^2$ and $x_2\sim-2/n$ reaches any $[a:b:c]$ with $c\ne0$: the destination is all of $\R\PP^2$. These are the fibres $X_{w_0\tau}$ of dimensions $0,1,2$ for the cells $\tau$ of dimensions $2,1,0$: Proposition~\ref{prop:unipotent}.
\end{example}
\begin{remark}[Answer to the cusp question]\label{rem:cuspanswer}
Thus, for $p^n$ with $p$ a principal unipotent (cusp element of a Hitchin-type group), the attracting and repelling flags coincide: $\xi^+=\xi^-=F$. Points transverse to $F$ are sent to $F_\theta$, and the non-transverse points have the Schubert-variety fibres of Proposition~\ref{prop:unipotent}(b). This coincidence is compatible with antipodality of the limit set, since antipodality compares distinct flags. In a transverse group, the convergence on the open cell is the parabolic fixed-point case of the convergence principle of Canary--Zhang--Zimmer. The calculation here describes the destinations of the unipotent sequence directly; comparison with a transverse-group theorem retains its own hypotheses. For a unipotent that is not principal, the same argument with the Jacobson--Morozov grading gives a type $\Theta=\{\alpha:\alpha(h)\ne0\}$ and a single flag of $G/P_\Theta$ playing both roles. The nonprincipal case is not treated here.
\end{remark}

%--------------------------------------------------------------
%--------------------------------------------------------------
\section{The dictionary with Kapovich--Leeb--Porti, all types}\label{sec:KLPfull}
For types $\Theta\subsetneq\Delta$, the residual Levi factor joins Bruhat cells into double-coset cells. We first check that destinations still follow the closure order, settling Remark~\ref{rem:Eorder}; this is the extra cell geometry needed to extend Proposition~\ref{prop:KLPdict}.

\begin{lemma}[Bruhat order on double cosets]\label{lem:doublecosets}\
For $I,J\subseteq\Delta$ and $D,D'\in W_I\backslash W/W_J$ let $d,d'$ be the elements of minimal length. The following are equivalent: (i) $d'\le d$ in the Bruhat order of $W$; (ii) some element of $D'$ is $\le$ some element of $D$ (equivalently, $\le\max D$); (iii) the $P_I$-orbit of $D'$ in $G/P_J$ lies in the closure of that of $D$. Moreover $d$ is minimal in $W_Id$ and in $dW_J$, and $\{D'':D''\le D\}=\{W_I\tau W_J:\ \tau\le d\}$.
\end{lemma}
\noindent(This is standard; see Bj\"orner--Brenti~\cite{BjornerBrenti}, Ch.~2, or Richardson--Springer~\cite{RichardsonSpringer1990}.)

\begin{proposition}[The cell graph is the order graph]\label{prop:Eorder}\
For every nonempty $\Theta\subseteq\Delta$, every $\theta\subseteq\Delta$ and every $\sigma\in\mathcal C_\Theta=W_\Theta\backslash W/W_\theta$, $E(\sigma)=\{\tau\in\mathcal C_\Theta:\tau\le\sigma\}$. This answers Remark~\ref{rem:Eorder}.
\end{proposition}
\begin{proof}
Take $\sigma$ to be the minimal representative. Then $\sigma^{-1}(\Phi^+\cap\Phi(L_\Theta))\subset\Phi^+$, so $\sigma^{-1}N_L\sigma\subset N\subset Q$ for $N_L=N\cap L_\Theta$, and $N_L\sigma Q=\sigma Q$. Since $N=N_\Theta N_L$, $N\sigma Q=N_\Theta\sigma Q$, and $\overline{N_\Theta\sigma Q}=\overline{B\sigma Q}=\bigcup_{\tau'\le\sigma}B\tau'Q$ (the Schubert variety of $\sigma W_\theta$). Each $B\tau'Q$ lies in the cell $O^+_{[\tau']}$ of its double coset, so $E(\sigma)=\{[\tau']:\tau'\le\sigma\}$, which is $\{\tau\le\sigma\}$ by Lemma~\ref{lem:doublecosets}.
\end{proof}

\paragraph{Setting.} Let $\tau\subseteq\Delta$ be invariant under $\iota=-w_0$ (so $w_0$ normalises $W_\tau$). A discrete $\Gamma<G$ is $\tau$-\emph{convergent} if every divergent sequence has a subsequence of type $\Theta\supseteq\tau$. Let $\Lambda_\tau\subset G/P_\tau$ be the closed set of images of the attracting flags of such subsequences (it contains the repelling ones, which are attracting flags of the inverse sequences). Put $\mathcal C_\tau=W_\tau\backslash W/W_\theta$; for $\Theta\supseteq\tau$ the projection $\pi:\mathcal C_\Theta\to\mathcal C_\tau$ is surjective and order preserving (Lemma~\ref{lem:doublecosets}), and left multiplication by $w_0$ is an order-reversing involution of $\mathcal C_\tau$. For $\xi=gP_\tau$ and $\mathrm{Th}\subseteq\mathcal C_\tau$ let $\mathrm{Th}(\xi)=\{x\in X:\ x\in gP_\tau\tau'Q\ \text{for some }[\tau']\in\mathrm{Th}\}$.

For $X=G/B$ ($W_\theta=1$), $\mathcal C_\tau=W_\tau\backslash W$ and $\mathrm{Th}$ is a left $W_\tau$-invariant subset of $W$: a thickening in the sense of Kapovich--Leeb--Porti. We write $\mathrm{Th}^\vee:=w_0(\mathcal C_\tau\setminus\mathrm{Th})$, so that $\mathcal C_\tau=\mathrm{Th}\sqcup w_0\mathrm{Th}^\vee$, which is the relation $W=\mathrm{Th}\sqcup w_0\mathrm{Th}^c$ of KLP with $\mathrm{Th}^c=\mathrm{Th}^\vee$. Fat, slim, balanced mean $\mathrm{Th}\supseteq,\subseteq,=\mathrm{Th}^\vee$. Since $w_0\mathrm{Th}^\vee=\mathcal C_\tau\setminus\mathrm{Th}$, fat means $\mathrm{Th}\cup w_0\mathrm{Th}=\mathcal C_\tau$, and balanced means $w_0\mathrm{Th}=\mathcal C_\tau\setminus\mathrm{Th}$.

\begin{theorem}[Dictionary, general case]\label{thm:KLPfull}
Assume $\Delta\ne\varnothing$. Let $\Gamma$ be $\tau$-convergent and $\mathrm{Th}\subseteq\mathcal C_\tau$ a Bruhat down-set.
\begin{enumerate}[label=\textup{(\roman*)},nosep]
\item For nonempty $\Theta\supseteq\tau$ put $\mathfrak T_\Theta=\pi^{-1}(\mathrm{Th})$, $\mathfrak S_\Theta=\pi^{-1}(w_0\mathrm{Th})$ \textup{(}complete the system by any vertex covers for $\Theta\not\supseteq\tau$, types which do not occur\textup{)}. Then $(\mathfrak S_\Theta,\mathfrak T_\Theta)_\Theta$ is a vertex cover of every such $\mathcal G_\Theta$ \emph{if and only if} $\mathrm{Th}$ is fat.
\item In that case $\Lambda_{\mathfrak B}(\Gamma)=\bigcup_{\xi\in\Lambda_\tau}\mathrm{Th}(\xi)$, which is closed, and $\Gamma$ acts properly discontinuously on
\[
\Omega_{\mathfrak B}=X\setminus\bigcup_{\xi\in\Lambda_\tau}\mathrm{Th}(\xi).
\]
For $X=G/B$ this is the fat-thickening properness domain of KLP \cite[Theorems~6.8 and~6.13]{KLP}. Their properness theorem also restricts to partial flag orbits. Here the ideal on $W$ is the pullback from the double-coset space, so it has right $W_\theta$-invariance.
\item The system is balanced in the sense of Proposition~\ref{prop:sublevel} \textup{(}$\mathfrak S_\Theta,\mathfrak T_\Theta$ complementary\textup{)} for all nonempty $\Theta\supseteq\tau$ iff $\mathrm{Th}$ is balanced.
\end{enumerate}
\end{theorem}
\begin{proof}
(i) Necessity follows from the loops for $\Theta=\Delta$ and surjectivity of $\pi$: (B) forces $\mathcal C_\tau\setminus w_0\mathrm{Th}\subseteq\mathrm{Th}$, the fat condition. Conversely, under this inclusion, $\sigma\notin\pi^{-1}(w_0\mathrm{Th})$ implies $\pi(\sigma)\in\mathrm{Th}$; for every destination label $\sigma'\le\sigma$, order preservation and the down-set property give $\pi(\sigma')\in\mathrm{Th}$. Proposition~\ref{prop:Eorder} now gives (B).
(ii) Let a subsequence have data $(\Theta,k,k',b)$ and projected flags $\xi^+=kP_\tau$, $\xi^-=k'^{-1}w_0P_\tau$ (the attracting flag of the inverse sequence, of type $\iota\Theta\supseteq\iota\tau=\tau$). Then $T=k\bigcup_{\sigma\in\mathfrak T_\Theta}P_\Theta\sigma Q=\mathrm{Th}(\xi^+)$, because $P_\Theta\sigma Q\subset P_\tau\sigma Q$ and the double coset of $\sigma$ in $\mathcal C_\tau$ is $\pi(\sigma)$. For the sources, $k'^{-1}P^-_\Theta\sigma Q=(k'^{-1}w_0)P_{\iota\Theta}\,w_0\sigma Q$, whose label in $\mathcal C_\tau$ is $w_0\pi(\sigma)$; if $\sigma\in\mathfrak S_\Theta$ then $w_0\pi(\sigma)\in\mathrm{Th}$, so $S=\mathrm{Th}(\xi^-)$. Every $\xi\in\Lambda_\tau$ occurs as some $\xi^+$ and as some $\xi^-$. Since $\mathrm{Th}$ is a down-set, $\mathrm{Th}(\xi)$ is a union of Schubert varieties, hence closed, and $\xi\mapsto\mathrm{Th}(\xi)$ is continuous in the Hausdorff topology, so the union over the compact $\Lambda_\tau$ is closed. Properness is Theorem~\ref{thm:balanced}.
(iii) is immediate from the definitions.
\end{proof}

Opposition invariance let us put the forward and backward flags in
the same flag manifold and remove one family of thickenings. If the
two flag types differ, we keep the two families apart. The same cell
argument then removes a complementary thickening at the source and
the chosen thickening at the target. Fatness enters when the two
families can be absorbed into one; the two-family construction itself
does not need it.

\begin{proposition}[Two thickenings for arbitrary flag types]\label{r3klp:nonself}
Allow an arbitrary $\tau\subseteq\Delta$ in the definition of
$\tau$-convergence, and put $\tau^*=\iota\tau$. For a $\tau$-convergent
discrete group $\Gamma$, let $\Lambda_\tau^+\subset G/P_\tau$ and
$\Lambda_\tau^-\subset G/P_{\tau^*}$ be the closed sets of forward
and backward projected limit flags. Thus a sequence of type
$\Theta\supseteq\tau$ has flags
$\xi^+=kP_\tau$ and $\xi^-=k'^{-1}w_0P_{\tau^*}$; the latter is an
attracting flag of the inverse sequence, of type
$\iota\Theta\supseteq\tau^*$.
Let $I\subset W$ be a Bruhat lower ideal with
$W_\tau I=I=IW_\theta$, and put
$I^\vee=w_0(W\setminus I)=W\setminus w_0I$.
For a flag $\xi$ of type $\rho$, write $\mathrm{Th}_J(\xi)$ for the
points of $X$ whose positions lie in $W_\rho\backslash J/W_\theta$.
Then every such sequence satisfies
\begin{equation}\label{r3klp:row-absorption}
 x\notin\mathrm{Th}_{I^\vee}(\xi^-)
 \quad\Longrightarrow\quad
 \D_{(g_m)}(x)\subseteq\mathrm{Th}_I(\xi^+).
\end{equation}
Consequently the closed invariant set
\begin{equation}\label{r3klp:two-thickenings}
 L_I=\bigcup_{\xi^-\in\Lambda_\tau^-}\mathrm{Th}_{I^\vee}(\xi^-)
       \ \cup\!
       \bigcup_{\xi^+\in\Lambda_\tau^+}\mathrm{Th}_I(\xi^+)
\end{equation}
is a vertex cover of $\D_\Gamma$, and $\Gamma$ acts properly on
$X\setminus L_I$. No fatness is needed for this two-set union.
For $\tau^*=\tau$ and fat $I$, it reduces to
$L_I=\bigcup_{\xi\in\Lambda_\tau}\mathrm{Th}_I(\xi)$.
\end{proposition}
\begin{proof}
The ideal $I^\vee$ is left $W_{\tau^*}$-invariant and right
$W_\theta$-invariant, so both thickenings have the required types.
Use the cell argument of Theorem~\ref{thm:KLPfull}: if the source cell
has representative $\sigma$, its position relative to $\xi^-$ is
represented by $w_0\sigma$. Exclusion from $I^\vee$ means $\sigma\in I$;
every destination cell has a representative $\sigma'\le\sigma$
(Proposition~\ref{prop:Eorder}), hence lies in $I$.
This proves \emph{all}-destination absorption
\eqref{r3klp:row-absorption}. Extracting a Cartan sequence from any witnessing
sequence proves the vertex-cover assertion. Closedness and invariance
follow from compactness and equivariance of thickenings as in
Theorem~\ref{thm:KLPfull}; properness follows from
Theorem~\ref{thm:relational}. In the self-opposed case inversion gives
$\Lambda_\tau^- =\Lambda_\tau^+=\Lambda_\tau$, and fatness gives
$I^\vee\subseteq I$.
\end{proof}

The two-family properness conclusion is that of KLP~\cite[Theorem~6.8]{KLP}, expressed here through the exact destination fibres.

\subsection{Absorption of every destination}\label{r2bd:absorption}
The implication in Proposition~\ref{r3klp:nonself} concerns every point
of the destination fibre. This is exactly what properness asks for.
For a closed invariant set $L\subset X$, Theorem~\ref{thm:relational}
can be written in the fibre form
\begin{equation}\label{r2bd:all-destinations}
 \Gamma\curvearrowright X\setminus L\text{ is proper}
 \quad\Longleftrightarrow\quad
 \D_\Gamma(x)\subset L\quad\text{for every }x\notin L.
\end{equation}
Indeed both sides say that $\D_\Gamma$ has no pair in
$(X\setminus L)^2$. A source outside the removed set must therefore
have \emph{all} its destinations inside. A single destination in a
small flag limit set tells us much less, as the rank-one picture
already shows.

\begin{example}[The existential test already fails in rank one]
\label{r2bd:existential}
Let $\Gamma<\PSL(2,\C)$ be a classical Schottky group of rank two,
with Cantor limit set $\Lambda\subset X=\C\PP^1$.
For $W=\{e,s\}$ the ideal $\mathrm{Th}=\{e\}$ is balanced and
$\mathrm{Th}(\Lambda)=\Lambda$. But
\[
 \{x\in X:\D_\Gamma(x)\cap\Lambda\ne\varnothing\}=X.
\]
To see the quantifiers directly, choose a loxodromic
$h(z)=\lambda z$, $|\lambda|>1$, with fixed points
$a^-=0$, $a^+=\infty$ in $\Lambda$. Then
\[
 \D_{(h^n)}=(X\times\{a^+\})\cup(\{a^-\}\times X).
\]
Off $0$ there is locally uniform convergence to $\infty$.
Every finite $y$ is reached from $0$ by
$z_n=\lambda^{-n}y$; the destination $\infty$ is reached by
$z_n=\lambda^{-\lfloor n/2\rfloor}$.
Thus every source has a destination in $\Lambda$, while the removed
thickening is the proper subset $\Lambda$. Formula
\eqref{r2bd:all-destinations}, rather than existential intersection,
is the properness test.
\end{example}

With the complementary-thickening convention of Section~\ref{sec:compare}, Theorem~\ref{thm:KLPfull} reads KLP properness directly in the fibres. Their construction applies to arbitrary discrete groups and restricts to partial flag orbits \cite[\S6.4, Propositions~6.15--6.16]{KLP}. In the double-coset description the thickening is itself right $W_\theta$-invariant. These properness statements assert neither nonemptiness nor the largest-domain property. The latter holds exactly when no destination pair joins two non-recurrent points (Theorem~\ref{thm:relational}); in the projective Borel--Anosov case it also follows from Stecker's classification, as explained in Section~\ref{sec:compare}.

\chapter{The adjoint action}
\label{ch:adjoint}

Every flag manifold of $G$ is a different space, and each seems to need its own treatment. But there is one projective space attached to $G$ with no choices at all: the projective space of its own Lie algebra. A semisimple group with trivial centre sits inside $\mathrm{GL}(\g)$ as $\mathrm{Aut}(\g)_0$ (Section~\ref{sec:setting-sum}), so the adjoint representation is just the inclusion $\mathrm{Ad}:G\hookrightarrow\PGL(\g)$, and a discrete $\Gamma<G$ acts on $\PP(\g)$ like a subgroup of $\PSL(n+1,\C)$ acts on $\C\PP^n$. For a simple noncompact $G$ whose highest restricted root space is a line, $\PP(\g)$ contains a closed $G$-orbit, the \emph{adjoint variety} (the projectivised minimal nilpotent orbit), which is a flag manifold of $G$. For higher root multiplicity the same flag is represented by a subspace of the Lie algebra, and its closed orbit lies in a Grassmannian.

\medskip  The restricted roots give the singular values of $\mathrm{Ad}(g_m)$. The root spaces in Cartan coordinates identify the middle and top blocks, and the limiting compact factors position them.
\begin{restated}{Theorem~\ref{thm:adjoint-spectrum}}Let $\g$ be simple and noncompact, and let $g_m=k_me^{\mu_m}k'_m$ have type $\Theta$. Then $\mathrm{Ad}(g_m)$ has singular values $e^{\lambda(\mu_m)}$ ($\lambda$ a restricted root or $0$). In Cartan coordinates its \emph{middle block} is the Levi subalgebra $\g_{\Theta}^0$ of $P_\Theta$ and its \emph{top block} is the top graded piece $\g^{\mathrm{top}}_\Theta$ of the $\Theta$-grading. The first pseudo-projective limit is $[\mathrm{Ad}(k)\,\mathrm{Ad}(b)\,\Pi_{\mathrm{top}}\,\mathrm{Ad}(k')]$. When the highest restricted root space is a line, the limit has rank one, with image in the adjoint variety, exactly when $\Theta\supseteq\Theta_{\mathrm{ad}}$.\end{restated}

The idea is as follows.  The Cartan decomposition of $g_m$ is already a singular value decomposition of $\mathrm{Ad}(g_m)$ for the form $B_\vartheta$, with singular values $e^{\lambda(\mu_m)}$. The exponents that stay bounded give the Levi subalgebra, and the largest exponents give the top graded piece. The middle block identifies $\Theta$ through its labelled root spaces in Cartan coordinates. The numerical singular values alone can coincide for different types. After extraction, the ordered adjoint blocks record which root exponents have bounded differences and how the resulting blocks are ordered. Different numerical rates can give the same blocks (Section~\ref{sec:ratetypes}). For the adjoint variety, the condition $\Theta\supseteq\Theta_{\mathrm{ad}}$ identifies the rank-one first limits when the highest restricted root space is a line. The first target space of $\mathrm{Ad}(g_m)$ that crosses the middle dimension contributes to the \emph{critical wall} used in Chapters~\ref{ch:parabolic} and~\ref{ch:psl}. Its parabolic description does not settle the question of recurrence: a matrix may preserve a limit flag without being able to return to itself.

\medskip When $\g$ is simple and noncompact, we write $\beta_{\max}$ for the highest restricted root, the maximum of the irreducible system $\Phi$ for the dominance order, and $d_{\max}=\dim\g_{\beta_{\max}}$. For a semisimple algebra these roots are taken separately in its noncompact simple factors, as in Remark~\ref{rem:semisimple}; compact factors have no restricted roots. The letter $\theta$ keeps its meaning of the index of a flag manifold.

%--------------------------------------------------------------
\section{The spectrum of \texorpdfstring{$\mathrm{Ad}(g_m)$}{Ad(gm)}}\label{sec:adspec}
Let $B_\vartheta(X,Y)=-B(X,\vartheta Y)$ be the positive definite form on $\g$ built from the Killing form $B$ and the Cartan involution $\vartheta$. The root space decomposition $\g=\g_0\oplus\bigoplus_{\lambda\in\Phi}\g_\lambda$ is $B_\vartheta$-orthogonal, $\mathrm{Ad}(K)$ is orthogonal, and $\mathrm{Ad}(e^{\mu})$ ($\mu\in\ak$) acts on $\g_\lambda$ by $e^{\lambda(\mu)}$ and trivially on $\g_0=\mathfrak m\oplus\ak$. Hence
\[
\mathrm{Ad}(g_m)=\mathrm{Ad}(k_m)\,\mathrm{Ad}(e^{\mu_m})\,\mathrm{Ad}(k'_m)
\]
is a singular value decomposition for $B_\vartheta$: the singular values are $e^{\lambda(\mu_m)}$, $\lambda\in\Phi\cup\{0\}$, with multiplicity $\dim\g_\lambda$ ($\dim\g_0$ for $\lambda=0$). Their logarithms are symmetric about $0$; on the chamber, their numerical ordering is compatible with dominance, while the order of incomparable roots depends on the Cartan parameter.

Fix a type $\Theta$ and write $\mathrm{coord}_\Theta(\lambda)=(n_\alpha(\lambda))_{\alpha\in\Theta}$ for the $\Theta$-coordinates of a root $\lambda=\sum n_\alpha\alpha$. The \emph{$\Theta$-grading} is $\g=\bigoplus_{c\in\mathbb Z^\Theta}\g^{(c)}$, $\g^{(c)}=\bigoplus_{\mathrm{coord}_\Theta\lambda=c}\g_\lambda$ (with $\g_0$ in degree $0$). Then $\g^{(0)}=\g_\Theta^0=\mathrm{Lie}(L_\Theta)$ is the Levi factor of $P_\Theta$, $\bigoplus_{c>0}\g^{(c)}=\mathrm{Lie}(N_\Theta)$, and the \emph{top piece} is $\g^{\mathrm{top}}_\Theta=\g^{(c_{\max})}$ with $c_{\max}=\mathrm{coord}_\Theta(\beta_{\max})$.

\begin{theorem}[Spectrum and blocks]\label{thm:adjoint-spectrum}
Let $\g$ be simple and let $(g_m)$ be a divergent sequence with Cartan data $(\Theta,k,k',b)$.
\begin{enumerate}[label=\textup{(\alph*)},nosep]
\item Two roots (or $0$) $\lambda,\lambda'$ with $\mathrm{coord}_\Theta\lambda=\mathrm{coord}_\Theta\lambda'$ have a bounded gap $|\lambda(\mu_m)-\lambda'(\mu_m)|$. In Cartan coordinates, the middle block of $\mathrm{Ad}(g_m)$ (the singular values whose logarithms stay bounded) is exactly $\g^{(0)}_\Theta=\mathrm{Lie}(L_\Theta)$, of dimension $\dim L_\Theta$. The remaining blocks lie in the positive part $\mathrm{Lie}(N_\Theta)$ or the negative part $\mathrm{Lie}(N^-_\Theta)$, which are exchanged by $g\mapsto g^{-1}$.
\item The top Cartan block of $\mathrm{Ad}(g_m)$ is $\g^{\mathrm{top}}_\Theta$: a root $\lambda$ has $\beta_{\max}(\mu_m)-\lambda(\mu_m)$ bounded iff $\mathrm{coord}_\Theta\lambda=c_{\max}$.
\item $\Theta=\{\alpha\in\Delta:\ \g_\alpha\not\subset\g_\Theta^{0}\}$: the type of $(g_m)$ is read off from the middle block in Cartan coordinates.
\item The first pseudo-projective limit of $(\mathrm{Ad}(g_m))$ is
\[
\tau^+=\bigl[\mathrm{Ad}(k)\,\mathrm{Ad}(b)\,\Pi_{\mathrm{top}}\,\mathrm{Ad}(k')\bigr],\qquad
\mathrm{Im}\,\tau^+=\PP\bigl(\mathrm{Ad}(k)\g^{\mathrm{top}}_\Theta\bigr),\quad
\ker\tau^+=\PP\bigl(\mathrm{Ad}(k')^{-1}(\g^{\mathrm{top}}_\Theta)^{\perp}\bigr),
\]
where $\Pi_{\mathrm{top}}$ is the $B_\vartheta$-orthogonal projection onto $\g^{\mathrm{top}}_\Theta$ and $\perp$ denotes the orthogonal complement with respect to $B_\vartheta$.
\end{enumerate}
\end{theorem}
\begin{proof}
(a) By Lemma~\ref{lem:split}, $\lambda(\mu_m)=\lambda(\nu_m)+\lambda(\eta_m)$ with $\eta_m$ convergent and $\lambda(\nu_m)=\sum_{\beta\in\Theta}n_\beta(\lambda)\beta(\nu_m)$, where $\beta(\nu_m)=\beta(\mu_m)-\beta(\eta_m)\to\infty$ for $\beta\in\Theta$. If the $\Theta$-coordinates of $\lambda$ and $\lambda'$ agree, $(\lambda-\lambda')(\nu_m)=0$ and the gap is $(\lambda-\lambda')(\eta_m)$, bounded. For a positive root, $\lambda(\mu_m)\to\infty$ iff some $n_\beta(\lambda)>0$ with $\beta\in\Theta$ (all terms $n_\alpha\alpha(\mu_m)$ are $\ge0$); so the bounded exponents are the roots with $\mathrm{coord}_\Theta=0$ and $0$, and for negative roots symmetrically.
(b) For every root $\lambda$ one has $\beta_{\max}-\lambda=\sum c_\alpha\alpha$ with $c_\alpha\ge0$, so $(\beta_{\max}-\lambda)(\mu_m)=\sum c_\alpha\alpha(\mu_m)$ is bounded iff $c_\alpha=0$ for every $\alpha\in\Theta$, i.e.\ iff $\mathrm{coord}_\Theta\lambda=\mathrm{coord}_\Theta\beta_{\max}$.
(c) By (a), $\g_\alpha$ ($\alpha\in\Delta$) lies in the middle block iff $\alpha(\mu_m)$ is bounded.
(d) Normalise by $e^{-\beta_{\max}(\mu_m)}$. On $\g^{\mathrm{top}}_\Theta$, $\mathrm{Ad}(e^{\mu_m})=\mathrm{Ad}(z_m)\mathrm{Ad}(b_m)$ with $z_m=e^{\nu_m}$ acting by the scalar $e^{\beta_{\max}(\nu_m)}$ (all $\lambda$ in the top piece have the same $\Theta$-coordinates, hence the same value on $\nu_m\in\ak_\Theta$), and $\mathrm{Ad}(b_m)\to\mathrm{Ad}(b)$ preserves $\g^{\mathrm{top}}_\Theta$ because $b\in L_\Theta$. On every other block the normalised action tends to $0$. Passing to the limit with $k_m\to k$, $k'_m\to k'$ gives the formula.
\end{proof}
\begin{remark}[Relative growth in the adjoint blocks]\label{rem:speeds}
Part (a) fixes the middle block. The full adjoint flag can require a further extraction on which the differences of root exponents converge in the extended real line. The positive part $\mathrm{Lie}(N_\Theta)$ can split into several blocks of $\mathrm{Ad}(g_m)$, according to how the divergent numbers $\beta(\nu_m)$, $\beta\in\Theta$, compare: two roots with different $\Theta$-coordinates $c\ne c'$ lie in the same block iff $\sum_{\beta\in\Theta}(c_\beta-c'_\beta)\beta(\nu_m)$ is bounded, which depends on relative growth and is invisible to the flag destinations of $G$ for fixed $(\Theta,k,k',b)$ (Section~\ref{sec:setting-sum}). For example, for $G=\PSL(2,\C)$ the Cartan type is $\Theta=\Delta$. The holomorphic adjoint blocks have dimensions $(1,1,1)$, and the underlying real adjoint blocks have dimensions $(2,2,2)$. For $\Theta=\Delta$ in higher rank the adjoint type depends on the sequence through the heights of the roots along $\mu_m$: if $\mu_m=m\mu_0$ with $\mu_0$ regular rational, the blocks are the level sets of $\lambda\mapsto\lambda(\mu_0)$.
\end{remark}

\subsection*{Examples}
For the complex groups in the examples below, we use the holomorphic adjoint module, complex root-space dimensions and complex projective spaces. The underlying real adjoint action instead uses real multiplicities and real projectivisation, as do the real forms considered below.

\begin{example}[$\PSL(2,\C)$]\label{ex:adsl2}
$\mathrm{Ad}:\PSL(2,\C)\to\mathrm{SO}(3,\C)=\PSL(2,\C)<\PSL(3,\C)$ is the principal image $N=2$ of Section~\ref{sec:principal} (the representation $\mathrm{Sym}^2$), and the adjoint variety is the conic $\{[X]\in\PP(\mathfrak{sl}_2):\ X\ \text{nilpotent}\}=\{\text{double points}\}\cong\C\PP^1$. For a divergent sequence with attractor $a$ and repeller $b$, the first limit has image line $\mathrm{Ad}(k)\g_{\beta_{\max}}=\C\,e_a$ and projective kernel $\PP(\{X:\ B(X,e_b^-)=0\})=$ the tangent line to the conic at $b$. This recovers Proposition~\ref{prop:sym2}: $\Lambda_{\Myr}(\mathrm{Ad}\,\Gamma)=\bigcup_{b\in\Lambda}T_b$, and the conic meets each $T_b$ only at $b$, so $\Lambda_{\Myr}\cap X_{\mathrm{ad}}=\Lambda$.
\end{example}
\begin{example}[$\PSL(n+1,\C)$]\label{ex:adpsl}
Here $\mathrm{Ad}$ is the action on $\mathfrak{sl}_{n+1}$ by conjugation, and $\Phi=\{\varepsilon_i-\varepsilon_j\}$ with $\beta_{\max}=\varepsilon_1-\varepsilon_{n+1}=\alpha_1+\dots+\alpha_n$. Singular values: $\sigma_i/\sigma_j$ and $1$. For a type with blocks $b_1,\dots,b_k$ (Section~\ref{sec:ext}) the middle block is $\mathfrak s(\mathfrak{gl}_{b_1}\oplus\dots\oplus\mathfrak{gl}_{b_k})$ and the top block is $\mathrm{Hom}(E_k,E_1)$ (a $b_1b_k$-dimensional space): the matrix units sending the slowest block to the fastest. For $k=2$ the top block is the whole $\mathrm{Hom}(E_2,E_1)=\mathrm{Lie}(N_\Theta)$, an abelian Lie algebra. The adjoint variety is $X_{\mathrm{ad}}=\{(p,H):p\in H\}\subset\PP^n\times\PP^{n*}$, the set of rank-one nilpotents $v\otimes\varphi$, $\varphi(v)=0$; $\Theta_{\mathrm{ad}}=\{\alpha_1,\alpha_n\}$.
\end{example}

%--------------------------------------------------------------
\begin{example}[Rank one: $\mathrm{PU}(n,1)$ and $\mathrm{SO}(n,1)$]\label{ex:adrank1}
For $\g=\mathfrak{su}(n,1)$ with $n\ge2$, the restricted roots are $\pm\alpha$ (multiplicity $2n-2$), $\pm2\alpha$ (multiplicity $1$) and $\g_0=\mathfrak m\oplus\ak$. Thus $\beta_{\max}=2\alpha$, $d_{\max}=1$, $\Theta_{\mathrm{ad}}=\Delta=\{\alpha\}$ and the adjoint variety is the whole boundary $\partial\mathbb H^n_\C=S^{2n-1}$, embedded in $\PP(\g)=\R\PP^{n^2+2n-1}$ as the projectivised minimal nilpotent orbit (the centre of the Heisenberg algebra at each boundary point). The Cartan blocks of $\mathrm{Ad}(g_m)$ are $e^{2\alpha(\mu_m)}$ (dimension $1$), $e^{\alpha(\mu_m)}$ (dimension $2n-2$), the middle block $\g_0$ and the symmetric ones. For $\g=\mathfrak{so}(n,1)$ the roots are $\pm\alpha$ with multiplicity $n-1$, so $d_{\max}=n-1$ and the positive part consists of one block at scale $e^{\alpha(\mu_m)}$: the limit $\tau^+$ has rank $n-1$ and its image is the projective subspace $\PP(\mathrm{Ad}(k)\g_\alpha)$ of dimension $n-2$; the kernel has codimension $n-1$. In both cases, Theorem~\ref{thm:adjvar} gives $\Lambda_{\Myr}(\mathrm{Ad}\,\Gamma)=\bigcup_{b\in\Lambda}K_b$ with $K_b$ the $B$-orthogonal complement of $\mathrm{Ad}(h)\g_{-\beta_{\max}}$, and for $\mathrm{PU}(n,1)$ one has $\Lambda=\Lambda_{CG}$, so the adjoint action detects exactly the Chen--Greenberg set.
\end{example}

\begin{remark}[Semisimple groups]\label{rem:semisimple}
If $\g=\g_1\oplus\dots\oplus\g_r$ is semisimple, $G=\prod G_i$ (we assume as in Section~\ref{sec:setting-sum} that $G$ has trivial centre), $\mathrm{Ad}=\bigoplus\mathrm{Ad}_i$ and $\Phi=\bigsqcup\Phi_i$. A sequence $g_m=(g_m^{(i)})$ has type $\Theta=\bigsqcup\Theta_i$, where $\Theta_i\subseteq\Delta_i$ is the type of the $i$-th component (possibly empty), and Theorem~\ref{thm:adjoint-spectrum} holds factor by factor: the middle block is $\bigoplus_i\g_{\Theta_i}^0$ (including every factor where $\Theta_i=\varnothing$, which is bounded), and the top block of a noncompact $i$-th factor is $\g^{\mathrm{top}}_{\Theta_i}$. A compact factor contributes only a bounded block; the highest-root exponents below refer to noncompact factors. What changes is the interleaving: the exponents $\beta^{(i)}_{\max}(\mu^{(i)}_m)$ diverge in the escaping factors, where their rates can vary independently, so the Cartan blocks of $\mathrm{Ad}(g_m)$ in $\PSL(\g)$ are unions of blocks of the factors that have bounded mutual gap. In particular the type in $\PSL(\g)$ of a sequence in a product is not determined by the types of its components; this is the first place where the adjoint action sees genuinely more than the product of the flag manifolds of the factors (cf.\ Remark~\ref{rem:speeds}). The restricted roots may also be non-reduced only inside a factor.
\end{remark}

%--------------------------------------------------------------
\section{The adjoint variety}\label{sec:advar}
Suppose $\g$ is simple and noncompact. Let
\[
\Theta_{\mathrm{ad}}=\{\alpha\in\Delta:\ \beta_{\max}-\alpha\in\Phi\cup\{0\}\}.
\]
Here $0$ is a weight of $\mathrm{Ad}$, not a root. These are the simple roots linked to $\beta_{\max}$, equivalently those with $\langle\beta_{\max},\alpha^\vee\rangle>0$; in the reduced case they are the nodes adjacent to the extra node of the extended Dynkin diagram. For $A_n$ they are $\alpha_1,\alpha_n$, including $A_1$, where the two names denote the same root.

The minimal parabolic preserves $\g_{\max}:=\g_{\beta_{\max}}$. For $\alpha\notin\Theta_{\mathrm{ad}}$, the bracket $[\g_{-\alpha},\g_{\max}]$ vanishes because its weight is neither a root nor $0$. For $\alpha\in\Theta_{\mathrm{ad}}$, it does not vanish. Indeed, choose a normalised $\alpha$-triple $(E_\alpha,F_\alpha,H_\alpha)$. Since $[E_\alpha,Z]=0$ for $Z\in\g_{\max}$, the $\mathfrak{sl}_2$ identity gives
\[
 [E_\alpha,[F_\alpha,Z]]=\beta_{\max}(H_\alpha)Z,
 \qquad \beta_{\max}(H_\alpha)>0.
\]
Thus a negative simple root space belongs to the stabiliser precisely when its root is outside $\Theta_{\mathrm{ad}}$. The stabiliser, which contains the minimal parabolic, is the parabolic $P_{\Theta_{\mathrm{ad}}}$; its Levi has simple roots $\Delta\setminus\Theta_{\mathrm{ad}}$. In reduced $A_1$ the lowering step is $[F_\alpha,E_\alpha]=-H_\alpha\in\g_0$, so $\Theta_{\mathrm{ad}}=\Delta$ and the stabiliser is $B$, as it should be. Hence the orbit map
\[
\iota_{\mathrm{ad}}:\ G/P_{\Theta_{\mathrm{ad}}}\longrightarrow\mathrm{Gr}_{d_{\max}}(\g),\qquad gP\mapsto\mathrm{Ad}(g)\g_{\max},
\]
is a closed equivariant embedding. When $d_{\max}=1$ it takes values in $\PP(\g)$: this is the \emph{adjoint variety} $X_{\mathrm{ad}}$ (the projectivised minimal nilpotent orbit). When $d_{\max}>1$ the points are projective subspaces and $X_{\mathrm{ad}}$ is a family of projective subspaces of dimension $d_{\max}-1$.

\begin{theorem}[The adjoint action and the adjoint variety]\label{thm:adjvar}
Let $\g$ be simple and let $(g_m)$ have Cartan data $(\Theta,k,k',b)$.
\begin{enumerate}[label=\textup{(\alph*)},nosep]
\item $\Theta\supseteq\Theta_{\mathrm{ad}}$ iff $\g^{\mathrm{top}}_\Theta=\g_{\max}$. In that case the first pseudo-projective limit of $\mathrm{Ad}(g_m)$ is
$\tau^+=[\mathrm{Ad}(k)\,\Pi_{\g_{\max}}\mathrm{Ad}(k')]$ of rank $d_{\max}$, with image $\PP(\mathrm{Ad}(k)\g_{\max})=\PP\bigl(\iota_{\mathrm{ad}}(kP_{\Theta_{\mathrm{ad}}})\bigr)$ and kernel the projective subspace $K=\{[X]\in\PP(\g):\ B(X,\mathrm{Ad}(k')^{-1}\g_{-\max})=0\}$ of codimension $d_{\max}$. For $d_{\max}=1$,
$\sigma_1/\sigma_2(\mathrm{Ad}\,g_m)=e^{\min_{\alpha\in\Theta_{\mathrm{ad}}}\alpha(\mu_m)}$ up to bounded factors.
\item If $\Gamma<G$ is discrete and $P_{\Theta_{\mathrm{ad}}}$-divergent, then $\mathrm{Ad}(\Gamma)$ is discrete and
\[
\Lambda_{\Myr}(\mathrm{Ad}\,\Gamma)=\bigcup_{b\in\Lambda_{\Theta_{\mathrm{ad}}}(\Gamma)}K_b,\qquad K_b=\{[X]\in\PP(\g):\ B(X,\g^{\,b}_{-\max})=0\},
\]
where, for $b=hP^-\!\in G/P^-_{\Theta_{\mathrm{ad}}}$, $\g^{\,b}_{-\max}=\mathrm{Ad}(h)\g_{-\max}$ is the lowest root space at $b$. The union of the images of the first pseudo-projective limits is
$\bigcup_{x\in\Lambda_{\Theta_{\mathrm{ad}}}(\Gamma)}\iota_{\mathrm{ad}}(x)$, identifying each Grassmannian point with its projective subspace. When $d_{\max}=1$ this is $\Lambda_{\CoG}(\mathrm{Ad}\,\Gamma)$; when $d_{\max}>1$ these limits have rank $d_{\max}$ and the rank-one Conze--Guivarc'h set is empty.
\item If moreover $d_{\max}=1$, then
\[
\Eq(\mathrm{Ad}\,\Gamma)\cap X_{\mathrm{ad}}=\{x\in X_{\mathrm{ad}}:\ x\ \text{is opposite to every }b\in\Lambda_{\Theta_{\mathrm{ad}}}(\Gamma)\},
\]
where $x=\mathrm{Ad}(g)[e_{\max}]$ is opposite to $b$ iff $B(\mathrm{Ad}(g)e_{\max},\g^{\,b}_{-\max})\ne0$, i.e.\ iff $(x,b)$ lies in the open $G$-orbit of $G/P_{\Theta_{\mathrm{ad}}}\times G/P_{\Theta_{\mathrm{ad}}}^-$.
\end{enumerate}
\end{theorem}
\begin{proof}
(a) Every weight other than $\beta_{\max}$ is reached from the highest restricted root space by lowering operators at simple restricted roots; the zero-weight space is included. The negative restricted root algebra is generated by the simple negative restricted root spaces; in a nonreduced system this includes the generation of the double-root spaces by brackets. Since the positive root algebra annihilates $\g_{\max}$ and $\g_0$ preserves it, the triangular decomposition and PBW express its generated submodule as the span of such lowering words. This nonzero submodule is all of the adjoint module because $\g$ is simple. Along a nonzero lowering word, its first step has weight $\beta_{\max}-\alpha\in\Phi\cup\{0\}$, hence $\alpha\in\Theta_{\mathrm{ad}}$. If the last weight has the same $\Theta$-coordinates as $\beta_{\max}$, every lowering step lies outside $\Theta$. This is impossible when $\Theta\supseteq\Theta_{\mathrm{ad}}$. Conversely, if $\alpha\in\Theta_{\mathrm{ad}}\setminus\Theta$, its nonzero first lowering step has the same $\Theta$-coordinates as $\beta_{\max}$ and lies outside $\g_{\max}$; when $\beta_{\max}=\alpha$, this step lies in $\g_0$. This proves the equivalence, including reduced $A_1$.

The formula for $\tau^+$ is Theorem~\ref{thm:adjoint-spectrum}(d): the factor $b=\exp\eta$ acts by the scalar $e^{\beta_{\max}(\eta)}$ on $\g_{\max}$, even when $d_{\max}>1$, so it disappears projectively. The kernel formula follows because $B_\vartheta(Y,Z)=-B(Y,\vartheta Z)$ and $\vartheta\g_{\max}=\g_{-\max}$. If $d_{\max}=1$, every other weight is at most $\beta_{\max}-\alpha$ for one of the allowed first steps, and each such first step occurs. Thus the second exponent is $\max_{\alpha\in\Theta_{\mathrm{ad}}}(\beta_{\max}-\alpha)(\mu_m)$, including the exponent $0$ in reduced $A_1$, and the stated ratio follows.
(b) By (a), $P_{\Theta_{\mathrm{ad}}}$-divergence means $\Theta\supseteq\Theta_{\mathrm{ad}}$ for every divergent sequence. The kernel of the first pseudo-projective limit is $K_b$ with $b=k'^{-1}P^-$ the repelling point. The identity $-w_0\beta_{\max}=\beta_{\max}$ and the invariance of $\Phi\cup\{0\}$ under $-w_0$ show that $\Theta_{\mathrm{ad}}$ is opposition invariant, so these $b$ exhaust $\Lambda_{\Theta_{\mathrm{ad}}}$. The opposite parabolic is the full stabiliser of $\g_{-\max}$, which also makes $K_b$ independent of the representative $h$. The adjoint image is discrete in $\mathrm{GL}(\g)$ because $\mathrm{Ad}:G\to\mathrm{Aut}(\g)_0$ is an isomorphism onto a closed subgroup. It lies in $\mathrm{SL}(\g)$, since its determinant character is trivial on the connected semisimple group $G$. The map $\mathrm{SL}(\g)\to\PSL(\g)$ is a proper finite covering, so the projective image is discrete as well. Proposition~\ref{prop:kernel} now gives $\Lambda_{\Myr}$ as the union of these kernels. The union is the projection of the compact incidence set $\{(b,[X]):[X]\in K_b\}$ over the compact limit set. The first-limit images are precisely the subspaces stated in (b). For $d_{\max}=1$ they give the Conze--Guivarc'h set. For $d_{\max}>1$ there are no rank-one limits.
(c) A point $x=[\mathrm{Ad}(g)e_{\max}]\in X_{\mathrm{ad}}$ lies in $K_b$ iff $B(\mathrm{Ad}(g)e_{\max},\mathrm{Ad}(h)e_{-\max})=0$. Writing $c(u)=B(\mathrm{Ad}(u)e_{\max},e_{-\max})$, this is $c(h^{-1}g)$, which is $\ne0$ exactly when $h^{-1}gP$ lies in the big cell $N^-P/P$ (for $u=n^-p$ one finds $c(u)\ne0$ because $\mathrm{Ad}(n^-)e_{\max}\in e_{\max}+\bigoplus_{\lambda<\beta_{\max}}\g_\lambda$, and for $u$ in a smaller Bruhat cell the weight of $\mathrm{Ad}(u)e_{\max}$ differs from $\beta_{\max}$ in the $\g_{\max}$-component), i.e.\ exactly when $x$ and $b$ are opposite. The equicontinuity region is the complement of $\Lambda_{\Myr}$ by Proposition~\ref{prop:kernel}.
\end{proof}
For $d_{\max}=1$, Theorem~\ref{thm:adjvar}(b),(c) is the adjoint counterpart of Proposition~\ref{prop:Tits-comparison}: images give the embedded limit flags, kernels the Myrberg wall, and their complement the transverse flags. Proposition~\ref{prop:highest-weight} in Chapter~\ref{ch:parabolic} extends both comparisons to representations whose highest restricted weight space is a line.

For the restriction comparison we take the highest restricted root space to be a line, so that the adjoint variety lies in the adjoint projective space.

\begin{corollary}[A restriction principle]\label{cor:restriction}
Let a countable group $\Gamma$ act by homeomorphisms on a compact metrisable space $X$, and let $Y\subset X$ be a closed $\Gamma$-invariant subset. Then $\D_\Gamma^Y\subseteq\D^X_\Gamma\cap(Y\times Y)$; hence $\mathscr R^Y_\Gamma\subseteq\mathscr R^X_\Gamma\cap Y$, every domain of proper discontinuity of $X$ restricts to one in $Y$, and $\Lambda^Y_{\mathrm{exc}}\subseteq\Lambda^X_{\mathrm{exc}}\cap Y$. In particular
\[
\Omega_{\Myr}(\mathrm{Ad}\,\Gamma)\cap X_{\mathrm{ad}}\ \subseteq\ \Omega_{\Myr}^{X_{\mathrm{ad}}}(\Gamma).
\]
\end{corollary}
\begin{proof}
For each sequence $(\gamma_m)$ of pairwise distinct elements and each $y\in Y$, allowing sources in $X$ instead of $Y$ gives
\[
 \D^Y_{(\gamma_m)}(y)\subseteq\D^X_{(\gamma_m)}(y).
\]
Taking unions over sequences gives the inclusion of group relations, hence the recurrence and properness assertions by Theorem~\ref{thm:relational}. The exceptional-source comparison is made before taking this union: if $y$ has two distinct destinations along one such sequence in $Y$, the same two points are destinations along that sequence in $X$. Thus $E^Y_{(\gamma_m)}\subseteq E^X_{(\gamma_m)}\cap Y$, and taking unions gives the asserted inclusion for $\Lambda_{\mathrm{exc}}$. Taking closures gives the Myrberg-domain inclusion.
\end{proof}
For the adjoint embedding by the highest restricted root line, the exceptional-source inclusion is an equality for each extracted sequence. In Cartan coordinates the kernel of the first ambient limit is the sum of weight spaces below the maximal $A_\Theta$-weight. The highest root line is outside that kernel, and its $P_\Theta^-$-orbit stays outside: $L_\Theta$ preserves the top weight space and $N_\Theta^-$ acts trivially on the top quotient. This orbit is the open source cell, which is exactly the singleton-fibre locus in $X_{\mathrm{ad}}$ by the dimension formula of Lemma~\ref{lem:cells} and Theorem~\ref{thm:dest}. Conversely a point outside the ambient kernel has one ambient destination by Proposition~\ref{prop:kernel}, so its nonempty intrinsic fibre is also a singleton. The union of ambient kernels is closed by Proposition~\ref{prop:kernel}. Taking unions therefore gives $\Lambda_{\Myr}^{X_{\mathrm{ad}}}(\Gamma)=\Lambda_{\Myr}(\mathrm{Ad}\,\Gamma)\cap X_{\mathrm{ad}}$.

The relation itself can still have a strict inclusion. In $\PSL(3,\C)$ with $g_m=\diag(e^{4m},e^m,e^{-5m})$, the ambient destination of $[E_{12}]$ contains $[E_{23}]$: its adjoint weight $6m$ is greater than the source weight $3m$. On the incidence flag variety, the point component of every source tending to $[E_{12}]$ tends to $[e_1]$ and its image tends to $[e_1]$, whereas the point component of $[E_{23}]$ is $[e_2]$. Under the hypotheses of Theorem~\ref{thm:adjvar}(c), the ambient equicontinuity region meets $X_{\mathrm{ad}}$ in the set of points opposite to \emph{all} of $\Lambda_{\Theta_{\mathrm{ad}}}$. This trace can be empty even when $X_{\mathrm{ad}}$ has a nonempty balanced domain, as in the full-flag example of Section~\ref{sec:lost}.

%--------------------------------------------------------------
\section{The critical wall of the adjoint action}\label{sec:adcrit}
The spectrum of $\mathrm{Ad}$ is symmetric about $0$, and the bounded exponents form the single middle block $\g^0_\Theta$. The \emph{critical level} of $\mathrm{Ad}(g_m)$ in $\PSL(\g)$ is the first space of the target flag that reaches the middle dimension. This cut lies at the end of the middle block. For a real adjoint module we use the paired critical system of Section~\ref{sec:critsys} with real projective spaces and real Grassmannians, taking every dimension over $\R$. The construction and its cover proof use the singular-value block formula and compact incidence, which hold over both fields.

\begin{theorem}[The critical wall of $\mathrm{Ad}$]\label{thm:adcrit}
Let $(g_m)$ have Cartan data $(\Theta,k,k',b)$, and put $N=\dim\g$ and $\ell=\dim\g^0_\Theta$. The first space of the target flag of $\mathrm{Ad}(g_m)$ of dimension $\ge\lceil N/2\rceil$, for any further extraction giving the full adjoint flag, is
\[
\mathrm{Ad}(k)\,\mathfrak p_\Theta=\mathrm{Lie}\bigl(kP_\Theta k^{-1}\bigr),\qquad\mathfrak p_\Theta=\g^0_\Theta\oplus\mathfrak n_\Theta,
\]
the Lie algebra of the parabolic subgroup fixing the attracting point $x^+=kP_\Theta\in G/P_\Theta$; its dimension is $(N+\ell)/2$. For the inverse sequence the corresponding space is $\mathrm{Ad}(k')^{-1}\mathfrak p_\Theta^-=\mathrm{Lie}(\mathrm{Stab}\,x^-)$, $x^-$ the repelling point. Consequently, let $\Gamma<G$ be infinite and discrete. If every Cartan extraction of a divergent sequence in $\Gamma$ has the same type $\Theta=\iota\Theta$ (for example, if $\Gamma$ is $P_\Delta$-divergent and $\Theta=\Delta$), then
\[
\Lambda_{\mathrm{crit}}(\mathrm{Ad}\,\Gamma)=\bigcup_{x\in\Lambda_\Theta(\Gamma)}\PP(\mathfrak p_x),\qquad\mathfrak p_x=\mathrm{Lie}(\mathrm{Stab}_G(x)),
\]
the union of the projectivised parabolic subalgebras of the limit points; equivalently $\Omega_{\mathrm{crit}}(\mathrm{Ad}\,\Gamma)=\{[X]\in\PP(\g):\ X\notin\mathfrak p_x\ \ \forall x\in\Lambda_\Theta\}$.
\end{theorem}
\begin{proof}
Applying Theorem~\ref{thm:adjoint-spectrum}(a) factor by factor, the positive blocks fill $\mathfrak n_\Theta$, followed by the middle block $\g^0_\Theta$ and the negative blocks in $\mathfrak n^-_\Theta$. The two nilradicals have dimension $(N-\ell)/2<N/2$, while $\dim\mathfrak p_\Theta=(N+\ell)/2\ge\lceil N/2\rceil$. Hence the first target space crossing the middle dimension is $\mathrm{Ad}(k)\mathfrak p_\Theta$, the Lie algebra of the stabiliser $kP_\Theta k^{-1}$. The gap after the middle block diverges, and the corresponding target singular subspace converges along the original sequence because $k_m\to k$. Inverting negates the exponents and reverses the blocks, giving $\mathrm{Ad}(k')^{-1}\mathfrak p_\Theta^-$. Under the common-type hypothesis, the forward and inverse middle crossings are proper, since the positive and negative blocks are both nonempty. Thus every extracted orientation is admissible and the common cut of Section~\ref{sec:critsys} has vector dimension $p_\Gamma=(N+\ell)/2$. It selects precisely the parabolic subalgebra at each attractor. The complementary cut in the inverse flag has dimension $N-p_\Gamma=(N-\ell)/2$ and selects the nilradical inside the parabolic subalgebra at the repeller. Each repeller also occurs as an attractor by inversion, so this complementary family adds no points to the union of projectivised parabolic subalgebras. The attractors exhaust $\Lambda_\Theta$, and compactness of $\Lambda_\Theta$ and continuity of $x\mapsto\mathfrak p_x$ make that union closed.
\end{proof}

\begin{example}[$\PSL(2,\C)$]\label{ex:adcritsl2}
$\g=\mathfrak{sl}_2$, $N=3$, $\ell=1$: the critical space is $\mathfrak p_x$ = the Borel subalgebra at $x$, a line of $\PP^2$, the tangent line at $x$ to the conic $X_{\mathrm{ad}}$. So $\Lambda_{\mathrm{crit}}(\mathrm{Ad}\,\Gamma)=\bigcup_{x\in\Lambda}T_x=\Lambda_{\Myr}$, as in Example~\ref{ex:adsl2} and Theorem~\ref{thm:C} with $N=2$.
\end{example}
\begin{example}[Hitchin representations in $\mathrm{SL}(3,\R)$]\label{ex:adcrithitchin}
Let $\Gamma$ be a closed-surface Hitchin image in $\mathrm{SL}(3,\R)$.
Here $\g=\mathfrak{sl}_3(\R)$, $N=8$, $\Theta=\Delta$ and
$\ell=\dim\g_0=2$. Its flag limit curve is
$\widehat C=\{(p,\ell_p):p\in\partial\Omega_c\}$, where
$\partial\Omega_c$ is strictly convex and $C^1$, and $\ell_p$ is its
tangent line (see~\cite[\S1, Theorems~1.3--1.4]{r2appLabourie2006}).
Thus
\[
 W:=\Lambda_{\mathrm{crit}}(\mathrm{Ad}\,\Gamma)
   =\bigcup_{F\in\widehat C}\PP(\mathfrak b_F)
   \subset\R\PP^7,
 \qquad \dim\mathfrak b_F=5.
\]
The flag curve is rectifiable: on each compact convex graph arc
$y=f(x)$, $f'$ is bounded and monotone, so both the point and its
tangent line have bounded variation. These arcs cover the compact
curve. Smooth local trivialisations of the projective Borel bundle
then cover its restriction to $\widehat C$ by Lipschitz images of
subsets of $\R\times\R^4$. Projecting this incidence onto $W$ gives
$\dim_H W\le5$.

For the reverse inequality choose $F_0\in\widehat C$ and
$X_0\in\mathfrak b_{F_0}$ with three distinct real eigenvalues.
Near $(F_0,[X_0])$, the projection of the universal incidence
$\{(F,[X]):X\in\mathfrak b_F\}$ onto $\PP(\g)$ is a local
diffeomorphism: the labelled eigenspaces give its smooth local inverse.
Its restriction over a small arc of $\widehat C$ therefore maps a
space homeomorphic to that arc times a four-dimensional ball
homeomorphically into $W$. Hence $\dim_{\mathrm{top}}W\ge5$.
Consequently $\dim_{\mathrm{top}}W=\dim_H W=5<7$.
The incidence is compact, so $W$ is closed and
$\Omega_{\mathrm{crit}}=\R\PP^7\setminus W$ is open and dense:
it consists of the projective classes of matrices preserving no limit
flag.
\end{example}

The critical wall is therefore easy to describe in terms of limit flags. It does not follow that every matrix in one of these parabolic subalgebras returns to itself. In rank one the tangent-line picture suggests this, whereas in higher rank a matrix with simple spectrum has several invariant flags, and a return must account for them as well. The next example makes this obstruction explicit.

\begin{proposition}[An obstruction to adjoint critical recurrence]
\label{prop:adrec-obstruction}
Let $\Gamma<\operatorname{PSL}(3,\mathbb R)$ be Zariski dense and
Borel--Anosov, and let
$\Lambda\subset\operatorname{Fl}(1,2;\mathbb R^3)$ be its flag limit set.
Choose $F=(L\subset P)\in\Lambda$ and a flag
$F'=(L'\subset P')\notin\Lambda$ opposite to $F$. In the direct sum
\[
 \mathbb R^3=L\oplus(P\cap P')\oplus L'
\]
let $X$ act by $1,2,-3$, respectively. Then
\[
 [X]\in\mathbb P(\mathfrak b_F)
       \subset\Lambda_{\mathrm{crit}}(\operatorname{Ad}\Gamma),
 \qquad
 [X]\notin\mathscr R_{\operatorname{Ad}\Gamma}.
\]
The nonrecurrence assertion holds both in
$\mathbb P(\mathfrak{sl}_3(\mathbb R))$ and in
$\mathbb P(\mathfrak{sl}_3(\mathbb C))$. In particular,
Zariski density and Borel--Anosovness do not imply that the adjoint
critical wall is the recurrent set.
\end{proposition}

\begin{proof}
We use two features of the limit flags: the points of $\Lambda$ are
pairwise opposite, and every divergent sequence in $\Gamma$, after
extraction, has divergent singular-value gaps at both simple roots;
its target flag and the target flag of its inverse belong to $\Lambda$.
Both come from the Borel--Anosov hypothesis; see the flag and divergence formulation in \cite{KLP,GuichardWienhard}.

We can always choose the flag $F'$. Indeed, choose
$H\in\Lambda\setminus\{F\}$. A small nonconstant perturbation of the
plane of $H$, with its line fixed, remains opposite to $F$. It cannot
belong to $\Lambda$, since it is distinct from $H$ and is not opposite
to $H$.

The endomorphism $X$ preserves $F$, so $X\in\mathfrak b_F$.
Theorem~\ref{thm:adcrit} therefore puts $[X]$ in the critical wall.

Conjugate the chosen splitting to the coordinate splitting, so that
$X=D=\operatorname{diag}(1,2,-3)$. The argument below works over
$\mathbb K=\mathbb R$ and over $\mathbb K=\mathbb C$. Suppose that
$[D]$ were recurrent. There would be distinct $\gamma_j\in\Gamma$,
matrices $Y_j\in\mathfrak{sl}_3(\mathbb K)$ tending to $D$, and
nonzero scalars $c_j\in\mathbb K$ such that
\[
 c_j\gamma_jY_j\gamma_j^{-1}\longrightarrow D.
\]
Taking traces of the second and third powers, and using
$\operatorname{tr}(D^2)=14$ and $\operatorname{tr}(D^3)=-18$, gives
\[
 c_j^2\longrightarrow1,\qquad c_j^3\longrightarrow1,
 \qquad\text{hence }c_j\longrightarrow1.
\]
Consequently $Z_j:=\gamma_jY_j\gamma_j^{-1}\to D$.

The eigenvalues of $D$ are simple. Label the eigenvalues of $Y_j$ by
$\lambda_{i,j}\to D_{ii}$. The same labelling applies to $Z_j$, since
$Y_j$ and $Z_j$ are conjugate. Local diagonalisation at a matrix
with simple spectrum supplies $S_j,T_j\to I$ such that
\[
 Y_j=S_jD_jS_j^{-1},\qquad
 Z_j=T_jD_jT_j^{-1},\qquad
 D_j=\operatorname{diag}(\lambda_{1,j},\lambda_{2,j},\lambda_{3,j}).
\]
Thus $T_j^{-1}\gamma_jS_j$ commutes with $D_j$ and is diagonal:
\[
 \gamma_j=T_j\operatorname{diag}(d_{1,j},d_{2,j},d_{3,j})S_j^{-1}.
\]
After extraction, fix an ordering $\pi$ of the moduli $|d_{i,j}|$.
Multiplication by $S_j^{-1},T_j\to I$ changes each singular value
by a factor tending to $1$. Borel divergence therefore gives
\[
 \frac{|d_{\pi(1),j}|}{|d_{\pi(2),j}|}\longrightarrow\infty,
 \qquad
 \frac{|d_{\pi(2),j}|}{|d_{\pi(3),j}|}\longrightarrow\infty.
\]
The normalised first and second exterior powers show that the target
flag of this sequence is
\[
 E_\pi=
 (\mathbb K e_{\pi(1)}\subset
             \mathbb K e_{\pi(1)}+\mathbb K e_{\pi(2)}),
\]
and that the target flag of the inverse sequence is $E_{\pi^{\rm rev}}$.
Both are real flags, both belong to $\Lambda$, and they are opposite
and invariant under $D$.

There are exactly six $D$-invariant full flags. Among them, the
unique flag opposite to the chosen $F$ is $F'$. Since $F\in\Lambda$
and $\Lambda$ is pairwise opposite, every other $D$-invariant flag
in $\Lambda$ would have to be $F'$, which was excluded. Thus $F$
is the only $D$-invariant flag in $\Lambda$. This contradicts the
presence of the two distinct flags $E_\pi,E_{\pi^{\rm rev}}$ and
proves nonrecurrence.
\end{proof}

The critical recurrence theorem does not force a return here, because its fullness condition concerns the flag module of the ambient projective space. The critical flag lines may span an irreducible representation of the original group and still occupy only a proper part of that ambient module. We can see the difference by counting dimensions.

\begin{remark}[Fullness in the ambient flag module]
\label{rem:ad-fullness-obstruction}
For $G=\operatorname{PSL}(3,\mathbb R)$, the fullness condition
of Theorem~\ref{thm:crit-recurrence} asks for more than the adjoint
critical flags can span. Work over $\mathbb C$ and put
$V=\mathfrak{sl}_3(\mathbb C)$, $a=3$, $p=5$. The critical flag line
\[
 v_F=\bigwedge^3\mathfrak n_F\otimes
                         \bigwedge^5\mathfrak b_F
\]
is a highest-weight line of weight $4\rho=4\omega_1+4\omega_2$
for $\operatorname{SL}_3$. Its $\operatorname{SL}_3$-span is the
irreducible module of dimension
\[
 \frac{(4+1)(4+1)(4+4+2)}2=125.
\]
Zariski density implies that the critical flag lines of $\Gamma$
span this module: their nonzero span is $\Gamma$-invariant, hence
$G$-invariant, and is contained in this irreducible module.
On the other hand, $W^{\mathrm{fl}}_{3,5}(V)$ is the
$\operatorname{SL}_8$ Cartan component of highest weight
$\omega_3+\omega_5$. The Weyl dimension formula, or the partition
$(2,2,2,1,1,0,0,0)$, gives
\[
 \dim W^{\mathrm{fl}}_{3,5}(V)
 =\prod_{1\le i<j\le8}
       \frac{\lambda_i-\lambda_j+j-i}{j-i}=2352.
\]
Thus the generated critical module is strongly irreducible but fills
only a proper part of the ambient flag module. Irreducibility of this
$G$-module leaves the ambient fullness condition unfulfilled.
\end{remark}

Coalescing opposite flags also raises a question about the middle block needed for a regular seed. Their parabolic subalgebras meet in a Cartan subalgebra before coalescence; in the limit this subalgebra cannot remain complementary to the nilradical.

\begin{lemma}[Coalescing Borel flags do not give a complementary middle block]
\label{lem:ad-seed-obstruction}
Let $F_j^+,F_j^-$ be opposite full flags in $\mathbb R^3$ with
$F_j^\pm\to F$. Put
$C_j=\mathfrak b_{F_j^+}\cap\mathfrak b_{F_j^-}$.
No limit $C$ of the $2$-planes $C_j$ is a complement to
$\mathfrak n_F$ in $\mathfrak b_F$.
\end{lemma}

\begin{proof}
Every $C_j$ is a split Cartan subalgebra. If a limit $C$ were a
complement, its projection to the diagonal quotient
$\mathfrak b_F/\mathfrak n_F$ would be an isomorphism. It would
contain an element $H$ whose diagonal entries are pairwise distinct.
This upper-triangular matrix has simple real spectrum. Choose
$H_j\in C_j$ with $H_j\to H$. For large $j$, the matrices $H_j$
also have simple real spectrum, and their six invariant flags are
the six orderings of their three eigenlines. These six flags converge
to the six distinct invariant flags of $H$. Both $F_j^+$ and $F_j^-$
are invariant under $H_j$ and tend to the same flag $F$.
They must therefore have the same ordering for large $j$, which
would make them equal, contrary to opposition.
\end{proof}

\begin{remark}[Real and complex adjoint spaces]
The field matters here as well. For a real simple Lie algebra, irreducibility of its real adjoint
module does not imply irreducibility after complexification. The
complex projective action on
$\mathfrak g_{\mathbb R}\otimes_{\mathbb R}\mathbb C$ is strongly
irreducible for a real Zariski-dense subgroup exactly when this
complex Lie algebra is simple. In particular, the underlying real
Lie algebra of a complex simple Lie algebra has complexification
equal to a sum of two simple ideals. Thus the adjoint module of a complex simple group stays irreducible when used
directly over the complex numbers, whereas complexifying its
underlying real adjoint module gives the two summands just described.
\end{remark}

%--------------------------------------------------------------
\section{Rate types}\label{sec:ratetypes}
Theorem~\ref{thm:dest} says that on a flag manifold only the type $\Theta$ and the datum $(k,k',b)$ matter. The middle adjoint block and the critical wall were determined by the same data. Away from the middle, however, Remark~\ref{rem:speeds} and Example~\ref{ex:exponents} show that other representations see more: they compare the speeds of roots which all diverge. We now record these comparisons, together with the finite limits that remain when two weights stay a bounded distance apart.

\paragraph{Definition.} Let $P\subset\ak^*$ be a lattice containing the simple restricted roots and the restricted weights under consideration. For a representation which descends only projectively to $G$, take a linear lift on a finite covering group and enlarge $P$ by finite index when necessary. Let $(g_m)$ be a divergent sequence with Cartan projection $\mu_m\in\ak^+$. Since $P$ is countable and the compact factors lie in $K$, one common extraction gives convergence of $k_m,k'_m$ and of $\nu(\mu_m)$ in $[-\infty,+\infty]$ for every $\nu\in P$. Put
\[
C_0=\{\nu\in P:\ \nu(\mu_m)\text{ converges in }\R\},\qquad \nu>0\ \text{ if }\nu(\mu_m)\to+\infty.
\]
The \emph{rate type} of $(g_m)$ is the triple $(C_0,\le,\chi)$, where $\le$ is the order on $P/C_0$ defined by $\nu\le\nu'$ iff $(\nu'-\nu)(\mu_m)\to+\infty$ or $\nu'-\nu\in C_0$, and $\chi:C_0\to\R$ is the limit, $\chi(\nu)=\lim\nu(\mu_m)$.

These conventions are compatible with finite-index enlargements of the lattice. If $P\subset P'$, the bounded subgroup becomes
\[
 C'_0=\{\nu\in P':\ q\nu\in C_0\text{ for some integer }q\ge1\},
 \qquad \chi'(\nu)=\frac{\chi(q\nu)}q\quad(\nu\in C'_0).
\]
The value of $\chi'$ is independent of the choice of $q$, and the order on $P'/C'_0$ is extended by positive scaling: a nonzero class is positive exactly when a positive integral multiple in $P/C_0$ is positive. Thus weights of a projective representation and their differences can be used in Theorem~\ref{thm:ratetypes} after this enlargement.

\begin{theorem}[Rate types]\label{thm:ratetypes}
\begin{enumerate}[label=\textup{(\alph*)},nosep]
\item $C_0$ is a saturated subgroup of $P$ ($q\nu\in C_0$, $q\in\mathbb Z\setminus\{0\}$ $\Rightarrow\nu\in C_0$), $\chi$ is a homomorphism, and $\le$ is a total order on $P/C_0$ compatible with addition, for which every simple root is $\ge0$.
\item $\Theta=\{\alpha\in\Delta:\alpha\notin C_0\}$. The normalised Levi datum of Lemma~\ref{lem:split} is $b=\exp\eta$, where $\eta$ is the unique element of $\operatorname{span}_{\mathbb R}\{\alpha^\vee:\alpha\in\Delta\setminus\Theta\}$ satisfying $\alpha(\eta)=\chi(\alpha)$ for every $\alpha\in\Delta\setminus\Theta$ (and $\eta=0$ when $\Theta=\Delta$). So $(\Theta,b)$ is a function of the rate type, but the rate type contains more: the order on $P/C_0$, and the bounded relations $\nu\in C_0$ among the divergent roots.
\item Every triple satisfying the conditions of \textup{(a)}, together with $C_0\ne P$ and $\chi(\alpha)\ge0$ for every $\alpha\in\Delta\cap C_0$, is the rate type of some divergent sequence in $\ak^+$ (hence in $G$). The ordered quotients $P/C_0$ are classified, up to equivalence, by finite sequences of real vectors (Robbiano's~\cite{Robbiano1985} description of orders on $\mathbb Z^n$).
\item The destination relation on any flag manifold $X_\theta$ is a function of $(\Theta,k,k',b)$ alone, hence of the rate type only through $(\Theta,b)$. For an irreducible $\varrho:G\to\PGL(V)$, the block structure is the partition of $\mathrm{Wt}(V)$ into classes modulo $C_0$, listed in decreasing order for $\le$. This partition can have a single block. For a block $B$ choose $\nu_B\in B$ and put
\[
E_B=\bigoplus_{\nu\in B}V_\nu,\qquad D_B|_{V_\nu}=e^{\chi(\nu-\nu_B)}\operatorname{Id}_{V_\nu}.
\]
The graded maps on these blocks are induced by $\varrho(k)D_B\varrho(k')$, up to one scalar on each block. Thus $\D$ on $\PP(V)$ is determined by the ordered partition, the values of $\chi$ on the bounded differences of weights, and the limiting compact factors. In particular, different values of $\chi$ can give different projective destination relations even when the ordered partition and $(\Theta,k,k',b)$ agree. With the limiting compact factors fixed, the projective destination relation recovers the divergence signs and the finite limits of the differences in $\mathrm{Wt}(V)-\mathrm{Wt}(V)$.
\end{enumerate}
\end{theorem}
\begin{proof}
(a) A sum of convergent sequences converges, so $C_0$ is a subgroup and $\chi$ additive; if $q\nu\in C_0$ for a fixed integer $q\ne0$, then $\nu(\mu_m)=(q\nu)(\mu_m)/q$ converges as $m\to\infty$, so $\nu\in C_0$. If $\nu\notin C_0$ then $\nu(\mu_m)\to\pm\infty$ by the extraction, so $\le$ is total, and compatibility with addition is clear ($\nu\le\nu'$, $\nu''$ arbitrary: the differences do not change). Simple roots are $\ge0$ because $\mu_m\in\ak^+$. (b) $\alpha(\mu_m)\ge0$ is bounded or tends to $+\infty$ by definition of $\Theta$ in Section~\ref{sec:thmA}; for $\alpha\notin\Theta$, $\alpha(\mu_m)=\alpha(\eta_m)\to\chi(\alpha)$. The Cartan matrix of the Levi root subsystem is invertible, so these values determine a unique $\eta$ in the stated coroot span. This is the limit selected in Lemma~\ref{lem:split}, and $b=\exp\eta$. (c) By the Hahn--Robbiano embedding the ordered group $P/C_0$ (torsion-free, finitely generated) embeds order-preservingly in $(\R^s,\mathrm{lex})$, $\nu\mapsto(\nu(u_1),\dots,\nu(u_s))$ with $u_j\in\ak$. Since $C_0$ is saturated in the finite-rank lattice $P$, a $\mathbb Z$-basis of $C_0$ extends to one of $P$. Extend $\chi$ linearly to $P\otimes_{\mathbb Z}\mathbb R=\ak^*$, and let $\eta_0\in\ak$ be the corresponding vector, so $\nu(\eta_0)=\chi(\nu)$ for all $\nu\in C_0$. This is a new choice, which need not equal the normalised Levi vector $\eta$ of~(b). Put $\mu_m=\sum_jm^{s+1-j}u_j+\eta_0$. Then $\nu(\mu_m)\to+\infty$ exactly when the first non-zero coordinate of $\nu$ is positive. For a simple root outside $C_0$, the first non-zero coordinate is positive. For $\alpha\in\Delta\cap C_0$ all those coordinates vanish and $\alpha(\mu_m)=\chi(\alpha)\ge0$. Thus $\mu_m\in\ak^+$ for large $m$; since $C_0\ne P$, the sequence is divergent. (d) The first sentence is Theorem~\ref{thm:dest}. We have $\varrho(g_m)=\varrho(k_m)\varrho(e^{\mu_m})\varrho(k'_m)$, and $\varrho(e^{\mu_m})$ acts on $V_\nu$ by $e^{\nu(\mu_m)}$. Two weights belong to the same block iff $\nu-\nu'\in C_0$, and the blocks are ordered by $\le$. On the block $E_B$, division by $e^{\nu_B(\mu_m)}$ gives the limit $D_B$, because $e^{(\nu-\nu_B)(\mu_m)}\to e^{\chi(\nu-\nu_B)}$. With more than one block, Proposition~\ref{prop:graded} gives the destination relation over $\R$ or $\C$. Conversely, the fibre dimensions and the quotient maps recover the ordered blocks and the projective graded maps, as in the proof of Theorem~\ref{thm:graphcomp}. With one block, normalisation by $e^{\nu_B(\mu_m)}$ gives an invertible limit $\varrho(k)D_B\varrho(k')$, and the destination relation is its ordinary projective graph. In either case the order recovers the divergence signs. With the compact factors fixed, the positive diagonal ratios within each block recover the finite limits of weight differences.
\end{proof}
\begin{example}[$\PSL(3)$]
Let $G=\PSL(3,\R)$ and $\Theta=\Delta=\{\alpha_1,\alpha_2\}$. The rate type decides whether $\alpha_1(\mu_m)-\alpha_2(\mu_m)$ stays bounded ($\alpha_1-\alpha_2\in C_0$, as for the principal $\PSL(2)$), tends to $+\infty$, or to $-\infty$. With the limiting compact factors $k,k'$ fixed, the destinations on $\R\PP^2$ and on the other flag manifolds are the same in the three cases. On $V=\R^3$ with weights $e_1,e_2,e_3$ the differences are $\alpha_1,\alpha_2,\alpha_1+\alpha_2$ and $\alpha_1-\alpha_2$ is \emph{not} a difference of two weights: $\R^3$ does not see the question. On $V=\g=\mathfrak{sl}_3$ (weights the roots and $0$) the difference $\alpha_1-\alpha_2$ is the difference of the weights $\alpha_1,\alpha_2$: the block structure of $\mathrm{Ad}(g_m)$ separates the singular values $e^{\alpha_1(\mu_m)}$ and $e^{\alpha_2(\mu_m)}$ or groups them in one block, according to the rate type. This is the content of Remark~\ref{rem:speeds}.
\end{example}
\begin{remark}
The ``full set of rate types'' is the set of triples of Theorem~\ref{thm:ratetypes}; for each $\Theta$ with $|\Theta|\ge2$, choose $u\in\mathfrak a_\Theta$ with positive rational values on $\Theta$. Varying the ratio of two of these values in $\mu_m=m\,u$ gives infinitely many bounded subgroups $C_0=\{\nu\in P:\nu(u)=0\}$. For $|\Theta|=1$ the bounded subgroup and the order are determined by $\Theta$, leaving only $\chi$ to specify. A representation $V$ resolves the part of the order lying in $\mathrm{Wt}(V)-\mathrm{Wt}(V)$; the family $V_{N\lambda}$, for a compatible regular dominant complex highest weight and $N\to\infty$, resolves the order on the restricted root lattice. Every restricted root-lattice vector is the restriction of a difference of two nonnegative sums of complex simple roots. For large $N$ each of the corresponding two weight subtractions from the complex highest weight $N\lambda$ occurs, by the highest-weight presentation of $V_{N\lambda}$. Positive integral multiples then determine the order on $P$. Two questions remain: which rate types matter for recurrence (the critical type $(a,p)$ of Chapter~\ref{ch:psl} is such a case, in $\PSL(n+1,\C)$), and how to express the Cartan module $W^{\mathrm{fl}}_{a,p}$ in the language of $P/C_0$ for general $G$.
\end{remark}

\begin{example}[A bounded weight difference remains visible]\label{ex:bounded-rate-difference}
Even when the ordered blocks stay fixed, their graded maps may turn. In $\PSL(3,\R)$ let
\[
 \mu_m^{(c)}=\diag(m+2c/3,-c/3,-m-c/3),\qquad g_m^{(c)}=e^{\mu_m^{(c)}}.
\]
For every fixed real $c$, $\Theta=\Delta$, $k=k'=1$ and the Levi datum in the original group is trivial. The ordered adjoint weight blocks are the same, but
$\chi(\alpha_1-\alpha_2)=c$. On $\PP(\mathfrak{sl}_3)$ the block formula gives
\[
 \D_{(\Ad g_m^{(c)})}([E_{12}+E_{23}])
 =\PP\operatorname{Span}(E_{13},e^cE_{12}+E_{23}).
\]
Indeed, the two coefficients in the common block have a fixed ratio $e^c$, the faster $E_{13}$ coefficient may be varied by a vanishing perturbation of the source, and all slower coefficients vanish after normalisation at their common rate. Thus finite weight differences can turn a destination even when the ordered partition and the intrinsic flag-manifold data are fixed.
\end{example}

\section{What the adjoint action adds}\label{sec:adadd}
The adjoint action has given us a \emph{canonical} projective space, with no choice of representation. Its ordered blocks reflect relative growth (Remark~\ref{rem:speeds}) and the interleaving of factors in the semisimple case (Remark~\ref{rem:semisimple}). At the middle dimension it recovers a union of Lie algebras of parabolic subgroups of the limit flags (Theorem~\ref{thm:adcrit}). This is the critical wall; its geometric simplicity does not make every point recurrent, even for Zariski-dense Borel--Anosov groups (Proposition~\ref{prop:adrec-obstruction}).

For a simple noncompact group, the same picture contains $X_{\mathrm{ad}}=G/P_{\Theta_{\mathrm{ad}}}$, the closed orbit of the highest restricted root space. When that space is a line, this is the adjoint variety in $\PP(\g)$, and Corollary~\ref{cor:restriction} compares its domains with those of the ambient projective space. For higher multiplicity the orbit parametrises projective subspaces in a Grassmannian.

\chapter{Parabolic subgroups, flag varieties and the adjoint action}
\label{ch:parabolic}

The middle level of the adjoint target flag has recovered a parabolic subalgebra from a divergent sequence. We now turn this around: a flag can itself be represented by the Lie algebra of its stabiliser. Let $G$ be semisimple with trivial centre and $P\subset G$ a parabolic subgroup. The quotient $G/P$ is a compact variety (a complex projective variety when $G$ is complex), $G$ acts on it, and a discrete $\Gamma<G$ gives an action on $G/P$ exactly as a subgroup of $\PSL(n+1,\C)$ acts on $\C\PP^n$. We shall represent $G/P$ by projective constructions arising from the adjoint representation, and ask what happens to its destinations and limit sets when the sources are allowed to leave the flag variety.

\medskip  Since $G$ has trivial centre, it acts on $G/P$ through $\mathrm{Ad}:G\to\mathrm{Aut}(\g)$, and a parabolic subgroup is the normaliser of its Lie algebra. Hence
\[
j_P:\ G/P\longrightarrow\mathrm{Gr}_{\dim\mathfrak p}(\g),\qquad gP\longmapsto\mathrm{Ad}(g)\mathfrak p,
\]
is an equivariant embedding of $G/P$ into a Grassmannian of $\g$. The Killing form recovers the nilradical by $\mathfrak p^{\perp_B}=\mathfrak n_P$: it pairs opposite root spaces, so the orthogonal complement consists exactly of the positive roots outside the Levi. Thus one first takes $W\mapsto W^{\perp_B}$ and then the Pl\"ucker embedding of the resulting $\dim\mathfrak n_P$-plane. Their composite gives an equivariant embedding of $G/P$ in the projective space of an irreducible representation $V_P\subset\wedge^{\dim\mathfrak n_P}\g$ (the anticanonical embedding). So the dynamics on $G/P$ can be compared with the dynamics of $\mathrm{Ad}\,\Gamma$ on projective spaces, where Chapter~\ref{ch:psl} applies.

\medskip We make the comparison in three steps, beginning with the limit flags themselves:

\emph{Theorem~\ref{thm:Quint-adjoint}.} If $\Gamma$ is infinite and discrete and every Cartan extraction of a divergent sequence in $\Gamma$ has the same type $\Theta=\iota\Theta$, then $j_{P_\Theta}(\Lambda_\Theta(\Gamma))$ is the family of critical target subspaces of $\mathrm{Ad}\,\Gamma$ in $\mathrm{Gr}(\g)$. Quint's limit flags parametrise this family, and the union of its projectivisations is the critical wall in $\PP(\g)$.

\emph{Theorems~\ref{thm:traces} and~\ref{thm:trace-rank1}.} In the rank-one comparison, opposition computes the ambient Myrberg trace. For Kulkarni's set we must also compare the allowed sources, and their difference gives an explicit defect. In rank one this difference disappears for highest-weight embeddings with a highest restricted weight line, including the adjoint embedding when its highest restricted root space is a line.

\emph{Theorems~\ref{thm:forget} and~\ref{thm:embed}.} Forgetting flags commutes with destinations. Embedding a flag manifold in projective space allows more moving sources; the ambient equicontinuity region cuts out the transverse domain on the flag manifold.

\medskip
We write $G$ for a real semisimple group as in Chapter~\ref{ch:cartan}. When $G$ is complex, we use the holomorphic adjoint action, holomorphic representations, and complex lines and dimensions. For a real form, the adjoint model uses real lines and dimensions, and intrinsic witnesses stay in the real flag orbit in $G_\C/P_\C$. In a highest-weight embedding, the ambient projective space is over the field of the representation. For a complex representation of a real form, ambient witnesses can move throughout that complex projective space.

%--------------------------------------------------------------
\section{Parabolic subgroups and parabolic subalgebras}\label{sec:parabolics}
We use the notation of Section~\ref{sec:setting-sum}: $P_\Theta=L_\Theta N_\Theta$ and $X_\theta=G/P_\theta$. The parabolic subgroups containing $B=P_\Delta$ correspond to subsets $\theta\subseteq\Delta$ ($P_\Theta\subseteq P_\theta$ iff $\Theta\supseteq\theta$), and the natural projections are
\[
\pi_{\Theta\theta}:\ G/P_\Theta\longrightarrow G/P_\theta\qquad(\Theta\supseteq\theta),
\]
fibre bundles with compact fibre $P_\theta/P_\Theta$, a flag manifold of the Levi factor of $P_\theta$.

\begin{example}[$\PSL(n+1,\C)$]\label{ex:parabolics-psl}
A parabolic subgroup is the stabiliser of a flag $V_{i_1}\subset\dots\subset V_{i_s}$ and $\theta=\{\alpha_{i_1},\dots,\alpha_{i_s}\}$. For $\theta=\{\alpha_1\}$ one gets $\C\PP^n$, for $\{\alpha_n\}$ the dual projective space, for $\{\alpha_i\}$ the Grassmannian $\mathrm{Gr}_i(\C^{n+1})$, for $\Delta$ the complete flags. The projection $\mathrm{Fl}(\C^{n+1})\to\C\PP^n$ forgets everything but the line.
\end{example}

\paragraph{Groups inside parabolic subgroups.} A subgroup $\Gamma<P_\theta$ fixes the point $o=eP_\theta\in X_\theta$. This is the \emph{reducible} case, and it is where the failure of a largest domain is typical.
\begin{proposition}\label{prop:reducible}
Let $\Gamma<G$ be infinite and discrete. (a) Every point of $X_\theta$ with infinite stabiliser in $\Gamma$ is recurrent; hence $L_0\subseteq\mathscr R_\Gamma$ (notation of Section~\ref{sec:KulMyr}). (b) If $\Gamma<P_\theta$ then $o\in\mathscr R_\Gamma$, and so is every point with a finite $\Gamma$-orbit.
\end{proposition}
\begin{proof}
(a) If $\gamma_m\in\mathrm{Stab}_\Gamma(x)$ are distinct, $x_m=x$ gives $\gamma_mx_m=x$, so $x\in\D_\Gamma(x)$. (b) $\Gamma$ is infinite and fixes $o$, so the point has an infinite stabiliser; a point with finite orbit has a finite-index subgroup of $\Gamma$ in its stabiliser.
\end{proof}
For example, the toral group $\Gamma_M<\PSL(3,\C)$ of Section~\ref{sec:toral} lies in the stabiliser of the line at infinity, every point of which has an infinite stabiliser (the translation subgroup), so $\ell_\infty\subset\mathscr R_\Gamma$: this is the mechanism behind the recurrent set found there.

%--------------------------------------------------------------
\section{The action factors through the adjoint representation}\label{sec:factors-ad}
Recall that $G=\mathrm{Aut}(\g)_0$, so $\mathrm{Ad}:G\to\mathrm{GL}(\g)$ is the inclusion (Section~\ref{sec:setting-sum}); any semisimple group acts on its flag varieties through $G/Z(G)=\mathrm{Ad}(G)$. The adjoint representation does more than factor the action: it contains a canonical copy of $G/P$.

\begin{proposition}[Canonical embeddings]\label{prop:canonical}
Let $P$ be a parabolic subgroup of $G$ with Lie algebra $\mathfrak p$, nilradical $\mathfrak n$, $d=\dim\mathfrak n$, $k=\dim\mathfrak p$.
\begin{enumerate}[label=\textup{(\alph*)},nosep]
\item $P=N_G(\mathfrak p)$ (a parabolic subgroup is self-normalising). Hence $j_P:G/P\to\mathrm{Gr}_k(\g)$, $gP\mapsto\mathrm{Ad}(g)\mathfrak p$, is a $G$-equivariant embedding onto the orbit of $\mathfrak p$, which is closed because $G/P$ is compact.
\item The line $\wedge^d\mathfrak n\subset\wedge^d\g$ is stabilised by $P$, so $\kappa_P:G/P\to\PP(\wedge^d\g)$, $gP\mapsto[\wedge^d\mathrm{Ad}(g)\mathfrak n]$ is an equivariant map; it is an embedding, and its image spans an irreducible subrepresentation $V_P\subseteq\wedge^d\g$ with highest weight $\chi_P=\sum_{\lambda\in\Phi(\mathfrak n)}m_\lambda\lambda$, the sum of the roots of $\mathfrak n$ with multiplicities over the field of the representation. (This is the anticanonical embedding: the dual of $\wedge^d\mathfrak n$ is the fibre of $K_{G/P}^{-1}$ at the base point.)
\item The support of $\chi_P$ is exactly $\theta$, where $P=P_\theta$ (in the notation of Section~\ref{sec:setting-sum}).
\end{enumerate}
\end{proposition}
\begin{proof}
(a) By Chevalley's theorem a parabolic subgroup is the normaliser of its Lie algebra (equivalently, the stabiliser of the corresponding flag in any faithful representation). The orbit map $G/P\to\mathrm{Gr}_k(\g)$ is injective and an immersion (its differential at the base point is $\g/\mathfrak p\to\mathrm{Hom}(\mathfrak p,\g/\mathfrak p)$, $X\mapsto(Y\mapsto[X,Y])$, whose kernel is the Lie algebra normaliser modulo $\mathfrak p$, hence zero since $N_\g(\mathfrak p)=\mathfrak p$); a compact injective immersion is an embedding.
(b) The line $\wedge^d\mathfrak n$ is extremal. Choose a chamber element central in the Levi and positive on the roots of $\mathfrak n$. Its positive eigenspaces in $\mathfrak g$ are precisely $\mathfrak n$. Since their dimensions add up to $d$, their top exterior power is the unique maximal eigenline in $\wedge^d\g$ for this element. It is $P_{\min}$-stable, so it is a highest weight line; in the completely reducible representation $\wedge^d\g$ a highest weight vector with one-dimensional weight space generates an irreducible submodule. The embedding assertion follows from (a), since $\mathfrak p=\mathfrak n^{\perp_B}$ and the Killing-orthogonal map and Pl\"ucker map are embeddings. The tangent space at the base point is $\g/\mathfrak p\cong\mathfrak n^*$ as a $P$-module, by the Killing pairing. Hence the fibre of $K^{-1}$ is $\wedge^{\mathrm{top}}(\g/\mathfrak p)\cong(\wedge^d\mathfrak n)^*$, while the tautological line of the projective embedding is $\wedge^d\mathfrak n$. Thus $\kappa_P^*\mathcal O(1)=K^{-1}_{G/P}$; for complex $G$, its global holomorphic sections are $V_P^*$.
(c) $\chi_P$ is the determinant character of $\mathrm{Ad}|_{\mathfrak n}$, trivial on the commutator of $P$. For $\alpha\notin\theta$, $s_\alpha$ permutes the roots of $\mathfrak n$ with their multiplicities, so $\langle\chi_P,\alpha^\vee\rangle=0$. For $\alpha\in\theta$, write $\chi_P$ as the sum of all positive restricted roots with multiplicities minus the corresponding sum for the Levi. The first sum has positive pairing with $\alpha^\vee$: reflection in $\alpha$ reverses only the positive roots proportional to $\alpha$ and permutes the others. The Levi sum has nonpositive pairing, since it is a nonnegative combination of simple roots different from $\alpha$. Thus $\langle\chi_P,\alpha^\vee\rangle>0$, also for a non-reduced restricted system.
\end{proof}
\begin{example}[Projective space]\label{ex:veronese}
$G=\PSL(n+1,\C)$, $P=P_{\alpha_1}$, $G/P=\C\PP^n$. Then $K^{-1}_{\C\PP^n}=\mathcal O(n+1)$ and $\kappa_P$ is the Veronese embedding $\nu_{n+1}:\C\PP^n\to\PP(\mathrm{Sym}^{n+1}\C^{n+1})$, with $V_P=\mathrm{Sym}^{n+1}\C^{n+1}$ (the irreducible representation of highest weight $(n+1)\omega_1$). The identity representation, whose line bundle is $\mathcal O(1)$, is a different equivariant embedding of the same variety. Both give the same action on $\C\PP^n$ and the same ambient equicontinuity trace. Indeed, the symmetric power of a nonzero normalised linear limit is a nonzero limit in the symmetric-power representation. A pure power belongs to its kernel exactly when the original vector belongs to the kernel of the linear limit. Proposition~\ref{prop:kernel} therefore gives equality of the traces. This argument applies in every Veronese degree. For $n=1$, $G=\PSL(2,\C)$: $K^{-1}=\mathcal O(2)$, $V_P=\mathrm{Sym}^2=\mathfrak{sl}_2$, and $\kappa_P:\C\PP^1\to\PP(\mathfrak{sl}_2)$ is the conic of Example~\ref{ex:adsl2}: the canonical embedding is the adjoint variety.
\end{example}
\begin{remark}
For a simple group with $d_{\max}=1$, the anticanonical embedding of $G/P_{\Theta_{\mathrm{ad}}}$ agrees with its adjoint projective realisation exactly when its nilradical is a line. Indeed, $\chi_P$ contains $\beta_{\max}$ once, together with the positive roots of the remaining nilradical, so $\chi_P=\beta_{\max}$ precisely in that case. This gives the $\mathfrak{sl}_2$ conic above. When the nilradical is larger, the same abstract flag manifold has two different projective realisations: the adjoint one and the anticanonical one. For $d_{\max}>1$ the adjoint orbit lies instead in a Grassmannian; the canonical projective embeddings, for this flag manifold and the others, still live in the exterior powers of $\g$.
\end{remark}

%--------------------------------------------------------------
\section{Quint's limit set as a limit set of the adjoint action}\label{sec:Quint-adjoint}
The middle adjoint block in Cartan coordinates determines the type. The limiting parabolic subalgebra at the middle cut determines the attracting flag.

\begin{theorem}[Limit points are critical subspaces of $\mathrm{Ad}$]\label{thm:Quint-adjoint}
Let $(g_m)$ be a divergent sequence with Cartan data
$(\Theta,k,k',b)$ and put $\ell=\dim\g^0_\Theta$.
For the opposite parabolic use the orbit embedding
\[
 j_{P_\Theta^-}:G/P_\Theta^-\longrightarrow
 \mathrm{Gr}_{(\dim\g+\ell)/2}(\g),\qquad
 gP_\Theta^-\longmapsto\mathrm{Ad}(g)\mathfrak p_\Theta^-.
\]
Then:
\begin{enumerate}[label=\textup{(\alph*)},nosep]
\item the type $\Theta$ consists of the simple roots $\alpha$ for which $\g_\alpha$ is outside the middle block in Cartan coordinates, before applying the compact factors (Theorem~\ref{thm:adjoint-spectrum}(c));
\item for any further extraction giving the full adjoint flag, the first target space of $\mathrm{Ad}(g_m)$ of dimension
$\ge\lceil\dim\g/2\rceil$ is
$\mathrm{Ad}(k)\mathfrak p_\Theta=j_{P_\Theta}(x^+)$,
where $x^+=kP_\Theta\in G/P_\Theta$ is the attracting flag
(Theorem~\ref{thm:adcrit}); for the inverse sequence the corresponding
space is
\[
 \mathrm{Ad}(k')^{-1}\mathfrak p_\Theta^-
 =j_{P_\Theta^-}(x^-),\qquad
 x^-=k'^{-1}P_\Theta^-\in G/P_\Theta^-;
\]
\item consequently, if $\Gamma<G$ is infinite and discrete and every Cartan extraction of a divergent sequence in $\Gamma$ has the same type $\Theta=\iota\Theta$, then
\[
j_{P_\Theta}\bigl(\Lambda_\Theta(\Gamma)\bigr)=\Lambda^{\mathrm{Gr}}_{\mathrm{crit}}(\mathrm{Ad}\,\Gamma),
\]
the closure in $\mathrm{Gr}_{(\dim\g+\ell)/2}(\g)$ of the set of critical target spaces of divergent sequences of $\mathrm{Ad}\,\Gamma$ (here $\ell=\dim\g^0_\Theta$). The critical limit set of $\mathrm{Ad}\,\Gamma$ in $\PP(\g)$ is the union $\bigcup_{x\in\Lambda_\Theta}\PP(\mathfrak p_x)$ of the projectivisations of these subspaces.
\end{enumerate}
\end{theorem}
\begin{proof}
(a) and (b) follow from Theorems~\ref{thm:adjoint-spectrum}
and~\ref{thm:adcrit}. The two orbit maps have the indicated domains;
$N_G(\mathfrak p_\Theta)=P_\Theta$ and
$N_G(\mathfrak p_\Theta^-)=P_\Theta^-$ show that each critical space
determines its flag. For (c), Theorem~\ref{thm:Quint}(a) identifies
$\Lambda_\Theta(\Gamma)$ with the attracting flags of sequences of
type containing $\Theta$, here exactly $\Theta$. By (b), their images
under $j_{P_\Theta}$ are precisely the critical target spaces.
Compactness and the embedding property give the asserted closure
identity; the projective union is Theorem~\ref{thm:adcrit}.
\end{proof}
\begin{corollary}[Independence of the representation]\label{cor:indep}
Quint's limit set and the transverse domain $\Omega^{\mathrm{trans}}_\theta(\Gamma)$ depend only on $\mathrm{Ad}\,\Gamma<\mathrm{Aut}(\g)$ in the present centre-free setting. More generally, discrete subgroups of semisimple groups with the same adjoint group and the same adjoint image have the same destination relations and domains of proper discontinuity on corresponding flag varieties, provided the adjoint maps have finite kernels on these subgroups. Finite fibres preserve sequences of distinct group elements after extraction, and a finite kernel preserves properness; an infinite central kernel would fix every flag and need not preserve either conclusion. In $\PSL(n+1,\C)$, $\Gamma$ and its adjoint image determine each other, and $j_{P_\Theta}\Lambda_\Theta(\Gamma)$ is the Grassmannian critical limit set of that image under the hypotheses of Theorem~\ref{thm:Quint-adjoint}(c).
\end{corollary}
\begin{example}[$\PSL(2,\C)$]\label{ex:Quint-sl2}
$\g=\mathfrak{sl}_2$, $N=3$, $\ell=1$, $\Theta=\Delta$. The critical target space is the Borel subalgebra $\mathfrak b_x$ (dimension $2$), a point of $\mathrm{Gr}_2(\mathfrak{sl}_2)=\PP(\mathfrak{sl}_2^*)$, and $j(x)=\mathfrak b_x$ is the dual conic. Thus the classical limit set $\Lambda\subset\C\PP^1$ parametrises the critical target spaces of $\mathrm{Sym}^2\Gamma$ in the dual projective plane. In the primal plane $\PP(\mathfrak{sl}_2)$, their projectivisations are the tangent lines to the nilpotent conic at the corresponding points, as in Proposition~\ref{prop:sym2}. The Killing form identifies the primal conic with the dual conic.
\end{example}

%--------------------------------------------------------------
\section{Limit sets of \texorpdfstring{$G/P$}{G/P} as traces of limit sets of \texorpdfstring{$\mathrm{Ad}$}{Ad}}\label{sec:traces}
Let $X=\kappa(G/P)\subset Y=\PP(V)$ be a canonical or highest-weight embedding (Section~\ref{sec:factors-ad}), so $\Gamma$ acts on $X$ and, through $\rho$, on $Y$. By Corollary~\ref{cor:restriction} one always has $\D^X\subseteq\D^Y\cap(X\times X)$. We ask how each limit set of $X$ is related to the corresponding limit set of $Y$, for instance to the traces $\Lambda^Y\cap X$. There are two questions here: whether a destination ends on the flag variety, and whether it can be reached with sources which stay there. In rank one, opposition computes the ambient Myrberg trace for the highest-weight embeddings considered below. It agrees with the intrinsic set, and the Kulkarni traces agree as well. In higher rank the comparison of allowed sources leaves an explicit obstruction.

We begin with the limits for which the source does not move.

\begin{lemma}[Fixed-source limits from $L_0$]\label{lem:K1}
Let $Y$ be a compact space on which $\Gamma$ acts continuously.  If
$x\in L_0^Y$, every accumulation point of the fixed orbit $\Gamma x$
belongs to $L_0^Y$.
\end{lemma}
\begin{proof}
The set $L_0^Y$ is closed and $\Gamma$-invariant.  Thus
$\gamma x\in L_0^Y$ for every $\gamma\in\Gamma$, and closedness gives
the assertion.
\end{proof}

Once the sources move, $\D_\Gamma^Y(x)$ can reach much further.
For a cyclic diagonal group with strictly decreasing eigenvalue moduli,
the slowest eigenpoint belongs to $L_0$, but vanishing perturbations of
that point have positive-power destinations filling the whole projective
space. The Kulkarni wall is only the union of the two extreme
hyperplanes, as computed in Proposition~\ref{r2app:cyclic-ambient-open}.

\begin{theorem}[Traces]\label{thm:traces}
Let a countable group $\Gamma$ act by homeomorphisms on a compact
metrisable space $Y$, and let $X\subset Y$ be closed and $\Gamma$-invariant.
Write $B^Z=L_0^Z\cup L_1^Z$, for $Z=X,Y$.
\begin{enumerate}[label=\textup{(\alph*)},nosep]
\item \emph{(Destinations, Myrberg and recurrence.)}
$\D_\Gamma^X\subseteq\D_\Gamma^Y\cap(X\times X)$,
$\Lambda^X_{\mathrm{exc}}\subseteq\Lambda^Y_{\mathrm{exc}}\cap X$,
and $\mathscr R^X\subseteq\mathscr R^Y\cap X$.  Every domain of proper
discontinuity of $Y$ restricts to one of $X$.
\item \emph{(The sets $L_0$ and $L_1$.)}
$L_0^X\subseteq L_0^Y\cap X$ and
$L_1^X\subseteq B^Y\cap X$.  In particular,
$B^X\subseteq B^Y\cap X$.
\item \emph{(The set $L_2$.)} Put
\[
 E=(X\cap B^Y)\setminus B^X.
\]
Then
\begin{align*}
 \Lambda_{\Kul}^X
 &\subseteq (\Lambda_{\Kul}^Y\cap X)\cup
 \bigl(\overline{\D_\Gamma^X(E)}\cap X\bigr)\\
 &\subseteq (\Lambda_{\Kul}^Y\cap X)\cup
 \bigl(\overline{\D_\Gamma^Y(E)}\cap X\bigr).
\end{align*}
In particular, $E=\varnothing$ implies
$\Lambda_{\Kul}^X\subseteq\Lambda_{\Kul}^Y\cap X$.
\end{enumerate}
\end{theorem}
\begin{proof}
Part \textup{(a)} is Corollary~\ref{cor:restriction}.  For
\textup{(b)}, the stabilisers of a point of $X$ are the same in the two
actions.  Taking closures, first in $X$ and then in $Y$, gives the
$L_0$ inclusion.  Suppose $y=\lim\gamma_mx$ for a fixed source
$x\in X\setminus L_0^X$ and pairwise distinct $\gamma_m$.  If
$x\notin L_0^Y$, this limit belongs to $L_1^Y$ by definition.  If
$x\in L_0^Y$, Lemma~\ref{lem:K1} puts it in $L_0^Y$.  Since $B^Y$
is closed, taking the closure in the definition of $L_1^X$ proves
$L_1^X\subseteq B^Y\cap X$.

For \textup{(c)}, first take a limit from a compact source set,
before taking the closure in the definition:
$z=\lim\gamma_mk_m$, where $k_m\in K\subset X\setminus B^X$ and $K$
is compact.  After extraction $k_m\to k\in K$.  If $k\notin B^Y$,
then $z\in\D_\Gamma^Y(k)\subseteq L_2^Y$, by
Proposition~\ref{prop:L2}(a).  If $k\in B^Y$, then $k\in E$ and
$z\in\D_\Gamma^X(E)$.  Taking the closure of all these limits gives
\[
 L_2^X\subseteq L_2^Y\cup\overline{\D_\Gamma^X(E)}.
\]
Together with \textup{(b)} this proves the first inclusion.  The second
follows from \textup{(a)}.  Closedness of $X$ ensures that all intrinsic
limits remain in $X$.
\end{proof}

The sources excluded only by the ambient action split according to whether they lie in its stabiliser wall or in its fixed-orbit wall:
\[
 E=\bigl((X\cap L_0^Y)\setminus B^X\bigr)
 \ \cup\
 \bigl((X\cap L_1^Y)\setminus(L_0^Y\cup B^X)\bigr).
\]
The first part comes from infinite stabilisers, but its destinations
are still formed with moving sources: Lemma~\ref{lem:K1} controls only
the fixed orbit. In rank one both parts disappear, and we can compute
the traces directly.

\begin{theorem}[Rank one: all limit sets of $G/P$ are traces]
\label{thm:trace-rank1}
Let $G$ have real rank one, let
$X=G/P=\partial_\infty\mathbb H$, and let $\Gamma<G$ be non-elementary
and discrete.  Let $\rho$ be a nontrivial irreducible representation with a highest restricted
weight line and an equivariant embedding $\iota_\rho:X\to Y=\PP(V)$,
with discrete projective image and finite kernel on $\Gamma$.  Identifying
$X$ with its image, one has
\[
 \Lambda_{\Kul}^X=\Lambda_{\Myr}^X=\Lambda(\Gamma)
 =\Lambda_{\Kul}^Y\cap X
 =\Lambda_{\Myr}^Y\cap X
 =\Lambda_{\CoG}^Y\cap X,
 \qquad
 \Omega^{\mathrm{trans}}=X\setminus\Lambda(\Gamma).
\]
In particular, this applies to the canonical highest-weight embedding,
and to the adjoint embedding in the module of the noncompact simple
factor when its highest restricted root space is a line.
\end{theorem}
\begin{proof}
In rank one, two boundary points are opposite precisely when they are
distinct.  Proposition~\ref{prop:highest-weight} therefore gives
\[
 \Lambda_{\Myr}^Y\cap X=\Lambda(\Gamma),
 \qquad
 \Lambda_{\CoG}^Y=\iota_\rho\Lambda(\Gamma).
\]
The endpoints of loxodromic elements are dense in the classical limit
set.  Each such endpoint, viewed in $Y$, has infinite stabiliser, since
it is fixed by all powers of its loxodromic element.  Consequently
\[
 \Lambda(\Gamma)\subseteq L_0^Y\cap X
 \subseteq\Lambda_{\Kul}^Y\cap X
 \subseteq\Lambda_{\Myr}^Y\cap X=\Lambda(\Gamma).
\]
The comparison $\Lambda_{\Kul}^Y\subseteq\Lambda_{\Myr}^Y$ is
Proposition~\ref{prop:kulkarni-myrberg}.  For the intrinsic sets, the exceptional sources of an extracted sequence
are its repeller alone (Corollary~\ref{cor:rank-one}), so their union is
$\Lambda(\Gamma)$. The loxodromic endpoints also give
$\Lambda(\Gamma)\subseteq L_0^X$. The rank-one relation in
Proposition~\ref{prop:D-rank1} gives $\mathscr R_\Gamma^X=\Lambda(\Gamma)$,
so $L_0^X=\Lambda(\Gamma)$. Every destination of a source outside
$\Lambda(\Gamma)$ lies in $\Lambda(\Gamma)$. This puts $L_1^X$ there
and, by Proposition~\ref{prop:L2}(a), puts $L_2^X$ there as well.
This proves the displayed chain.  It also shows
$X\cap B^Y\subseteq\Lambda(\Gamma)=B^X$, so the source defect $E$
is empty.
\end{proof}
When the highest restricted root space is a line, the same adjoint trace
equalities hold in the full adjoint module even with compact factors. Projection to the noncompact factor is proper and
thus gives a discrete image of $\Gamma$ with finite fibres. Along an
escaping sequence, the adjoint operators on the compact summands stay
bounded and vanish after normalisation at the highest restricted root
scale. The rank-one image and kernel trace are therefore the same as
for the noncompact factor. Loxodromic endpoints still have infinite
stabilisers in the full projective action, so the preceding inclusions
give the Kulkarni trace there as well.

\begin{example}[The trace on a rational normal curve]
For $G=\PSL(2,\C)$, take a non-elementary discrete subgroup with limit set $\Lambda\ne\C\PP^1$. Put $Y=\PP(\mathrm{Sym}^N\C^2)$ for $N\ge2$ and let $X$ be the rational normal curve. Theorem~\ref{thm:C}(b) then gives $\Lambda^Y_{\Kul}=\bigcup_{\xi\in\Lambda}H_\xi$ and $H_\xi\cap X=\{N\xi\}$, so $\Lambda^Y_{\Kul}\cap X=\{N\xi:\xi\in\Lambda\}=\Lambda$, as Theorem~\ref{thm:trace-rank1} predicts. The ambient set is much larger than its trace: the hyperplane $H_\xi$ is the part of the ambient limit set consisting of divisors with a root at $\xi$ and never meets the curve outside $N\xi$.
\end{example}

\paragraph{Higher rank.}
In higher rank, Theorem~\ref{thm:traces}(c) isolates the sources excluded only by the ambient action:
$E=(X\cap B^Y)\setminus B^X$.
They may be limits of ambient infinite-stabiliser points or ambient fixed-orbit limits. Outside $B^Y$, the compact-source calculation already contributes to $L_2^Y$; the remaining intrinsic limits lie in $\overline{\D_\Gamma^X(E)}$. This gives the Kulkarni trace inclusion when $E=\varnothing$. For Myrberg's set, under the hypotheses of Proposition~\ref{prop:highest-weight} the ambient trace and the intrinsic set both equal the complement of the transverse locus (Theorem~\ref{thm:embed}). The equality of exceptional sets requires only the coincidence of the singleton-fibre loci. The fibres themselves can differ, as in the adjoint example following Corollary~\ref{cor:restriction}.

\paragraph{The critical trace.} For $Y=\PP(\g)$ and $X=X_{\mathrm{ad}}$, with the highest restricted root space a line, the trace of the critical set can be written explicitly. By Theorem~\ref{thm:adcrit}, if every divergent sequence has type $\Theta=\iota\Theta$ then $\Lambda^Y_{\mathrm{crit}}\cap X_{\mathrm{ad}}=\{[e]\in X_{\mathrm{ad}}:\ e\in\mathfrak p_x\ \text{for some }x\in\Lambda_\Theta\}$: the points of the adjoint variety whose root vector lies in the parabolic subalgebra of a limit point. By Corollary~\ref{cor:restriction} its complement $\Omega^Y_{\mathrm{crit}}\cap X_{\mathrm{ad}}$ is a domain of proper discontinuity in $X_{\mathrm{ad}}$. To identify this as a largest domain, a return argument would have to reach every point of the removed trace. The non-recurrent witness in Proposition~\ref{prop:adrec-obstruction} is regular semisimple in $\mathfrak{sl}_3$, while $X_{\mathrm{ad}}$ parametrises rank-one nilpotents there. That ambient obstruction leaves the intrinsic returns on this trace to be checked separately.

\begin{question}\label{q:trace-higher}
For which $(G,P,\Gamma)$ and equivariant highest-weight embeddings does
$\Lambda_{\Kul}^X=\Lambda_{\Kul}^Y\cap X$ hold?
Which hypotheses ensure $X\cap B^Y=B^X$, or, more generally,
$\overline{\D_\Gamma^X(E)}\subseteq\Lambda_{\Kul}^Y$?
These are sufficient conditions for the inclusion in Theorem~\ref{thm:traces}(c). Equality also requires the reverse trace inclusion. The rank-one case under the stated representation hypotheses is Theorem~\ref{thm:trace-rank1}.
\end{question}

%--------------------------------------------------------------
\section{Forgetting flags}\label{sec:forget}
Embedding allowed the source to move outside the variety. Forgetting part of a flag is a different operation: nearby partial flags can be lifted to nearby full flags, so destinations pass through the projection exactly. Let $\Gamma<G$ be discrete. For each $\theta\subseteq\Delta$ write $\D^\theta=\D^{X_\theta}_\Gamma$ for the destination relation of $\Gamma$ on $X_\theta$.

\begin{theorem}[Forgetful maps]\label{thm:forget}
Let $\Theta\supseteq\theta$ and $\pi=\pi_{\Theta\theta}$. Then for every $x\in X_\Theta$,
\[
\D^\theta_\Gamma(\pi x)=\pi\bigl(\D^\Theta_\Gamma(x)\bigr).
\]
Consequently: (a) $\D^\theta=(\pi\times\pi)(\D^\Theta)$ as a subset of $X_\theta\times X_\theta$; (b) $\pi(\mathscr R^\Theta)\subseteq\mathscr R^\theta$; (c) if $U\subseteq X_\theta$ is a $\Gamma$-invariant open set on which $\Gamma$ acts properly discontinuously, then so does $\Gamma$ on $\pi^{-1}U\subseteq X_\Theta$; and (d) $\pi^{-1}(\Lambda^\theta_{\mathrm{exc}})\subseteq\Lambda^\Theta_{\mathrm{exc}}$, hence $\Omega_{\Myr}^\Theta\subseteq\pi^{-1}\Omega_{\Myr}^\theta$.
\end{theorem}
\begin{proof}
``$\supseteq$'': if $y=\lim\gamma_mx_m$ with $x_m\to x$ then $\pi y=\lim\gamma_m\pi x_m$ and $\pi x_m\to\pi x$. ``$\subseteq$'': let $y'=\lim\gamma_mz_m$ with $z_m\to\pi x$. Since $\pi$ is a locally trivial fibre bundle, we can lift $z_m$ to $x_m\to x$ with $\pi x_m=z_m$; by compactness of $X_\Theta$ a subsequence of $\gamma_mx_m$ converges to some $w$, and $\pi w=\lim\gamma_mz_m=y'$.
(a) is a restatement. (b) If $x\in\D^\Theta(x)$ then $\pi x\in\pi\D^\Theta(x)=\D^\theta(\pi x)$. (c) $\pi$ is proper (compact fibres), so a compact $K\subset\pi^{-1}U$ has compact $\pi K\subset U$, and $\{\gamma:\gamma K\cap K\ne\varnothing\}\subseteq\{\gamma:\gamma\pi K\cap\pi K\neq\varnothing\}$ is finite. (d) Fix one extracted sequence $(\gamma_m)$ for which $x$ has two distinct destinations. Lift the two realising source sequences through any prescribed $\tilde x\in\pi^{-1}(x)$, and extract limiting lifted targets. Their projections are distinct. Hausdorff convergence of the graphs along this same sequence (Theorem~\ref{thm:dest}(d)) places both targets in its whole-sequence destination at $\tilde x$, so $\tilde x$ is exceptional for that sequence. Take unions over sequences and then closures; the locally trivial projection is open, so inverse image commutes with this closure.
\end{proof}
\begin{corollary}[The full flag manifold knows everything]\label{cor:fullflag}
For every $\theta$ and every $x\in X_\theta$ and any lift $\tilde x\in G/B$, $\D^\theta_\Gamma(x)=\pi_\theta\bigl(\D^\Delta_\Gamma(\tilde x)\bigr)$.
\end{corollary}
Thus the relations of Chapters~\ref{ch:cartan}--\ref{ch:balanced} on any flag manifold are determined by the one on $G/B$. For a single sequence, Theorem~A combines the Bruhat cell of $\tilde x$ with the type $\Theta$ and Levi class $b$ (Section~\ref{sec:thmA}); the corollary allows any convenient lift of $x$.

\begin{remark}[Pullback of balanced systems]\label{rem:pullback-balanced}
A balanced system $\mathfrak B=(\mathfrak S,\mathfrak T)$ on $X_\theta$ (Chapter~\ref{ch:balanced}) pulls back to a balanced system $\pi^*\mathfrak B$ on $X_\Theta$ whose domain is $\pi^{-1}\Omega_{\mathfrak B}$ (this is (c) above). A domain upstairs descends precisely when its final closed wall is saturated for the fibres of $\pi$. Indeed, $\pi$ is open, and compact sets downstairs have compact full preimages, so properness descends on a saturated domain. Saturation is tested after taking the union over sequences and its closure. The cell graphs on partial flag manifolds are obtained by projecting source and target cell labels from the cell graph on $G/B$ for the same Cartan type.
\end{remark}

%--------------------------------------------------------------
\section{Embedding in projective space}\label{sec:embedding}
Forgetting part of a flag preserved every destination. We now return to embeddings, where the ambient space supplies extra directions in which a source can move. Let $\varnothing\ne\theta\subseteq\Delta$ and let $\rho:G\to\mathrm{GL}(V)$ be a nontrivial irreducible representation (over $\R$ or $\C$) with a highest restricted weight line $\ell_\chi$ of support $\theta$, so $P_\theta$ is the stabiliser of $\ell_\chi$. When only the projective representation descends to $G$, the linear formulas are read on a finite covering group; their projectivisations are independent of the lift. We choose a Cartan-compatible inner product (Hermitian over $\C$), for which $\rho(K)$ is isometric and $\rho(\exp\ak)$ is self-adjoint. Then the orbit map defines the \emph{highest-weight embedding}
\[
\iota_\rho:\ X_\theta=G/P_\theta\hookrightarrow\PP(V),\qquad gP_\theta\mapsto\rho(g)\ell_\chi,
\]
a closed equivariant embedding (the closed orbit in $\PP(V)$). The Tits representation $\rho_\theta$ of Chapter~\ref{ch:limit}, with highest weight $\chi=\sum_{\alpha\in\theta}c_\alpha\omega_\alpha$ for the chosen positive integers $c_\alpha$, is one choice; for a simple noncompact group, the adjoint representation (Section~\ref{sec:advar}) is another when the highest restricted root space is a line, with $\theta=\Theta_{\mathrm{ad}}$.

\begin{proposition}[Highest-weight comparison]\label{prop:highest-weight}
Let $\rho:G\to\mathrm{GL}(V)$ be nontrivial and irreducible with highest restricted weight $\chi$ whose weight space $V_\chi$ is a line $\ell_\chi$, let $\varnothing\ne\theta=\{\alpha\in\Delta:\ \chi-\alpha\ \text{is a weight}\}$ (so that $P_\theta=\mathrm{Stab}(\ell_\chi)$), and let $\Gamma<G$ be $P_\theta$-divergent with $\iota\theta=\theta$. Then the projective restriction $\overline\rho|_\Gamma$ has finite kernel and closed discrete image; in particular $\rho(\Gamma)$ is discrete, and:
\begin{enumerate}[label=\textup{(\roman*)},nosep]
\item for every sequence with Cartan data $(\Theta,k,k',b_{\mathrm L})$ and $\Theta\supseteq\theta$, $\sigma_2/\sigma_1(\rho g_m)=e^{-\min_{\alpha\in\theta}\alpha(\mu_m)}$ up to bounded factors, so the first pseudo-projective limit is $\tau^+=[\rho(k)\Pi_\chi\rho(k')]$ of rank one with image $\iota_\rho(kP_\theta)$ and kernel $K_{b}=\{[v]\in\PP(V):\langle v,\rho(k')^{-1}v_{\chi}\rangle=0\}$, $b=k'^{-1}P^-_\theta$;
\item $\Lambda_{\CoG}(\rho\Gamma)=\iota_\rho\Lambda_\theta$, $\Lambda_{\Myr}(\rho\Gamma)=\bigcup_{b\in\Lambda_\theta}K_b$, and $\iota_\rho^{-1}\Eq(\rho\Gamma)=\{x\in X_\theta:\ x\ \text{opposite to every }b\in\Lambda_\theta\}$.
\end{enumerate}
\end{proposition}
\begin{proof}
(i) The lowering argument of Theorem~\ref{thm:adjvar}(a) applies to $\ell_\chi$: every other weight is reached by a lowering word whose first step is $\chi-\alpha$ with $\alpha\in\theta$. Thus every weight other than $\chi$ is at most one of these first-step weights, and each first step occurs. The second exponent is therefore $\chi(\mu_m)-\min_{\alpha\in\theta}\alpha(\mu_m)$ up to a bounded error; $\theta$-divergence makes the gap divergent. Normalising $\rho(g_m)=\rho(k_m)\rho(e^{\mu_m})\rho(k'_m)$ by $e^{\chi(\mu_m)}$ leaves only the orthogonal projection onto $\ell_\chi$. The limiting compact factors give $\tau^+=[\rho(k)\Pi_\chi\rho(k')]$, with the stated image and kernel.
(ii) The divergent gap also proves discreteness. The function $q([A])=\sigma_2(A)/\sigma_1(A)$ is continuous and strictly positive on $\PGL(V)$, and is unchanged by scalar multiplication. For every sequence of distinct $\gamma_m\in\Gamma$, $\theta$-divergence and (i), after Cartan extraction, give $q(\overline\rho(\gamma_m))\to0$. If a compact subset of $\PGL(V)$ had infinitely many preimages in $\Gamma$, a sequence of such preimages would contradict the positive minimum of $q$ on that compact set. Thus compact preimages are finite, which proves closed discreteness of the projective image and finiteness of its fibres, including its kernel. This proves discreteness of the linear image as well, even if $\rho$ is trivial on a whole simple factor of $G$.

By opposition invariance, the repelling flags exhaust $\Lambda_\theta$. The image formula in (i), finite fibres and Proposition~\ref{prop:kernel} give the Conze--Guivarc'h and Myrberg identities, as in Theorem~\ref{thm:adjvar}(b); the kernel union is closed because it varies continuously over compact $\Lambda_\theta$. For equicontinuity, the highest-weight coefficient test in the proof of Theorem~\ref{thm:adjvar}(c) applies to $c(u)=\langle\rho(u)v_\chi,v_\chi\rangle$: it is nonzero on the big cell $N^-_\theta P_\theta/P_\theta$ and vanishes on the other Bruhat cells. Thus $\iota_\rho(x)\notin K_b$ exactly when $x$ is opposite to $b$. Proposition~\ref{prop:kernel} identifies the equicontinuity region with the complement of these kernels.
\end{proof}
This contains Proposition~\ref{prop:Tits-comparison} (take for $\rho$ the chosen Tits representation, $\chi=\sum_{\alpha\in\theta}c_\alpha\omega_\alpha$) and, when $d_{\max}=1$, Theorem~\ref{thm:adjvar}(b),(c) (take $\rho=\mathrm{Ad}$ and $\theta=\Theta_{\mathrm{ad}}$).

\begin{theorem}[Passage to projective space]\label{thm:embed}
Let $\Gamma<G$ be discrete and $P_\theta$-divergent with $\theta=\iota\theta$ opposition invariant, $\rho$ as above, and $\Lambda_\theta=\Lambda_\theta(\Gamma)\subset X_\theta$ the Quint limit set.
\begin{enumerate}[label=\textup{(\alph*)},nosep]
\item $\rho(\Gamma)$ is discrete, and every first pseudo-projective limit of a divergent sequence has rank one, with image in $\iota_\rho(X_\theta)$. In particular, $\Lambda_{\CoG}(\rho\Gamma)=\iota_\rho(\Lambda_\theta)$.
\item Identifying $X_\theta$ with its image for the relation, $\D^\theta_\Gamma\subseteq\D^{\PP(V)}_{\rho\Gamma}\cap\iota_\rho(X_\theta)^2$. The exceptional sets satisfy $\Lambda_{\mathrm{exc}}^{X_\theta}=\iota_\rho^{-1}\Lambda_{\mathrm{exc}}^{\PP(V)}$, and $\iota_\rho^{-1}\Omega_{\Myr}(\rho\Gamma)=\Omega_{\Myr}^{X_\theta}(\Gamma)$.
\item $\iota_\rho^{-1}\Eq(\rho\Gamma)=\Omega^{\mathrm{trans}}_\theta(\Gamma):=\{x\in X_\theta:\ x\ \text{is opposite to every }b\in\Lambda_\theta\}$.
\end{enumerate}
In particular $\Omega^{\mathrm{trans}}_\theta(\Gamma)=\Omega_{\Myr}^{X_\theta}(\Gamma)$. It is the same for the canonical embedding $\kappa_P$ of Proposition~\ref{prop:canonical} and for any Tits-type embedding with support $\theta$. In rank one, for $\theta=\{\alpha\}$, this is the classical domain $\partial_\infty\mathbb H\setminus\Lambda$.
\end{theorem}
\begin{proof}
(a),(c) are Proposition~\ref{prop:highest-weight}. Its finite projective fibres identify escaping sequences after extraction and preserve properness: the kernel fixes $X_\theta$, and compact transporters in $\Gamma$ are finite-fibre preimages of those in the image. Corollary~\ref{cor:restriction} gives the relation inclusion. To compare exceptional sources, fix an extracted sequence of type $\Theta$. In adapted coordinates the open source cell is $P_\Theta^-P_\theta/P_\theta$, and the intrinsic singleton fibres occur exactly there by Theorem~\ref{thm:dest} and the cell dimensions. Its highest-weight image lies outside the kernel of the first ambient limit: the highest weight line survives on the top quotient, which is preserved by the Levi and on which $N_\Theta^-$ acts trivially. Conversely, outside this kernel a source has one ambient destination, so its intrinsic fibre is a singleton. Thus the exceptional sets agree on $X_\theta$ for the sequence. By the compact kernel incidence of Proposition~\ref{prop:kernel}, the ambient exceptional union is closed, and taking the union proves the equalities in (b).
\end{proof}
\begin{example}[The rational normal curve]\label{ex:rnc}
Let $\Gamma<\PSL(2,\C)$ be non-elementary and discrete with proper classical limit set $\Lambda\subsetneq\C\PP^1$. For $N\ge1$ let $\rho_N=\mathrm{Sym}^N$ and let $\iota:\C\PP^1\to\PP^N$ be the rational normal curve. Here $\Omega^{\mathrm{trans}}=\C\PP^1\setminus\Lambda=\Omega_{\Myr}^{\C\PP^1}$ (equality), and the ambient largest domain of Theorem~\ref{thm:C}, $\Omega_{\max}(\rho_N\Gamma)=\{D:\mathrm{mult}_\xi D<\lceil N/2\rceil\ \forall\xi\in\Lambda\}$, meets the curve exactly in $\{N\xi:\xi\notin\Lambda\}$, the same set. But for $N\ge3$ it also contains divisors off the curve that meet $\Lambda$ (e.g.\ $\xi+E$ with $\xi\in\Lambda$ and $E$ an effective divisor of degree $N-1$ supported outside $\Lambda$), which no point of the curve sees: the ambient dynamics gives a bigger domain around the curve while on the curve it is the classical $\C\PP^1\setminus\Lambda$.
\end{example}

\paragraph{Comparison with the flag manifold.} The ambient space supplies domains around points off $\iota_\rho(X_\theta)$, as in the preceding example, and the projective questions of critical level, strong irreducibility and proximality. The trichotomy of Chapter~\ref{ch:limit} applies to $\rho\Gamma$; Section~\ref{sec:lost} asks what its first pseudo-projective limits leave unseen on the flag manifold.

%--------------------------------------------------------------
\section{The model case: \texorpdfstring{$\PSL(n+1,\C)$}{PSL(n+1,C)} and \texorpdfstring{$\C\PP^n$}{complex projective n-space}}\label{sec:psl-case}
For $G=\PSL(n+1,\C)$ and $P=P_{\alpha_1}$ the variety $G/P$ is $\C\PP^n$ and there are three equivariant projective realisations of the same action:
\begin{itemize}[nosep]
\item the \emph{identity} projective action on $\PP(V)$, $V=\C^{n+1}$, induced by the standard representation of $\mathrm{SL}(n+1,\C)$ of highest weight $\omega_1$, with line bundle $\mathcal O(1)$;
\item the \emph{canonical} (anticanonical) one, $V_P=\mathrm{Sym}^{n+1}V\subset\wedge^{n}\mathfrak{sl}(V)$, line bundle $\mathcal O(n+1)$, Example~\ref{ex:veronese};
\item the Grassmannian embedding $j_P:\C\PP^n\to\mathrm{Gr}_{n^2+n}(\mathfrak{sl}(V))$, $[v]\mapsto\mathfrak p_{[v]}=\{X:Xv\in\C v\}$.
\end{itemize}
All three give the same intrinsic relation $\D_\Gamma$ on $\C\PP^n$, since each is an equivariant embedding of the same compact space. The ambient equicontinuity traces require a separate comparison: Proposition~\ref{prop:highest-weight} applies under its stated support hypothesis, while the symmetric-power argument of Example~\ref{ex:veronese} identifies the identity and Veronese traces. No such equality is asserted here for the full Grassmannian ambient space. What differs is what the ambient space sees outside the variety: in the identity representation the ambient space \emph{is} the variety, so Chapters~\ref{ch:cartan}--\ref{ch:limit} specialise to the block form of Theorem~A (Theorem~\ref{thm:blockform}). This is what makes $\C\PP^n$ the model for the block calculations. For a general $G/P$, a representation into $\PSL(n+1,\C)$ gives ambient calculations on $\C\PP^n$ through the Pl\"ucker projectivisation of $j_P$ or through $\kappa_P$; the restriction theorem states what passes back to the flag variety. The other flag manifolds reach $\C\PP^n$-like spaces by the forgetful maps of Section~\ref{sec:forget} and the fundamental exterior-power embeddings:
\[
\mathrm{Fl}(\C^{n+1})\xrightarrow{\ \pi\ }\mathrm{Gr}_k(\C^{n+1})\xrightarrow{\ \text{Pl\"ucker}\ }\PP(\wedge^k\C^{n+1}),
\]
with $\wedge^k$ the fundamental representation of highest weight $\omega_k$ and $\theta=\{\alpha_k\}$.

\begin{proposition}[Grassmannians and $\C\PP^n$]\label{prop:grass-proj}
Let $(g_m)$ be a divergent sequence in $\PSL(n+1,\C)$ with extracted Cartan blocks $E_1,\dots,E_s$ and let $1\le k\le n$. The Cartan exponents of $\wedge^kg_m$ are the ordered sums $\lambda_{i_1}+\dots+\lambda_{i_k}$. To read its full Cartan type, extract further if necessary so that the adjacent gaps in these ordered sums converge in $[0,\infty]$. A rank-one first limit occurs iff $k$ is a block boundary, $k=r_j$, in which case the image is the point $[G_j]\in\mathrm{Gr}_k$ and the kernel is the projective hyperplane in $\PP(\wedge^k\C^{n+1})$ whose Grassmannian trace consists of the $k$-planes not transverse to the source $(n+1-k)$-plane $F_{j+1}$.
\end{proposition}
\begin{proof}
The two largest exponents of $\wedge^kg_m$ differ by $\lambda_k-\lambda_{k+1}$, so a rank-one limit occurs exactly at a block boundary. At such a boundary, the highest-weight normalisation in the proof of Proposition~\ref{prop:highest-weight}(i), with $\chi=\omega_k$, leaves the image $[\wedge^kG_j]$. The kernel is the projective hyperplane on which the top exterior coefficient vanishes; on a decomposable $k$-vector this is precisely failure of transversality to $F_{j+1}$.
\end{proof}
For an individual extracted sequence, Theorem~\ref{thm:graphcomp} recovers the full exterior datum from its projective destination relation; Theorem~\ref{thm:dest} then gives the intrinsic Grassmannian relation. Conversely, over $\C$, Theorems~\ref{r2gr:brion} and~\ref{r2gr:destinations} identify the reduced graph families on $\C\PP^n$ and $\mathrm{Gr}_k$ with the same wonderful compactification of $\PGL(n+1,\C)$. A reduced graph is determined by its support, so either entire sequence relation determines that common datum and the other relation. This uses the whole intrinsic relation, rather than only the rank-one limits or a trace of an ambient projective relation.

%--------------------------------------------------------------
%--------------------------------------------------------------
\section{A principal image on a quadric}\label{sec:q3}
On a quadric we can compute what an ambient critical wall removes from a flag variety in two ways: by the Schubert chain of Chapter~\ref{ch:balanced} and by the equation of the variety. Let $\iota:\PSL(2,\C)\to\mathrm{PSp}(4,\C)\cong\mathrm{SO}(5,\C)$ be the principal homomorphism. Its lift from $\mathrm{SL}(2,\C)$ acts by $\mathrm{Sym}^3$ on $\C^4$; the descended representation on $\C^5$ is $\mathrm{Sym}^4$. Let $X=\mathrm{LG}(2,4)=Q^3\subset\PP^4$ be the Lagrangian Grassmannian, which is the quadric $\{I=0\}$ of binary quartics for the quadratic invariant $I$. Let $\xi=(\xi^1\subset\xi^2):\PP^1\to G/B$ be the Frenet (osculating) flag curve; $\xi^2(x)\in Q^3$ lies on the rational normal curve $\mathcal C=\{\ell^4\}$. For $x\in\PP^1$ put
\[
L_x=\{L\in\mathrm{LG}(2,4):\ \xi^1(x)\subset L\}\ \cong\ \PP^1,
\]
the line of $Q^3$ through $\xi^2(x)$ tangent to $\mathcal C$ (in the quartic model, $\{\ell^3m\}$ for $\ell=x$).
\begin{proposition}[The largest domain on $Q^3$]\label{prop:q3}
Let $\Gamma<\PSL(2,\C)$ be non-elementary and discrete with limit set $\Lambda\subsetneq\PP^1$, acting on $Q^3$ through $\iota$. Then
\[
\mathscr R_\Gamma=\Lambda_{\mathfrak B}(\Gamma)=\bigcup_{x\in\Lambda}L_x,
\]
the balanced system being the sublevel system $\mathrm{Th}=\{0,1\}$ (the point and the line of the Schubert chain relative to a flag), and $\Omega=Q^3\setminus\bigcup_{x\in\Lambda}L_x$ is the \emph{largest} domain of proper discontinuity. It is exactly the restriction to $Q^3$ of the largest domain in $\PP^4$ of Theorem~E. Scheme-theoretically, the plane of quartics divisible by $\ell^2$ cuts the quadric in a double line; its reduced support is $L_\ell$.
\end{proposition}
\begin{proof}
Every divergent sequence of $\iota(\Gamma)$ has type $\Delta$, with data $\xi^+=\xi(a)$, $\xi^-=\xi(b)$ where $a,b\in\Lambda$ are the attracting and repelling points of the sequence in $\PP^1$ (as in Proposition~\ref{prop:unipotent}); every pair $(a,b)\in\Lambda\times\Lambda$, including $a=b$, occurs for a non-elementary $\Gamma$ (take loxodromics whose ordered fixed points converge to $(a,b)$, then take powers large enough for their Cartan gaps and endpoint data to realise this limit). The Weyl group of $C_2$ has order $8$ and $W_\theta$ order $2$, so there are four Schubert cells, of complex dimensions $0,1,2,3$, forming a chain for the Bruhat order, and $w_0$ reverses it. By Proposition~\ref{prop:KLPdict}(ii), $\mathrm{Th}=\{0,1\}$ satisfies the fat condition (the complement of $w_0\mathrm{Th}=\{2,3\}$ is $\{0,1\}\subset\mathrm{Th}$), and $\mathrm{Th}(\xi)=L_x$ for $\xi=\xi(x)$ (the cells of dimension $\le1$ relative to $\xi$ form the unique Schubert line $L_x=\{L\supset\xi^1\}$ for that flag, tangent to the rational normal curve at $\xi^2(x)$). Hence $\Lambda_{\mathfrak B}=\bigcup_{x\in\Lambda}L_x$, closed, and $\mathscr R_\Gamma\subset\Lambda_{\mathfrak B}$ by Corollary~\ref{cor:trichotomy}. Conversely, let $L\in L_x$, $x\in\Lambda$, and take a sequence with $\xi^+=\xi^-=\xi(x)$. By Proposition~\ref{prop:unipotent}(b) (whose proof only uses $k'^{-1}w_0B=kB$), $\D(L)=X_{w_0\tau}(\xi(x))$ where $\tau\in\{0,1\}$ is the position of $L$; this is the Schubert variety of dimension $3$ or $2$ (all of $Q^3$, or the Maslov cycle $\{L':L'\cap\xi^2(x)\ne0\}$), and it contains $L_x\ni L$. So $(L,L)\in\D_\Gamma$, $L\in\mathscr R_\Gamma$. Thus $\mathscr R_\Gamma=\Lambda_{\mathfrak B}$ is a vertex cover, and Corollary~\ref{cor:trichotomy}(i) gives the largest domain. For the comparison, use plain coefficients $f=ax^4+bx^3y+cx^2y^2+dxy^3+ey^4$, for which the invariant is $I=12ae-3bd+c^2$ up to a nonzero scalar. With $\ell=x$ and $q=Ax^2+Bxy+Cy^2$ one obtains $I(\ell^2q)=C^2$. Thus
\[
 \bigl(Q^3\cap\PP(\ell^2\operatorname{Sym}^2\C^2)\bigr)_{\mathrm{red}}
 =\PP(\ell^3\operatorname{Sym}^1\C^2)=L_\ell.
\]
The scheme cut is a double line, and the reduced equality holds for every $\ell$ by equivariance. Taking the union over $\Lambda$ identifies the two removed supports and hence the two open domains. The two supports agree independently of the size of the classical limit set.
\end{proof}
\begin{remark}
The point limit curve, its tangent-line family, and the family of Maslov cycles are different levels of the flag geometry. On $Q^3$, the critical removed set is the tangent-line family $\bigcup_{x\in\Lambda}L_x$; the exceptional wall comes from the larger Maslov cycles. The ambient middle-multiplicity wall restricts to the same critical support on this quadric, although its scheme cut has multiplicity two. The support determines the removed set, while the double line records how the ambient plane meets the quadric. Here the largest domains agree. For a strict restriction of destination relations, we turn to the Segre example in Section~\ref{r2gr:restriction}.
\end{remark}

\begin{remark}[Cominuscule flag manifolds and Demazure's list]\label{rem:demazure}\
When $P_\theta$ is cominuscule for the acting group $G$, its nilradical $N_\theta$ is abelian and the open cell has the corresponding vector-space coordinates. The standard examples include projective spaces, Grassmannians, quadrics, Lagrangian and spinor Grassmannians, and the two exceptional ones, with their cominuscule presentations. This gives useful coordinates for possible generalisations of Chapter~\ref{ch:psl}; we suggest Grassmannians, then $\mathrm{LG}(n,2n)$ and quadrics. For a general sequence, however, the group in Theorem~\ref{thm:dest} is $N_\Theta$, indexed by its Cartan type, and need not be abelian. The destination fibres remain unipotent-orbit closures and need not be projective subspaces in a projective embedding. For a complex simple adjoint group $G$ and a maximal parabolic $P$, Demazure's classification gives $\mathrm{Aut}^0(G/P)=G$ outside a short list of exceptions. Examples of larger connected automorphism groups are those of $\mathrm{Sp}(2n)/P_1=\PP^{2n-1}$ (automorphism group $\PSL(2n)$) and $G_2/P_1=Q^5\subset\PP^6$ (automorphism group $\mathrm{SO}(7)$). The cominuscule condition concerns the pair $(G,P_\theta)$: for instance the line parabolic of $\mathrm{Sp}(2n)$ has nonabelian nilradical for $n\ge2$, while the line parabolic of $\PSL(2n)$ has abelian nilradical. For $\Gamma\subset G$ acting on such an $X$ both groups compute the same relation $\D_\Gamma$, since it is intrinsic; a consistency check is Table~\ref{tab:sp4}, whose $\R\PP^3$ sequences agree with Proposition~\ref{prop:proj} for $\PSL(4)$. What is \emph{not} automatic is that the Cartan type of $g\in G$ computed in the larger group is the image of its type in $G$; we have checked this only for $\mathrm{Sp}(2n)\subset\PSL(2n)$, where the gaps of $\PSL(2n)$ are the symmetrised gaps of $\mathrm{Sp}(2n)$, and not for $G_2\subset\mathrm{SO}(7)$.
\end{remark}

\section{What the passage to projective space loses}\label{sec:lost}
Theorem~\ref{thm:embed}(c) identifies the equicontinuity region of $\rho\Gamma$, restricted to $X_\theta$, with the set of flags opposite to \emph{every} limit flag. For rank one this is the right domain. In higher rank it may be empty, because incidence conditions force flags to meet the limit set.

\begin{example}[Hitchin representations in $\mathrm{SL}(3,\R)$]\label{ex:hitchin-empty}
Let $\Gamma<\mathrm{SL}(3,\R)$ be the image of a closed surface group under a Hitchin (= convex projective) representation, $\Omega_c\subset\R\PP^2$ the properly convex domain it divides, and $C=\partial\Omega_c$ the limit curve, a strictly convex $C^1$ curve with tangent line $\ell_b$ at $b\in C$. The limit set in the full flag manifold $\mathrm{Fl}=\{(p,\ell):p\in\ell\}$ is $\Lambda=\{(b,\ell_b)\}$. A flag $(p,\ell)$ is opposite to $(b,\ell_b)$ iff $p\notin\ell_b$ and $b\notin\ell$. Hence $(p,\ell)\in\Omega^{\mathrm{trans}}$ iff $\ell\cap C=\varnothing$ and $p\notin\bigcup_b\ell_b$. But $\ell$ misses $C$ only if $\ell$ lies outside $\overline{\Omega_c}$, and then $p\in\ell$ lies outside $\overline{\Omega_c}$ and every such point lies on a tangent line to $C$. So $\Omega^{\mathrm{trans}}=\varnothing$: the trace of the equicontinuity region of the Tits image of $\Gamma$ on the full flag manifold is empty. Nevertheless $\Gamma$ acts properly discontinuously and cocompactly on a nonempty open subset of $\mathrm{Fl}$ (Guichard--Wienhard; Kapovich--Leeb--Porti, balanced thickenings). The same $\Gamma$ in $\mathrm{SL}(3,\R)$ acts on $\R\PP^2$ with the domain $\Omega_c$.
\end{example}
In the full-flag setting of this example, the transverse domain corresponds to the fat thickening $\mathrm{Th}_{\mathrm{trans}}=W\setminus\{w_0\}$: a point is removed iff it is not opposite to some limit flag. This thickening gives a vertex cover by Theorem~\ref{thm:KLPfull}; no additional condition on the Levi datum is required. In type $A_2$, the balanced thickening $\{e,s_1,s_2\}$ removes fewer Bruhat positions and can therefore give a larger nonempty domain. Recovering those domains asks which destination pairs must be excluded, not only whether a source is transverse to every limit flag. What a first pseudo-projective limit records depends on the representation $\rho_\theta$: a regular highest-weight representation can encode the full exterior boundary (Theorem~\ref{thm:extG}). Thus the information lost by a first pseudo-projective limit is determined by the representation in which it is taken; this is why Chapters~\ref{ch:cartan}--\ref{ch:balanced} work with the full relation, and why $X_{\mathrm{ext}}$ (Section~\ref{sec:ext}) records all gaps at once.
\paragraph{The balanced domain in $A_2$.} For $W=S_3$ the thickening $\mathrm{Th}=\{e,s_1,s_2\}$ (the elements of length $\le1$) is balanced: $w_0=s_1s_2s_1$ and $w_0\{s_1s_2,s_2s_1,w_0\}=\{s_1,s_2,e\}=\mathrm{Th}$. In the Hitchin example the corresponding domain is, up to the convention for the order,
\[
\Omega_{\mathrm{bal}}=\{(p,\ell)\in\mathrm{Fl}:\ p\notin C,\ \ell\ \text{is not tangent to }C\},
\]
the complement of the flags in relative position $\le s_i$ to some limit flag (they share the point or the line). It is nonempty, and indeed open and dense, in contrast with $\Omega^{\mathrm{trans}}=\varnothing$, and Kapovich--Leeb--Porti \cite[Theorem~7.12]{KLP} give a cocompact proper action of the Anosov group $\Gamma$ on it. In the notation of Chapter~\ref{ch:balanced}, $\Omega^{\mathrm{trans}}$ comes from the fat thickening $W\setminus\{w_0\}$, whose complementary thickening is $\{e\}$. The balanced ideal $\{e,s_1,s_2\}$ is properly contained in it and its vertex cover gives $\Omega_{\mathrm{bal}}$. Thus $\Omega^{\mathrm{trans}}\subseteq\Omega_{\mathrm{bal}}$, with strict inclusion in this example.

\chapter{The case of \texorpdfstring{$\PSL(n+1,\C)$}{PSL(n+1,C)}}
\label{ch:psl}

The projective questions studied here arise from the theory of complex Kleinian groups. A discrete subgroup $\Gamma<\PSL(n+1,\C)$ is called a \emph{complex Kleinian group} if it acts properly discontinuously on some nonempty $\Gamma$-invariant open subset of $\C\PP^n$. This is the notion introduced by Seade--Verjovsky~\cite{SeadeVerjovsky2001}; see also~\cite{SeadeVerjovsky2002,CNS}. The existence of such an open set is part of the definition. It does not require that the equicontinuity region be nonempty, that it coincide with the Kulkarni region, or that a largest domain exist. We continue to consider arbitrary discrete projective groups, including those with no nonempty domain of proper discontinuity.

In dimension two, Navarrete's work on complex hyperbolic and cyclic projective groups~\cite{Navarrete2006,Navarrete2008}, and the work of Barrera--Cano--Navarrete on equicontinuity, Kulkarni's set and configurations of projective lines~\cite{BCNLimit2011,BCNFour2011,BCNLines2016}, established the distinctions that motivate this chapter. Cano--Seade~\cite{CanoSeade2014} studied the geometry of projective groups with compact quotient components, and Cano--Navarrete--Seade~\cite{CNS} developed the general framework. The present comparison continues that programme: it retains the asymptotic information of all realised Cartan types and the admissibility of the moving sources, then asks which parts of the resulting boundary force recurrence. The recent flag description in dimension two~\cite{BCNSFlag2026} is a direct point of comparison, with its divergence hypotheses kept explicit.

The preceding chapters treat flag manifolds of semisimple groups, with simplicity or effectiveness hypotheses where specified. We now specialise to $G=\PSL(n+1,\C)$ acting on $\C\PP^n$. We compare the limit sets of a discrete group $\Gamma$, determine how they are ordered, and ask when a single natural one gives the largest domain of proper discontinuity. The projective case is special because the whole geometry can be written in terms of subspaces of $V=\C^{n+1}$: a divergent sequence has an ordered block decomposition of $V$ with a flag of kernels and a flag of images (Section~\ref{sec:blockform}), the exterior powers turn Cartan data into highest-weight lines, and the flags give the projective spaces that enter the limit-set calculations. For Kulkarni's set, which sources are admissible also depends on the fixed-source set.

\medskip We begin by bringing together the limit sets of Chapters~\ref{ch:cartan}--\ref{ch:balanced}. Their projective subspaces let us compare the sets directly, but the comparison must remember which sources are allowed to move. Write $\operatorname{Lim}_{\rm SP}(\Gamma)\subset\PP(\mathrm{End}\,V)$ for the set of limits of sequences of distinct elements of $\Gamma$, $\Eq(\Gamma)$ for the equicontinuity region, $B_\Gamma$ for the fixed-source set of Kulkarni, $\D_\Gamma$ for the destination relation and $\mathscr R_\Gamma$ for the recurrent set. $\Lambda_{\mathrm Z}^{(k)}(\Gamma)$ is the closure of the union, over divergent sequences, of the target spaces $\PP(G_j)$ of projective dimension $\le k$ of their Cartan flags (Section~\ref{sec:flagLS}), and $\Lambda_{\mathrm{crit}}$, $\Omega_{\mathrm{crit}}$ are the set and the domain of the critical system of Section~\ref{sec:critsys}: admissible oriented sequences supply the common global level $q_\Gamma$ and its paired complementary level $n-q_\Gamma-1$.

\paragraph{Rank-one images are recurrent.}
Let $[T]$ be a rank-one pseudo-projective limit of distinct elements
$\gamma_m\in\Gamma$, and put $a=\PP(\operatorname{Im}T)$.
For every fixed $x\notin\PP(\ker T)$ one has $\gamma_mx\to a$.
Thus
\[
 \bigl(\PP(V)\setminus\PP(\ker T)\bigr)\times\{a\}
 \subseteq\D_\Gamma.
\]
The first factor is dense, and $\D_\Gamma$ is closed by
Lemma~\ref{lem:Dbasic}; hence $(a,a)\in\D_\Gamma$.
Choose a fixed source outside $\PP(\ker T)$. If it lies in $L_0$, its limit $a$ also lies in $L_0$, by closedness and invariance; otherwise $a$ is a generator of $L_1$. Taking closures and applying Proposition~\ref{prop:kulkarni-myrberg} gives
\[
 \Lambda_{\CoG}(\Gamma)\subseteq B_\Gamma\subseteq\mathscr R_\Gamma.
\]
No irreducibility or proximal-element hypothesis is needed.

\begin{theorem}[The flag of limit sets]\label{thm:flag-of-limit-sets}
Let $\Gamma<\PSL(n+1,\C)$ be an infinite discrete group. The following statements compare the closed $\Gamma$-invariant limit sets and their domains in $\C\PP^n$. The Cartan flags supply the flag levels and the critical wall, and the fixed-source set $B_\Gamma$ enters the formula for Kulkarni's set.
\begin{enumerate}[label=\textup{(\alph*)}]
\item \emph{The flag.}
\[
\Lambda_{\CoG}(\Gamma)=\Lambda_{\mathrm Z}^{(0)}(\Gamma)\subseteq\Lambda_{\mathrm Z}^{(1)}(\Gamma)\subseteq\cdots\subseteq\Lambda_{\mathrm Z}^{(n-1)}(\Gamma)=\Lambda_{\Myr}(\Gamma),
\]
where $\Lambda_{\Myr}(\Gamma)=\C\PP^n\setminus\Eq(\Gamma)=\bigcup_{\tau\in\operatorname{Lim}_{\rm SP}(\Gamma)}\PP(\ker\tau)$, and $\Gamma$ acts properly discontinuously on $\Eq(\Gamma)$.
\item \emph{The bottom.} If $\Gamma$ is strongly irreducible and contains a proximal element, $\Lambda_{\CoG}(\Gamma)$ is the closure of the attracting fixed points of proximal elements and is the unique non-empty minimal closed invariant set.
\item \emph{The critical level.} $\Lambda_{\CoG}(\Gamma)\subseteq\mathscr R_\Gamma\subseteq\Lambda_{\mathrm{crit}}(\Gamma)$, and $\Lambda_{\mathrm{crit}}$ is a closed invariant vertex cover, so $\Gamma$ acts properly discontinuously on $\Omega_{\mathrm{crit}}$. Exactly one of three cases occurs: $\Lambda_{\mathrm{crit}}=\mathscr R_\Gamma$, and $\Omega_{\mathrm{crit}}$ is the largest domain; or $\mathscr R_\Gamma\subsetneq\Lambda_{\mathrm{crit}}$ and $\mathscr R_\Gamma$ is a vertex cover, in which case the largest domain $\C\PP^n\setminus\mathscr R_\Gamma$ strictly contains $\Omega_{\mathrm{crit}}$; or $\mathscr R_\Gamma$ is not a vertex cover, and no largest domain exists.
\item \emph{The Kulkarni level.} $\Gamma$ acts properly discontinuously on $\Omega_{\Kul}(\Gamma)$; $\Eq(\Gamma)\subseteq\Omega_{\Kul}(\Gamma)$, that is $\Lambda_{\Kul}(\Gamma)\subseteq\Lambda_{\Myr}(\Gamma)$; and $\Lambda_{\Kul}(\Gamma)=B_\Gamma\cup\overline{\D_\Gamma(\C\PP^n\setminus B_\Gamma)}$. The equality $\Lambda_{\Kul}=\Lambda_{\Myr}$ holds exactly when $\Omega_{\Kul}=\Eq$. In general $\Lambda_{\Kul}(\Gamma)$ is not a level of the filtration in \textup{(a)}.
\item \emph{Position of the levels.} Without any irreducibility or largest-domain hypothesis,
\[
\Lambda_{\CoG}(\Gamma)\subseteq B_\Gamma\subseteq\mathscr R_\Gamma\subseteq\Lambda_{\Kul}(\Gamma)\subseteq\Lambda_{\Myr}(\Gamma).
\]
If $\Omega_{\mathrm{crit}}$ is the largest domain, then
$\Lambda_{\mathrm{crit}}=\mathscr R_\Gamma$, so the critical set can
replace the recurrent set in this chain. Without that property the
critical set need not lie below $\Lambda_{\Myr}$.
\item \emph{Largest domains.} An invariant open set is a domain of proper discontinuity exactly when it contains no pair of $\D_\Gamma$. The union of all such domains is $\C\PP^n\setminus\mathscr R_\Gamma$. A largest domain exists exactly when no two non-recurrent points are related, and otherwise there are at least two incomparable maximal domains.
\end{enumerate}
\end{theorem}
Most of the proof is already in Chapters~\ref{ch:limit} and~\ref{ch:balanced}; we give the comparison in Section~\ref{sec:flagLS}. The relational criterion and the compact-source formula apply on general flag manifolds (Theorem~\ref{thm:relational} and Proposition~\ref{prop:L2}). Their Cartan-level filtration is Theorem~\ref{thm:L}, with its stated hypotheses. The projective kernel description in (a) is specific to the projective setting; the unconditional inclusion $\Lambda_{\Kul}\subseteq\Lambda_{\Myr}$ in (d) holds for the general flag actions of Section~\ref{sec:KulMyr}. The critical wall records a different property: which of its subspaces consist of recurrent points.

\medskip By (c), the critical domain is the largest domain precisely when $\Lambda_{\mathrm{crit}}=\mathscr R_\Gamma$, that is, when there is a return, with moving sources, to each point of the critical wall. A largest domain may also exist and strictly contain the critical domain, as in the second case of (c). The second theorem gives conditions under which these returns can be constructed:
\begin{restated}{Theorem~\ref{thm:crit-recurrence}}If $\Gamma<\PSL(n+1,\C)$ is discrete with self-opposed uniform critical type, full critical flag module and strongly irreducible critical action, and either $p-a=1$ or the critical flags carry a regular seed, then $\Lambda_{\mathrm{crit}}=\mathscr R_\Gamma$: the critical domain is the largest domain of proper discontinuity.\end{restated}
The return construction in Section~\ref{sec:crit-proj} asks us to understand both the span of the critical flags and the geometry of their middle spaces. For lattice images, Borel density can establish the module hypotheses (Corollary~\ref{r3bor:lattice-largest}). For Borel--Anosov images the complete boundary flag gives another direct return construction in every dimension (Proposition~\ref{r3an:borel-return}). We then examine the limits of the criterion: strong irreducibility can hold with no largest domain at all (Section~\ref{sec:cp3}), and the critical wall itself can be too large (Section~\ref{r4voice:limits}).

\medskip The spaces $\PP(G_j)$ of the Cartan flags of all divergent sequences give the levels of the flag, ordered by dimension. For full-flag sequences the two ends are the Conze--Guivarc'h points ($\dim 0$) and the Myrberg hyperplanes ($\dim n-1$, equivalently the kernels of the pseudo-projective limits); in general the kernel spaces need not be hyperplanes. Between the ends, the critical system pairs a middle target space with its complementary inverse space so that every dynamical pair meets the resulting wall. We prove recurrence by constructing returns on one carrier and propagating them through its orbit closure (Section~\ref{sec:crit-proj}). The examples will show why the distinction from Kulkarni's set matters. For the symmetric power of degree $N\ge3$ of a non-elementary Kleinian group with proper classical limit set, the largest domain is strictly larger than the Kulkarni domain. The difference is a multiplicity threshold (Theorem~E, Section~\ref{sec:principal}).

In projective dimension two, Cano--Cano Cordero--Loeza~\cite{CCL2026} prove the existence of a largest domain outside two explicitly described exceptional families, assuming that a nonempty invariant domain of proper discontinuity exists. The general recurrence criterion here should be read alongside that theorem; it is not a new assertion of existence in all of those two-dimensional cases.

The section has four parts. The first reads Theorem~A in block form and compares the limit sets of Theorem~\ref{thm:flag-of-limit-sets}, including Kulkarni's set through its admissible sources. The second defines the critical system and proves the recurrence theorem, with its applications to lattices, cusped Hitchin groups and Borel--Anosov groups, and its limits. The third treats principal images of Kleinian groups and unitary groups. The fourth constructs families with prescribed behaviour and studies the complex geometry of their equicontinuity regions.

%--------------------------------------------------------------
\sectionpart{Limit sets in projective space}

\section{Theorem A in block form}\label{sec:blockform}
Let $V=\C^{d}$, $d=n+1$. A divergent sequence $(g_m)$ in $\PSL(V)$ has, after extraction, Cartan blocks $E_1,\dots,E_k$ of dimensions $b_1,\dots,b_k$ (Section~\ref{sec:ext}), source flag $F_1\supset F_2\supset\cdots\supset F_k$ (vector dimensions $d,\,d-b_1,\ldots$) and target flag $G_1\subset\cdots\subset G_k=V$ with $\dim G_j=b_1+\dots+b_j$. We first read the destination of one source from these two flags. Proposition~\ref{prop:graded} becomes the following form of Theorem~A.

\begin{theorem}[Theorem A for $\PSL(V)$]\label{thm:blockform}
Let $x=[\xi]$ with $\xi\in F_j\setminus F_{j+1}$. Then
\[
\D_{(g_m)}(x)=\PP\bigl(G_{j-1}\oplus\C\,v_x\bigr),\qquad v_x=k\,D_j\,\mathrm{pr}_{E_j}(k'\xi),
\]
a projective space of dimension $b_1+\dots+b_{j-1}$ (a point if $j=1$). Equivalently: $\D(x)$ is the preimage, under $V\supset G_j\to G_j/G_{j-1}$, of the line $A_j(\C\,q_j\xi)$ in the $j$-th graded piece of the target flag, where $A_j:F_j/F_{j+1}\to G_j/G_{j-1}$ is the \emph{incidence arrow}. The inverse sequence has the transposed relation.
\end{theorem}
Thus: a source in the $j$-th stratum is sent onto the projective span of everything the target flag has accumulated before level $j$, together with one more direction determined by the $j$-th graded piece. Two special cases are used constantly. For $j=1$ (the generic stratum) the destination is the single point $[k D_1\mathrm{pr}_{E_1}k'\xi]$ in $\PP(G_1)$. For the deepest stratum $j=k$ the destination is the projective span $\PP(G_{k-1}\oplus\C v_x)$, a projective space of dimension $d-b_k$; when $k=2$ and $b_2=1$ this is the pseudo-projective picture (kernel $\PP(F_2)$, image $\PP(G_1)$).

\begin{corollary}\label{cor:blockexc}
For a divergent sequence with extracted Cartan data, $E_{(g_m)}=\PP(F_2)$ and
$\dim_\C\PP(F_2)=d-b_1-1$. The Myrberg set is the closure of the union of these exceptional sets over escaping sequences with extracted Cartan data, equivalently the closed union of the projectivised kernels of all pseudo-projective boundary limits. If every such limit has rank one, this is a union of hyperplanes. The converse does not follow from that description: other group sequences can enlarge the union without changing its expression by hyperplanes.
\end{corollary}
\begin{proof}
The generic stratum has a singleton destination; all other strata have positive-dimensional destinations by Theorem~\ref{thm:blockform}. The source space $F_2$ has vector dimension $d-b_1$. The remaining statements follow from the kernel description in Proposition~\ref{prop:kernel}.
\end{proof}

%--------------------------------------------------------------
\section{The flag of limit sets}\label{sec:flagLS}
Letting the divergent sequence vary through the group, the target spaces of its flags give an increasing family of limit sets. The general comparison $B_\Gamma\subseteq\mathscr R_\Gamma\subseteq\Lambda_{\Kul}\subseteq\Lambda_{\Myr}$ of Section~\ref{sec:KulMyr} applies here; the rank-one observation preceding Theorem~\ref{thm:flag-of-limit-sets} adds $\Lambda_{\CoG}\subseteq B_\Gamma$. In projective dimension two, Barrera--Cano--Navarrete--Seade
\cite{BCNSFlag2026} study pseudo-projective limits and their flag geometry. Here we use all realised limiting Cartan flags; any comparison with that work must respect the divergence hypotheses imposed on the sequences. We recall four sets, all defined in Chapter~\ref{ch:limit}, and add one.

\begin{itemize}
\item $\Lambda_{\CoG}(\Gamma)$, the Conze--Guivarc'h set: the closure of the union of the images of rank-one pseudo-projective limits. If $\Gamma$ is strongly irreducible and contains a proximal element, it is the unique minimal closed $\Gamma$-invariant set (Conze--Guivarc'h).
\item $\Lambda_{\Myr}(\Gamma)=\Lambda_{\mathrm{exc}}(\Gamma)$, the complement of the equicontinuity region (Proposition~\ref{prop:kernel}).
\item $\Lambda_{\Kul}(\Gamma)$, Kulkarni's set (Definition in Section~\ref{sec:KulMyr}).
\item $\Lambda_{\mathrm{crit}}(\Gamma)$, the critical set of Section~\ref{sec:critsys} below.
\end{itemize}
In between lie the \emph{Cartan-flag levels} $\Lambda_{\mathrm Z}^{(k)}(\Gamma)$: the closure of the union of the projectivised spaces $\PP(G_{j})$ of the target flags, over divergent sequences, of projective dimension $\dim_\C\PP(G_j)\le k$; these sets collect points from the Cartan target spaces without imposing the compact-source condition. (This is the level $2k$ of Theorem~\ref{thm:L}, whose dimensions are real.) The theorem of the section, stated above, is Theorem~\ref{thm:flag-of-limit-sets}; here is its proof.

\begin{proof}[Proof of Theorem~\ref{thm:flag-of-limit-sets}]
(a) $\Lambda_{\mathrm Z}^{(k)}$ is closed and invariant and increasing in $k$ by definition. A target space $\PP(G_j)$ is a point precisely when $j=1$ and the first block has vector dimension $b_1=1$. By Proposition~\ref{prop:SP}, $b_1$ is the rank of the first pseudo-projective limit and $G_1$ is its image. Every first limit can be retained in a Cartan extraction, so $\Lambda_{\mathrm Z}^{(0)}$ is the closure of the rank-one images, namely $\Lambda_{\CoG}$. The last equality $\Lambda_{\mathrm Z}^{(n-1)}=\Lambda_{\Myr}$ is Proposition~\ref{prop:ends}(b): the target spaces of dimension $<n$ are nested and have union the proper projective subspace $\PP(G_{k-1})$, which for $(\gamma_m)$ is $\PP(F_2)$ of $(\gamma_m^{-1})$, and $\Lambda_{\mathrm{exc}}=\bigcup\PP(F_2)=\bigcup_\tau\PP(\ker\tau)$ by Proposition~\ref{prop:kernel}. The equality $\Lambda_{\Myr}=\C\PP^n\setminus\Eq$ follows from Proposition~\ref{prop:kernel}, and properness on $\Eq=\Omega_{\Myr}$ is Theorem~\ref{thm:exc}.
(b) Conze--Guivarc'h and Abels--Margulis--Soifer~\cite{ConzeGuivarch,AbelsMargulisSoifer}.
(c) Proposition~\ref{r3crit:cover} gives the closed invariant vertex cover and properness; Theorem~\ref{thm:relational} gives its trichotomy and $\mathscr R_\Gamma\subseteq\Lambda_{\mathrm{crit}}$. Together with the rank-one observation above, this proves $\Lambda_{\CoG}\subseteq\mathscr R_\Gamma\subseteq\Lambda_{\mathrm{crit}}$.
(d) Proposition~\ref{prop:L2} gives properness on $\Omega_{\Kul}$ and the formula for $\Lambda_{\Kul}$. The general inclusion $\Lambda_{\Kul}\subseteq\Lambda_{\Myr}$ is proved in Section~\ref{sec:KulMyr}; the projective identity $\Lambda_{\Myr}=\C\PP^n\setminus\Eq$ gives $\Eq\subseteq\Omega_{\Kul}$. The equality criterion is a restatement of this inclusion, and the last sentence is Example~\ref{ex:cp5}.
(e) The rank-one observation above gives $\Lambda_{\CoG}\subseteq B_\Gamma$. The general chain $B_\Gamma\subseteq\mathscr R_\Gamma\subseteq\Lambda_{\Kul}\subseteq\Lambda_{\Myr}$ is proved in Section~\ref{sec:KulMyr}. If $\Omega_{\mathrm{crit}}$ is the largest domain, Theorem~\ref{thm:relational} gives $\Lambda_{\mathrm{crit}}=\mathscr R_\Gamma$. Without that property, $\Lambda_{\mathrm{crit}}$ may exceed $\Lambda_{\Myr}$, as shown by Proposition~\ref{prop:oversized} and its arithmetic realisation, Proposition~\ref{r2arith:existence}.
(f) Theorem~\ref{thm:relational}(i)--(iv).
\end{proof}

The largest-domain hypothesis is essential here: without it the critical wall need not lie below the Kulkarni set (Proposition~\ref{r2arith:existence}). Even the Kulkarni set itself can fall strictly between consecutive Cartan-flag levels. We begin with a cyclic example, where both flags can be written down and the source of the difference can be seen.

\begin{example}[A cyclic group in $\C\PP^5$ whose Kulkarni set is not a Cartan level]\label{ex:cp5}
Let $V=E\oplus F$, $E=\langle e_1,e_2,h\rangle$, $F=\langle f_1,f_2,f_3\rangle$, and
\[
g=\mathrm{diag}\bigl(4e^{N_2},\,4,\,\tfrac14e^{N_3}\bigr)\quad\text{on }\langle e_1,e_2\rangle\oplus\langle h\rangle\oplus F,
\]
with $N_2e_2=e_1$ and $N_3f_3=f_2$, $N_3f_2=f_1$. For $\Gamma=\langle[g]\rangle$ the singular scales of $g^m$ are, up to constants, $4^mm,\ 4^m,\ 4^m/m,\ 4^{-m}m^2,\ 4^{-m},\ 4^{-m}m^{-2}$, so every adjacent ratio diverges and the type is the full flag. The target flags of $g^{m}$ and $g^{-m}$ are
\[
\langle e_1\rangle\subset\langle e_1,h\rangle\subset E\subset E{+}\langle f_1\rangle\subset A\subset V,\qquad
\langle f_1\rangle\subset\langle f_1,f_2\rangle\subset F\subset C\subset F{+}\langle e_1,h\rangle\subset V,
\]
with $A=E+\langle f_1,f_2\rangle$ and $C=F+\langle e_1\rangle$. One finds
\[
\Lambda_{\Myr}=\PP(A)\cup\PP(F{+}\langle e_1,h\rangle),\qquad
\Lambda_{\Kul}=\PP(A)\cup\PP(C),
\]
so that, with $\Lambda_{\mathrm Z}^{(k)}$ the union of the target spaces of projective dimension $\le k$ (here $\Lambda_{\mathrm Z}^{(3)}=\PP(E{+}\langle f_1\rangle)\cup\PP(C)$ and $\Lambda_{\mathrm Z}^{(4)}=\Lambda_{\Myr}$),
\[
\Lambda_{\mathrm Z}^{(3)}\ \subsetneq\ \Lambda_{\Kul}\ \subsetneq\ \Lambda_{\mathrm Z}^{(4)}.
\]
The fixed-source set is $S=\PP\langle e_1,h\rangle\cup\{[f_1]\}$. Indeed, these are the projective points fixed by nonzero powers of $g$, and every fixed-source accumulation point is among them by the Jordan form of $g$. For positive powers the deepest source stratum is $\{[f_1]\}$ and is excluded. The fifth source stratum contains the admissible source $[f_2]$, whose destination is $\PP(A)$. For negative powers the fifth source stratum, the first whose destinations reach outside $\PP(C)$, is $\PP\langle e_1,h\rangle\setminus\{[e_1]\}$. This stratum and the sixth, $\{[e_1]\}$, both lie in $S$. The fourth source stratum contains the admissible source $[e_2]$, whose destination is $\PP(C)$. The block formula and the compact-source condition give the stated $\Lambda_{\Kul}$. In particular $[e_2+f_2]\in\Lambda_{\Kul}\setminus\Lambda_{\mathrm Z}^{(3)}$ and $[h+f_3]\in\Lambda_{\mathrm Z}^{(4)}\setminus\Lambda_{\Kul}$.
\end{example}

%--------------------------------------------------------------
\section{Kulkarni's set in a refined carrier system}
\label{r2kul:section}

In the last example the two flags are known, yet one of their spaces contributes only part of the Kulkarni set. The reason lies in the sources: Kulkarni's construction excludes
the fixed-source set before taking destinations. We keep that restriction
and the graded map determining each destination. The resulting carriers
recover the Kulkarni set even when dimensions alone cannot pick it out.

Throughout this section, $V=\C^d$, $d\ge2$, $X=\PP(V)$ and
$\Gamma<\PSL(V)$ is discrete. Write
\[
 S=S_\Gamma=B_\Gamma=L_0(\Gamma)\cup L_1(\Gamma),\qquad
 C=\Lambda_{\Myr}(\Gamma)=X\setminus\Eq(\Gamma).
\]
For $A\subset X$, $\D_\Gamma(A)$ denotes the set of destinations of
sources in $A$, with the moving-source convention of
Definition~\ref{def:Dgamma}. Proposition~\ref{prop:L2} gives
\begin{equation}\label{r2kul:compact-source}
 L_2(\Gamma)=\overline{\D_\Gamma(X\setminus S)},\qquad
 \Lambda_{\Kul}(\Gamma)=S\cup\overline{\D_\Gamma(X\setminus S)}.
\end{equation}
The order matters: we exclude $S$ before taking the closure. Indeed, a
realisation $x_m\to x\notin S$ eventually lies in a compact neighbourhood
of $x$ disjoint from $S$; conversely, sources in a compact subset of
$X\setminus S$ have subsequential limits there.

\subsection{Exact comparison with the Myrberg set}

Let $B_{\rm ext}\Gamma$ be the set of limits of escaping sequences of
$\Gamma$ in the exterior compactification of Section~\ref{sec:ext}.
For $\xi\in B_{\rm ext}\Gamma$, its graded datum is
\[
 F_1=V\supsetneq\cdots\supsetneq F_{k+1}=0,\qquad
 0=G_0\subsetneq\cdots\subsetneq G_k=V,\qquad
 [A_j]:F_j/F_{j+1}\longrightarrow G_j/G_{j-1}.
\]
Here $k\ge2$. For $x=[v]\in\PP(F_j)\setminus\PP(F_{j+1})$, put
\begin{equation}\label{r2kul:carrier}
 W_\xi(x)=\pi_j^{-1}\bigl(A_j(\C q_jv)\bigr)\subset G_j,
 \qquad \D_\xi(x)=\PP(W_\xi(x)),
\end{equation}
where $q_j:F_j\to F_j/F_{j+1}$ and
$\pi_j:G_j\to G_j/G_{j-1}$ are the quotient maps. This is independent of
the representative of $[A_j]$, and uses no complement of $G_{j-1}$.
Theorem~\ref{thm:ext-boundary} identifies it with the whole-sequence
destination for every realising sequence.

\begin{lemma}\label{r2kul:source-restriction}
Every fixed-source accumulation point belongs to $S$. For each
$\xi\in B_{\rm ext}\Gamma$ one has $\PP(F_k)\subset S$. Consequently,
if $x\notin S$ belongs to stratum $j$, then $j<k$ and
\[
 W_\xi(x)\subset G_j\subset G_{k-1},
\]
where $G_{k-1}$ is the kernel of the leading pseudo-projective limit of
the inverse sequence.
\end{lemma}
\begin{proof}
If a fixed source lies in $L_0$, its orbit and all its limits lie in
$L_0$, since $L_0$ is closed and invariant. Otherwise every accumulation
point belongs to $L_1$ by definition. The leading inverse limit has image
$\PP(F_k)$. Every point of that image is reached from some fixed source
outside its kernel, so the first assertion puts the whole image in $S$.
Thus an admissible source cannot lie in the deepest stratum. Formula~\eqref{r2kul:carrier}
and the inverse flag description give the remaining assertions.
\end{proof}

\begin{proposition}[Exact equality criteria]\label{r2kul:exact-criteria}
One has
\[
 \Lambda_{\CoG}(\Gamma)\subseteq S\subseteq\mathscr R_\Gamma
 \subseteq\Lambda_{\Kul}(\Gamma)\subseteq C,\qquad
 \Eq(\Gamma)\subseteq\Omega_{\Kul}(\Gamma).
\]
Moreover, the following are equivalent:
\begin{enumerate}[label=\textup{(\roman*)},nosep]
\item $\Lambda_{\Kul}(\Gamma)=\Lambda_{\Myr}(\Gamma)$;
\item $C\setminus S\subseteq\overline{\D_\Gamma(X\setminus S)}$;
\item $\Omega_{\Kul}(\Gamma)=\Eq(\Gamma)$.
\end{enumerate}
These equalities hold if $C\subseteq\mathscr R_\Gamma$, and also if
$\Eq(\Gamma)$ is the largest invariant open set of proper discontinuity.
\end{proposition}
\begin{proof}
First $S\subseteq C$. Indeed, properness on $\Eq(\Gamma)$ excludes
infinite stabilisers there, hence $L_0\subset C$. If
$\gamma_mx\to y\in\Eq(\Gamma)$, invariance and openness put $x$ in
$\Eq(\Gamma)$; properness then excludes accumulation of this fixed
orbit inside that open set. This proves $L_1\subset C$.

By Proposition~\ref{prop:kernel}, $C$ is the union of the projectivised
kernels of pseudo-projective boundary limits. This union is closed:
the boundary limit set is compact, and the incidence $Tv=0$ is closed
in its product with $X$. If $(x,y)\in\D_\Gamma$ and $x\notin S$,
extract a realising sequence with exterior datum $\xi$.
Lemma~\ref{r2kul:source-restriction} puts $y$ in an inverse boundary
kernel. Thus $\overline{\D_\Gamma(X\setminus S)}\subseteq C$,
and~\eqref{r2kul:compact-source} proves the comparison. The inclusion $\mathscr R_\Gamma\subseteq\Lambda_{\Kul}$ is
Proposition~\ref{prop:L2}(c). The rank-one fixed-source observation gives
$\Lambda_{\CoG}\subseteq S$, and $S\subseteq\mathscr R_\Gamma$ is the
general fixed-source recurrence statement of Section~\ref{sec:KulMyr}.

Since $S\subset C$, formula~\eqref{r2kul:compact-source} makes
\textup{(i)} equivalent to \textup{(ii)}; taking complements gives
\textup{(iii)}. If $C\subseteq\mathscr R_\Gamma$, the inclusions force
equality. If $\Eq(\Gamma)$ is largest, properness on
$\Omega_{\Kul}(\Gamma)$ puts that open set in $\Eq(\Gamma)$, and the
universal reverse inclusion gives equality.
\end{proof}

The equality criterion compares two known sets. We may instead be given
a candidate wall and want to recognise it from the behaviour of compact
sets in its complement. The following test makes that possible.

\begin{proposition}[A compact-source recognition criterion]
\label{r3kul:quasiminimality}
Let $L\subset X=\PP(V)$ be closed and $\Gamma$-invariant. Suppose that
for every compact $K\subset X\setminus L$, every limit of
$\gamma_jx_j$, with $x_j\in K$ and $\gamma_j\in\Gamma$ pairwise
distinct, belongs to $S\cap L$. Then $\Lambda_{\Kul}(\Gamma)\subseteq L$.
If, moreover, $L$ is the closure of a union of points of $S$ and
destination carriers $\PP(W_\xi(x))$ with $x\notin S$, then
$L=\Lambda_{\Kul}(\Gamma)$.
\end{proposition}
\begin{proof}
The singleton test excludes infinite stabilisers outside $L$, so
$L_0\subset L$. Fixed-orbit limits from $L$ stay in $L$ by invariance;
those from its complement lie in $L$ by the same test. Thus $S\subset L$.
If $\gamma_jx_j\to y\notin L$ generated $L_2$, with
$x_j$ in a compact subset of $X\setminus S$, take $x_j\to x\notin S$.
The tail images lie in a compact neighbourhood of $y$ disjoint from
$L$. Applying the hypothesis to $\gamma_j^{-1}$ forces $x\in S$,
a contradiction. Closedness gives $\Lambda_{\Kul}\subset L$.
Conversely, every stated generator of $L$ belongs to $\Lambda_{\Kul}$
by~\eqref{r2kul:compact-source} and~\eqref{r2kul:carrier}; taking
closures proves equality. Compare~\cite[\S3.3.1]{CNS}.
\end{proof}

There is a useful shortcut when the Kulkarni set contains hyperplanes.
Their defining covectors can force every boundary kernel into the same
closed invariant set.

\begin{corollary}\label{r2kul:hyperplane-comparison}
If $\Lambda_{\Kul}(\Gamma)$ contains $d$ hyperplanes with linearly
independent defining covectors, then
$\Lambda_{\Kul}(\Gamma)=\Lambda_{\Myr}(\Gamma)$ and
$\Omega_{\Kul}(\Gamma)=\Eq(\Gamma)$. If $\Gamma$ is irreducible, one
hyperplane contained in $\Lambda_{\Kul}(\Gamma)$ suffices.
\end{corollary}
\begin{proof}
Apply Lemma~\ref{lem:myrberg-hyperplane-propagation} to the closed
invariant set $\Lambda_{\Kul}(\Gamma)$, and use
Proposition~\ref{r2kul:exact-criteria} for the other inclusion. Under
irreducibility, the span of the defining covectors of the orbit of one
hyperplane is the whole dual space; a finite subcollection is a basis.
\end{proof}

\subsection{Carriers and their lower ideals}

A destination contains all the projective subspaces lying in its carrier.
We can remember this containment before taking the union of their points.
This is the reason for using lower ideals.

Let $\mathscr P_\Gamma$ be the partially ordered set, under inclusion,
of nonzero vector subspaces consisting of all lines in $V$, all source
and target Cartan flag spaces of escaping sequences in $\Gamma$, all
carriers~\eqref{r2kul:carrier}, and all pseudo-projective boundary
kernels. A \emph{refined level} is a lower ideal
$\mathcal I\subseteq\mathscr P_\Gamma$, with shadow
\[
 \operatorname{Sh}(\mathcal I)
 =\overline{\bigcup_{W\in\mathcal I}\PP(W)}.
\]
For a collection $\mathcal A\subseteq\mathscr P_\Gamma$, write
$\downarrow\mathcal A=\{U\in\mathscr P_\Gamma:
U\subseteq W\text{ for some }W\in\mathcal A\}$. Lower ideals form a
complete distributive lattice under union and intersection. They do not
close a family under vector-space sums: two points in a limit set need
not bring the projective line joining them into that set.

For $x\in X$, write $\ell_x$ for its representing line, and let
$\operatorname{Lim}_{\rm SP}(\Gamma)$ denote the singular
pseudo-projective limits of escaping sequences. Define
\begin{align*}
 \mathcal I_{\CoG}
 &=\downarrow\{\operatorname{Im}T:
       [T]\in\operatorname{Lim}_{\rm SP}(\Gamma),\ \operatorname{rank}T=1\},\\
 \mathcal I_{\rm fix}&=\downarrow\{\ell_x:x\in S\},\\
 \mathcal I_{\Kul}&=\downarrow\bigl(\{\ell_x:x\in S\}
       \cup\{W_\xi(x):\xi\in B_{\rm ext}\Gamma,\ x\notin S\}\bigr),\\
 \mathcal I_{\Myr}&=\downarrow\{\ker T:
       [T]\in\operatorname{Lim}_{\rm SP}(\Gamma)\}.
\end{align*}
The inverse orientation is already included, since inverses of escaping
sequences in $\Gamma$ are again escaping sequences in $\Gamma$.

\begin{theorem}[The refined Kulkarni level]\label{r2kul:ideal-level}
These ideals satisfy
\[
 \mathcal I_{\CoG}\subseteq\mathcal I_{\rm fix}
 \subseteq\mathcal I_{\Kul}\subseteq\mathcal I_{\Myr},
\]
and their shadows are, respectively,
\[
 \Lambda_{\CoG}(\Gamma),\qquad S,\qquad
 \Lambda_{\Kul}(\Gamma),\qquad\Lambda_{\Myr}(\Gamma).
\]
\end{theorem}
\begin{proof}
A rank-one image is a fixed-source limit, so its line belongs to
$\mathcal I_{\rm fix}$ by Lemma~\ref{r2kul:source-restriction}.
The second inclusion is definitional. Each admissible carrier lies in
an inverse boundary kernel by the same lemma. Every point of $S$ lies in
some boundary kernel by Proposition~\ref{r2kul:exact-criteria} and the
closed kernel-union description of $C$. This proves the last inclusion
at the level of ideals, not merely at the level of shadows.

The Conze--Guivarc'h shadow is its rank-one definition; the fixed-source
shadow is $S$, which is closed. The kernel shadow is $C$ by
Proposition~\ref{prop:kernel}. Finally, the exterior destination formula
and~\eqref{r2kul:compact-source} give
\begin{equation}\label{r2kul:L2-carriers}
 L_2(\Gamma)=\overline{\bigcup_{\xi\in B_{\rm ext}\Gamma}
                  \ \bigcup_{x\notin S}\PP(W_\xi(x))}.
\end{equation}
Taking lower closure adds only subspaces of the generators and hence
leaves their projective union unchanged. This proves the Kulkarni shadow
identity.
\end{proof}

The four ideals end at the kernel ideal, but the lattice continues.
For example, it contains a level with shadow $X$. Before taking closures,
the point union preserves intersections of ideals. A point belongs to
both unions exactly when its representing line belongs to both ideals.
After taking closures this equality need not hold. At the other end,
strong irreducibility and the existence of a proximal element give
minimality of the rank-one shadow, as recalled in Section~\ref{sec:flagLS}.

\subsection{Typed Grassmannian families}

We can also gather carriers of the same dimension in a Grassmannian.
Two numbers enter the destination formula: the dimension accumulated in
the target before the active block, and the dimension of that block.
For a source in stratum $j$, set
\[
 r=\dim G_{j-1},\qquad s=\dim(G_j/G_{j-1}).
\]
The destination has projective dimension $r$, and its active quotient has
vector dimension $s$. An admissible source has $j<k$, so its type belongs
to the finite set
\[
 \mathcal T_d=\{(r,s):r\ge0,\ s\ge1,\ r+s<d\}.
\]
For $(r,s)\in\mathcal T_d$, define
\begin{equation}\label{r2kul:typed-family}
 \mathcal C_{r,s}(\Gamma)=
 \overline{\left\{W_\xi(x):
 \begin{array}{l}
 \xi\in B_{\rm ext}\Gamma,\quad
 x\in\PP(F_j)\setminus\PP(F_{j+1}),\quad x\notin S,\\
 \dim G_{j-1}=r,\quad\dim(G_j/G_{j-1})=s
 \end{array}\right\}}
 \ \subseteq\operatorname{Gr}(r+1,V).
\end{equation}
The closure is taken in this fixed Grassmannian. The triple $(\xi,x,j)$ comes from one source and one exterior datum.
Choosing its entries from separate marginal limit sets would lose this
condition; likewise $A_j$ must be the realised graded map, rather than
an arbitrary map with the same source and target dimensions.

For a closed Grassmannian family $\mathcal C$, write
\[
 |\mathcal C|=\{y\in X:\ y\in\PP(W)\text{ for some }W\in\mathcal C\}.
\]
The empty family has empty shadow. We shall use the following elementary
consequence of Grassmannian incidence.

\begin{lemma}\label{r2kul:incidence}
If $\mathcal A\subseteq\operatorname{Gr}(a,V)$, then
\[
 |\overline{\mathcal A}|
 =\overline{\bigcup_{W\in\mathcal A}\PP(W)}.
\]
In particular, the shadow of a closed Grassmannian family is closed.
\end{lemma}
\begin{proof}
The incidence $\{(W,y):y\in\PP(W)\}$ is closed in the product of two
compact spaces. Its projection over $\overline{\mathcal A}$ is therefore
closed and contains the right-hand side. Conversely, if $W_m\in\mathcal A$
tends to $W$ and $y\in\PP(W)$, orthogonal projections onto $W_m$ of a
nonzero representative of $y$ give, for large $m$, points
$y_m\in\PP(W_m)$ tending to $y$. Thus every point on the left belongs
to the right-hand closure.
\end{proof}

\begin{theorem}[The typed Kulkarni filtration]\label{r2kul:typed-levels}
Each $\mathcal C_{r,s}(\Gamma)$ is compact and $\Gamma$-invariant, and
each of its members is contained in a pseudo-projective boundary kernel.
Moreover,
\begin{equation}\label{r2kul:L2-types}
 L_2(\Gamma)=\bigcup_{(r,s)\in\mathcal T_d}|\mathcal C_{r,s}(\Gamma)|.
\end{equation}
Define, for $0\le q\le d-2$,
\begin{equation}\label{r2kul:Kq}
 \mathcal K_q(\Gamma)=S\cup
 \bigcup_{\substack{(r,s)\in\mathcal T_d\\r\le q}}
                  |\mathcal C_{r,s}(\Gamma)|.
\end{equation}
These are closed invariant sets satisfying
\begin{equation}\label{r2kul:chain}
 \Lambda_{\CoG}(\Gamma)\subseteq S=\mathcal K_0(\Gamma)
 \subseteq\cdots\subseteq\mathcal K_{d-2}(\Gamma)
 =\Lambda_{\Kul}(\Gamma)\subseteq\Lambda_{\Myr}(\Gamma).
\end{equation}
Each $\mathcal K_q$ is the shadow of the lower ideal generated by the
fixed-source lines and the actual admissible carriers with $r\le q$.
\end{theorem}
\begin{proof}
Compactness follows from the Grassmannian closure. If $g\in\Gamma$,
left multiplication sends a realising sequence $(\gamma_m)$ to
$(g\gamma_m)$, whose destination relation is $g\D_\xi$ on the target
factor. Hence its carrier at $x$ is $gW_\xi(x)$. Bounded multiplication
preserves the Cartan block dimensions (Lemma~\ref{lem:bdd}), so $r,s$
are unchanged. This proves invariance of the generating family and then
of its closure.

Let $W_m\to W$ be generating carriers of one fixed type. By
Lemma~\ref{r2kul:source-restriction}, choose boundary limits $[T_m]$ with
$W_m\subseteq\ker T_m$. Normalise $\|T_m\|=1$ and extract
$T_m\to T\ne0$. Compactness of the pseudo-projective boundary puts
$[T]$ in that boundary. Closed incidence gives $TW=0$, so the assertion
about kernels persists through the Grassmannian closure.

There are finitely many types. Apply Lemma~\ref{r2kul:incidence} to each
of them in~\eqref{r2kul:L2-carriers}; the resulting finite union is
closed, proving~\eqref{r2kul:L2-types}. Equations~\eqref{r2kul:Kq} and
\eqref{r2kul:compact-source} now give the terminal equality and
closedness; invariance and the remaining inclusions follow from the
preceding theorem. If $r=0$, the source lies in the first stratum. Its
carrier is the single ordinary fixed-source limit, hence lies in $S$.
Closedness of $S$ and Lemma~\ref{r2kul:incidence} show that all
$|\mathcal C_{0,s}|$ remain in $S$, proving $\mathcal K_0=S$.
Finally, taking the lower ideal of the actual generators does not change
their projective union, and taking its shadow supplies precisely the
closures already present in~\eqref{r2kul:Kq}.
\end{proof}

Thus $\mathcal K_q$ grows by admissible destination dimension, whereas
$\Lambda_{\mathrm Z}^{(q)}$ grows by flag-space dimension without
restricting the sources. The two indices $r,s$ tell us more than $q$
alone: they distinguish the space already accumulated in the target
from the active quotient. The smallest nontrivial projective example
shows how these distinctions arise.

\begin{example}\label{r2kul:diagonal}
For $d=3$ and
$\Gamma=\langle[\diag(e^2,1,e^{-2})]\rangle$, the two orientations give
\[
 \Lambda_{\CoG}=\{[e_1],[e_3]\},\qquad
 S=\mathcal K_0=\{[e_1],[e_2],[e_3]\}.
\]
The admissible sources of the middle positive stratum are
$[ae_2+be_3]$ with $ab\ne0$; their carrier is
$\langle e_1,e_2\rangle$. Negative powers give
$\langle e_2,e_3\rangle$. Therefore
\[
 \mathcal K_1=\Lambda_{\Kul}=\Lambda_{\Myr}
 =\PP\langle e_1,e_2\rangle\cup\PP\langle e_2,e_3\rangle.
\]
The level $S$ is strictly between the rank-one and Kulkarni sets, whereas
the scalar levels are the two extreme points and the union of the two
lines. Although the deepest source in this example has destination $X$,
it belongs to $S$ and contributes no admissible carrier.
\end{example}

\section{One recurrent hyperplane}
\label{r2kul:recurrent-hyperplanes}

In the cyclic example the Kulkarni and Myrberg sets agree, but their
common complement is not the largest domain. Suppose now that every
point of one hyperplane is recurrent, in the moving-point sense. The
realising sequence may vary from point to point; nevertheless the orbit
of this one hyperplane can determine the whole Myrberg set.

\begin{proposition}\label{r2kul:recurrent-hyperplane}
Suppose that a hyperplane $H\subseteq\mathscr R_\Gamma$ satisfies
$\bigcap_{g\in\Gamma}gH=\varnothing$. Then
\[
 \mathscr R_\Gamma=\Lambda_{\Myr}(\Gamma)
 =\Lambda_{\Kul}(\Gamma)
 =\overline{\bigcup_{g\in\Gamma}gH},
\]
and $\Eq(\Gamma)=\Omega_{\Kul}(\Gamma)$ is the largest invariant open
set of proper discontinuity. If $\Gamma$ is irreducible, one hyperplane
contained in $\mathscr R_\Gamma$ suffices.
\end{proposition}
\begin{proof}
Put $F=\overline{\bigcup_{g\in\Gamma}gH}$. Closedness and invariance of
$\mathscr R_\Gamma$ give $F\subseteq\mathscr R_\Gamma$, while
properness on $\Eq(\Gamma)$ gives $\mathscr R_\Gamma\subseteq C$.
The empty-intersection condition says exactly that the covectors defining
the orbit hyperplanes span $V^*$. A finite subcollection already spans,
so Lemma~\ref{lem:myrberg-hyperplane-propagation} gives $C\subseteq F$.
Thus $F=\mathscr R_\Gamma=C$. Proposition~\ref{r2kul:exact-criteria}
adds the Kulkarni equality. Every domain of proper discontinuity avoids
$\mathscr R_\Gamma$, so their common complement $\Eq(\Gamma)$ is largest.
Under irreducibility, the intersection of the vector hyperplanes in this
orbit is a proper invariant subspace and hence is zero.
\end{proof}

\begin{samepage}
\begin{proposition}[A rank-one kernel criterion]\label{r2kul:rank-one-largest}
Suppose that $\Gamma$ is irreducible and admits a rank-one
pseudo-projective boundary limit $[T_0]$. Put $H_0=\PP(\ker T_0)$.
The following are equivalent:
\begin{enumerate}[label=\textup{(\roman*)},nosep]
\item $\Eq(\Gamma)$ is the largest invariant open set of proper discontinuity;
\item $\mathscr R_\Gamma=\Lambda_{\Myr}(\Gamma)$;
\item some hyperplane is contained in $\mathscr R_\Gamma$;
\item $H_0\subseteq\mathscr R_\Gamma$.
\end{enumerate}
When these conditions hold,
\[
 \mathscr R_\Gamma=\Lambda_{\Myr}(\Gamma)
 =\Lambda_{\Kul}(\Gamma)
 =\overline{\bigcup_{g\in\Gamma}gH_0},\qquad
 \Omega_{\max}(\Gamma)=\Eq(\Gamma)=\Omega_{\Kul}(\Gamma).
\]
\end{proposition}
\begin{proof}
By Theorem~\ref{thm:relational}, the union of all domains is
$X\setminus\mathscr R_\Gamma$. Thus \textup{(i)} implies
\textup{(ii)}. Since $T_0$ has rank one, $H_0$ is a hyperplane in
$C$, so \textup{(ii)} implies \textup{(iv)}, which implies
\textup{(iii)}. Proposition~\ref{r2kul:recurrent-hyperplane} gives
\textup{(iii)}$\Rightarrow$\textup{(ii)}. Finally, \textup{(ii)} makes the domain $\Eq(\Gamma)=X\setminus\mathscr R_\Gamma$
contain every other such open set, proving \textup{(i)}. Applying the
same proposition to $H_0$ gives the orbit formula and the remaining
equalities.
\end{proof}
\end{samepage}

\sectionpart{The critical system and recurrence}

\section{The critical system}\label{sec:critsys}
A recurrent hyperplane can determine the largest domain, but in many
examples the recurrent carriers have smaller dimension, and the relevant
spaces lie in the middle of the Cartan flag. The destination formula suggests pairing a
target space there with a complementary space from the inverse sequence:
if a source avoids the second space, its destination lies in the first.
We choose the orientations for which this middle target is proper, and
keep the complementary inverse space paired with that choice.

\begin{example}[The critical system of $\PSL(n+1,\C)$]\label{ex:critsys}
Let $\Gamma<\PSL(V)$ be infinite and discrete, $d=n+1\ge2$, and put
$m=\lfloor n/2\rfloor$, $c=m+1=\lceil d/2\rceil$.
Write $0=r_0<r_1<\cdots<r_k=d$ for the endpoints of an escaping
Cartan sequence $h_\bullet$, and use its flags $G_j,F_j$ as in
Theorem~\ref{thm:blockform}, with $G_0=F_{k+1}=0$.
Its orientation is \emph{admissible} when $r_{k-1}\ge c$.
For an admissible orientation put
\[
 t_h=\min\{j:r_j\ge c\}<k,\qquad
 A_h=G_{t_h-1},\quad B_h=G_{t_h},\quad
 q_{\rm or}(h_\bullet)=r_{t_h}-1.
\]
At least one orientation is admissible: failure for both would give
$b_1,b_k\ge d-c+1$, contradicting $b_1+b_k\le d$.
Let $\mathfrak A_\Gamma$ consist of all admissible oriented extractions
of divergent sequences in $\Gamma$, and define
\begin{equation}\label{r3crit:levels}
 q_\Gamma=\min_{h_\bullet\in\mathfrak A_\Gamma}
          q_{\rm or}(h_\bullet),\qquad
 q_\Gamma^*=n-q_\Gamma-1,\qquad p_\Gamma=q_\Gamma+1.
\end{equation}
Thus $m\le q_\Gamma\le n-1$, and the minimum is attained. This is
equivalently the minimum, over unoriented sequences, of their smaller defined
orientation level. The first-middle crossing of an inadmissible sequence is
$V$; including those values would not change this scalar minimum.

For $\dim A<s+1\le\dim B$ set
\[
 \Sigma_s(A,B)=\{L\in\mathrm{Gr}_s(\PP(V)):
                   \PP(A)\subseteq L\subseteq\PP(B)\}.
\]
Here $\mathrm{Gr}_s(\PP(V))$ parametrises projective $s$-planes and
$\PP(0)=\varnothing$. For every $h_\bullet\in\mathfrak A_\Gamma$,
\begin{equation}\label{r3crit:absorption}
 \dim A_h<c\le p_\Gamma\le\dim B_h.
\end{equation}
In the increasing target flag of the inverse sequence, let
$A_h^*\subset B_h^*$ be the consecutive spaces satisfying
$\dim A_h^*<d-p_\Gamma\le\dim B_h^*$.
In the block calculation of Proposition~\ref{prop:ext-data}, the active
exterior image consists of wedges containing the lower space and using
only the upper space. Its decomposable lines therefore identify
$\Sigma_{q_\Gamma}(A_h,B_h)$ and
$\Sigma_{q_\Gamma^*}(A_h^*,B_h^*)$ with the decomposable loci of the
active target blocks at degrees $p_\Gamma$ and $d-p_\Gamma$ for
$h_\bullet$ and $h_\bullet^{-1}$ respectively. Both are nonempty.
In particular $q_\Gamma$ is the largest projective level at least $m$
active in every admissible first-middle crossing: these common levels
are exactly $m,\ldots,q_\Gamma$. Indeed,~\eqref{r3crit:absorption}
gives all those levels, while a \emph{primitive realiser}, with
$\dim B_h=q_\Gamma+1$, excludes every larger one and has the singleton
family $\{\PP(B_h)\}$ at level $q_\Gamma$.

The two families are indexed by the \emph{same admissible oriented sequences}:
\begin{align*}
 \mathcal Z_{q_\Gamma}(\Gamma)
 &=\overline{\bigcup_{h_\bullet\in\mathfrak A_\Gamma}
               \Sigma_{q_\Gamma}(A_h,B_h)},\\
 \mathcal Z_{q_\Gamma}^{\vee}(\Gamma)
 &=\overline{\bigcup_{h_\bullet\in\mathfrak A_\Gamma}
               \Sigma_{q_\Gamma^*}(A_h^*,B_h^*)}.
\end{align*}
The paired inverse family is retained even when the inverse orientation
is inadmissible; no level-$q_\Gamma$ target family is added from an
inadmissible orientation. For a family $\mathcal A$ of projective planes
write $|\mathcal A|=\bigcup_{L\in\mathcal A}L$. Define
\begin{equation}\label{r3crit:wall}
 \Lambda_{\mathrm{crit}}(\Gamma)
   =|\mathcal Z_{q_\Gamma}(\Gamma)|
      \cup|\mathcal Z_{q_\Gamma}^{\vee}(\Gamma)|,
 \qquad \Omega_{\mathrm{crit}}=\PP(V)\setminus\Lambda_{\mathrm{crit}}.
\end{equation}
Each nonempty $\Sigma_s(A,B)$ above has shadow $\PP(B)$: every vector
of $B$ lies in an $(s+1)$-space between $A$ and $B$. Compactness and
local triviality of incidence therefore give the useful formula
\begin{equation}\label{r3crit:shadows}
 \Lambda_{\mathrm{crit}}(\Gamma)
   =\overline{\bigcup_{h_\bullet\in\mathfrak A_\Gamma}
        \bigl(\PP(B_h)\cup\PP(B_h^*)\bigr)}.
\end{equation}
\end{example}

\begin{proposition}\label{r3crit:cover}
The paired critical wall~\eqref{r3crit:wall} is closed, invariant,
and a dynamical vertex cover. Hence the action on
$\Omega_{\mathrm{crit}}$ is properly discontinuous.
\end{proposition}
\begin{proof}
The two Grassmannian families are compact, so closedness follows from
compactness of incidence. For invariance fix $g\in\Gamma$.
Multiplication on either side by $g$ changes every singular value by
bounded factors; thus it preserves each divergent boundary gap and
keeps all ratios inside each block bounded. After extraction, the
target boundary flags of $gh_j$ are $gG_i$, and the slow source
boundary flags of $h_jg^{-1}$ are $gF_i$. One can read this directly
from the images and kernels of the rank-one exterior limits at the
block endpoints (Proposition~\ref{prop:ext-data}). Both modified sequences
retain the original ordered block dimensions and admissibility.
Consequently $g$ transports the target families and, using the second
modified sequence, the paired inverse families into their respective
unions. Taking closures and then applying $g^{-1}$ proves invariance.

Only closures in fixed Grassmannians are used here. For fixed endpoint
dimensions, $\Sigma_s(A_j,B_j)$ converges in the Hausdorff topology to
$\Sigma_s(A,B)$ when $(A_j,B_j)\to(A,B)$: local graph coordinates
identify these families with $\mathrm{Gr}_{s+1-\dim A}(B/A)$.
This does not assert that an arbitrary limit of sequences is realised with
the same block type. Such a diagonal realisation requires uniform
bounds on the internal singular ratios; the fixed multiplications
above preserve those bounds for the one sequence being transported.

Take a pair $(x,y)\in\D_\Gamma$ and extract its witnessing sequence.
If necessary invert it and interchange $x,y$, so the orientation is
admissible. Write $t=t_h$. Since $p_\Gamma\le r_t$,
\[
 \dim F_{t+1}=d-r_t\le d-p_\Gamma,
 \qquad F_{t+1}\subseteq B_h^*.
\]
If $x\notin\PP(B_h^*)$, then its source stratum is
$\PP(F_j\setminus F_{j+1})$ for some $j\le t$.
Theorem~\ref{thm:blockform} gives
$y\in\PP(G_j)\subseteq\PP(B_h)$.
Thus either $x$ or $y$ lies in~\eqref{r3crit:shadows}, including the
endpoint case $p_\Gamma=r_t$, where $B_h^*=F_{t+1}$.
The properness criterion of Theorem~\ref{thm:relational} gives the conclusion.
\end{proof}

\paragraph{The all-sequence variant.}
Instead of fixing an orientation, one may take
the closure of the first-middle target shadows
of \emph{all} escaping sequences, allowing the crossing to be $V$; denote
this larger wall by $\Lambda_{\mathrm{crit}}^{\mathrm{all}}$.
The same block argument makes it a vertex cover. Moreover
\begin{equation}\label{r3crit:comparison}
 \Lambda_{\mathrm{crit}}\subseteq
 \Lambda_{\mathrm{crit}}^{\mathrm{all}}.
\end{equation}
Indeed $B_h$ is one of its crossings, and
$d-p_\Gamma\le c$, so $B_h^*$ is contained in the inverse sequence's
first-middle crossing. The inclusion can be strict. For
$\Gamma=\langle[\operatorname{diag}(2,1,1)]\rangle$, put
$z=[e_1]$ and $H=\PP\langle e_2,e_3\rangle$.
Positive powers have type $(1,2)$ and are inadmissible; negative
powers have type $(2,1)$ and are admissible. Thus
\[
 q_\Gamma=1,\quad q_\Gamma^*=0,\qquad
 \mathcal Z_1(\Gamma)=\{H\},\quad
 \mathcal Z_1^\vee(\Gamma)=\{\{z\}\},\qquad
 \Lambda_{\mathrm{crit}}=H\cup\{z\},\quad
 \Lambda_{\mathrm{crit}}^{\mathrm{all}}=\C\PP^2.
\]
For $n=1$ the only escaping type is $(1,1)$, both orientations are
admissible, and $q_\Gamma=q_\Gamma^*=0$; the two singleton families
recover the classical limit set. More generally, a fully regular
sequence forces $q_\Gamma=m$. This scalar fact imposes no uniform type on
the other sequences and does not prevent a paired complementary shadow
from being all of $\PP(V)$ (Proposition~\ref{prop:oversized}).

Toledo-Acosta~\cite[Definition~3.4.1]{ToledoAcosta2019} introduced
a middle-block wall, a precursor of the all-sequence variant. The
cutoff $\lfloor d/2\rfloor$ in that definition requires adjustment in
odd vector dimension. For instance, when $d=3$, the powers of
$\operatorname{diag}(2,1,1/2)$ in the two orientations give only the
first and last coordinate points at that cutoff, leaving the fixed
middle coordinate point in the complement. Our paired wall uses
$\lceil d/2\rceil$, admissible orientations and complementary inverse
spaces. The preceding cyclic example distinguishes it from the
all-sequence wall with the same cutoff.

\section{When the critical domain is the largest}\label{sec:crit-proj}

The paired wall gives a proper action on its complement. By
Theorem~\ref{thm:relational}, this domain is the largest one once every
point of the wall returns. The construction has two parts: make one critical carrier recurrent,
then use the orbit of its flag to reach all the others. The exterior
powers control this orbit; the middle block controls the return.

Put $d=\dim_{\mathbb C}V$. Recurrence is always meant with moving
sources (Definition~\ref{def:Dgamma}); the recurrent set is closed and
invariant (Lemma~\ref{lem:Dbasic}), and every domain of proper
discontinuity avoids it (Theorem~\ref{thm:relational}).

Fix a Hermitian norm. For a singular-value decomposition
\[
 h=U\operatorname{diag}(\sigma_1,\ldots,\sigma_d)V^*,\qquad
 \sigma_1\ge\cdots\ge\sigma_d>0,
\]
and integers $1\le a<p<d$ with $a+p=d$, define the finite critical
flags by
\[
 \widehat A_h^+=U\langle e_1,\ldots,e_a\rangle,
 \quad \widehat B_h^+=U\langle e_1,\ldots,e_p\rangle,
\]
\[
 \widehat A_h^-=V\langle e_{p+1},\ldots,e_d\rangle,
 \quad \widehat B_h^-=V\langle e_{a+1},\ldots,e_d\rangle.
\]
The minus flag is the slow source flag; it is also the target flag
of $h^{-1}$, with the singular values put in decreasing order.

\paragraph{The hypotheses.}
We first ask that the middle crossing have the same dimensions along
every escaping sequence, including inverse sequences. Uniform self-opposed
critical type $(a,p)$ means that every divergent
Cartan extraction has its first middle crossing immediately after
dimension $a$ and at dimension $p$, with $a+p=d$. In singular values this
gives
\begin{equation}\label{eq:rec-critical-boundary-gaps}
 \frac{\sigma_a(h_j)}{\sigma_{a+1}(h_j)}\longrightarrow\infty,
 \qquad
 \frac{\sigma_p(h_j)}{\sigma_{p+1}(h_j)}\longrightarrow\infty,
 \qquad
 \sup_j\frac{\sigma_{a+1}(h_j)}{\sigma_p(h_j)}<\infty.
\end{equation}
These statements are made after extraction of all adjacent singular
ratios; the last assertion says that the crossing block does not split.
Let $\mathscr F$ be the closure of the target critical flags of such
sequences. The inverse sequence has the same middle crossing $(a,p)$, so its target flags, and
therefore all limiting slow source flags, also belong to $\mathscr F$.
Both orientations are admissible: reversal changes the consecutive
crossing $(a,p)$ to $(d-p,d-a)=(a,p)$. Hence $q_\Gamma=p-1$.
The paired complementary degree is $a$, whose inverse target space
is contained in the inverse critical $p$-space. Thus the
paired wall~\eqref{r3crit:wall}, and also the all-sequence variant in this
uniform case, is
\[
 \Lambda_{\mathrm{crit}}
   =\bigcup_{(A\subset B)\in\mathscr F}\mathbb P(B).
\]
This union is closed by compactness of the incidence correspondence.

Let $W=W_{a,p}^{\mathrm{fl}}$ be the irreducible Cartan component of
$\bigwedge^aV\otimes\bigwedge^pV$, and let
$\jmath:\operatorname{Fl}(a,p;V)\to\mathbb P(W)$ be its
Segre--Pl\"ucker embedding. Write $\tau$ for the induced projective
representation. The critical module is
\[
 W_{a,p}^{\mathrm{crit}}(\Gamma)
   =\operatorname{Span}_{\mathbb C}
         \{v_F:[v_F]=\jmath(F),\ F\in\mathscr F\}.
\]
It is \emph{full} if it equals $W$. The action is
\emph{critically strongly irreducible} if the actions on $V$ and on
this critical module are strongly irreducible. Here strong
irreducibility means that every finite-index subgroup acts irreducibly.

In the spectral construction using a regular seed, the common middle
spaces must approach a complement of the lower flag space. This gives
the middle directions needed to construct a return throughout the
upper flag space. We make this requirement explicit.

\begin{definition}[Regular critical seed]\label{def:rec-regular-seed}
A critical flag $F=(A\subset B)\in\mathscr F$ is a regular critical
seed if there are $s_j\in\Gamma$ with strict spectral modulus gaps at
$a,p$ such that the following hold. Their positive and slow negative
spectral flags
$F_j^\pm=(A_j^\pm\subset B_j^\pm)$ both tend to $F$. Their common
middle spectral blocks
\[
 C_j=B_j^-\cap B_j^+,
 \qquad V=A_j^-\oplus C_j\oplus A_j^+,
\]
have relatively compact cyclic projective action
$\overline{\langle[s_j|_{C_j}]\rangle}\subset\operatorname{PGL}(C_j)$.
Finally, after extraction,
\begin{equation}\label{eq:rec-seed-complement}
 C_j\longrightarrow C_\infty,
 \qquad B=A\oplus C_\infty.
\end{equation}
The complement condition is part of the hypothesis. Coalescence of the
two spectral flags and compactness of their middle actions alone do
not imply it.
\end{definition}

We begin with the orbit of the flags. Fullness lets the exterior
contractions act on the whole ambient Cartan component, and strong
irreducibility lets us move a flag away from their kernels. In the
one-dimensional middle case the same argument brings both spectral
flags to a prescribed critical flag.

\begin{lemma}[Critical contractions and coalescence]
\label{lem:rec-contractions}
Suppose that the critical module is full and its action is strongly
irreducible. Then $\tau(\Gamma)$ contains a proximal element and
\[
 \mathscr F=\overline{\Gamma F_0}\qquad(F_0\in\mathscr F).
\]
If $p=a+1$, then, for each realised critical flag $F$, there are
elements $s_j\in\Gamma$ with spectral gaps at $a,p$ whose positive
and slow negative spectral flags both tend to $F$.
\end{lemma}

\begin{proof}
On the ambient exterior tensor, the largest singular exponent of a
diagonal Cartan factor is
\[
 \sum_{i=1}^a\log\sigma_i+\sum_{i=1}^p\log\sigma_i.
\]
Its tensor line belongs to $W$ and is unique. A different $a$-subset
loses at least $\log(\sigma_a/\sigma_{a+1})$, and a different
$p$-subset loses at least $\log(\sigma_p/\sigma_{p+1})$.
Consequently
\[
 \frac{\sigma_2(\tau(h_j))}{\sigma_1(\tau(h_j))}
 \le
 \max\left\{\frac{\sigma_{a+1}(h_j)}{\sigma_a(h_j)},
             \frac{\sigma_{p+1}(h_j)}{\sigma_p(h_j)}\right\}
 \longrightarrow0.
\]
The compact Cartan factors act unitarily for the induced norm. Hence,
for a sequence realising $F$, normalised representatives converge, after
extraction, to $T=v\otimes\phi$ of rank one with $[v]=\jmath(F)$.
Fullness ensures that the module in use is the entire ambient Cartan
component.

Irreducibility supplies $g\in\Gamma$ with
$\phi(\tau(g)v)\ne0$. The normalised representatives of
$\tau(h_jg)$ converge to $T\tau(g)$, whose unique nonzero eigenvalue
is simple. Continuity of the roots of the characteristic polynomial
therefore makes $\tau(h_jg)$ proximal for large $j$.

The set $\mathscr F$ is invariant. Indeed, left multiplication by a
fixed $g$ multiplies the nonzero rank-one limits of
$\bigwedge^ah_j$ and $\bigwedge^ph_j$ by the corresponding exterior
powers of $g$, so their image flag becomes $gF$. The singular gaps
and bounded internal ratios are preserved by multiplication by a
fixed invertible matrix. Passing to closures gives invariance.

Fix $F_0\in\mathscr F$. For the rank-one limit $T$ just constructed,
irreducibility gives a $g$ such that $\jmath(gF_0)$ avoids
$\mathbb P(\ker T)$: the span of the orbit of any nonzero line is
the whole module. Rank-one convergence now gives $h_jgF_0\to F$.
Every realised flag, and then its closure, is therefore contained in
$\overline{\Gamma F_0}$. The reverse inclusion is invariance.

For the last assertion, choose $s$ such that $\tau(s)$ is proximal.
Let $\mu_1,\ldots,\mu_d$ be the eigenvalues of a lift of $s$, in
decreasing order of modulus. In an eigenbasis of its semisimple part,
the highest tensor vector is
\[
 v_0=(e_1\wedge\cdots\wedge e_a)\otimes
                      (e_1\wedge\cdots\wedge e_p).
\]
The two nonzero vectors obtained by the elementary lowering operators
$E_{a+1,a}$ and $E_{p+1,p}$ lie in $W$. Their eigenvalues equal that
of $v_0$ times $\mu_{a+1}/\mu_a$ and $\mu_{p+1}/\mu_p$,
respectively. Equality in either boundary modulus gap would thus
contradict proximality, including algebraic multiplicity. Hence
\[
 |\mu_a|>|\mu_{a+1}|>|\mu_{a+2}|\qquad(p=a+1).
\]
The corresponding top and bottom spectral flags $F_s^+,F_s^-$ are
opposite. They belong to $\mathscr F$: powers of the exterior
representations at $a,p$, divided by their simple dominant
eigenvalues, converge to rank-one spectral projectors. Their image
lines give the target spectral flags, and the inverse powers give the
slow negative flags.

We need one element moving \emph{both} flag lines outside a prescribed
hyperplane. This follows from strong irreducibility. Let $L$ be the
complex Zariski closure of $\tau(\Gamma)$. The subgroup
$\tau(\Gamma)\cap L^0$ has finite index and is Zariski dense in the
irreducible algebraic group $L^0$. Each orbit $L^0\ell$ spans $W$;
otherwise that finite-index subgroup would preserve a proper
subspace. For two lines $\ell^+,\ell^-$, the conditions
$g\ell^\pm\subset\ker T$ are therefore proper algebraic closed
conditions on $L^0$. The intersection of their open complements is
nonempty and meets $\tau(\Gamma)\cap L^0$.

Apply this to $\ell^\pm=\jmath(F_s^\pm)$. For a fixed resulting
$g\in\Gamma$, both $h_jgF_s^\pm$ tend to $F$ by rank-one
convergence. The conjugates $s_j=h_jgsg^{-1}h_j^{-1}$ have exactly
these spectral flags and retain the two strict spectral gaps.
\end{proof}

We can now construct the returns. When the middle quotient is a line,
its singular action is controlled by one scalar. In higher dimension
the regular seed gives the middle return and the limiting complement
needed to lift it to the carrier.

\begin{theorem}[Critical recurrence]\label{thm:crit-recurrence}
Let $\Gamma<\operatorname{PSL}(V)$ be infinite and discrete, of
self-opposed uniform critical type $(a,p)$, where
$1\le a<p<d$ and $a+p=d$. Suppose that its critical flag module is
full and its action is critically strongly irreducible. Assume either
that $p-a=1$, or that a critical flag is a regular critical seed in
the sense of Definition~\ref{def:rec-regular-seed}. Then
\[
 \Lambda_{\mathrm{crit}}=\mathscr R_\Gamma.
\]
Consequently $\Omega_{\mathrm{crit}}$ is the largest
$\Gamma$-invariant open set on which the action is proper. It may be
empty.
\end{theorem}

\begin{proof}
We first produce a recurrent critical carrier.

\smallskip\noindent\emph{The case $p-a=1$.}
Fix a realised target flag $F=(A\subset B)$. By
Lemma~\ref{lem:rec-contractions}, there are $s_j$ with opposite
spectral flags $F_j^\pm\to F$ and strict boundary spectral gaps.
For each fixed $j$, the finite singular critical flags of $s_j^M$
tend to its spectral flags as $M\to\infty$. To see both source and
target convergence, normalise powers of $\bigwedge^as_j$ and
$\bigwedge^ps_j$ by their simple dominant eigenvalues. The resulting
rank-one spectral projectors give the limiting top left singular
lines in the exterior powers. The orthogonal complements of the
projector kernels give the limiting top right singular lines there.
The Pl\"ucker embedding identifies these lines with the limiting
fast source planes in $V$, and their orthogonal complements in $V$
are the slow source spaces.
All complementary Jordan factors are polynomially bounded relative
to a strictly smaller exponential rate.

Choose $M_j$ sufficiently large that, for $\gamma_j=s_j^{M_j}$,
the finite singular flags $\widehat F_j^\pm$ are within $1/j$ of
$F_j^\pm$, the ratios $\sigma_a/\sigma_p$ and
$\sigma_p/\sigma_{p+1}$ exceed $j$, and $\gamma_j$ differs from all
previously selected elements. All these requirements are possible
for each fixed $s_j$. The ratios also show that the selected sequence
escapes compact subsets.

Put
\[
 C_j^\pm=\widehat B_j^\pm\cap(\widehat A_j^\pm)^\perp,
 \qquad C=B\cap A^\perp.
\]
These are lines tending to $C$. Choose unit vectors
$c_j^\pm\to c\in C$. The singular-value decomposition gives
\[
 \gamma_jc_j^-=\sigma_p(\gamma_j)\eta_jc_j^+,
 \qquad |\eta_j|=1.
\]
Choose a lift of $\gamma_j$ and put
$\lambda_j=\eta_j^{-1}/\sigma_p(\gamma_j)$. Then
$\lambda_j\gamma_jc_j^-=c_j^+$. With
$H_j^-=(\widehat B_j^-)^\perp$, the same decomposition gives
\[
 \lambda_j\gamma_jH_j^-=\widehat A_j^+,
 \quad
 \|\lambda_j\gamma_j|_{\widehat A_j^-}\|
       =\frac{\sigma_{p+1}}{\sigma_p}\longrightarrow0,
 \quad
 \|(\lambda_j\gamma_j|_{H_j^-})^{-1}\|
       =\frac{\sigma_p}{\sigma_a}\longrightarrow0.
\]
For $x=[u+zc]\in\mathbb P(B)\setminus\mathbb P(A)$ choose
$u_j^\pm\in\widehat A_j^\pm$ tending to $u$. The vectors
\[
 v_j=u_j^-+zc_j^-+
              (\lambda_j\gamma_j|_{H_j^-})^{-1}u_j^+
\]
converge to $u+zc$, and
$\lambda_j\gamma_jv_j\to u+zc$. Hence these points are recurrent.
Closedness of $\mathscr R_\Gamma$ includes $\mathbb P(A)$, so
$\mathbb P(B)\subset\mathscr R_\Gamma$.

\smallskip\noindent\emph{The regular-seed case.}
Let $F=(A\subset B)$ and $s_j$ satisfy
Definition~\ref{def:rec-regular-seed}. Choose identifications
$J_j:C_\infty\to C_j$ tending to inclusion. For each fixed $j$,
compactness of the cyclic projective middle action gives arbitrarily
large positive integers $N$ for which
$[J_j^{-1}s_j^NJ_j]$ is arbitrarily close to the identity. This holds
for a finite cyclic group by taking multiples of its order, and in
general by taking differences of sufficiently close powers in its
compact closure.

For completeness, projective compactness also controls the scalar
normalisation needed in the other two blocks. It implies that
$s_j|_{C_j}$ is diagonalisable with all eigenvalues of one modulus
$q_j>0$: unequal moduli or a nontrivial Jordan block would make the
condition numbers of its powers unbounded. Thus
$q_j^{-N}s_j^N|_{C_j}$ and its inverse are bounded for fixed $j$.
The strict spectral gaps imply
\[
 q_j^{-N}\|s_j^N|_{A_j^-}\|\longrightarrow0,
 \qquad
 q_j^N\|(s_j^N|_{A_j^+})^{-1}\|\longrightarrow0.
\]
Polynomial factors from Jordan blocks outside $C_j$ are absorbed by
the strict exponential gaps.

For each fixed $j$, take a sequence of middle return times
$N\to\infty$ and choose scalars $\lambda_{j,N}$ such that
$J_j^{-1}(\lambda_{j,N}s_j^N|_{C_j})J_j\to I$ along that sequence.
The modulus of its determinant is
$(|\lambda_{j,N}|q_j^N)^{p-a}$, so
$|\lambda_{j,N}|q_j^N\to1$. Thus
$\|\lambda_{j,N}s_j^N|_{A_j^-}\|\to0$ and
$\|(\lambda_{j,N}s_j^N|_{A_j^+})^{-1}\|\to0$ along these return times.
Choose a return time $N_j$ so that all three errors are smaller than
$1/j$, both singular
boundary ratios exceed $j$, and $\gamma_j=s_j^{N_j}$ is distinct from
all earlier choices. The strict spectral gaps ensure that these powers
escape. With $\lambda_j=\lambda_{j,N_j}$, this gives
\[
 J_j^{-1}(\lambda_j\gamma_j|_{C_j})J_j\longrightarrow I,
 \quad \|\lambda_j\gamma_j|_{A_j^-}\|\longrightarrow0,
 \quad \|(\lambda_j\gamma_j|_{A_j^+})^{-1}\|
                                        \longrightarrow0.
\]

For a nonzero $v=u+c\in A\oplus C_\infty=B$, choose
$u_j^\pm\in A_j^\pm$ tending to $u$ and set $c_j=J_jc$. Then
\[
 v_j=u_j^-+c_j+
              (\lambda_j\gamma_j|_{A_j^+})^{-1}u_j^+
       \longrightarrow v,
 \qquad
 \lambda_j\gamma_jv_j\longrightarrow v.
\]
Thus the entire carrier $\mathbb P(B)$ is recurrent. This argument
uses the complementary limit in
\eqref{eq:rec-seed-complement}; cyclic compactness alone would not
supply the decomposition $v=u+c$ in the limit.

\smallskip\noindent\emph{Propagation and properness.}
In either case a critical carrier $\mathbb P(B_0)$ is recurrent.
Lemma~\ref{lem:rec-contractions} gives
$\mathscr F=\overline{\Gamma F_0}$. Invariance and closedness of
$\mathscr R_\Gamma$, and continuity of the upper flag projection,
therefore give
$\Lambda_{\mathrm{crit}}\subset\mathscr R_\Gamma$.

We verify the reverse inclusion and properness directly. Suppose
$x_j\to x$ and $\gamma_jx_j\to y$ for a divergent sequence in $\Gamma$.
Extract its singular flags. If $x\notin\Lambda_{\mathrm{crit}}$,
then $x$ avoids its limiting slow-source $p$-plane $B^-$, which is
a target carrier of the inverse sequence. The fast source components of
unit lifts of $x_j$ are therefore bounded below by a positive
constant. Their images have norm at least that constant times
$\sigma_a(\gamma_j)$, whereas the remaining image components have
norm at most a fixed multiple of $\sigma_{a+1}(\gamma_j)$.
The boundary gap forces $y\in\mathbb P(A^+)$, inside the critical
wall. Hence no dynamical pair has both endpoints in the complement.
If a compact subset of the complement met infinitely many of its
translates, compactness would produce just such a pair. The action
on $\Omega_{\mathrm{crit}}$ is therefore proper. It avoids all
recurrent points, proving
$\mathscr R_\Gamma\subset\Lambda_{\mathrm{crit}}$.
Finally every open set of proper action avoids every recurrent point,
so it is contained in $\Omega_{\mathrm{crit}}$.
\end{proof}

\begin{proposition}[A regular seed with a wider middle block]
\label{r6seed:wide-example}
For every $d\ge4$ there is an infinite discrete, complex-Zariski-dense
$\Gamma<\operatorname{PSL}_d(\mathbb C)$ satisfying all the hypotheses
of Theorem~\ref{thm:crit-recurrence} with $(a,p)=(1,d-1)$ and a
regular critical seed. In particular the regular-seed branch is
non-vacuous when $p-a=d-2>1$, and its largest domain is nonempty.
\end{proposition}
\begin{proof}
Write $V=\mathbb Ce_1\oplus C\oplus\mathbb Ce_d$, with
$C=\langle e_2,\ldots,e_{d-1}\rangle$, and use the Hermitian form
\[
 J=\begin{pmatrix}0&0&1\\0&I_{d-2}&0\\1&0&0\end{pmatrix}.
\]
Let $u,s,w$ act identically on $C$ and, on
$\langle e_1,e_d\rangle$ in that order, by
\[
 u=\begin{pmatrix}1&i\\0&1\end{pmatrix},\qquad
 s=\begin{pmatrix}2&i\\-3i&2\end{pmatrix},\qquad
 w=\begin{pmatrix}0&i\\i&0\end{pmatrix}.
\]
For $z\in\{2e_k,2ie_k:2\le k\le d-1\}\subset C$, also put
\[
 n(z)=\begin{pmatrix}
 1&-z^*&-\|z\|^2/2\\0&I_{d-2}&z\\0&0&1
 \end{pmatrix}.
\]
These matrices have Gaussian-integer entries, determinant one, and
$g^*Jg=J$. Let $\Gamma$ be the projective group they generate.
It is infinite because $u$ has infinite order, and is discrete:
bounded determinant-one Gaussian-integer matrices form a finite set,
and the projective scalar ambiguity has finite fibre.

We check density explicitly. Powers of $n(2e_k)$ and $n(2ie_k)$
lie on the one-parameter unipotent subgroups with infinitesimal
matrices
\[
 X_k=-2E_{1k}+2E_{kd},\qquad
 Y_k=2iE_{1k}+2iE_{kd}.
\]
Their complex-Zariski closures contain those whole subgroups, since a
polynomial vanishing at every integer parameter vanishes identically.
The complex Lie algebra of the closure therefore contains
$E_{1k},E_{kd}$; conjugation by $w$ supplies $E_{dk},E_{k1}$.
Their brackets supply every off-diagonal matrix unit and all traceless
diagonal matrices. Hence the closure is $\operatorname{SL}_d(\mathbb C)$.
Every finite-index subgroup has the same closure, since this algebraic
group is connected. Thus the actions on $V$ and
$W_{1,d-1}^{\mathrm{fl}}\simeq\mathfrak{sl}_d(\mathbb C)$ are
strongly irreducible.

For uniform type, put $v=(e_1-e_d)/\sqrt2$, so $J=I-2vv^*$.
The identity $g^*Jg=J$ gives
\[
 g^*g-I=2(g^*vv^*g-vv^*),\qquad
 \operatorname{rank}(g^*g-I)\le2.
\]
Since $\det(g^*g)=1$, the singular values are
$e^t,1,\ldots,1,e^{-t}$ for some $t\ge0$.
Along an escaping sequence $t\to\infty$, giving uniform critical type
$(1,d-1)$. The nonzero critical flag span is invariant; the ambient
adjoint module is irreducible under the Zariski-dense group, so this
span is the whole $W_{1,d-1}^{\mathrm{fl}}$.

Finally $s$ has eigenvalues $2+\sqrt3$, $1$ on $C$, and
$2-\sqrt3$, with extremal eigenlines
$\mathbb C(e_1\mp i\sqrt3e_d)$. For $s_j=u^jsu^{-j}$ these lines
become
\[
 \mathbb C\bigl((1\pm j\sqrt3)e_1\mp i\sqrt3e_d\bigr)
 \longrightarrow A=\mathbb Ce_1.
\]
The attracting and repelling critical spectral flags thus both tend
to $F=(A\subset B)$ with $B=A\oplus C$; their common middle block
is exactly $C$, and every $s_j$ acts identically there. Powers of each
$s_j$ have finite singular critical flags tending to its spectral
flags, by the rank-one exterior spectral projectors at degrees
$1,d-1$. A diagonal choice of sufficiently large powers therefore
realises $F$ as a limiting target critical flag. Thus $F\in\mathscr F$
is a regular seed, with $C_\infty=C$ and $B=A\oplus C_\infty$.
All hypotheses of Theorem~\ref{thm:crit-recurrence} have now been checked.
The $J$-negative projective lines form a nonempty invariant ball. The
ambient unitary group acts properly on this ball, since its point
stabilisers are compact; hence so does its discrete subgroup $\Gamma$.
The largest domain supplied by the theorem contains this ball.
\end{proof}

%--------------------------------------------------------------
\section{Borel density and the critical module}
\label{r3bor:section}

For a lattice image, the ambient semisimple group shows that the critical
flags span the ambient module. Borel density
transfers invariant subspaces from the lattice to that group
\cite{r3bor:Borel1960}; ambient irreducibility then forces the nonzero
critical span to be full. We use a form for continuous complex
representations, so no algebraic matrix realisation of the ambient
group is required.

\begin{proposition}\label{r3bor:transfer}
Let $S$ be a connected semisimple Lie group without compact factors,
let $\widetilde\Gamma<S$ be a lattice, and let
$\tau:S\to\operatorname{GL}(E)$ be a continuous finite-dimensional
complex representation. Every complex subspace invariant under a
finite-index subgroup of $\widetilde\Gamma$ is $S$-invariant.
Consequently, if $E$ is irreducible as a complex $S$-module, then
$\tau(\widetilde\Gamma)$ is strongly irreducible. The transfer applies
equally to any continuous tensor or subquotient representation.
\end{proposition}

\begin{proof}
A finite-index subgroup $\widetilde\Gamma_0$ is again a lattice.
Let $U\subset E$ be $\widetilde\Gamma_0$-invariant. The cases
$U=0,E$ are immediate. Otherwise use real underlying spaces, put
$r=\dim_{\mathbb R}U$, and choose a nonzero
$u\in\bigwedge^r_{\mathbb R}U$. The continuous orbit map
\[
 S/\widetilde\Gamma_0\longrightarrow
 \mathbb P_{\mathbb R}\!\left(\bigwedge^r_{\mathbb R}E\right),
 \qquad s\widetilde\Gamma_0\longmapsto
 [\bigwedge^r_{\mathbb R}\tau(s)u]
\]
pushes forward the invariant probability measure to an $S$-invariant probability measure whose
support contains $[u]$. For connected semisimple groups without compact
factors, the projective-measure form of Borel density says that this
support consists of $S$-fixed points
\cite[Lemma~3 and the Theorem on p.~211]{r3bor:Furstenberg1976}.
Thus $S$ preserves $U$, as required. Any finite-index subgroup of
$\tau(\widetilde\Gamma)$ has finite-index inverse image in
$\widetilde\Gamma$, so ambient irreducibility gives strong
irreducibility. Realification here imposes no
rationality assumption on the complex representation.
\end{proof}

When $S$ is given by a real algebraic group and $\tau$ is rational,
one can equivalently use real Zariski density of every finite-index
lattice subgroup: the stabiliser of a complex subspace is real
algebraic. This algebraic proof uses rationality; the preceding
continuous-representation proof does not.

Apply this transfer both to the original representation and to its
Cartan component. The remaining questions are geometric: the type of
the middle crossing and the existence of the return.

\begin{corollary}\label{r3bor:lattice-largest}
Let $S$ be a connected semisimple Lie group without compact factors,
let $\widetilde\Gamma<S$ be a lattice, and let
$\tau:S\to\operatorname{SL}(V)$ be a continuous complex representation.
Assume that $\Gamma=[\tau(\widetilde\Gamma)]<\operatorname{PSL}(V)$
is infinite and discrete. For integers $1\le a<p<d=\dim_{\mathbb C}V$
with $a+p=d$, suppose that
\begin{enumerate}[label=\textup{(\roman*)},nosep]
\item $\Gamma$ has self-opposed uniform critical type $(a,p)$;
\item $V$ and the ambient critical flag module
      $W_{a,p}^{\mathrm{fl}}$ are irreducible complex $S$-modules;
\item either $p-a=1$, or the critical flag set contains a regular
      critical seed as in Definition~\ref{def:rec-regular-seed}.
\end{enumerate}
Then $W_{a,p}^{\mathrm{crit}}(\Gamma)=W_{a,p}^{\mathrm{fl}}$, the
action is critically strongly irreducible, and
$\Lambda_{\mathrm{crit}}(\Gamma)=\mathscr R_\Gamma$.
Hence $\Omega_{\mathrm{crit}}(\Gamma)$ is the largest invariant open
set of proper discontinuity, possibly empty.
\end{corollary}

\begin{proof}
The generated critical module is a nonzero
$\widetilde\Gamma$-invariant subspace of $W_{a,p}^{\mathrm{fl}}$.
Proposition~\ref{r3bor:transfer} makes it $S$-invariant, and ambient
irreducibility makes it the whole module. If $\Gamma_0<\Gamma$ has
finite index, its inverse image in $\widetilde\Gamma$ has the same
finite index. The proposition applied to this inverse image on both
$V$ and $W_{a,p}^{\mathrm{fl}}$ gives irreducibility on both. Scalar
lifts preserve the same subspaces, so $\Gamma$ is critically strongly
irreducible. Theorem~\ref{thm:crit-recurrence}, with hypotheses
\textup{(i)} and \textup{(iii)}, now gives the conclusions.
\end{proof}

Here we used irreducibility of the ambient flag module. A generated
critical submodule may itself be irreducible and still be a proper
summand of that ambient module. The corollary leaves uniform critical
type and the return hypothesis to be checked, and assumes separately
that the projective lattice image is discrete.

\section{Cusped Hitchin groups and a nonempty largest domain}
\label{r11:cusped-section}

The one-dimensional-middle case also applies when parabolic elements
prevent traditional Borel--Anosov dynamics. We use the cusped Hitchin
representations of Canary--Zhang--Zimmer. Their flag map is defined on
the classical limit set of a geometrically finite Fuchsian group;
the singular gaps diverge even along parabolic powers.

\begin{theorem}\label{r11:cusped-largest}
Let $\Delta<\mathrm{PSL}(2,\R)$ be torsion-free, non-elementary and
geometrically finite. Let $d=2m+1\ge5$, and let
$\rho:\Delta\to\mathrm{SL}_d(\R)$ be a cusped Hitchin representation
whose image is real Zariski dense. Write
$\xi:\Lambda(\Delta)\to\mathrm{Fl}(\R^d)$ for its positive flag map,
$\xi^q_\C(x)=\xi^q(x)\otimes_\R\C$, and
$\Gamma=[\rho(\Delta)]<\mathrm{PSL}_d(\C)$. Then
\[
 \mathscr R_\Gamma=\Lambda_{\mathrm{crit}}(\Gamma)
 =\bigcup_{x\in\Lambda(\Delta)}\mathbb P(\xi^{m+1}_\C(x)),
 \qquad q_\Gamma=m.
\]
Its complement is nonempty and is the largest domain of proper
discontinuity. If $\Delta$ contains a parabolic element, $\rho$ is
not Borel--Anosov in the traditional word-hyperbolic sense.
\end{theorem}

\begin{proof}
Canary--Zhang--Zimmer prove that $\rho$ is cusped $P_k$--Anosov for
every $k$, and that parabolic elements have regular unipotent images
\cite[Theorem~1.4]{r11CZZCusped}. Their gap estimate and Cartan-flag
convergence \cite[Theorem~1.1(1) and Corollary~4.2(1)]{r11CZZCusped}
give, along every escaping sequence $\delta_j$,
\[
 \frac{\sigma_k(\rho(\delta_j))}
 {\sigma_{k+1}(\rho(\delta_j))}\longrightarrow\infty
 \quad(1\le k<d),
\]
and limiting target flags precisely $\xi_\C(\Lambda(\Delta))$.
Indeed an escaping sequence in $\Delta$ has
$d_{\mathbb H^2}(o,\delta_j o)\to\infty$; after extraction it
converges to a classical limit point. Conversely, powers of a
loxodromic element realise its attracting flag, and such fixed points
are dense in $\Lambda(\Delta)$. The gap estimate also makes the
projective image discrete: a compact set of projective matrices has
bounded singular ratios, so its preimage in $\Delta$ is finite.
Thus the critical type is uniformly $(a,p)=(m,m+1)$.

Real Zariski density in $\mathrm{SL}_d(\R)$ implies complex Zariski
density in $\mathrm{SL}_d(\C)$. Every finite-index subgroup has the
same complex closure, since $\mathrm{SL}_d(\C)$ is connected.
The critical flags form a nonempty invariant compact set. Their
Segre--Pl\"ucker lines therefore span a nonzero invariant subspace
of the irreducible ambient Cartan component
$W^{\mathrm{fl}}_{m,m+1}$, hence the whole component. The actions on
this module and on $\C^d$ are strongly irreducible.
Theorem~\ref{thm:crit-recurrence}, with $p-a=1$, proves the displayed
equality, and that the complement is the largest domain, without a seed assumption.

For nonemptiness choose distinct $x,y\in\Lambda(\Delta)$ and nonzero
real vectors $u\in\xi^1(x)$, $v\in\xi^1(y)$. The point $[u+iv]$
belongs to none of the removed planes. If
$u+iv\in\xi^{m+1}_\C(z)$, reality of that plane gives
$u,v\in\xi^{m+1}(z)\subset\xi^{d-1}(z)$. Transversality of the
flag curve forces both $x=z$ and $y=z$, a contradiction.
Finally, a parabolic element has infinite order and unipotent image,
whereas every infinite-order element in a traditional Borel--Anosov
representation has a strict spectral gap. This proves the last assertion.
\end{proof}

These groups are outside the rank-one unitary family as well. A real
Zariski-dense subgroup of $\mathrm{SL}_d(\R)$, $d\ge3$, cannot preserve
a nonzero Hermitian form on $\C^d$: invariance would extend to
$\mathrm{SL}_d(\R)$, and its positive diagonal matrices force every
matrix coefficient of the form to vanish. The same argument applies
after a complex change of basis. Thus, in the presence of a cusp, this
application is covered by neither the unitary largest-domain theorem
nor Stecker's traditional Borel--Anosov classification. KLP's
regular-group construction already gives properness of this domain
\cite[Theorem~6.13]{KLP}. The cusped
flag and gap estimates are the inputs from the existing theory; the
critical recurrence criterion supplies the largest-domain conclusion.

\begin{proposition}\label{r11:cusped-existence}
There are real Zariski-dense cusped Hitchin representations of the
three-punctured sphere group into $\mathrm{SL}_5(\R)$, with regular
unipotent peripheral monodromies. Their projective images satisfy
Theorem~\ref{r11:cusped-largest}.
\end{proposition}
\begin{proof}
Equip the three-punctured sphere $S$ with its complete finite-area
hyperbolic metric, and let $\Delta<\mathrm{PSL}(2,\R)$ be its holonomy.
Thus $\Delta\cong\pi_1(S)=F_2$ and $\chi(S)=-1$.
Since $\mathrm{PGL}_5(\R)\cong\mathrm{SL}_5(\R)$,
the positive decorated local systems of Fock--Goncharov
form a cell $\mathcal A^+_{\mathrm{SL}_5,S}$ of dimension $24$.
On the positive locus, forgetting decorations has fibres
$H(\R_{>0})^3\cong(\R_{>0})^{12}$, where $H$ is the rank-four
Cartan torus, and its image has unipotent peripheral monodromies
\cite[(1.12), (2.9)--(2.10)]{r11FGPositive}. The projection is positive:
in their coordinates, $p^*X_i=\prod_j\Delta_j^{\varepsilon_{ij}}$
\cite[Proposition~9.1, (9.6)]{r11FGPositive}. Its positive image
therefore has dimension $24-3\cdot4=12$.

The positive peripheral framing is injective, so no peripheral element
can have trivial image. The representations therefore preserve
parabolic type. They are cusped Hitchin by
\cite[Theorem~9.2]{r11CZZCusped}. They are irreducible and every
peripheral image is regular unipotent
\cite[Theorem~1.4]{r11CZZCusped}. Such an image fixes a unique complete
flag. Forgetting the framing thus preserves the dimension $12$ of
the family of conjugacy classes.

In dimension five, the Zariski closure of a Hitchin image of a
Fuchsian lattice is either $\mathrm{SL}_5(\R)$ or is contained in
a conjugate of $\mathrm{SO}(3,2)$
\cite[Theorem~11.10, following Sambarino]{r11CZZEntropy}.
The irreducible $F_2$-characters in the latter locus have dimension
at most $2\dim\mathrm{SO}(3,2)-\dim\mathrm{SO}(3,2)=10$:
the space of pairs has dimension $20$, and irreducibility makes
their centraliser finite. Indeed, the commuting real division algebra
is $\R$, $\C$ or $\mathbb H$ by Schur's lemma; the latter two require
even real dimension, so here the centraliser in $\mathrm{SL}_5(\R)$
consists of scalars with $t^5=1$. Allowing the preserved form to vary does
not increase the dimension after conjugation in $\mathrm{SL}_5$.
Equivalently, the incidence of a pair and a preserved form up to
scalar has dimension at most $20+14$, and its irreducible
$\mathrm{SL}_5$-orbits have dimension $24$.
The positive twelve-dimensional family cannot lie in this locus.
A representation outside it has the required Zariski density.
\end{proof}

\section{Borel--Anosov critical return}
\label{r3an:section}

For a Borel--Anosov group, the largest projective domain follows from
Stecker's minimal-fat-ideal classification, as explained in Section~\ref{sec:compare}. Here the returns are constructed directly. Conjugating one group element along powers of another makes its target
Cartan flag and inverse target flag coalesce. The destination formula
returns every point of a middle carrier, in both parities, without the
module or seed hypotheses of Theorem~\ref{thm:crit-recurrence}.

\begin{proposition}\label{r3an:borel-return}
Let $\Delta$ be a non-elementary word-hyperbolic group and let
$\rho:\Delta\to\operatorname{PSL}(V)$ be Borel--Anosov, where
$d=\dim_{\mathbb C}V\ge2$. Put $\Gamma=\rho(\Delta)$ and
$p=\lceil d/2\rceil$, and write $\xi$ for the complete boundary flag
map. For the paired critical system of Example~\ref{ex:critsys},
\[
 q_\Gamma=p-1,\qquad
 \Lambda_{\mathrm{crit}}(\Gamma)
 =\bigcup_{z\in\partial\Delta}\mathbb P(\xi^p(z))
 =\mathscr R_\Gamma.
\]
Consequently $\Omega_{\mathrm{crit}}(\Gamma)$ is the largest
$\Gamma$-invariant open set of proper discontinuity. It may be empty;
no cocompactness is asserted.
\end{proposition}

\begin{proof}
The Borel--Anosov properties
\cite[Theorem~1.3 and Theorem~5.3(3)]{R63GGKW} give a discrete image
with finite kernel, divergence of every adjacent Cartan gap along
escaping sequences, and limiting target complete flags precisely
$\xi(\partial\Delta)$. The map $\xi$ is continuous, equivariant,
transverse and dynamics-preserving. In particular an infinite-order
element has distinct eigenvalue moduli and its two eigenflags are the
boundary flags at its fixed points. All Cartan blocks are therefore
lines. The critical level is $p-1$, and the two paired shadows have
vector dimensions $p$ and $d-p$. Since $d-p\le p$, their union is the
displayed union of $p$-planes. It is closed by compactness of the
boundary and incidence. Properness of its complement, proved for the
critical system in Proposition~\ref{r3crit:cover}, gives
$\mathscr R_\Gamma\subseteq\Lambda_{\mathrm{crit}}(\Gamma)$.

For the reverse inclusion choose an infinite-order $a\in\Delta$, put
$z=a^+$ and $b=a^-$, and choose $h\in\Delta$ with $hb\ne b$.
Non-elementarity guarantees such a choice. In an eigenbasis of a lift
of $\rho(a)$ write
\[
 A=\operatorname{diag}(\lambda_1,\ldots,\lambda_d),\qquad
 |\lambda_1|>\cdots>|\lambda_d|,
\]
so $\xi^q(z)=\langle e_1,\ldots,e_q\rangle$ and
$\xi^q(b)=\langle e_{d-q+1},\ldots,e_d\rangle$.
Let $H$ lift $\rho(h)$ and set $M_m=A^mHA^{-m}$. For
$I_q=\{1,\ldots,q\}$ and $J_q=\{d-q+1,\ldots,d\}$, the
$(I_q,J_q)$ coefficient of $\bigwedge^qM_m$ is
\[
 \left(\frac{\lambda_1\cdots\lambda_q}
 {\lambda_{d-q+1}\cdots\lambda_d}\right)^m
 \det H_{I_q,J_q}\qquad(1\le q<d).
\]
The minor is nonzero by transversality of $\xi^q(hb)$ to
$\xi^{d-q}(b)$. Its exponential factor has strictly larger modulus
than every other coefficient's factor. Division by that factor thus
gives a nonzero rank-one limit with image $\bigwedge^q\xi^q(z)$.
The same holds for $M_m^{-1}=A^mH^{-1}A^{-m}$, since $h^{-1}b\ne b$.
These limits show that the projective sequence escapes, all its Cartan
blocks are lines, and both its target flag and its inverse target flag
are $\xi(z)$. After taking distinct terms and extracting Cartan data,
its decreasing source and increasing target flags are consequently
\[
 F_j=\xi^{d-j+1}(z),\qquad G_j=\xi^j(z)
 \quad(1\le j\le d).
\]
Thus actual Cartan flags coalesce; coalescence of eigenflags alone
would not establish this assertion.

If $x\in\mathbb P(\xi^p(z))$ belongs to source stratum
$\mathbb P(F_j)\setminus\mathbb P(F_{j+1})$, then
$j\ge d-p+1\ge p$. Since the graded blocks are one-dimensional,
Proposition~\ref{prop:graded} gives
\[
 \D_{(M_m)}(x)=\mathbb P(G_j)
 \supseteq\mathbb P(\xi^p(z))\ni x.
\]
The destination formula realises this inclusion by moving sources, so
every point of this critical carrier is recurrent. The inequality
includes both $d=2r$, $p=r$, and $d=2r+1$, $p=r+1$.
Attracting fixed points of infinite-order elements are dense in
$\partial\Delta$. Continuity of $\xi$ and closedness of
$\mathscr R_\Gamma$ therefore make every critical carrier recurrent.
Equality follows. Every domain of proper discontinuity avoids the recurrent
set, so the complement is the largest domain.
\end{proof}

\section{Limits of the recurrence criterion}
\label{r4voice:limits}

The two return proofs make clear why we must distinguish spectral
flags, singular flags and the common middle space. A regular seed is
one way to obtain the required return, but the one-dimensional middle
argument can succeed even when the seed's limiting complement fails,
as a Hitchin--Schottky group shows.

\begin{example}[A strict improvement of the seed hypothesis]
\label{r2bd:strict-seed}
The one-dimensional-middle case of
Theorem~\ref{thm:crit-recurrence} applies to Zariski-dense groups
which have no regular critical seed. Take a strongly dense Hitchin
representation $\rho:\pi_1(S)\to\mathrm{SL}_5(\R)$ of a closed
surface group; such representations exist by
\cite[Theorem~1]{r2bd:LRW2025}. Restrict to a rank-two
convex-cocompact free subgroup $H<\pi_1(S)$, obtained by taking
sufficiently large powers of two independent elements. Let $C$ be
its boundary Cantor set and $\Gamma=\rho(H)$, acting on $\C\PP^4$.
The representation and its contragredient have Frenet limit flags,
and the Anosov simple-root gaps diverge along every escaping sequence
\cite[Theorems~1.4 and~4.1]{r2appLabourie2006}.
Thus $\Gamma$ has self-opposed uniform critical type $(2,3)$ and
critical flags $(\xi^2(x)_\C\subset\xi^3(x)_\C)$, $x\in C$.

Strong density makes $\rho(H)$ Zariski dense in
$\mathrm{SL}_5(\R)$, and its complex Zariski closure is
$\mathrm{SL}_5(\C)$. Every finite-index subgroup has the same
closure, which is connected. Its actions on $\C^5$ and on the
irreducible Cartan component $W_{2,3}^{\mathrm{fl}}$ are therefore
strongly irreducible. The nonzero invariant span of the critical flag
lines is the full module. The critical complement is nonempty:
choose distinct $x,y\in C$ and nonzero real
$u\in\xi^1(x)$, $v\in\xi^1(y)$. If
$u+iv\in\xi^3(t)_\C$, both real vectors would lie in $\xi^3(t)$,
contrary to Frenet transversality. Theorem~\ref{thm:crit-recurrence}
therefore gives a nonempty largest domain.

There is no regular critical seed. If spectral critical flags at
$x_j,y_j\in C$, $x_j\ne y_j$, both converge to the flag at $x$,
the boundary embedding gives $x_j,y_j\to x$. Put
$\zeta^r(t)=(\xi^{5-r}(t))^\perp$. Applying the Frenet confluence
property to the dual Hitchin representation gives
\[
 \zeta^2(x_j)\oplus\zeta^2(y_j)\longrightarrow\zeta^4(x).
\]
Taking annihilators, their common middle lines satisfy
\[
 \xi^3(x_j)\cap\xi^3(y_j)
       \longrightarrow(\zeta^4(x))^\perp
       =\xi^1(x)\subset\xi^2(x).
\]
The same holds after complexification, excluding the complement
condition in Definition~\ref{def:rec-regular-seed}. Hence removing
that condition when $p-a=1$ allows the criterion to cover more groups.
The direct Borel--Anosov return proof above applies here as well.
Stecker's classification already gives the largest domain; the return
calculation reaches every point of its critical wall even when the
seed's limiting complement fails.
\end{example}

To recognise recurrence without choosing a spectral construction,
we can work at a fixed precision. Compactness makes a finite family
of moving-source returns suffice at each precision, although the
family may grow as the error decreases.

\begin{proposition}[Finite return witnesses at each precision]
\label{r2bd:finite-returns}
Let a countable group $\Gamma$ act on a compact metric space $(X,d)$
by homeomorphisms, and let $\ell:\Gamma\to\mathbb N$ have finite
sublevels. For $g\in\Gamma$ set
\[
 r_g(x)=\min_{z\in X}\max\{d(z,x),d(gz,x)\}.
\]
A nonempty compact set $K\subset X$ is contained in
$\mathscr R_\Gamma$ if and only if for every $n\ge1$ there is
$L(n)>n$ such that
\begin{equation}\label{r2bd:finite-return-test}
 \max_{x\in K}\ \min_{n<\ell(g)\le L(n)}r_g(x)<2^{-n}.
\end{equation}
Here the minimum over an empty set is $+\infty$.
If $K$ is also invariant and the action on $X\setminus K$ is proper,
this test is equivalent to $K=\mathscr R_\Gamma$, and then
$X\setminus K$ is the largest domain of proper discontinuity.
\end{proposition}
\begin{proof}
Each $r_g$ is $1$-Lipschitz, by the triangle inequality followed by
minimisation. Suppose $K\subset\mathscr R_\Gamma$ and fix $n$.
For each $x\in K$ choose a moving-source return $(g_x,z_x)$ with
$\ell(g_x)>n$ and both errors less than $2^{-n-1}$.
The inequalities with the larger bound $2^{-n}$ persist in an open
neighbourhood of $x$. A finite subcover gives finitely many $g_x$;
let $L(n)$ be their largest length. The minimum in
\eqref{r2bd:finite-return-test} is continuous and strictly below
$2^{-n}$ on compact $K$, so its maximum has that property.
Conversely, for each fixed $x\in K$ choose minimising $g_n,z_n$ in
\eqref{r2bd:finite-return-test}. Then $\ell(g_n)>n$ and
$z_n,g_nz_n\to x$. Finite sublevels imply that $g_n$ leaves every
finite subset; extracting pairwise distinct elements witnesses
$x\in\mathscr R_\Gamma$. Finally, properness off $K$ gives
$\mathscr R_\Gamma\subset K$, and equality leaves no dynamical
pair in the complement because that action is proper.
Theorem~\ref{thm:relational} identifies it as the largest domain.
\end{proof}
The finite family depends on $n$, and the argument gives no bound on
$L(n)$. Thus the test characterises recurrence of a compact critical
wall whose complement is already proper, but a check at finitely many
precisions gives no stopping criterion for all the others.

\paragraph{Other return constructions.}
Principal representations and Nori groups have direct recurrence proofs below. In the unitary case we use the independent largest-domain theorem of Cano--Liu--L\'opez~\cite{CLL}, which does not require the regular-seed criterion. Compactness of a unitary middle action does not by itself check the limiting complement. Likewise, the adjoint obstruction of Proposition~\ref{prop:adrec-obstruction} asks us to examine the ambient module as well as the middle-space geometry before applying the recurrence theorem.

\paragraph{A critical wall larger than the recurrent set.}
The Hitchin example required a different return construction. A more
serious difficulty appears when the critical type varies: a
complementary inverse crossing can be the whole vector space, even
though the group has a nonempty equicontinuity region.

\begin{proposition}\label{prop:oversized}\
Let $\Gamma<\PGL_{16}(\C)$ be discrete with critical level $q_\Gamma=7$ and a divergent sequence with Cartan blocks of dimensions $(10,1,4,1)$. Then $\Lambda_{\mathrm{crit}}=\C\PP^{15}$. If the equicontinuity region is nonempty, then $\mathscr R_\Gamma\subsetneq\Lambda_{\mathrm{crit}}$.
\end{proposition}
\begin{proof}
The middle cut is vector dimension $8$. The direct sequence crosses it first at the space of dimension $10$; the inverse sequence has dimensions $(1,4,1,10)$, cumulative $1,5,6,16$, and crosses the complementary cut only at all of $\C^{16}$. The complementary family of an admissible sequence is retained, so the paired critical wall is all of $\C\PP^{15}$. On the equicontinuity region the action is properly discontinuous, so no point there is recurrent; hence the recurrent set is strictly smaller than the wall.
\end{proof}
\paragraph{An arithmetic realisation.}
We now realise the hypotheses of Proposition~\ref{prop:oversized} in
one Zariski-dense group: the large Cartan block responsible for the wall,
and an open region where the
action is equicontinuous. The combination theorem of
Danciger--Gu\'eritaud--Kassel lets us keep the arithmetic subgroup
producing the prescribed sequence while enlarging it to a Zariski-dense
group.

\begin{proposition}\label{r2arith:existence}
There is a discrete, complex-Zariski-dense subgroup
$G<\PGL_{16}(\C)$ with $\Eq(G)\ne\varnothing$, critical level
$q_G=7$, and a divergent Cartan extraction with block dimensions
$(10,1,4,1)$. Moreover, $G$ contains a semisimple element whose
sixteen eigenvalue moduli are distinct. Thus
\[
 \mathscr R_G\subsetneq\Lambda_{\mathrm{crit}}(G)=\C\PP^{15}.
\]
\end{proposition}

\begin{proof}
\emph{The subgroup to be preserved.}
Let
\[
 W=\operatorname{Sym}^2\R^5\oplus\R,
 \qquad \tau=\operatorname{Sym}^2\oplus\mathbf1,
 \qquad \Gamma_0=[\tau(\mathrm{SL}_5(\mathbb Z))]<\PGL(W).
\]
Here $\dim W=16$ and $\det\tau(g)=(\det g)^6=1$.
This projective representation is faithful on $\mathrm{SL}_5(\R)$.
Indeed, a projectively scalar $\tau(g)$ is the identity on its trivial
summand, so $\operatorname{Sym}^2g=I$. Then every $gv$ is equal to
$v$ or $-v$; linearity gives $g=\pm I$, and determinant one in odd
dimension gives $g=I$. In the usual integral symmetric-tensor basis,
$\tau(\mathrm{SL}_5(\mathbb Z))$ consists of integral determinant-one
matrices, so it is discrete. Its projective image is discrete as well:
projective convergence to an invertible matrix permits bounded
determinant-one representatives, and a bounded set of integral matrices
is finite.

The hypotheses of the combination theorem are verified in
\cite[Proposition~4.3, Theorem~1.6 and \S4.6]{r2arith:DGK2025}.
The positive definite symmetric matrices give an invariant properly
convex open subset of $\PP(\operatorname{Sym}^2\R^5)$.
Proposition~4.3 of \cite{r2arith:DGK2025} thickens it to a $\Gamma_0$-invariant properly convex
open subset $\Omega_0$ of $\PP(W)$ whose full orbital limit set is
contained in $\PP(\operatorname{Sym}^2\R^5)$. Its convex core is
therefore a proper subset of $\Omega_0$. Section~4.6 of \cite{r2arith:DGK2025} chooses a
Zariski-dense discrete subgroup $\Gamma_1<\PGL(W)$ preserving a
properly convex open set with a proper convex core.
Theorem~1.6 provides $b\in\PGL(W)$ such that
\[
 G=\langle\Gamma_0,b\Gamma_1b^{-1}\rangle<\PGL(W)
\]
is discrete and preserves a nonempty properly convex open subset
$\Omega$ of $\PP(W)$. It is real-Zariski-dense because it contains
$b\Gamma_1b^{-1}$. This application retains the given $\Gamma_0$
unchanged. No cocompactness of $\Gamma_0$ is needed; Corollary~1.9
explicitly allows a nonuniform prescribed lattice.

\emph{A sequence with bounded internal gaps.}
Put $K=\mathbb Q(\zeta_{11}+\zeta_{11}^{-1})$.
The cyclotomic field has degree ten, and complex conjugation fixes
$K$ with index two; its five embeddings are real. Thus $K$ is
totally real of degree five. In an integral basis of $\mathcal O_K$,
multiplication by a totally positive unit $v$ gives
$M(v)\in\mathrm{SL}_5(\mathbb Z)$: its determinant is the norm of
$v$, which is positive and equal to $\pm1$, hence equal to one.
The real embedding map simultaneously diagonalises these matrices,
with positive eigenvalues $\sigma_i(v)$, $1\le i\le5$.

By the full-lattice form of Dirichlet's unit theorem
\cite[Theorem~5.9]{r2arith:Milne2020}, the logarithms of totally
positive units form a full lattice in
\[
 H=\{s\in\R^5:\textstyle\sum_i s_i=0\}.
\]
Indeed, they contain the squares of all units, whose logarithms are
twice the full unit lattice. A bounded fundamental domain gives units
$v_n$ and a constant $C$ such that, for
$\omega=(1,1,1,1,-4)$,
\begin{equation}\label{r2arith:unit-ray}
 |\log\sigma_i(v_n)-n\omega_i|\le C\qquad(1\le i\le5).
\end{equation}
These logarithms escape every bounded set; take a subsequence of
distinct units. On $W$, the logarithms of the eigenvalues of
$\tau(M(v_n))$ are grouped as
\begin{equation}\label{r2arith:scales}
 \underbrace{2n+O(1),\ldots,2n+O(1)}_{10},\quad
 0,\quad
 \underbrace{-3n+O(1),\ldots,-3n+O(1)}_{4},\quad
 -8n+O(1),
\end{equation}
with uniform error bounds. These matrices are conjugate to their
positive diagonal forms by one fixed real matrix $S$.
For each singular value, multiplication by $S$ and $S^{-1}$ changes
its logarithm by at most
$\log(\|S\|\|S^{-1}\|)$, by the min--max inequalities.
Consequently every adjacent singular ratio inside each of the four
displayed groups stays bounded, while all three ratios between groups
diverge. Extracting the bounded ratios and the compact singular frames
therefore gives the Cartan block partition $(10,1,4,1)$ for a sequence in
this same group $G$.

\emph{The critical level.}
The same unit lattice supplies a spectrally regular element. Exclude
from $H$ the finitely many hyperplanes
\[
 s_i+s_j=s_k+s_l\quad(\{i,j\}\ne\{k,l\}),
 \qquad s_i+s_j=0\quad(i\le j).
\]
They are proper: distinct pair weights restrict distinctly to $H$,
and no pair weight restricts to zero. A full lattice cannot be
contained in a finite union of proper hyperplanes. To see this, write
the product of their defining linear forms in a lattice basis; a
nonzero polynomial cannot vanish on every point of $\mathbb Z^4$,
by induction on the number of variables. Choose a totally positive
unit whose logarithm avoids these hyperplanes. Its image under $\tau$
is semisimple with sixteen distinct eigenvalue moduli, including the
extra eigenvalue one. Powers of this element have all adjacent
singular ratios divergent, since the diagonalising matrix is fixed.
Their first middle crossing has vector dimension eight and projective
dimension seven. Every middle crossing has projective dimension at
least seven, so the definition of the critical level gives $q_G=7$.

\emph{Complexification and equicontinuity.}
The inclusion $\PGL_{16}(\R)\hookrightarrow\PGL_{16}(\C)$ is a
closed topological embedding, so $G$ remains discrete. It is
complex-Zariski-dense: the real and imaginary parts of a complex
polynomial vanishing on $G$ vanish on the real ambient group by real
Zariski density, and the real form is Zariski-dense in its
complexification.

We verify complex equicontinuity for the whole group. Let $C_0$ be
the closed pointed cone over $\Omega$, with nonempty interior. Choose
a Euclidean unit functional $\varphi\in\operatorname{int}C_0^*$.
Compactness of the section $\{y\in C_0:\varphi(y)=1\}$ gives
$M>0$ such that
\[
 \|y\|\le M\varphi(y)\qquad(y\in C_0).
\]
Fix a real unit vector $x\in\operatorname{int}C_0$ and
$\varepsilon>0$ with $x\pm\varepsilon u\in C_0$ for every real
unit vector $u$. Each element of $G$ has a real lift $A$ preserving
$C_0$. Hence
\[
 2\varepsilon\|Au\|
 \le\|A(x+\varepsilon u)\|+\|A(x-\varepsilon u)\|
 \le2M\varphi(Ax)\le2M\|Ax\|.
\]
The real and complex Euclidean operator norms agree, so
$\|Ax\|\ge(\varepsilon/M)\|A\|$.
For $z\in W\otimes\C$ with $\|z-x\|<\varepsilon/(2M)$,
\begin{equation}\label{r2arith:kernel-avoidance}
 \|Az\|/\|A\|\ge\varepsilon/(2M)>0.
\end{equation}
Every normalised matrix limit is therefore defined throughout the
corresponding complex projective neighbourhood of $[x]$.
Compactness of the unit sphere of matrices then gives local uniform
convergence of a subsequence of every sequence of projective maps
there. This neighbourhood lies in $\Eq(G)$.

For clarity, it also contains no recurrent point. Otherwise choose
distinct $g_j$ and $z_j\to z$ with $g_jz_j\to z$, and extract
unit-norm lifts $A_j\to A$ and unit-norm inverse lifts $B_j\to B$.
Discreteness makes the projective sequence divergent, so $BA=0$.
Inequality~\eqref{r2arith:kernel-avoidance}, applied to both sequences,
gives $Az\ne0$ and $Bz\ne0$, whereas recurrence gives
$[Az]=[z]$. This contradicts $BAz=0$.

Finally the inverse of the sequence~\eqref{r2arith:scales} has block
dimensions $(1,4,1,10)$ and cumulative dimensions $1,5,6,16$.
The forward sequence $(10,1,4,1)$ is admissible and $q_G=7$, so its paired inverse degree is $16-(7+1)=8$. This degree crosses the inverse flag only at the whole vector space. The inverse orientation is not itself admissible; it is retained only as that complementary pair. Thus
Proposition~\ref{prop:oversized} gives
$\Lambda_{\mathrm{crit}}(G)=\C\PP^{15}$. The nonempty neighbourhood
just constructed contains no recurrent point, proving strictness.
\end{proof}

\begin{remark}\label{r2arith:scope}
The regular element above fixes the minimal critical level. The
group has sequences with different middle crossings, so the
uniform-type hypothesis of Theorem~\ref{thm:crit-recurrence} fails.
\end{remark}

The construction keeps a regular element while allowing the middle
crossings to vary. We have used subspaces and exterior powers to
understand this distinction. For a general flag manifold, which
features of the projective return construction have an intrinsic
replacement?

\begin{question}\label{q:tits}
For a general simple $G$ and flag manifold, what replaces the four
hypotheses of Theorem~\ref{thm:crit-recurrence}? A possible direction
is to use balanced systems supported on self-opposed types
$\iota\Theta=\Theta$, and to check recurrence in the Tits
representations for the corresponding fundamental weights.
What intrinsic condition replaces a one-dimensional middle quotient?
In type $A$, the middle quotient $B/A$ has dimension $r=p-a$;
the block group acting on it is $\mathrm{GL}_r$, whose derived factor
$\mathrm{SL}_r$ has semisimple rank $r-1$.
Thus $p-a=1$ gives a toral middle block (reductive rank one and
semisimple rank zero), while semisimple rank one corresponds to
$p-a=2$. Whether these structures provide the required return
construction, including a regular seed when needed, is part of the
question.
\end{question}

\section{Strong irreducibility does not force a largest domain}\label{sec:cp3}
A wall that is too large leaves open the possibility of a larger domain. The other obstruction is a dynamical pair whose two endpoints are non-recurrent. Such a pair rules out every largest
domain. Alderete--Cabrera--Cano--M\'endez construct Schottky examples
for which $\Omega_{\Kul}$ is not the largest domain of proper
discontinuity \cite{AldereteCabreraCanoMendez}. The example below gives a strongly irreducible one.

Theorem~\ref{thm:crit-recurrence} gives a sufficient criterion involving the critical type, the module and the middle block. Strong irreducibility, the algebraic analogue of ``non-elementary'', does not suffice by itself, already in $\C\PP^3$.

\begin{proposition}\label{prop:cp3-strong-irr-no-largest}
There is a faithful, discrete, strongly irreducible subgroup $G_0<\PSL(4,\C)$, isomorphic to $F_2\times F_2$, and points $x,y\in\C\PP^3\setminus\mathscr R_{G_0}$ with $(x,y)\in\D_{G_0}$. Consequently $G_0$ has no largest domain of proper discontinuity, and there are at least two incomparable maximal ones (Theorem~\ref{thm:relational}(iii)--(iv)).
\end{proposition}
\begin{proof}
\emph{The group.} Let $T=\begin{psmallmatrix}1&10\\0&1\end{psmallmatrix}$, $R=\begin{psmallmatrix}21&-40\\10&-19\end{psmallmatrix}$, $\Gamma_0=\langle T,R\rangle<\mathrm{SL}(2,\mathbb Z)$, and let $G_0$ be the image of $\Gamma_0\times\Gamma_0$ acting on $\PP(M_2(\C))=\C\PP^3$ by $(A,B):[M]\mapsto[AMB^{-1}]$. With $D_0=\{|z|<\tfrac14\}$, $D_2=\{|z-2|<\tfrac14\}$, $U_T=\{|\Re z|>5\}$ in the affine coordinate of $\C\PP^1$ one checks $T^k(D_0\cup D_2)\subset U_T$ and $R^k(D_0\cup U_T)\subset D_2$ for $k\ne0$; ping-pong gives $\Gamma_0\cong F_2$, $-I\notin\Gamma_0$, and $D_0$ wanders. In particular no escaping sequence in $\Gamma_0$ has a diagonal rank-one normalised limit: a limit proportional to $E_{22}$ would force $A_m(0)\to 0$, and the case $E_{11}$ follows by applying this to inverses (for $\det=1$ the adjugate is the inverse).

\emph{Discreteness, faithfulness.} The operator $M\mapsto AMB^{-1}$ has integral entries and determinant $1$; a bounded set of integral matrices is finite, so the image is projectively discrete (the covering $\mathrm{SL}_4\to\PSL_4$ has finite fibre). If $AMB^{-1}=cM$ for all $M$ then $M=I$ gives $A=cB$ and $B$ commutes with everything, so $A,B$ are scalar, hence trivial.

\emph{Strong irreducibility.} A finite-index subgroup contains powers of $(T,I),(R,I),(I,T),(I,R)$, and $T=I+10E_{12}$, $R=I+10N$ with $N=\begin{psmallmatrix}2&-4\\1&-2\end{psmallmatrix}$, $N^2=0$. An invariant subspace of $M_2(\C)$ is stable under left and right multiplication by $E_{12}$ and $N$; with $I$ these generate $M_2(\C)$ (use $E_{12}N=E_{11}-2E_{12}$, $E_{22}=I-E_{11}$, $E_{21}=N-2E_{11}+4E_{12}+2E_{22}$). A nonzero subspace stable under all left and right multiplications is everything ($E_{pi}ME_{jq}\in\C^\times E_{pq}$). So every finite-index subgroup is irreducible, and a finite invariant union of proper subspaces would give a finite-index stabiliser of one of its members.

\emph{The dynamical pair.} Fix $0<\theta<\pi/4$, $a=e^{i\theta}$, $b=e^{-i\theta}$, $t_m=10m$. Then
\[
T^m\begin{pmatrix}a&0\\(b-a)/t_m&b\end{pmatrix}T^{-m}=\begin{pmatrix}b&0\\(b-a)/t_m&a\end{pmatrix},
\]
so $(T^m,T^m)$ sends sources tending to $x=[\mathrm{diag}(a,b)]$ to points tending to $y=[\mathrm{diag}(b,a)]$: $(x,y)\in\D_{G_0}$.

\emph{No returns at $x$ or $y$.} Let $M=\mathrm{diag}(a,b)$ (the case $\mathrm{diag}(b,a)$ is symmetric). If $x\in\mathscr R_{G_0}$ there are $M_m,N_m\to M$ (normalised to $\det=1$ by a local square root) and distinct $(A_m,B_m)$ with $A_mM_m=\varepsilon_mN_mB_m$, $\varepsilon_m=\pm1$. The norms of $A_m,B_m$ are comparable and tend to infinity; after extraction $A_m/\|A_m\|\to D$, $B_m/\|B_m\|\to E$, $\varepsilon_m\|B_m\|/\|A_m\|\to r\in\R^\times$, with $D,E$ real rank one and $DM=rME$, i.e.\ $D=rMEM^{-1}$. Conjugation by $M$ multiplies the $(1,2)$ and $(2,1)$ entries of $E$ by $e^{2i\theta}$, $e^{-2i\theta}$, which are non-real while $D,E,r$ are real: both off-diagonal entries of $E$ (hence of $D$) vanish. So $D,E$ would be diagonal rank-one limits in $\Gamma_0$, which was excluded. Hence $x,y\notin\mathscr R_{G_0}$, and Theorem~\ref{thm:relational}(iii)--(iv) gives the conclusion.
\end{proof}
\begin{remark}
The same obstruction occurs for the solvable toral group of Section~\ref{sec:toral} and for a diagonal real cyclic group in Kapovich--Leeb--Porti (Example~5.5). For $G_0$ the obstruction is a dynamical edge $(x,y)$ between two non-recurrent points for a strongly irreducible group. The uniform-type hypothesis of Theorem~\ref{thm:crit-recurrence} already fails here: the sequences $(T^m,I)$ and $(T^m,T^m)$ have Cartan block dimensions $(2,2)$ and $(1,2,1)$, respectively. The theorem therefore makes no largest-domain assertion for this group.

Toledo-Acosta \cite[Theorem~3.4.6]{ToledoAcosta2019} states a
smallest-limit-set conclusion under strong irreducibility.
If this means least among all closed invariant sets whose complement
has a proper action, the preceding example rules it out:
Theorem~\ref{thm:relational} puts $x$ and $y$ in such domains,
while their relation prevents a largest one.
A narrower class of limit sets would need to be specified.
Strong irreducibility alone therefore cannot replace the hypotheses
of Theorem~\ref{thm:crit-recurrence}.
\end{remark}

%--------------------------------------------------------------
\sectionpart{Principal and unitary images}

\section{Principal images of \texorpdfstring{$\PSL(2,\C)$}{PSL(2,C)}}\label{sec:principal}
In a principal representation, a projective point is a divisor on the
sphere. We can therefore read the destination relation by following its
roots: some stay away from the repeller, while the others may move towards
it. How many roots must coalesce for the divisor to return to itself?
The equality of the Kulkarni and Myrberg sets below uses a nonempty classical
ordinary set. This hypothesis excludes uniform lattices, for which the
two sets already differ in degree three, as the full-limit-set
calculation below shows.
For non-elementary convex-cocompact groups the largest projective domain also follows from Stecker's classification in Section~\ref{sec:compare}. The calculation below gives the recurrence wall for every non-elementary discrete group, with no finite-generation or convex-cocompactness assumption. When the classical limit set is proper, it also computes Kulkarni's set using sources in a compact set disjoint from $L_0\cup L_1$. The two walls are then separated by their multiplicity thresholds.
Fix an integer $N\ge1$ and let $\rho_N:\PSL(2,\C)\to\PGL(\mathrm{Sym}^N\C^2)$ be the irreducible representation, so that $\PP(\mathrm{Sym}^N\C^2)=\mathrm{Sym}^N(\C\PP^1)$ is the space of effective divisors of degree $N$ on $\C\PP^1$ (factor a binary form). Let $\Gamma<\PSL(2,\C)$ be non-elementary and discrete, $G=\rho_N(\Gamma)$, $\Lambda=\Lambda(\Gamma)\subset\C\PP^1$, $t=\lceil N/2\rceil$.
For a source Cartan parameter $s$ the extreme singular values of $\rho_N$ have ratio $e^{2Ns}$ in the $\mathrm{SU}(2)$-invariant norm. This ratio tends to infinity on every escaping sequence, so the projective representation is proper and its image of $\Gamma$ is discrete.

\begin{theorem}[Principal images of $\PSL(2,\C)$]\label{thm:C}
\begin{enumerate}[label=\textup{(\alph*)},nosep]
\item The critical set is
$\Lambda_{\mathrm{crit}}(G)=C_\Lambda:=\{D:\ \mathrm{mult}_\xi D\ge t\ \text{for some }\xi\in\Lambda\}=\mathscr R_G$. So $\Omega_{\max}=\{D:\ \mathrm{mult}_\xi D<t\ \forall\xi\in\Lambda\}$ is the largest domain of proper discontinuity. No assumption of convex cocompactness, finite generation or absence of parabolics is made.
\item If $\Lambda\ne\C\PP^1$, then $\Lambda_{\Myr}(G)=\Lambda_{\Kul}(G)=\bigcup_{\xi\in\Lambda}H_\xi=\{D:\ \mathrm{supp}D\cap\Lambda\ne\varnothing\}$ with $H_\xi=\{D:\mathrm{mult}_\xi D\ge1\}$, and $\Omega_{\Kul}=\Eq=\mathrm{Sym}^N(\C\PP^1\setminus\Lambda)$.
\item If $\Lambda\ne\C\PP^1$, the three regions $\Eq\subseteq\Omega_{\Kul}\subseteq\Omega_{\max}$ agree for $N=1,2$, whereas $\Eq=\Omega_{\Kul}\subsetneq\Omega_{\max}$ for $N\ge3$.
\end{enumerate}
\end{theorem}
For $\Lambda=\C\PP^1$ the complement of $C_\Lambda$ is exactly the set of \emph{stable} divisors of geometric invariant theory~\cite{ThomasGIT} (all multiplicities $<N/2$; semistability allows $=N/2$). For a proper $\Lambda$ the condition is imposed only over $\Lambda$, and this relative locus is the right object for the discrete group.

\begin{proof}[Proof of (a)]
\emph{Properness on $\Omega_{\max}$.} Let $\gamma_j\in\Gamma$ be distinct and $D_j\to D$, $\gamma_jD_j\to D'$. By the classical convergence principle there are $a,b\in\Lambda$ with $\gamma_j\to a$ locally uniformly on $\C\PP^1\setminus\{b\}$. Write divisors as unordered $N$-tuples. If $k=\mathrm{mult}_bD$ then at least $N-k$ of the points of $D_j$ stay away from $b$, so $\mathrm{mult}_aD'\ge N-k$. If $D\notin C_\Lambda$ then $k\le t-1$ and $N-k\ge N-t+1\ge t$, so $D'\in C_\Lambda$. Thus no pair of $\D_G$ has both ends outside $C_\Lambda$, and Theorem~\ref{thm:relational}(i) gives properness on $\Omega_{\max}$. ($C_\Lambda$ is compact: image of $\Lambda\times\mathrm{Sym}^{N-t}$ under $(\xi,E)\mapsto t\xi+E$.)

\emph{Recurrence of $C_\Lambda$.} Let $D=k\xi+E$ with $\xi\in\Lambda$, $k=\mathrm{mult}_\xi D\ge t$, $\xi\notin\mathrm{supp}E$, $\deg E=N-k$. By the coalescence lemma (Section~\ref{sec:coalesce}) choose loxodromic $s_j\in\Gamma$ with attracting and repelling points $a_j,b_j\to\xi$; for large $j$ neither lies in $\mathrm{supp}E$. For each fixed $j$ choose a power $\gamma_j=s_j^{M_j}$, distinct from the earlier ones, with $\gamma_jE$ within $1/j$ of $(N-k)a_j$ and $\gamma_j^{-1}E$ within $1/j$ of $(N-k)b_j$. Since $2k-N\ge0$ the divisor
\[
D_j=E+\gamma_j^{-1}E+(2k-N)a_j
\]
is effective of degree $N$, and $\gamma_ja_j=a_j$ gives $D_j\to E+k\xi=D$ and $\gamma_jD_j=\gamma_jE+E+(2k-N)a_j\to D$. So $D\in\D_G(D)$, i.e.\ $D\in\mathscr R_G$. Hence $C_\Lambda\subseteq\mathscr R_G$; with properness on the complement, $\mathscr R_G=C_\Lambda$.

\emph{Identification with $\Lambda_{\mathrm{crit}}$.} The weights of $\mathrm{Sym}^N$ are one-dimensional and ordered; the target flag of a divergent sequence with attractor $a$ is the osculating flag of the rational normal curve at $a$, $\mathcal O_k(a)=u^{N-k}\mathrm{Sym}^k\C^2$, whose critical member is $\PP(\mathcal O_m(a))$ with $m=\lfloor N/2\rfloor$. This is exactly $\{D:\mathrm{mult}_aD\ge N-m=t\}$. Since attractors of divergent sequences are dense in $\Lambda$ (loxodromic fixed points), $\Lambda_{\mathrm{crit}}=C_\Lambda$.
\end{proof}
\begin{proof}[Proof of (b) and (c)]
The first normalised projective limit of a sequence with classical attractor $a$ and repeller $b$ has image $[u_a^N]$ and kernel $H_b$. Thus Theorem~\ref{thm:dest} and the projective kernel criterion give
\[
 \Lambda_{\Myr}=\bigcup_{b\in\Lambda}H_b,
 \qquad \Eq=\operatorname{Sym}^N(\C\PP^1\setminus\Lambda).
\]
The union is closed: it is the image of the compact set
$\Lambda\times\operatorname{Sym}^{N-1}(\C\PP^1)$ under $(b,E)\mapsto b+E$.

For the Kulkarni assertion put $S=L_0\cup L_1$. We first prove
$S\subseteq\operatorname{Sym}^N(\Lambda)$. An infinite stabiliser of a divisor has an infinite subgroup fixing each of its finitely many support points. Along an escaping sequence in that subgroup, every such point must be the classical attractor or repeller, and hence belongs to $\Lambda$. Taking closures treats $L_0$. For a fixed-source orbit limit, every root other than the repeller tends to the attractor. If the repeller itself is a root, its images remain in the closed invariant set $\Lambda$. Thus all roots of the limit lie in $\Lambda$, which also treats $L_1$.

The general inclusion $\Lambda_{\Kul}\subseteq\Lambda_{\Myr}$ from Section~\ref{sec:KulMyr} gives one direction. For the other, write $D=k\xi+E$ with $\xi\in\Lambda$, $k\ge1$ and $\xi\notin\operatorname{supp}E$. Choose $z_0\notin\Lambda$ and, by Lemma~\ref{lem:principal-coalescing-fixed-pairs}, loxodromics $s_j$ with both fixed points $a_j,b_j$ tending to $\xi$. For large $j$ their fixed points avoid $\operatorname{supp}E$. Choose powers $\gamma_j=s_j^{M_j}$, escaping every compact subset and distinct from all preceding choices, such that
\[
 \gamma_jz_0\longrightarrow\xi,
 \qquad \gamma_j^{-1}E\longrightarrow(N-k)\xi.
\]
Both requirements follow from the north--south dynamics of each fixed $s_j$, followed by diagonal selection. Set
\[
 D_j=kz_0+\gamma_j^{-1}E,
 \qquad D_\infty=kz_0+(N-k)\xi.
\]
The set $K=\{D_\infty,D_1,D_2,\ldots\}$ is compact and disjoint from $S$, since each of its members has a root at $z_0\notin\Lambda$. Moreover
$\gamma_jD_j=k\gamma_jz_0+E\to D$. Consequently $D\in L_2$, as required by the compact-source definition; the admissibility of these sources is part of the argument.

For $N=1,2$ the threshold $t=\lceil N/2\rceil$ is $1$, so excluding a root over $\Lambda$ and excluding a root of multiplicity at least $t$ are the same condition. For $N\ge3$, $t\ge2$ and the divisor $\xi+(N-1)z_0$, with $\xi\in\Lambda$ and $z_0\notin\Lambda$, belongs to $\Omega_{\max}$ but not to $\Omega_{\Kul}$. This proves (c).
\end{proof}

The two multiplicity conditions also lead to different complex geometry.
Allowing a simple root over the limit set enlarges the proper-action
domain enough to admit compact curves.
\begin{corollary}[The largest domain and complex geometry]\label{cor:principal-separation}
Under the hypotheses of Theorem~\ref{thm:C}(b), if $N\ge3$, the largest domain $\Omega_{\max}$ is connected, contains a projective line, and is neither Stein nor Kobayashi hyperbolic. In contrast, every connected component of $\Eq=\Omega_{\Kul}$ is complete Kobayashi hyperbolic, hyperbolically embedded, taut, and Stein.
\end{corollary}
\begin{proof}
The locus of divisors with distinct roots is connected and dense in $\operatorname{Sym}^N(\C\PP^1)$. It is contained in $\Omega_{\max}$ because $t\ge2$. Every component of this open set meets that locus, so $\Omega_{\max}$ is connected.
Choose independent linear factors $u,v$ whose factor points lie outside $\Lambda$. The line $L=\PP\operatorname{Span}(u^N,v^N)$ lies in $\Omega_{\max}$. Its endpoints have their only root outside $\Lambda$; every other member $\alpha u^N+\beta v^N$ has $N$ distinct factors. Holomorphic functions on $\Omega_{\max}$ are constant on the compact curve $L$, so they do not separate its points and the domain is not Stein. The inclusion $\C\PP^1\hookrightarrow\Omega_{\max}$ also excludes Kobayashi hyperbolicity. For the final assertion, choose $2N+1$ distinct points $\xi_i\in\Lambda$. The hyperplanes $H_{\xi_i}$ are in general position: their defining covectors are evaluations of binary forms of degree $N$, and evaluation at any $N+1$ distinct points is independent by polynomial interpolation. They lie in $\Lambda_{\Myr}$ by Theorem~\ref{thm:C}(b); the equicontinuity region contains $Nz_0$ for $z_0\notin\Lambda$. Theorem~\ref{thm:limit-hyperplanes-complete-hyperbolicity} therefore applies.
\end{proof}
This distinguishes two questions already at the level of explicit divisors: proper discontinuity is governed by the middle multiplicity, whereas equicontinuity excludes every root over the classical limit set.

\subsection{Coalescing fixed points}\label{sec:coalesce}
Both moving-source constructions above use loxodromics whose two endpoints
approach the same prescribed point. We record why these are available at
every classical limit point, including parabolic fixed points.
\begin{lemma}\label{lem:principal-coalescing-fixed-pairs}
Let $\Gamma<\PSL(2,\C)$ be non-elementary and discrete. For every $\xi\in\Lambda(\Gamma)$ there are loxodromic $s_j\in\Gamma$ with $s_j^+\to\xi$ and $s_j^-\to\xi$.
\end{lemma}
\begin{proof}
Pick distinct $h_j\in\Gamma$ with $h_jo\to\xi$ in the ball model. After extraction the classical convergence principle gives a repeller $\eta$ such that $h_j\to\xi$ uniformly on compact subsets of $\C\PP^1\setminus\{\eta\}$. Choose a loxodromic $\ell_0\in\Gamma$. Non-elementarity makes the orbit of $\eta$ infinite, so there is $g\in\Gamma$ with $g^{-1}\eta\notin\{\ell_0^+,\ell_0^-\}$. Thus $\ell=g\ell_0g^{-1}$ has both fixed points different from $\eta$. The loxodromics $s_j=h_j\ell h_j^{-1}$ have fixed points $h_j\ell^\pm$, and both tend to $\xi$.
\end{proof}

%--------------------------------------------------------------
\section{Complex hyperbolic lattices and \texorpdfstring{$\mathrm{PU}(1,s)$}{PU(1,s)}}\label{sec:unitary}
For the standard rank-one unitary action, the middle wall and the
hyperplane wall coincide. The geometry of the Hermitian form makes this
visible: both are made of tangent hyperplanes to the boundary of the ball.
Let $h$ have one negative and $s$ positive directions on $V=\C^{s+1}$, $s\ge2$, with $\mathbb H^s_\C=\{[v]:h(v,v)<0\}$ and $S=\partial\mathbb H^s_\C$. For an isotropic line $\ell$, write $T_\ell=\PP(\ell^{\perp_h})$.

\begin{corollary}[$\mathrm{PU}(1,s)$]\label{cor:PU1s}
Let $\Gamma<\mathrm{PU}(1,s)$ be non-elementary and discrete, and let $\Lambda_0\subset S$ be its Chen--Greenberg limit set. Then
\[
 \Lambda_{\mathrm{crit}}(\Gamma)=\mathscr R_\Gamma
   =\Lambda_{\Kul}(\Gamma)=\Lambda_{\Myr}(\Gamma)
   =\bigcup_{\ell\in\Lambda_0}T_\ell.
\]
Its complement is $\Eq(\Gamma)$ and is the largest domain of proper discontinuity. If $\Lambda_0=S$, in particular for a lattice, this complement is the ball $\mathbb H^s_\C$.
\end{corollary}
\begin{proof}
The rank-one Cartan decomposition has singular blocks of dimensions $(1,s-1,1)$ and exponents $(t,0,-t)$. Its fast target line is isotropic and the sum of the first two target blocks is its $h$-orthogonal hyperplane. Consequently both the critical target carrier and the exceptional slow-source carrier are tangent hyperplanes. The attractors and repellers of divergent sequences range over $\Lambda_0$; compactness of that set makes the union of these hyperplanes closed. The block formula therefore gives
\[
 \Lambda_{\mathrm{crit}}=\Lambda_{\Myr}
       =\bigcup_{\ell\in\Lambda_0}T_\ell.
\]
Cano--Liu--L\'opez~\cite[arXiv:1506.08113v5, Theorems~0.1 and~5.1]{CLL} prove that this union is $\Lambda_{\Kul}$ and that its complement is $\Eq$ and the largest domain for every non-elementary discrete $\Gamma$. The relational criterion, Theorem~\ref{thm:relational}, then identifies the union of tangent hyperplanes with $\mathscr R_\Gamma$. No Zariski-density hypothesis is needed. If $\Lambda_0=S$, every point outside the ball lies on a tangent hyperplane, while no point of the ball does, proving the last assertion.
\end{proof}
Here the published unitary theorem shows that the domain is largest, and the Cartan
description identifies the wall; the hypotheses
of Theorem~\ref{thm:crit-recurrence} are not needed. In dimension two we can also construct
the destinations directly, as in Proposition~\ref{prop:PU21}.

\section{Higher unitary rank and absence of an ordinary region}
\label{r2app:unitary-higher-rank}

When the unitary rank increases, the ball domain disappears.
For a lattice in the standard projective representation of
$\mathrm{PU}(r,s)$, with $r,s\ge2$, every projective point is recurrent.
To see the difference from Corollary~\ref{cor:PU1s}, we look at point
stabilisers. They now contain escaping elements, and finite covolume
turns those elements into returns. The following form of Poincar\'e
recurrence carries out the passage from the homogeneous quotient to
moving sources.

\begin{lemma}[Finite volume and a noncompact stabiliser]
\label{r2app:lattice-stabilizer}
Let $L$ be a unimodular Lie group, let $\Gamma<L$ be a lattice, and let
$L$ act continuously on a compact metrisable space $X$.  Suppose $p\in X$
is fixed by an element $a\in L$ such that $a^n$ leaves every compact subset
of $L$ as $n\to+\infty$.  Then
\[
 Lp\subseteq\mathscr R_\Gamma.
\]
In particular, $p$ is recurrent with moving sources.  The conclusion does
not assert that the orbit of $p$ itself returns to $p$.
\end{lemma}

\begin{proof}
Right translation $T(\Gamma s)=\Gamma sa$ preserves the finite,
full-support invariant measure on $M=\Gamma\backslash L$.
Poincar\'e recurrence~\cite[\S8]{r2appPoincare1890}, applied to a
countable basis and followed by diagonal selection, gives a conull,
hence dense set of cosets with
\[
 T^{n_j}(\Gamma s)\longrightarrow\Gamma s,
 \qquad n_j\longrightarrow\infty.
\]
Their inverse image is dense in $L$, since the quotient map is open.

For a recurrent coset choose representatives of the returns so that
\[
 b_j:=\gamma_j s a^{n_j}\longrightarrow s,
 \qquad \gamma_j\in\Gamma.
\]
Such representatives exist in a local quotient chart.  Since $a^{n_j}p=p$,
\[
 \gamma_j(sp)=b_jp\longrightarrow sp.
\]
Moreover $\gamma_j$ leaves every compact subset of $L$: otherwise the
identity $a^{n_j}=s^{-1}\gamma_j^{-1}b_j$ contradicts the hypothesis on
$a$.  After extraction the $\gamma_j$ are pairwise distinct.  Hence $sp$
is recurrent, even with a fixed source, for a dense set of $s\in L$.
Continuity makes these points dense in $Lp$.  The moving-source recurrent
set is closed, as observed in Section~\ref{sec:crit-proj}; therefore it
contains the whole orbit $Lp$.
\end{proof}

\begin{proposition}[Lattices in higher unitary rank]
\label{r2app:unitary-full-recurrence}
Let $h$ be a nondegenerate Hermitian form of signature $(r,s)$ on $V$,
where $2\le r\le s$, and let $\Gamma<\mathrm{PU}(r,s)$ be a lattice in the
standard projective representation.  Then
\[
 \mathscr R_\Gamma=\Lambda_{\mathrm{crit}}(\Gamma)
 =\Lambda_{\Kul}(\Gamma)=\Lambda_{\Myr}(\Gamma)=\PP(V).
\]
Consequently $\Eq(\Gamma)=\varnothing$ and there is no nonempty open set
of proper discontinuity in $\PP(V)$.  The same absence of a nonempty
domain of proper discontinuity holds for every projective group containing $\Gamma$.
\end{proposition}

\begin{proof}
Fix a projective point $p=[v]$.  If $h(v,v)\ne0$, the orthogonal complement
$v^{\perp_h}$ has signature $(r-1,s)$ or $(r,s-1)$.  Both entries are
positive.  Choose a nondegenerate plane of signature $(1,1)$ inside this
complement.  Its special unitary group, extended by the identity on the
orthogonal complement of the plane, fixes $v$ and contains a one-parameter
subgroup with eigenvalues $e^t,e^{-t}$.  Its positive integral powers
escape compact subsets of $\mathrm{PU}(r,s)$.

If $h(v,v)=0$, choose an isotropic vector $w$ with $h(v,w)=1$.  In the Witt
splitting
\[
 V=\C v\oplus C\oplus\C w,
 \qquad C=(\C v+\C w)^{\perp_h},
\]
the maps
\[
 a_t=\diag(e^t,\operatorname{Id}_C,e^{-t})
\]
preserve $h$ and fix $[v]$.  Again $a_1^n$ escapes every compact set.
Lemma~\ref{r2app:lattice-stabilizer}, applied separately to each $p$, gives
$\mathscr R_\Gamma=\PP(V)$.

The projective image of $\Gamma$ is discrete; equivalently one may take
$\Gamma$ in the projective group from the outset.  The critical complement
is a domain of proper discontinuity by Proposition~\ref{r3crit:cover}, and the Kulkarni and
Myrberg complements are such domains by Proposition~\ref{prop:L2} and
Theorem~\ref{thm:exc}.  Every one of them avoids $\mathscr R_\Gamma$;
hence all three walls are $\PP(V)$.  The equicontinuity formula follows
from Proposition~\ref{prop:kernel}.  Finally, the recurrent witnesses from
$\Gamma$ remain witnesses for every larger group, which proves the last
assertion without a discreteness assumption on that larger group.
\end{proof}

Total recurrence answers the question about domains. We can still ask
where the Cartan flags first cross the middle dimension. The same
finite-volume returns realise the smallest possible crossing.
\begin{proposition}[The middle crossing for unitary lattices]
\label{r2app:unitary-crossing}
Let $1\le r\le s$ and let $\Gamma<\mathrm{PU}(r,s)$ be a lattice in its
standard representation on $V$.  Among the first target spaces of
projective dimension at least $\lfloor(r+s-1)/2\rfloor$ in the limiting
Cartan flags of divergent sequences of $\Gamma$, the smallest realised
projective dimension is $s-1$.
\end{proposition}

\begin{proof}
For a compatible positive-definite Hermitian norm the ordered singular
values of a unitary-indefinite transformation, up to a common scalar, are
\[
 e^{t_1},\ldots,e^{t_r},
 \underbrace{1,\ldots,1}_{s-r},
 e^{-t_r},\ldots,e^{-t_1},
 \qquad t_1\ge\cdots\ge t_r\ge0.
\]
After extraction, exactly $q$ of the $t_i$ tend to infinity and the
others remain bounded, where $1\le q\le r$.  Every bounded coordinate,
its reciprocal, and the $s-r$ unit coordinates belong to one central
Cartan block.  The spaces before that block have vector dimension at
most $q$; the space including it has vector dimension $r+s-q\ge s$.
If $r<s$, the middle vector threshold is
$\lceil(r+s)/2\rceil>r$, so the first crossing is this central block and
has projective dimension at least $s-1$.  If $r=s$, the threshold is $r$.
It is reached among the positive blocks only when $q=r$, at dimension
$r$; otherwise the central block again crosses it at dimension at least
$r$.  The same lower bound follows in either case, including the case of
an absent central block.

To realise equality, take a Cartan one-parameter subgroup $a_t$ with
$t_i=c_i t$, where $c_1>\cdots>c_r>0$.  Right translation by $a_1$ has a
recurrent coset in the finite-volume quotient, by the measure argument
in Lemma~\ref{r2app:lattice-stabilizer}.  Its returns give
$\gamma_j=b_j a_{-n_j}s^{-1}\in\Gamma$, where $b_j\to s$ and
$n_j\to\infty$.  The norm inequalities for left and right multiplication
by fixed bounded invertible matrices give
\[
 \log\sigma_i(\gamma_j)
   =\log\sigma_i(a_{-n_j})+O(1)
\]
for every $i$.  Thus the realised blocks have dimensions
$(1,\ldots,1,s-r,1,\ldots,1)$, with $r$ one-dimensional blocks on either
side and the zero block omitted.  The middle crossing has vector dimension
$s$ if $r<s$ and $r=s$ otherwise, hence projective dimension $s-1$.
\end{proof}

In rank one the stabiliser of a line in the ball is compact, so this
return construction is unavailable there. When both indices are at least
two, every nonisotropic orthogonal complement contains a $(1,1)$-plane,
and the obstruction fills projective space without a condition on a
limiting middle Levi action.

Finite covolume is one way to obtain total recurrence. A dense family of
infinite stabilisers gives another, even in a purely parabolic group.
\begin{example}[A purely parabolic group with total recurrence]
\label{r2app:parabolic-plane}
The same obstruction occurs already on $\C\PP^2$, without any loxodromic
elements.  Put
\[
 \Gamma_0=\left\{
 u(a,b)=\begin{pmatrix}1&a&b\\0&1&0\\0&0&1\end{pmatrix}:
 (a,b)\in\mathbb Z[i]^2\right\}<\PSL(3,\C).
\]
This is a discrete group isomorphic to $\mathbb Z^4$; every nonidentity
element is unipotent.  It satisfies
\[
 L_0(\Gamma_0)=\mathscr R_{\Gamma_0}
 =\Lambda_{\Kul}(\Gamma_0)=\Lambda_{\Myr}(\Gamma_0)=\C\PP^2,
 \qquad \Eq(\Gamma_0)=\varnothing.
\]
Indeed,
$u(a,b)[x:y:z]=[x+ay+bz:y:z]$.  For each nonzero Gaussian pair $(a,b)$,
every point of the line $ay+bz=0$ is fixed by all powers of $u(a,b)$ and
has infinite stabiliser.  These lines have dense union: on the chart
$y\ne0$ their parameters $z/y$ range over the dense field $\mathbb Q(i)$.
The closure in the definition of $L_0$ therefore gives $L_0=\C\PP^2$.
Points with infinite stabiliser are recurrent, so closedness also gives
total moving-source recurrence.  The remaining conclusions follow as in
Proposition~\ref{r2app:unitary-full-recurrence}.  This example belongs to
the setting of the classification of purely parabolic groups by
Barrera--Cano--Navarrete--Seade~\cite{r2appBCNSParabolic}; the calculation
above proves the claims needed here directly.
\end{example}

\section{Relative stability of principal images}
\label{r2app:principal-perturbations}

The strict separation in the principal examples of Theorem~\ref{thm:C}
survives suitable perturbations. For classical Schottky
groups it does: the moving roots and the rational normal curve remain
available as we deform the original M\"obius transformations. This is an
open family in the source group. An arbitrary ambient projective
deformation need not preserve the curve on which the calculation rests.

For $N\ge3$ write
\[
 X_N=\operatorname{Sym}^N(\C\PP^1),\qquad
 m=\lfloor N/2\rfloor,\qquad t=\lceil N/2\rceil.
\]
For a classical limit set $\Lambda$ put
\[
 H_\Lambda=\{D:\operatorname{supp}D\cap\Lambda\ne\varnothing\},
 \qquad
 C_\Lambda=\{D:\operatorname{mult}_\xi D\ge t
                       \text{ for some }\xi\in\Lambda\}.
\]
The first wall detects a root, and the second detects a middle
multiplicity.

\begin{theorem}[An open Schottky family in the source group]
\label{r2app:schottky-stability}
For every integer $k\ge2$ there is a nonempty Euclidean-open set
$\mathcal U\subset\PSL(2,\C)^k$ such that the group
$\Gamma_{\mathbf a}=\langle a_1,\ldots,a_k\rangle$ is free, discrete,
non-elementary, and has nonempty ordinary set for every
$\mathbf a\in\mathcal U$.  For its principal image
$G_{\mathbf a}=\rho_N(\Gamma_{\mathbf a})$, with classical limit set
$\Lambda_{\mathbf a}$, one has
\begin{align*}
 \Lambda_{\mathrm{crit}}(G_{\mathbf a})
 &=\mathscr R_{G_{\mathbf a}}=C_{\Lambda_{\mathbf a}}
 \subsetneq H_{\Lambda_{\mathbf a}}
 =\Lambda_{\Kul}(G_{\mathbf a})=\Lambda_{\Myr}(G_{\mathbf a}),\\
 \Eq(G_{\mathbf a})&=\Omega_{\Kul}(G_{\mathbf a})
 \subsetneq\Omega_{\mathrm{crit}}(G_{\mathbf a})
 =\Omega_{\max}(G_{\mathbf a}).
\end{align*}
The critical carrier of every divergent sequence has projective dimension
$m$.  There is a continuous choice of points
$D_{\mathbf a}\in H_{\Lambda_{\mathbf a}}\setminus C_{\Lambda_{\mathbf a}}$.
All these conclusions are preserved by simultaneous projective conjugation.
\end{theorem}

\begin{proof}
Choose pairwise disjoint closed round discs $D_i^\pm$ with nonempty open
complement $F$. Large powers of loxodromics with fixed points in the
respective interiors give
\begin{equation}\label{r2app:strict-pingpong}
 a_i(\C\PP^1\setminus\operatorname{int}D_i^-)
       \Subset\operatorname{int}D_i^+,
 \qquad
 a_i^{-1}(\C\PP^1\setminus\operatorname{int}D_i^+)
       \Subset\operatorname{int}D_i^-.
\end{equation}
Compactness preserves these strict inclusions on a neighbourhood
$\mathcal U$ of the marked tuple, with fixed points in their discs.
Every nonempty reduced word sends $F$ into its first-letter disc,
proving freeness and, by the positive distance from a fixed $z\in F$ to
the discs, discreteness. The disjoint fixed-point pairs give
non-elementarity. The translates of $F$ are pairwise disjoint and their
union is an open set on which the action is proper: a compact subset meets only finitely many
translates, and only finitely many group elements relate their indices.
Thus $F$ lies in the ordinary set throughout $\mathcal U$.

For the $\mathrm{SU}(2)$-invariant norm on $\operatorname{Sym}^N\C^2$,
a source sequence with Cartan parameter $s_j\to\infty$ has logarithmic
singular values, up to common normalisation,
\[
 Ns_j,(N-2)s_j,\ldots,-Ns_j.
\]
All consecutive gaps diverge. The ratio $e^{2Ns_j}$ of the extreme
singular values makes $\rho_N$ proper, and its faithful action on the
rational normal curve makes it injective. Hence the image is discrete
and escaping source and image sequences correspond. Each has a complete
osculating flag with first-middle projective dimension $m$.
Theorem~\ref{thm:C} gives the walls and domains.

The attracting fixed point $\xi_{\mathbf a}=a_1^+$ varies continuously
and belongs to $\Lambda_{\mathbf a}$.  The fixed point $z\in F$ stays
outside that limit set.  Hence
\[
 D_{\mathbf a}=\xi_{\mathbf a}+(N-1)z
\]
lies on $H_{\Lambda_{\mathbf a}}$ and has multiplicity only one over the
limit set, whereas $t\ge2$.  It is the asserted continuous witness to
strict separation.  Finally, conjugation transports destinations,
recurrent points, stabilisers and compact-source limits, proving the last
assertion.
\end{proof}

\begin{remark}[Deforming the source group]
\label{r2app:relative-openness}
The neighbourhood $\mathcal U$ is open in
$\operatorname{Hom}(F_k,\PSL(2,\C))$.  Its simultaneous-conjugation
saturation is open there and has open image in the corresponding source
character space with the quotient topology.  After composition with
$\rho_N$, the theorem describes this principal family and its projective
conjugates.  It gives no assertion of openness in
$\operatorname{Hom}(F_k,\PGL(N+1,\C))$.  Such an ambient perturbation can
lose the invariant rational normal curve on which the multiplicity proof
rests.  In particular it does not follow that equality of equicontinuity,
Kulkarni and largest proper-action regions is an open property for general
projective representations; even within the displayed family the first
two regions are strictly smaller than the third.
\end{remark}

For a single generator we can allow arbitrarily small ambient
perturbations and continue its eigenlines. The resulting cyclic family
also separates two questions: the critical wall remains smaller than the
Kulkarni wall, yet there is no largest domain. The cyclic calculations continue the dimension-two study of Navarrete~\cite{Navarrete2008} and the higher-dimensional analysis of Cano--Loeza--Ucan-Puc~\cite{CLU}; the compact-source argument below
identifies precisely the wall that persists.

\begin{proposition}[Ambient openness for a cyclic group]
\label{r2app:cyclic-ambient-open}
Let $d\ge4$, let $m=\lfloor(d-1)/2\rfloor$, and let
$g=[A]\in\PGL(d,\C)$ have an eigenbasis $e_1,\ldots,e_d$ with
\[
 |\lambda_1|>|\lambda_2|>\cdots>|\lambda_d|>0.
\]
Put $E_j=\langle e_1,\ldots,e_j\rangle$,
$E'_j=\langle e_{d-j+1},\ldots,e_d\rangle$, and $\Gamma=\langle g\rangle$.
Then $\Gamma$ is infinite and discrete, and
\begin{align*}
 \Lambda_{\mathrm{crit}}(\Gamma)
 &=\PP(E_{m+1})\cup\PP(E'_{m+1}),\\
 \Lambda_{\Kul}(\Gamma)=\Lambda_{\Myr}(\Gamma)
 &=\PP(E_{d-1})\cup\PP(E'_{d-1}),\\
 \mathscr R_\Gamma&=\{[e_1],\ldots,[e_d]\}.
\end{align*}
The critical wall is strictly smaller than the other two walls.  These
formulas and this strict separation hold throughout an open neighbourhood
of $g$ in $\operatorname{Hom}(\mathbb Z,\PGL(d,\C))=\PGL(d,\C)$, with
the eigenlines continued locally.  None of these groups has a largest
domain of proper discontinuity.
\end{proposition}

\begin{proof}
In eigenvector coordinates, positive powers have all consecutive
singular-value gaps divergent and complete fast flags $E_j$; negative
powers give $E'_j$.  Conjugation by one fixed basis change affects
logarithmic singular values by a bounded amount.  The normalised image in
$\bigwedge^jV$ tends to the line
$\C(e_1\wedge\cdots\wedge e_j)$ for positive powers, which verifies the
asserted target flags also in the original coordinates.  The powers therefore
escape compact sets as the absolute value of the exponent tends to
infinity, proving discreteness.  The critical cut is at projective
dimension $m$, giving the first formula by the critical-system definition;
the paired source carriers are contained in these same two spaces.

Let $S=\{[e_1],\ldots,[e_d]\}$.  Distinct eigenvalue moduli imply that
a point fixed by a nonzero power has a single nonzero coordinate.  Every
fixed-source orbit limit is also an eigenpoint.  Thus $L_0\cup L_1=S$.
For positive powers, a source whose first nonzero coordinate has index
$j$ has destination $\PP(E_j)$, by Theorem~\ref{thm:blockform}.  A compact
set disjoint from $S$ has no limiting source $[e_d]$, so all its
positive-power cluster points lie in $\PP(E_{d-1})$.  Negative powers give
$\PP(E'_{d-1})$, proving one inclusion for the Kulkarni wall.

For the other inclusion, let
$y=[e_{d-1}+\sum_{i<d-1}c_i e_i]$ and put
\[
 x_n=\left[e_d+e_{d-1}
       +\sum_{i<d-1}c_i(\lambda_{d-1}/\lambda_i)^n e_i\right].
\]
Then $x_n\to[e_{d-1}+e_d]\notin S$ and $g^nx_n\to y$.  The convergent
sequence of these sources, together with its limit, is a compact set
disjoint from $S$.  Hence $y\in L_2$.  Such points are dense in
$\PP(E_{d-1})$; repeat the argument for the inverse to obtain the other
hyperplane.  The normalised limits of positive and negative powers have
kernels exactly these two hyperplanes, which also proves the Myrberg
formula.  Since $m+1\le d-2$, the critical spaces have smaller dimension
and strict inclusion follows.

For recurrence, a self-destination along positive powers must have first
nonzero index $j$ and also belong to $\PP(E_j)$; it is therefore $[e_j]$.
Negative powers give the same result, and eigenpoints have infinite
stabiliser.  Nevertheless the two nonrecurrent points
\[
 x=[e_2+e_d],\qquad y=[e_1+e_2]
\]
are dynamically related: the sources
$[e_2+e_d+(\lambda_2/\lambda_1)^n e_1]$ tend to $x$ and their images under
$g^n$ tend to $y$.  Theorem~\ref{thm:relational} thus rules out a largest
domain.  Finally, simplicity of the spectrum and the strict inequalities
of moduli persist on an ambient neighbourhood, and its eigenlines vary
continuously.  Applying the same proof to each member gives all the
openness assertions.
\end{proof}

\section{Parabolics and full classical limit sets}
\label{r2app:first-kind-boundary}

The Schottky family keeps the classical limit set proper and has no parabolics. Both restrictions can be removed. The coalescing-points argument still works in the presence of parabolics;
replacing it by a conical-points argument would lose that generality.
When the classical limit set fills the sphere, it is instead the
compact-source calculation of the Kulkarni wall that changes. Strong
irreducibility alone does not decide either question.

\begin{example}[Strong irreducibility does not make every limit point conical]
\label{r2app:modular-safeguard}
The modular group $\PSL(2,\mathbb Z)$ is discrete, strongly irreducible on
$\C^2$, and contains the proximal matrix
$\bigl(\begin{smallmatrix}2&1\\1&1\end{smallmatrix}\bigr)$, but its
parabolic fixed point $\infty$ is not conical.  To check strong
irreducibility, note that a finite invariant set of lines would be fixed pointwise
by a positive power of $T(z)=z+1$ and therefore be contained in
$\{\infty\}$.  The element $S(z)=-1/z$ rules out that possibility.

For the last assertion, a conical sequence at $\infty$ would have
repeller $\infty$ and would keep the images of $\infty$ and $0$
separated in the limit.  Its integral determinant-one lifts have the form
$\gamma_j=\bigl(\begin{smallmatrix}a_j&b_j\\c_j&d_j\end{smallmatrix}\bigr)$.
The repeller condition says that the normalised first column tends to
zero and the normalised second column does not.  Thus the second-column
norm tends to infinity.  The first column is a nonzero integral vector,
so its norm is at least one.  Their chordal separation is
\[
 d_{\mathrm{ch}}(\gamma_j\infty,\gamma_j0)
 =\frac{1}{\sqrt{|a_j|^2+|c_j|^2}\sqrt{|b_j|^2+|d_j|^2}}
 \longrightarrow0,
\]
a contradiction.  Thus proximal endpoints cannot be promoted to conical
behaviour at every point by taking their closure.  The recurrence proof
of Theorem~\ref{thm:C} uses coalescing loxodromic endpoints and applies
also in the presence of parabolics.
\end{example}

For a full classical limit set, the recurrence proof still applies, but
we can no longer choose the root outside that set which made the
Kulkarni sources admissible. The third symmetric power shows exactly
what can happen at this boundary.
\begin{proposition}[A first-kind boundary for the third symmetric power]
\label{r2app:first-kind-sym3}
Let $\Gamma<\PSL(2,\C)$ be non-elementary and discrete, with
$\Lambda(\Gamma)=\C\PP^1$.  Let $G=\rho_3(\Gamma)$ act on
$X_3=\operatorname{Sym}^3(\C\PP^1)$, and let
\[
 \Delta=\{2a+b:a,b\in\C\PP^1\}
\]
be the discriminant of divisors with a multiple point.  Then
\[
 \mathscr R_G=\Lambda_{\mathrm{crit}}(G)=\Lambda_{\Kul}(G)=\Delta,
 \qquad \Lambda_{\Myr}(G)=X_3,
 \qquad \Eq(G)=\varnothing.
\]
The nonempty open set $X_3\setminus\Delta$ is the largest domain of proper
discontinuity and equals $\Omega_{\Kul}(G)$.  In particular the proper-limit-set
hypothesis in Theorem~\ref{thm:C}(b) cannot be removed.
\end{proposition}

\begin{proof}
The critical multiplicity threshold is $t=2$.  Theorem~\ref{thm:C}(a),
which does not assume a proper classical limit set, gives
$\Lambda_{\mathrm{crit}}=\mathscr R_G=\Delta$, and its complement is the
largest domain.  The first projective limits have the root hyperplanes
$H_b=\{D:b\in\operatorname{supp}D\}$ as kernels, for $b$ ranging over
$\C\PP^1$.  Their union is $X_3$, giving
$\Lambda_{\Myr}=X_3$ and $\Eq=\varnothing$.

We compute the Kulkarni wall independently.  An infinite stabiliser of a
divisor has a finite-index subgroup fixing its finite support pointwise.
A nonidentity M\"obius map cannot fix three distinct points.  Therefore
$L_0(G)\subseteq\Delta$.

The ordered pairs of endpoints of loxodromic elements are dense in the
ordered distinct pairs of the classical limit set.  Here is the short
argument needed in this case.  Given disjoint open sets $U,V$ meeting the
limit set, Lemma~\ref{lem:principal-coalescing-fixed-pairs} supplies
loxodromics $s,h$ with both endpoints in $U,V$, respectively. Divide
lifts of $s^n$ and $h^{-n}$ by the $n$th powers of the dominant eigenvalues
of fixed lifts of $s$ and $h^{-1}$, respectively. The limits are rank-one
maps with images $s^+,h^-$ and kernels $s^-,h^+$. Their product is nonzero
because $h^-\ne s^-$;
its trace is nonzero because $s^+\ne h^+$.  Hence normalised lifts of
$s^nh^{-n}$ tend to a rank-one map with a simple nonzero eigenvalue, and
the dominant eigenline tends to $s^+$.  Applying the same calculation to
$h^ns^{-n}$ gives the other endpoint tending to $h^+$.  In particular
$s^nh^{-n}$ is loxodromic for large $n$, with ordered endpoints in $U,V$.

For each loxodromic element with endpoints $a,b$, the divisor $2a+b$ is
fixed by all its powers.  Density of endpoint pairs, followed by closure,
gives every $2a+b$, including $a=b$, in $L_0(G)$.  Thus $L_0(G)=\Delta$.
For a fixed divisor outside $\Delta$, an escaping sequence has a classical
repeller containing at most one of its three distinct roots.  At least
two image roots therefore tend to its attractor, so $L_1(G)\subseteq\Delta$.

Finally, the space of ordered triples of distinct points of $\C\PP^1$
is $\PSL(2,\C)$, on which $\Gamma$ acts properly by left translation.
The unordered triple space $X_3\setminus\Delta$ is its quotient by a
finite permutation group, and the action there is still proper.  A compact
source set in $X_3\setminus\Delta$ cannot have escaping translates
accumulate in that same open set: adjoining a compact neighbourhood of the
putative limit would contradict properness.  Since
$L_0\cup L_1=\Delta$, this proves $L_2(G)\subseteq\Delta$, and hence
$\Lambda_{\Kul}(G)=\Delta$.
\end{proof}

A uniform lattice in $\PSL(2,\C)$ satisfies the hypothesis of the last
proposition: its orbit in hyperbolic three-space is a bounded-distance
net, and comparison with every geodesic ray gives the full sphere as its
classical limit set.  The discriminant complement is the space of triples
of distinct points.  It is therefore possible for a projective action to
have a nonempty largest domain of proper discontinuity while its equicontinuity
region is empty.  The equality of the two regions in the unitary rank-one
case and the equality of the Kulkarni and equicontinuity regions for proper
principal limit sets have different hypotheses and mechanisms.

\sectionpart{Families and complex geometry}

\section{Families with prescribed behaviour}\label{sec:families}
In the following families one can \emph{prescribe} the geometry: first a critical level with returns on every carrier, then
arbitrarily many limiting planes at every Cartan level. Both constructions
use Schottky ping-pong in $\C\PP^n$.

\paragraph{Nori groups.}
We use Nori's projective ping-pong construction~\cite{Nori}.
Let $d=2r$ and choose $2g$ complex subspaces $P_a^\pm\subset\C^{2r}$, $a=1,\dots,g$, $g\ge2$, of dimension $r$, any two of which are complementary. For a small parameter $\alpha$ put $\Lambda_\alpha=(1-\alpha)/\alpha$ and let $\tilde g_a$ be the linear map equal to $\Lambda_\alpha$ on $P_a^+$ and to $1$ on $P_a^-$, so that in the graph chart over $P_a^+$ the action on $r$-planes is $X\mapsto\Lambda_\alpha^{-1}X$. Let $G_\alpha=\langle g_1,\dots,g_g\rangle<\PSL(2r,\C)$.

Contraction is visible in these graph coordinates. To identify the
limiting planes with actual Cartan targets, we must translate it into a
singular-value gap. The following estimate does so uniformly over reduced
words.
\begin{lemma}[Pl\"ucker collapse and the middle singular gap]
\label{r3nor:collapse-gap}
Fix $r\ge1$ and a Hermitian space $V$ of dimension $2r$, and give
$\mathrm{Gr}_r(V)$ its Pl\"ucker angular Fubini--Study distance.  For each
nonempty open set $U\subset\mathrm{Gr}_r(V)$ there is a constant $C_U$
such that, for every $h\in\operatorname{GL}(V)$,
\[
 \operatorname{diam}(hU)\le\delta
 \quad\Longrightarrow\quad
 \frac{\sigma_{r+1}(h)}{\sigma_r(h)}\le C_U\delta.
\]
More generally, it suffices that the images of a fixed basis of
$\wedge^rV$ consisting of Pl\"ucker vectors in $U$ have diameter at most
$\delta$, with a constant depending on that basis.  If a finite symmetric
alphabet $\mathcal A$ is nonempty and has relatively compact nonempty
open ping-pong domains $U_a$
with pairwise disjoint closures and, for every $b\ne a^{-1}$,
\[
 a\overline{U_b}\subset U_a,\qquad
 \operatorname{Lip}(a|_{U_b})\le\rho,\qquad 0<\rho<1,
\]
then one constant $C$ gives
\[
 \frac{\sigma_{r+1}(w)}{\sigma_r(w)}\le C\rho^{|w|}
 \tag{$*$}\label{r3nor:word-gap}
\]
for every nonempty reduced word $w$. Write $\partial\mathcal A$
for the space of infinite reduced words in the alphabet, in which
consecutive letters are never inverse. The limits of the fast target
$r$-planes, over words with lengths tending to infinity, are precisely
$\xi(\partial\mathcal A)$; the slow source $r$-planes have the same set
of limits by applying this assertion to inverse words.
\end{lemma}
\begin{proof}
Put $W=\wedge^rV$.  Pl\"ucker vectors in any nonempty open subset of the
Grassmannian span $W$: a linear form vanishing there vanishes on the
whole Grassmannian by the identity theorem in graph coordinates, and
the coordinate decomposable wedges span $W$.  Choose a unit-vector basis
$v_1,\ldots,v_D$ from $U$, with dual basis $\varphi_1,\ldots,\varphi_D$,
and put $A=\wedge^rh$.  Let $Q$ be orthogonal projection onto
$\C Av_1$.  The diameter assumption gives
\[
 \|(I-Q)Av_i\|\le\delta\|Av_i\|,\qquad
 \|A-QA\|\le
 \delta\Bigl(\sum_{i=1}^D\|\varphi_i\|\Bigr)\|A\|.
\]
Since $QA$ has rank one, the min--max principle bounds
$\sigma_2(A)/\sigma_1(A)$ by the last relative error.  The singular
values of an exterior power give the exact identity
$\sigma_2(\wedge^rh)/\sigma_1(\wedge^rh)
=\sigma_{r+1}(h)/\sigma_r(h)$.  This proves the first assertion with
$C_U=\sum_i\|\varphi_i\|$.

For an alphabet consisting of a single self-inverse letter, every
nonempty reduced word has length one. The estimate then holds with
$C=\rho^{-1}$, and both limit sets are empty. Otherwise, choose such a
basis in each $U_b$. If $w=a_1\cdots a_m$ is reduced, choose
$b\ne a_m^{-1}$.  Its basis images have diameter at most
$D_0\rho^m$, where $D_0=\max_b\operatorname{diam}(U_b)$; taking the
maximum of the finitely many basis constants proves \eqref{r3nor:word-gap}.
Applying the estimate to $A=\wedge^rw$, a top unit left singular
vector $u$ satisfies
$\|(I-Q)u\|\le C\rho^m$.  Thus the top singular Pl\"ucker line is
within $O(\rho^m)$ of $[Av_1]$, a Pl\"ucker point in $wU_b$.
The compact sets
$a_1\cdots a_{k-1}\overline{U_{a_k}}$ are nested and have diameter at
most $D_0\rho^{k-1}$, defining $\xi(a_1a_2\cdots)$.
Every sequence of words of unbounded length has an extracted common
prefix of each length $k$; the preceding estimates identify its singular
limit with that nested intersection.  Conversely the successive prefixes
of each infinite word realise its value of $\xi$.  The fastest $r$-plane
of $w^{-1}$ is the slow source $r$-plane of $w$, proving the final claim.
\end{proof}

The estimate locates the middle cut without requiring that either half
remain a single Cartan block. We can now follow both the limiting planes
and the points which return on them.
\begin{proposition}[Nori groups]\label{prop:nori}
There is $\alpha_0>0$ such that for $0<\alpha<\alpha_0$ the group $G_\alpha$ is free and discrete, has a singular-value gap at index $r$ along every escaping sequence, with critical level $r-1$ (projective dimension of the critical carriers), and
\[
\Lambda_{\mathrm{crit}}(G_\alpha)=\bigcup_{\omega\in\partial F_g}\PP(\xi(\omega))=\mathscr R_{G_\alpha},
\]
where $\xi:\partial F_g\to\mathrm{Gr}_r(\C^{2r})$ is the equivariant boundary map. Hence $\Omega_{\mathrm{crit}}$ is the largest domain of proper discontinuity; it is nonempty and has compact quotient.
\end{proposition}
\begin{proof}
\emph{Ping-pong and the gap.}
Complementarity gives disjoint relatively compact neighbourhoods $U_a$
of the signed expanding planes, with
$\bigcup_{b\ne a^{-1}}\overline{U_b}$ in a bounded graph chart over
$P_a^+$ with complement $P_a^-$. The action
$X\mapsto\Lambda_\alpha^{-1}X$ gives the contraction bounds of
Lemma~\ref{r3nor:collapse-gap}. Classical Nori ping-pong gives freeness
and discreteness~\cite[Proposition~2.9 and \S3.1]{MO2016}. Here the
space of infinite reduced words is the boundary $\partial F_g$ of the
free group. The lemma gives the exponential middle singular gap and
identifies all fast target and slow source $r$-plane limits with
$\xi(\partial F_g)$. Thus the critical carriers are exactly
$\PP(\xi(\omega))$, giving the first equality.

\emph{Recurrence.} Fix $\omega\in\partial F_g$ and $L=\xi(\omega)$. Choose prefixes $w_j\to\omega$ and a signed generator $a$ whose two fixed boundary points avoid the repeller of $(w_j)$; after extraction one fixed choice of $a$ suffices. The convergence action on $\partial F_g$ and continuity of $\xi$ give spectral planes
\[
 U_j=w_jP_a^+\longrightarrow L,
 \qquad V_j=w_jP_a^-\longrightarrow L
\]
for $s_j=w_jg_aw_j^{-1}$. Fix $[v]\in\PP(L)$ and unit representatives $u_j\in U_j$, $v_j\in V_j$ with $u_j,v_j\to v$. For a sufficiently large power $\gamma_j=s_j^{N_j}$ let $R_j=\Lambda_\alpha^{N_j}\to\infty$. A projective lift acts by $R_j$ on $U_j$ and by $1$ on $V_j$. Hence
\[
 z_j=[v_j+R_j^{-1/2}u_j]\longrightarrow[v],
 \qquad \gamma_jz_j=[R_j^{-1/2}v_j+u_j]\longrightarrow[v].
\]
The powers can be chosen to escape a fixed compact exhaustion and to be pairwise distinct. Thus every point of $\PP(L)$ is recurrent. The two spectral planes approach the same plane, and the chosen source balances their two expansion rates to produce the return.

\emph{Largestness.} Proposition~\ref{r3crit:cover} makes the critical
complement a domain of proper discontinuity. The recurrence just proved
identifies it with $X\setminus\mathscr R_\Gamma$, which is the largest
domain by Theorem~\ref{thm:relational}(iii).

\emph{The classical domain.}
Put $\lambda=\Lambda_\alpha$, $\varepsilon=\lambda^{-1/2}$ and
$X=\PP(\C^{2r})$. With norms on each complementary pair, shared by inverse letters, set
\[
 C_s=\{[p+q]:p\in P_s,\ q\in P_{s^{-1}},\
                  \|q\|\le\varepsilon\|p\|\},\qquad
 K_s=\{\operatorname{graph}A:\|A\|\le\varepsilon,\
                         A:P_s\to P_{s^{-1}}\}.
\]
For small $\alpha$ the tubes are disjoint and $K_s\subset U_s$.
Rank-one extension gives $C_s=\bigcup_{L\in K_s}\PP(L)$, while
$\lambda\varepsilon^2=1$ gives the exact pairing
\[
 s(X\setminus\operatorname{int}C_{s^{-1}})=C_s,\qquad
 s(X\setminus C_{s^{-1}})=\operatorname{int}C_s.
\]
The classical domain $\Omega_N=\bigcup_w wF$, where
$F=X\setminus\bigcup_s\operatorname{int}C_s$, is nonempty with compact
quotient~\cite[Proposition~2.9 and \S3.1]{MO2016}.
The exact pairing gives $sK_t\subset\operatorname{int}K_s$ for
$t\ne s^{-1}$. For a reduced infinite word $\omega=s_1s_2\cdots$,
the nested graph sets
$s_1\cdots s_{k-1}K_{s_k}$ contract to $\xi(\omega)$; hence
$\bigcap_k s_1\cdots s_{k-1}C_{s_k}=\PP(\xi(\omega))$.
Successive reduction by tube letters, with cancellation excluded by the
exact pairing, puts every point outside $\Omega_N$ in such an
intersection. Conversely these intersections avoid every translate of
$F$. Thus $\Omega_N=\Omega_{\mathrm{crit}}$.
\end{proof}
The fast and slow $r$-planes are separated by the middle cut; extra gaps
inside either half are allowed. Every first-middle target crossing is an
$r$-plane, so the two carrier families for each extracted sequence
are singletons and
$q_\Gamma=q_\Gamma^*=r-1$. The return was constructed directly on these
planes. No intervening middle quotient or seed hypothesis is needed.

\paragraph{Exterior powers and generic Schottky groups.}
For the Nori groups, the middle cut suffices to determine the largest
domain. To fill the rest of the filtration, we use elements with a gap
at every index and follow their spectral planes through exterior powers.
A group $\Gamma<\PSL(V)$, $V=\C^{n+1}$, is \emph{exterior SIP} if for every $1\le k\le n$ the action of $\wedge^k\Gamma$ on $\PP(\wedge^kV)$ is strongly irreducible and contains a proximal element. An element $g$ is \emph{regular loxodromic} if its eigenvalue moduli are strictly ordered $|\lambda_1|>\dots>|\lambda_{n+1}|$.

\begin{lemma}\label{lem:regprox}
If $g$ is regular loxodromic, then $\wedge^kg$ is proximal for every $k$, with attracting line $[E_{\lambda_1}\wedge\cdots\wedge E_{\lambda_k}]$, which is decomposable. If $\Gamma$ is Zariski dense in $\PSL(V)$, every $\wedge^k\Gamma$ acts strongly irreducibly.
\end{lemma}
\begin{proof}
The eigenvalues of $\wedge^kg$ are the products $\lambda_{i_1}\cdots\lambda_{i_k}$ and the strict ordering makes $\lambda_1\cdots\lambda_k$ the unique one of largest modulus. For the second claim, $\wedge^kV$ is an irreducible algebraic representation; a finite union of proper subspaces invariant under $\wedge^k\Gamma$ is invariant under its Zariski closure, which is the image of the connected group $\mathrm{SL}(V)$; connectedness forbids permuting the finitely many maximal linear components, so each is invariant, contradicting irreducibility.
\end{proof}

The spectral planes are easy to conjugate. Since the filtration is
defined by singular flags, we first check that these planes occur as
Cartan limits, even for a nonnormal matrix.
\begin{lemma}[Spectral planes are Cartan limits]
\label{r3sch:spectral-cartan}
Let $g\in\operatorname{GL}(V)$, $\dim V=n+1$, have ordered eigenvalue
moduli $|\lambda_1|\ge\cdots\ge|\lambda_{n+1}|$ and a strict gap
$|\lambda_k|>|\lambda_{k+1}|$, $1\le k\le n$.
Let $E_k^+$ be the $k$-dimensional generalised spectral subspace
associated with the first $k$ eigenvalues, counted with algebraic
multiplicity.  Then
\[
 \frac{\sigma_{k+1}(g^m)}{\sigma_k(g^m)}\longrightarrow0,
 \qquad
 \text{the fast target $k$-plane of $g^m$}\longrightarrow E_k^+.
\]
Consequently, if a discrete projective group $\Gamma$ contains $[g]$,
then
\[
 \PP(E_k^+)\subset\Lambda_{\mathrm Z}^{(k-1)}(\Gamma).
\]
\end{lemma}
\begin{proof}
For $T=\wedge^kg$ the eigenvalue
$\mu=\lambda_1\cdots\lambda_k$ is algebraically simple and every
other eigenvalue has modulus at most
$|\mu|\,|\lambda_{k+1}/\lambda_k|<|\mu|$.
Jordan normal form therefore gives $\mu^{-m}T^m\to P$, the nonzero
rank-one spectral projector with image $\wedge^kE_k^+$.
This argument allows nonnormal $g$ and Jordan blocks.
Continuity of singular values gives
$\sigma_2(T^m)/\sigma_1(T^m)\to0$; its left singular line tends to
$\operatorname{Im}P$, since the component of a unit top left singular
vector orthogonal to $\operatorname{Im}P$ is bounded by
$\|\mu^{-m}T^m-P\|/\|\mu^{-m}T^m\|$.
These are exactly the ratio $\sigma_{k+1}(g^m)/\sigma_k(g^m)$ and
the Pl\"ucker line of its fast target $k$-plane.  The powers escape
projective compact sets by that gap.  After Cartan extraction, $E_k^+$
is an actual target-flag space of vector dimension $k$, so its
projectivisation belongs to level $k-1$ by the definition in
Section~\ref{sec:flagLS}.
\end{proof}

We may now choose conjugates whose planes are in general position at
every level at once. At each choice only finitely many proper algebraic
sets have to be avoided.
\begin{proposition}[Abundance at every Cartan level]
\label{r3sch:all-levels}
Suppose $\Gamma<\PSL(n+1,\C)$ is discrete and Zariski dense and
contains a regular loxodromic element $g$.  For every $N\ge1$ there
are $N$ conjugates $h_i g h_i^{-1}$ such that, simultaneously for
every $0\le j\le n-1$, their attracting projective $j$-planes lie
in $\Lambda_{\mathrm Z}^{(j)}(\Gamma)$ and are in Pl\"ucker general
position: every subfamily of at most $\binom{n+1}{j+1}$ of their Pl\"ucker vectors
is linearly independent.
\end{proposition}
\begin{proof}
For $k=j+1$ put $D_k=\dim\wedge^kV$ and let $E_k^+$ be the fast
spectral $k$-plane of $g$.  The Pl\"ucker Grassmannian is irreducible
(a homogeneous variety of the connected group $\PSL(V)$) and
nondegenerate because decomposable coordinate wedges span
$\wedge^kV$.  Inductively choose $h_i\in\Gamma$ so that
$[\wedge^k(h_iE_k^+)]$ avoids the spans of all subfamilies of at most
$D_k-1$ previously chosen points, for every $k$ simultaneously.
Each span cuts out a proper closed subset of the Grassmannian; its
inverse image in $\PSL(V)$ is proper and closed by transitivity.
The finite union of these forbidden sets is still proper, since
$\PSL(V)$ is irreducible, and the Zariski-dense group $\Gamma$
therefore has a point outside it.  This proves the general-position
assertions.  Each chosen plane is the spectral plane of a conjugate of
$g$, and Lemma~\ref{r3sch:spectral-cartan} places it in the asserted
scalar level, whose definition uses actual Cartan targets.
\end{proof}

To realise this abundance in an open Schottky family, we need Zariski
density as well as contraction. The adjoint action gives an open
condition which we can impose after taking the powers required by
ping-pong. For $n\ge1$ and $r\ge2$, put $G=\PSL(n+1,\C)$ and
$\mathfrak g=\mathfrak{sl}(n+1,\C)$, and let
$\mathcal O_{\mathrm{Ad}}\subset G^r$ be the set of $r$-tuples whose
adjoint images have no common nonzero proper invariant subspace of
$\mathfrak g$.

\begin{lemma}\label{lem:adopen}
For $n\ge1$ and $r\ge2$, $\mathcal O_{\mathrm{Ad}}$ is nonempty and Zariski open, and an infinite group generated by a tuple in it is Zariski dense in $G$.
\end{lemma}
\begin{proof}
For each $1\le k<\dim\mathfrak g$, the tuples preserving some
$k$-plane form a closed set: project the closed incidence subset of
$G^r\times\mathrm{Gr}_k(\mathfrak g)$ to $G^r$.  The Grassmannian is
projective, so this projection is closed.  Their finite union is the
complement of $\mathcal O_{\mathrm{Ad}}$.

For nonemptiness put $d=n+1$.  In $\mathrm{SL}(d,\C)$ choose
\[
 D=\diag\bigl(p_1,\ldots,p_{d-1},(p_1\cdots p_{d-1})^{-1}\bigr),
 \qquad N=uv^t,\quad
 u=(1,\ldots,1)^t,\quad v=(1,\ldots,1,-(d-1))^t,
\]
where the $p_i$ are distinct primes.  The powers of $D$ are Zariski dense
in the diagonal torus, since no nontrivial torus character takes the
value $1$ on $D$.  Also $N^2=0$, and the closure of the powers of $I+N$
contains $I+tN$ for every $t\in\C$.  The Lie algebra of the closure of
$\langle D,I+N\rangle$ therefore contains $N$ and is invariant under
the diagonal torus.  Every off-diagonal entry of $N$ is nonzero, so its
torus-weight components give every $E_{ij}$, $i\ne j$; their brackets
give the traceless diagonal matrices.  The closure is consequently
$\mathrm{SL}(d,\C)$.  Its projective pair acts irreducibly on
$\mathfrak g$, since $\mathfrak{sl}(d,\C)$ is simple.  Adding identity
generators gives a tuple in $\mathcal O_{\mathrm{Ad}}$ for every $r\ge2$.

Finally, if $A$ is the Zariski closure of an infinite $\Delta$ with
generating tuple in $\mathcal O_{\mathrm{Ad}}$, then
$\mathrm{Lie}(A)\ne0$ is invariant under the adjoint images of the
generators.  It equals $\mathfrak g$, so $A=G$ because $G$ is connected.
\end{proof}

\begin{theorem}[Open realisation of generic Schottky groups]\label{thm:open-realisation}
For all $n\ge2$ and $r\ge2$ there is a nonempty open set $\mathcal U\subset\mathrm{Hom}(F_r,\PSL(n+1,\C))$ such that every $\rho\in\mathcal U$ is faithful with free, discrete, Zariski-dense image $\Gamma_\rho$, regular loxodromic marked generators, exterior SIP, $\Eq(\Gamma_\rho)\neq\varnothing$, and $\Lambda_{\Myr}(\Gamma_\rho)$ contains $2n+1$ hyperplanes in general position, given by a fixed finite list of words throughout $\mathcal U$.
\end{theorem}
\begin{proof}
\emph{Ping-pong.} Choose regular loxodromic $a_{1,0},\dots,a_{r,0}$ in projective ping-pong position, with the attracting points of the signed alphabet distinct and outside the non-corresponding repelling hyperplanes (finitely many proper algebraic conditions). Choose disjoint compact convex neighbourhoods $D_s$ of the attracting points in affine charts, a coordinate ball $B_*$ disjoint from them and from the repelling hyperplanes, and $p_*\in B_*$. For a large common power $N$ the maps $a_{s,0}^N$ send $\overline{B_*}$ and the $D_t$ ($t\ne s^{-1}$) into $\mathrm{int}\,D_s$ with Lipschitz constants $<1$. Put $\rho_0(x_i)=a_{i,0}^N$. Ping-pong gives faithfulness and freeness; the positive distance from $p_*$ to $\bigcup D_s$ prevents nonidentity words converging to the identity, hence discreteness. All strict inclusions persist on a Euclidean-open neighbourhood $\mathcal U_0$, which we shrink so that the marked generators remain regular loxodromic. On it the family $\{\rho(w)\}$ is normal on $B_*$: from any sequence of nonempty reduced words, extract a subsequence with the same first letter $s$. All its images of $B_*$ lie in the same compact set $D_s$ in one affine chart, so the coordinate functions are bounded and Montel applies. A sequence containing infinitely many identity words is already constant along a subsequence. Thus $B_*\subset\Eq$.

\emph{Zariski density after taking powers.} We impose Zariski density on the \emph{final} generators, after the power $N$ has been fixed: $\mathcal O_{\mathrm{Ad}}$ is a nonempty Zariski-open set, hence dense in the Euclidean topology, so $\mathcal U_1=\mathcal U_0\cap\mathcal O_{\mathrm{Ad}}\ne\varnothing$. By Lemma~\ref{lem:adopen}, $\Gamma_\rho$ is Zariski dense; by Lemma~\ref{lem:regprox} it is exterior SIP.

\emph{Hyperplanes.} Fix $\rho_*\in\mathcal U_1$. Proposition~\ref{r3sch:all-levels}
applied to $\rho_*(x_1)$ gives $2n+1$ conjugates with attracting
hyperplanes in general position. Each hyperplane is the kernel of the
first pseudo-projective limit of the inverse powers, so lies in
$\Lambda_{\mathrm{exc}}=\Lambda_{\Myr}$. Write the conjugators as
$\rho_*(v_i)$ and fix the words $v_ix_1v_i^{-1}$. Their spectral
hyperplanes vary continuously, so general position persists on a smaller
open set $\mathcal U$.
\end{proof}

Continuity of finitely many spectral flags makes all-level abundance
open in the same family.
\begin{corollary}[Open realisation at every level]
\label{r3sch:open-all-levels}
For every $n\ge2$, $r\ge2$ and $N\ge1$, the open set in
Theorem~\ref{thm:open-realisation} has a nonempty open subset
$\mathcal U_N$ on which $N$ fixed words are regular loxodromic and,
for every $0\le j\le n-1$, their attracting projective $j$-planes
are in Pl\"ucker general position and are contained in
$\Lambda_{\mathrm Z}^{(j)}(\Gamma_\rho)$.
All the conclusions of that theorem, including its fixed-word
$2n+1$ Myrberg hyperplanes and common nonempty equicontinuity ball,
hold on $\mathcal U_N$.
\end{corollary}
\begin{proof}
Choose $\rho_*$ in that open set and apply
Proposition~\ref{r3sch:all-levels} to $g=\rho_*(x_1)$.
Write the selected $h_i$ as $\rho_*(v_i)$ and fix the words
$v_ix_1v_i^{-1}$.  Regular spectral flags vary continuously with
the marked tuple.  The finitely many nonzero minors expressing
Pl\"ucker general position remain nonzero on a nonempty open
neighbourhood $\mathcal U_N$ of $\rho_*$ inside the original open set.
Lemma~\ref{r3sch:spectral-cartan} gives level membership throughout
that neighbourhood.
\end{proof}

\paragraph{A real-open unitary family.}
The complex-open construction fills the filtration and leaves a common
equicontinuity ball. It does not identify all recurrent points. For an
open family with equality of the three regions, we return to the
rank-one unitary action: the tangent-hyperplane theorem supplies the
returns which general position alone does not give. This time openness
is taken inside the real unitary group.
Put $G=\PSL(n+1,\C)$ and $L=\mathrm{PU}(1,n)<G$.

\begin{lemma}\label{r3fam:real-form-density}
Every nonempty real-open subset of $L^r$ is complex Zariski dense in $G^r$.
\end{lemma}
\begin{proof}
Since $\mathfrak{su}(1,n)\otimes_{\R}\C=\mathfrak{sl}(n+1,\C)$,
exponential coordinates make $L^r$ maximally totally real in $G^r$.
A regular function vanishing on a real-open patch therefore vanishes on a
complex neighbourhood by the holomorphic identity theorem, and then on the
irreducible variety $G^r$.
\end{proof}

\begin{proposition}\label{r3fam:unitary-family}
For every $n\ge2$ and $r\ge2$ there is a nonempty real-open subset
$\mathcal V_r\subset L^r$ which is complex Zariski dense in $G^r$.
Every $\rho\in\mathcal V_r$ is faithful with free, discrete, strongly
irreducible, complex Zariski-dense, non-elementary image $\Gamma_\rho$, and
\[
 \Omega_{\Kul}(\Gamma_\rho)=\Eq(\Gamma_\rho)
      =\Omega_{\max}(\Gamma_\rho)\ne\varnothing.
\]
Every component of this common region is complete Kobayashi hyperbolic,
hyperbolically embedded, taut and Stein; its quotient by the component
stabiliser is a complete Kobayashi hyperbolic complex manifold.
\end{proposition}
\begin{proof}
Choose loxodromic elements of $L$ with $2r$ distinct signed fixed points.
Distinct isotropic lines for $h$ are nonorthogonal, so their
attracting points avoid every noncorresponding repelling hyperplane.
The strict ping-pong construction of Theorem~\ref{thm:open-realisation},
using only contraction at the signed proximal points, large common
powers and a ball $B\Subset\mathbb H^n_\C$, gives a
nonempty real-open neighbourhood $\mathcal V_r^0\subset L^r$ of faithful
free discrete representations with non-elementary image and
$B\subset\Eq(\Gamma_\rho)$.  After these powers have been fixed, set
$\mathcal V_r=\mathcal V_r^0\cap\mathcal O_{\mathrm{Ad}}$.
Lemma~\ref{r3fam:real-form-density} makes this intersection nonempty and
complex Zariski dense; Lemmas~\ref{lem:adopen} and~\ref{lem:regprox} give
full complex Zariski closure and strong irreducibility.

Corollary~\ref{cor:PU1s}, using Cano--Liu--L\'opez~\cite{CLL}, gives the
equality of the three regions.  The signed generators are proximal, so
Corollary~\ref{cor:proximal-equicontinuity-geometry} gives the component
geometry, and Corollary~\ref{cor:complete-hyperbolic-quotients} gives the
quotient assertion because $\Gamma_\rho$ is torsion-free.
\end{proof}

\paragraph{Marked density and the final Zariski topology.}
Both families are abundant in the complex Zariski topology, although one
is open in the complex group and the other only in its real form. We
first keep the generators marked, then pass to their image groups. This
distinguishes the density assertion from either kind of openness.
Write $X_r=\operatorname{Hom}(F_r,G)\cong G^r$, $r\ge2$, and let
$\mathcal K_r\subset X_r$ consist of representations whose image is
discrete, not virtually cyclic, and acts properly discontinuously on some
nonempty invariant open subset of $\C\PP^n$.
Give $\mathcal K_r$ its relative complex Zariski topology and put
\[
 \pi_r(\rho)=\rho(F_r),\qquad
 \mathcal D_{\mathrm{fg}}=\bigcup_{r\ge2}\pi_r(\mathcal K_r),\qquad
 \mathcal D_{\le q}=\bigcup_{2\le r\le q}\pi_r(\mathcal K_r).
\]
These are sets of image subgroups, without identifying conjugate groups.
Give $\mathcal D_{\mathrm{fg}}$ the final topology for all $\pi_r$:
$O$ is open exactly when each $\pi_r^{-1}(O)$ is relatively Zariski open.
Give $\mathcal D_{\le q}$ the subspace topology.

\begin{proposition}\label{r3gen:final-density}
For $n\ge2$ and each $r\ge2$, both the complex-open family $\mathcal U_r$ of
Theorem~\ref{thm:open-realisation} and the real-open family $\mathcal V_r$
of Proposition~\ref{r3fam:unitary-family} are Zariski dense in $X_r$ and
meet every nonempty relatively Zariski-open subset of $\mathcal K_r$.
More generally, if $\mathcal A_r\subset\mathcal K_r$ is relatively Zariski
dense for every $r\ge2$, then
$\bigcup_{r\ge2}\pi_r(\mathcal A_r)$ is dense in
$\mathcal D_{\mathrm{fg}}$, and
$\bigcup_{2\le r\le q}\pi_r(\mathcal A_r)$ is dense in
$\mathcal D_{\le q}$ for every $q\ge2$.
In particular, the image groups of each of the two families are dense in
these spaces, retaining all the corresponding geometric conclusions.
\end{proposition}
\begin{proof}
A nonempty Euclidean-open subset of the irreducible variety $X_r$ is
Zariski dense; the same density for $\mathcal V_r$ is
Lemma~\ref{r3fam:real-form-density}.  Since both families lie in
$\mathcal K_r$, ambient density implies their relative density there.

If $O\subset\mathcal D_{\mathrm{fg}}$ is nonempty and open, choose
$\rho\in\mathcal K_r$ with $\pi_r(\rho)\in O$; then the nonempty open
set $\pi_r^{-1}(O)$ meets $\mathcal A_r$.
For a nonempty subspace-open set $W=O\cap\mathcal D_{\le q}$, first choose
$\Gamma\in W$ and a lift $\rho\in\mathcal K_r$ with $r\le q$.
The same argument yields an element of $\pi_r(\mathcal A_r)\cap W$.
\end{proof}

The Nori family supplies returns on every critical carrier; the
complex Schottky family supplies all-level general position, fixed-word
Myrberg hyperplanes and a common equicontinuity ball. All-level abundance
does not identify its recurrent set: returns are a further condition, as
in Section~\ref{sec:cp3}. The real-open unitary family obtains equality of the
three regions from tangent geometry.

\section{Limiting hyperplanes and complex geometry}
\label{sec:limiting-hyperplanes-geometry}

The analytic geometry of ordinary regions is part of the complex Kleinian programme. Navarrete--Barrera~\cite{BNKobayashi2009} study the Kobayashi metric for complex hyperbolic groups acting on $\C\PP^2$. Barrera--Cano--Navarrete~\cite{BCNLimit2011} establish complete Kobayashi hyperbolicity and the domain-of-holomorphy property for equicontinuity components under their dimension-two hypotheses; see also the regular-set results of~\cite{CCL2026}. Here we use the limiting hyperplanes to formulate the corresponding sufficient conditions in arbitrary projective dimension.

The hyperplanes just constructed constrain holomorphic discs and prevent
them from approaching the boundary at finite Kobayashi distance. This
explains the analytic conclusions for the unitary family, and applies
equally to the equicontinuity regions of the complex Schottky family.
As the principal examples show, these properties need not extend to a
larger domain of proper discontinuity. Throughout this section
$G<\PSL(n+1,\C)$ is discrete, $n\ge1$, and
\[
 C=\Lambda_{\Myr}(G),\qquad \Eq(G)=\C\PP^n\setminus C.
\]
Here $C$ is the closure of the union of the kernels of pseudo-projective
boundary limits, as in Proposition~\ref{prop:kernel}; those limits may
have any rank. A finite family
$\PP(\ker\ell_i)$ of hyperplanes is in \emph{general position} if every
subfamily of at most $n+1$ defining covectors is linearly independent.

The Kobayashi pseudodistance $d_M$ is the infimum of sums of
Poincar\'e distances along chains of holomorphic discs. Holomorphic maps
decrease it; $M$ is \emph{hyperbolic} if it separates points and
\emph{complete hyperbolic} if it is complete. An open
$U\subset\C\PP^n$ is \emph{hyperbolically embedded} if distinct
$p,q\in\overline U$ have neighbourhoods $V_p,V_q$ with
\[
 \inf\{d_U(x,y):x\in U\cap V_p,\ y\in U\cap V_q\}>0.
\]
The manifold is \emph{taut} if every sequence of holomorphic discs has a
locally uniformly convergent subsequence into $M$ or a compactly divergent
one: images of compact subsets eventually avoid every compact subset of
$M$. For completeness we shall use boundary hyperplanes, and for
Steinness a plurisubharmonic exhaustion. We first propagate the finite
hyperplane family through the whole Myrberg wall.

\begin{lemma}[Propagation of limiting hyperplanes]
\label{lem:myrberg-hyperplane-propagation}
Let $F\subset\C\PP^n$ be closed and $G$-invariant.  Suppose that $F$
contains hyperplanes $H_i=\PP(\ker\ell_i)$ whose defining covectors span
$(\C^{n+1})^*$.  Every point of $\Lambda_{\Myr}(G)$ then lies on a
hyperplane contained in $F$.  In particular $\Lambda_{\Myr}(G)\subseteq F$.
\end{lemma}

\begin{proof}
Take normalised linear lifts $A_m$ of distinct elements $g_m\in G$ such
that $A_m\to T\ne0$ is a singular pseudo-projective limit.  Some
$\ell_iT$ is nonzero because the $\ell_i$ span the dual.  In the dual
projective space,
\[
 g_m^{-1}H_i=\PP\ker(\ell_iA_m)
       \longrightarrow K_T:=\PP\ker(\ell_iT).
\]
All the hyperplanes on the left lie in $F$, by invariance; hence
$K_T\subseteq F$, by closedness.  Moreover
$\PP(\ker T)\subseteq K_T$.

To pass to the closure in the definition of $\Lambda_{\Myr}$, let
$p_j\in\PP(\ker T_j)$ tend to $p$.  Choose $K_j\subseteq F$ as above
with $p_j\in K_j$.  Compactness of the dual projective space gives a
subsequence $K_j\to K$.  The incidence correspondence is closed, so
$p\in K$, and closedness of $F$ gives $K\subseteq F$.
\end{proof}

A finite arrangement in general position gives the initial hyperbolic
domain by Green's theorem. Additional limiting hyperplanes need not be in
general position with it. The following observation deals with any finite
number of them.
\begin{lemma}[Adding hyperplanes to a Green arrangement]
\label{lem:green-extra-hyperplanes}
Let $H_1,\ldots,H_{2n+1}$ be hyperplanes in general position in
$\C\PP^n$.  For any additional finite family of hyperplanes
$K_1,\ldots,K_s$, without a general-position hypothesis on them,
\[
 \C\PP^n\setminus
   \left(\bigcup_{i=1}^{2n+1}H_i\cup\bigcup_{j=1}^sK_j\right)
\]
is complete hyperbolic and hyperbolically embedded.
\end{lemma}

\begin{proof}
Discard repetitions.  For the original arrangement, a nonempty stratum
lying in the intersection of $r$ hyperplanes is an open subset of
$\C\PP^{n-r}$ obtained by deleting $2n+1-r$ hyperplanes in general
position.  If $n-r>0$, then
\[
 2n+1-r\ge2(n-r)+1,
\]
so Green's hyperplane theorem makes the stratum hyperbolic.  For
$n-r=0$ it is a point.  Thus no stratum admits a nonconstant entire map
from $\C$.  Adding the $K_j$ replaces each original stratum by its
intersections with selected $K_j$ and the complements of the others.
Every resulting stratum is contained in an original stratum, so it too
admits no nonconstant entire map.  Green's divisor criterion
\cite[Theorem~2]{Green}, in the form also recalled in
\cite[Theorem~4.1]{JavanpeykarLevin}, now applies to the enlarged divisor
and proves hyperbolic embeddedness. No normal-crossings assumption on the
enlarged arrangement is needed in Green's criterion.

For completeness, put
\[
 U_0=\C\PP^n\setminus\bigcup_{i=1}^{2n+1}H_i,
 \qquad U_j=U_{j-1}\setminus K_j\quad(1\le j\le s).
\]
Green's hyperplane theorem makes $U_0$ complete hyperbolic. Suppose
inductively that every component of $U_{j-1}$ is complete, and let
$(z_\nu)$ be a Kobayashi-Cauchy sequence in a component $\Omega$ of
$U_j$. Let $D$ be the component of $U_{j-1}$ containing $\Omega$.
Distance decrease under inclusion makes $(z_\nu)$ Cauchy in $D$, so
it converges there to a point $p$. If $p\in K_j$, then
$K_j\cap D$ is a smooth closed complex hypersurface in the hyperbolic
manifold $D$. The hypersurface-distance theorem
\cite[Theorem~2(A)]{IvashkovichRosay} places $p$ at infinite
Kobayashi distance from $D\setminus K_j$, contradicting the
boundedness of a Cauchy sequence in $\Omega$. Thus $p\notin K_j$.
A small coordinate ball about $p$ in $D\setminus K_j$ contains the
tail of the sequence and belongs to $\Omega$; its Kobayashi distance
bounds $d_\Omega$ from above and proves $d_\Omega(z_\nu,p)\to0$.
Induction proves completeness of every component of $U_s$.
Intersections among the added hyperplanes cause no difficulty: at each
step $K_j\cap D$ is a smooth hypersurface in the domain from which
the earlier hyperplanes have already been removed.
\end{proof}

To find the initial arrangement, it is enough to have one hyperplane and
strong irreducibility. Its orbit then supplies arbitrarily many
hyperplanes in general position.
\begin{proposition}[Amplification by strong irreducibility]
\label{prop:strong-irr-amplifies-hyperplane}
If $G$ is strongly irreducible, the orbit of every hyperplane contains
arbitrarily large finite families in general position.  Consequently a
closed $G$-invariant set containing one hyperplane contains such families.
\end{proposition}

\begin{proof}
The dual action is strongly irreducible: an invariant proper nonzero dual
subspace for a finite-index subgroup would give such a subspace in
$\C^{n+1}$ by taking its annihilator.  An orbit in $\PP((\C^{n+1})^*)$
cannot be contained in a finite union of proper projective subspaces.
Indeed, the Zariski closure of that orbit has finitely many irreducible
components, which $G$ permutes.  If it were contained in that union, each
component would lie in one of those proper subspaces.  The linear span
of any component would then be nonzero and proper and would be preserved
by its finite-index stabiliser in $G$, a contradiction.

Inductively choose orbit covectors avoiding the finitely many spans of
subfamilies of at most $n$ previously chosen covectors.  Each new covector
is independent of every such subfamily.  Thus any subfamily of at most
$n+1$ chosen covectors is independent, as required.
\end{proof}

We now pass from finite arrangements to the whole exceptional wall.
Propagation places a hyperplane through each of its points; these
hyperplanes give both the boundary-distance estimate and the exhaustion.
\begin{theorem}[Complex geometry of the exceptional domain]
\label{thm:limit-hyperplanes-complete-hyperbolicity}
Suppose that $\Eq(G)\ne\varnothing$ and that $\Lambda_{\Myr}(G)$
contains $2n+1$ hyperplanes in general position.  Then every connected
component of $\Eq(G)$ is complete Kobayashi hyperbolic, hyperbolically
embedded in $\C\PP^n$, taut, and Stein.  The hyperplane hypothesis is
automatic if $G$ is strongly irreducible and $\Lambda_{\Myr}(G)$ contains
one hyperplane.
\end{theorem}

\begin{proof}
Choose the stated hyperplanes $H_1,\ldots,H_{2n+1}\subseteq C$ and put
\[
 Y=\C\PP^n\setminus\bigcup_{i=1}^{2n+1}H_i.
\]
Green's hyperplane theorem~\cite[Theorem~1]{Green} says that $Y$ is
complete Kobayashi hyperbolic and hyperbolically embedded in
$\C\PP^n$.  Since the defining covectors span the dual,
Lemma~\ref{lem:myrberg-hyperplane-propagation}, applied to $F=C$, also
shows that every point $p\in C$ lies on a hyperplane $K_p\subseteq C$.
Thus $C$ is a union of complete projective hyperplanes.  The final
assertion follows from Proposition~\ref{prop:strong-irr-amplifies-hyperplane}.

For a component $\Omega$ of $\Eq(G)$, monotonicity gives
\[
 d_Y(x,y)\le d_\Omega(x,y)\qquad(x,y\in\Omega).
\]
It implies hyperbolicity and transfers the positive neighbourhood
separation defining hyperbolic embeddedness from $Y$ to $\Omega$.

For completeness, let $(z_j)$ be $d_\Omega$-Cauchy.  It is $d_Y$-Cauchy
and hence tends to a point $p\in Y$.  Suppose that $p\in C$.  Choose a
hyperplane $K_p\subseteq C$ through $p$ and let $Z$ be the component of
$Y\setminus K_p$ containing $\Omega$.  The set $K_p\cap Y$ is a closed
smooth complex hypersurface in the hyperbolic manifold $Y$.  The
hypersurface-distance theorem~\cite[Theorem~2(A)]{IvashkovichRosay}
therefore gives
\[
 d_Z(z_1,z_j)\longrightarrow+\infty
       \quad\hbox{whenever }z_j\longrightarrow p\in K_p\cap Y.
\]
This contradicts $d_Z\le d_\Omega$ and the boundedness of a Cauchy
sequence. The theorem's smoothness hypotheses hold because both the
complex structure and the hyperplane are smooth. Thus $p\notin C$.
A small coordinate ball about $p$ lies in $\Eq(G)$ and contains $z_j$
for all large $j$, so it belongs to $\Omega$.  The distance of that ball
dominates $d_\Omega$ and tends to zero as $z_j\to p$; hence
$d_\Omega(z_j,p)\to0$.

The Kobayashi distance of a hyperbolic manifold induces its topology
and is a length distance, so completeness and local compactness make its
closed bounded balls compact. For discs $f_j:\Delta\to\Omega$ which are
not compactly divergent, extract $\zeta_j$ in a fixed compact subset of
$\Delta$ with $f_j(\zeta_j)$ in a fixed compact subset of $\Omega$.
Distance decrease gives equicontinuity and bounds $f_j(K)$ in a compact
Kobayashi ball for each compact $K\subset\Delta$. Arzel\`a--Ascoli and
diagonal extraction give a holomorphic limit into $\Omega$, proving tautness.

Use $H_1$ as the hyperplane at infinity. Then $F=C\cap\C^n$ is a closed
nonempty union of affine hyperplanes, since $H_2\ne H_1$. For
$\delta(z)=\operatorname{dist}(z,F)$ and affine equations $L$ with
linear part of norm one,
\[
 -\log\delta(z)
  =\sup_{\{L=0\}\subseteq F}\bigl(-\log|L(z)|\bigr)
       \qquad(z\in\C^n\setminus F).
\]
The supremum of these pluriharmonic functions is continuous and locally
bounded above, hence plurisubharmonic. Thus
\[
 \phi(z)=\|z\|^2-\log\delta(z)
\]
is continuous and strictly plurisubharmonic on $\Omega$.
At finite boundary points $\delta$ tends to zero.  At infinity, fix one
affine hyperplane $\{L_0=0\}\subseteq F$ with normalised linear part.
For a constant $A>0$,
\[
 \delta(z)\le |L_0(z)|\le A(1+\|z\|),\qquad
 \phi(z)\ge\|z\|^2-\log A-\log(1+\|z\|)\longrightarrow+\infty.
\]
Thus $\phi$ exhausts $\Omega$. Richberg's theorem
\cite[Theorem~I.5.21]{Demailly} gives a smooth strictly plurisubharmonic
exhaustion $\phi\le\psi\le\phi+1/2$. The pseudoconvexity criterion and
Levi theorem~\cite[Theorems~I.7.2 and VIII.9.11(a)]{Demailly} make
$\Omega$ a domain of holomorphy, hence Stein
\cite[Theorem~I.6.11 and Definition~I.6.16]{Demailly}.
\end{proof}

The hyperplanes in the theorem can also arise from powers of a single
proximal element. Strong irreducibility then supplies their general
position, just as in the Schottky construction.
\begin{corollary}[A proximal sufficient condition]
\label{cor:proximal-equicontinuity-geometry}
If $G$ is strongly irreducible, contains a proximal element, and has
$\Eq(G)\ne\varnothing$, every connected component of $\Eq(G)$ has the
four properties in Theorem~\ref{thm:limit-hyperplanes-complete-hyperbolicity}.
\end{corollary}

\begin{proof}
For a proximal linear representative $A$, normalised positive powers
converge projectively to its rank-one spectral projection: the dominant
eigenvalue is simple and all other eigenvalues have smaller modulus.
These projective powers are distinct, and their limit is singular.
Its projectivised kernel is therefore a hyperplane in
$\Lambda_{\Myr}(G)$.  Apply the theorem and amplification proposition.
\end{proof}

If the initial hyperplanes already lie in the Kulkarni set, propagation
has a further consequence: it identifies that set with the whole
Myrberg wall.
\begin{corollary}[When the Kulkarni and Myrberg sets agree]
\label{cor:kulkarni-equicontinuity-geometry}
Assume $\Eq(G)\ne\varnothing$.  If $\Lambda_{\Kul}(G)$ contains $2n+1$
hyperplanes in general position, then
\[
 \Lambda_{\Kul}(G)=\Lambda_{\Myr}(G),\qquad
 \Omega_{\Kul}(G)=\Eq(G),
\]
and every connected component of this common region is complete
Kobayashi hyperbolic, hyperbolically embedded, taut, and Stein.  In
particular the same conclusion holds if $G$ is strongly irreducible and
$\Lambda_{\Kul}(G)$ contains a hyperplane.
\end{corollary}

\begin{proof}
The general inclusion $\Lambda_{\Kul}\subseteq\Lambda_{\Myr}$ is proved
in Section~\ref{sec:KulMyr} and recalled in
Theorem~\ref{thm:flag-of-limit-sets}(d).  Apply
Lemma~\ref{lem:myrberg-hyperplane-propagation} to
$F=\Lambda_{\Kul}$ for the reverse inclusion, and then apply
Theorem~\ref{thm:limit-hyperplanes-complete-hyperbolicity}.  In the strongly
irreducible case the amplification proposition supplies the required family.
\end{proof}

Finally, completeness descends to a quotient when the component
stabiliser acts freely. We use the distance formula for a covering to
keep track of Cauchy sequences downstairs.
\begin{corollary}[Complete hyperbolic quotients]
\label{cor:complete-hyperbolic-quotients}
Under the hypotheses of Theorem~\ref{thm:limit-hyperplanes-complete-hyperbolicity},
let $\Omega$ be a connected component of $\Eq(G)$ and let
$G_\Omega=\{g\in G:g\Omega=\Omega\}$.  If $G_\Omega$ acts freely on
$\Omega$, then $\Omega/G_\Omega$ is a complete Kobayashi hyperbolic complex
manifold.  Freeness holds in particular when $G$ is torsion-free.
\end{corollary}

\begin{proof}
Proper discontinuity follows from Theorem~\ref{thm:exc}.  Together with
freeness it makes the quotient projection a holomorphic covering.
Lifting chains of holomorphic discs gives the exact formula
\[
 d_{\Omega/G_\Omega}([x],[y])
     =\inf_{h\in G_\Omega}d_\Omega(x,hy).
\]
This infimum is positive for distinct orbits.  Otherwise, for some $h_j$,
$d_\Omega(x,h_jy)\to0$, hence $h_jy\to x$ in the manifold topology.
Proper discontinuity forbids infinitely many distinct $h_j$ with this
property; an element repeating along a subsequence would instead give
$x=hy$, also a contradiction.

For completeness, a Cauchy sequence downstairs has a subsequence whose
successive distances have finite sum.  Using the infimum formula, choose
successive lifts so that their distances upstairs are bounded by the
corresponding distances downstairs plus $2^{-j}$.  These lifts form a
Cauchy sequence and converge in $\Omega$.  Projection gives a limit of
the subsequence, and hence of the original Cauchy sequence.  Finally,
proper discontinuity makes all point stabilisers finite, so they are
trivial if $G$ is torsion-free.
\end{proof}

These conclusions concern components of $\Eq(G)$ and, when the two sets
agree, components of $\Omega_{\Kul}(G)$.  The principal examples of
Theorem~\ref{thm:C} distinguish these regions from the largest domain of
proper discontinuity: the latter can contain complete projective lines.

The exact destinations distinguish proper discontinuity, the existence of a largest domain, and equicontinuity. The principal examples answer these three questions and show how the boundaries can differ. The critical recurrence theorem gives another route to the largest domain. Its proof uses a regular seed when the middle block has dimension greater than one.

Over $\C$, the set of diagonal return points does not record the scheme structure of the diagonal intersection; that structure is the subject of Chapter~\ref{r2gr:chapter}.

\chapter{Complex graph geometry}\label{r2gr:chapter}

We now return to a single Cartan sequence and view its destination
relation as a limit of graphs. For a complex flag variety this graph
belongs to Brion's compactification. Its reduced boundary fibres and
Schubert components give an algebraic model of the complex destination
relation of Theorems~A and~B
\cite[Theorem~3 and Proposition~10]{r2grBrion2003}.

Since the boundary fibre is reduced, the destination relation determines
its scheme structure. We use this model to study how the graph meets the diagonal. The Cartan factors identify its boundary fibre and retain
the finite-gap Levi map. In projective space we use the block
resolution to calculate the component multidegrees and the cohomology
sheaves of the resolution defect of each block. The residual maps also give a return
operator under the quotient identifications stated below. These
calculations relate the size of an active block to its diagonal
intersection degree and describe the return scheme in those cases.

The flag variety in this section is complex. The real destination
formula remains Theorem~\ref{thm:dest}; the scheme and intersection
statements below use complex algebraic geometry. The single-sequence
statements require neither discreteness nor Anosov hypotheses.

\section{The intrinsic family of graphs}\label{r2gr:intrinsic}

Let $\mathbf G$ be a connected complex semisimple algebraic group acting
faithfully on $X=\mathbf G/P$. We take $\mathbf G$ of adjoint type. For
the graphs of the action this is no restriction: divide first by the
kernel, including any simple factors acting trivially. When we return
to a discrete group, however, its elements must keep their labels;
properness counts group elements, even when they induce the same map.

Choose opposite Borel subgroups $B,B^-$, with $B\subset P$, and put
$T=B\cap B^-$. Let $\Delta$ be the simple roots, $W$ the Weyl group,
and $I\subset\Delta$ the simple roots of the Levi subgroup of $P$.
Thus $I=\Delta\setminus\theta$ in the notation of Chapter~\ref{ch:cartan},
and $W_I=W_\theta$. Write $W^I$ for the minimal representatives of
$W/W_I$, and set
\[
 X_w=\overline{BwP/P},\qquad X_w^-=\overline{B^-wP/P}.
\]
For $w\in W^I$, their dimensions are $\ell(w)$ and
$\dim X-\ell(w)$, respectively. For an arbitrary index these symbols
mean the varieties attached to its right $W_I$-coset.

For $g\in\mathbf G$ let
\[
 \Gamma_g=\{(x,gx):x\in X\},\qquad
 \mathcal H_{X,\mathbf G}
 =\overline{\{[\Gamma_g]:g\in\mathbf G\}}
 \subset\operatorname{Hilb}(X\times X),
\]
where the Zariski closure has its reduced structure. All the graphs have
the same Hilbert polynomial. A polarisation can be chosen to place them in
one projective Hilbert scheme; the subschemes being parametrised are
intrinsic to $X\times X$. Denote the universal family restricted to
$\mathcal H_{X,\mathbf G}$ by $\mathcal Z$, and its fibre at $\xi$ by
$Z_\xi$.

\begin{theorem}[Brion]\label{r2gr:brion}
The variety $\mathcal H_{X,\mathbf G}$ is the wonderful compactification
of $\mathbf G$. The family $\mathcal Z\to\mathcal H_{X,\mathbf G}$ is
flat, and each $Z_\xi$ is reduced, Cohen--Macaulay and of pure dimension
$\dim X$. In characteristic zero the Hilbert--Chow map identifies this
compactification with the closure of the graph cycles in the Chow
variety.
\end{theorem}

These are the algebraic inputs from \cite[Theorems~3 and~12,
Proposition~10]{r2grBrion2003}. The comparison with the entire irreducible
component through the diagonal requires
$\mathbf G=\operatorname{Aut}^0(X)$; faithful transitivity alone gives the
stated orbit closure. For example, a transitive proper subgroup of the
connected automorphism group still has its own graph compactification.

The three properties of a boundary fibre have different roles.
Reducedness determines the scheme from its support and gives each
component coefficient one in the fundamental cycle. Pure dimension puts
all these components in dimension $\dim X$. Cohen--Macaulayness is a
further local algebraic property, useful for proper intersections.

\section{Flat limits and moving sources}\label{r2gr:flat}

Flatness lets us approach every limiting pair along the prescribed
sequence of graphs, as required for a destination.

\begin{lemma}\label{r2gr:flat-support}
Let $S$ be a compact complex analytic space, let $M$ be a complex
projective variety with a compatible metric, and let
$\mathcal Y\subset S\times M$ be a closed analytic family, flat over
$S$, with nonempty fibres. Then
\[
 S\longrightarrow\mathcal K(M),\qquad s\longmapsto|\mathcal Y_s|
\]
is continuous for the Hausdorff topology on nonempty compact subsets.
In particular, if $s_m\to s$ and $q\in|\mathcal Y_s|$, there are
$q_m\in|\mathcal Y_{s_m}|$ with $q_m\to q$.
\end{lemma}

\begin{proof}
A flat holomorphic map of complex analytic spaces is open
\cite[\S8.3, p.~68, note~13]{r2grDouady1966}. If $V$ is an open
neighbourhood of $q$ in $M$, the projection of
$\mathcal Y\cap(S\times V)$ is therefore an open neighbourhood of $s$.
It follows that every sufficiently late fibre meets $V$. Taking a
decreasing neighbourhood basis at $q$ and choosing the points in
successively smaller members gives $q_m\to q$ along the whole
sequence. The finitely many initial choices are immaterial.

Conversely, if $q_m\in|\mathcal Y_{s_m}|$ and $q_m\to q$, closedness of
$\mathcal Y$ gives $q\in|\mathcal Y_s|$. These statements give the two
Hausdorff inclusions uniformly. For the first, otherwise some
$q_m\in|\mathcal Y_{s_m}|$ stay a fixed positive distance from
$|\mathcal Y_s|$; compactness gives a contradictory convergent
subsequence. For the other, cover the compact set $|\mathcal Y_s|$ by
finitely many small balls centred at its points, and use the first
paragraph for those centres.
\end{proof}

\begin{theorem}\label{r2gr:destinations}
The support map
\[
 \mathcal H_{X,\mathbf G}\longrightarrow\mathcal K(X\times X),
 \qquad \xi\longmapsto|Z_\xi|
\]
is a homeomorphism onto the Hausdorff closure of the graphs of
$\mathbf G$. If $[\Gamma_{g_m}]\to\xi$, then, for each $x\in X$,
\begin{equation}\label{r2gr:fibre}
 \D^X_{(g_m)}(x)=\{y:(x,y)\in|Z_\xi|\}.
\end{equation}
Every pair in $|Z_\xi|$ is approached by
$(x_m,g_mx_m)$ along the whole given sequence.
\end{theorem}

\begin{proof}
Apply Lemma~\ref{r2gr:flat-support} to the analytification of the
universal family. It gives continuity of the support map.
Theorem~\ref{r2gr:brion} says that its fibres are reduced: two points of
the Hilbert scheme with the same support therefore define the same
closed subscheme. The support map is injective. Its domain is compact
and its target is Hausdorff, so it is a homeomorphism onto its image.
The open set of graph points is dense also in the complex topology;
the image is precisely the Hausdorff closure of these graphs.

If $x_m\to x$ and $g_mx_m\to y$, closedness of the universal family
gives $(x,y)\in|Z_\xi|$. In the other direction, the final assertion of
Lemma~\ref{r2gr:flat-support}, applied to the point $(x,y)$, gives
$(x_m,y_m)\in\Gamma_{g_m}$ converging to $(x,y)$. Since $y_m=g_mx_m$,
these are the required sources and destinations.
\end{proof}

The \emph{family of entire graphs} is flat in $\xi$, whereas the first
projection $Z_\xi\to X$ need not be flat. Destination dimensions may
therefore jump; specialising only the image of a fixed source would
lose most of a boundary graph.

Transposition on $X\times X$ takes $\Gamma_g$ to $\Gamma_{g^{-1}}$ and
preserves the graph compactification. Thus it takes the scheme limit to
its transpose, as well as transposing the destination relation. More
generally, for $a,c\in\mathbf G$,
\begin{equation}\label{r2gr:action}
 (a,c)\Gamma_g=\Gamma_{cga^{-1}}.
\end{equation}
This elementary convention fixes the placement of both compact factors
in the Cartan limit below.

\section{Schubert components of a Cartan limit}\label{r2gr:schubert-section}

Choose a maximal compact subgroup $K$ and a maximal real split torus
$A=\exp\ak$ compatibly with $(B,T)$: we require $A\subset T$ and
choose the positive chamber $\ak^+$ determined by $B$.
We use the corresponding Cartan notation of Chapter~\ref{ch:cartan}.
After extraction write
\[
 g_m=k_m\exp(H_m)k'_m,\qquad k_m\to k,\quad k'_m\to k',
 \quad H_m\in\ak^+,
\]
with every simple root having a limit in $[0,\infty]$. Let
\[
 \Theta=\{\alpha:\alpha(H_m)\to\infty\},\qquad
 J=\Delta\setminus\Theta.
\]
Thus $J$ records the roots whose gaps remain finite. Put
$\mathbf L_J=L_\Theta$, the standard Levi subgroup with root system
generated by $J$, and let $w_J$ be the longest element of $W_J$.
Use the splitting of Lemma~\ref{lem:split}:
\[
 H_m=\nu_m+\eta_m,\qquad z_m=\exp(\nu_m),\qquad
 b_m=\exp(\eta_m)\longrightarrow b.
\]
Here $z_m$ is central in $\mathbf L_J$,
$\alpha(\nu_m)=0$ for $\alpha\in J$, and the other simple roots of
$\nu_m$ tend to $+\infty$. The element $b\in\mathbf L_J$ retains the
finite gaps; the set $J$ does not determine it.

Let ${}^JW^I$ denote the minimal representatives of
$W_J\backslash W/W_I$. Define closed subvarieties of $X\times X$ by
\begin{equation}\label{r2gr:CJ}
 C_{J,w}=\diag(\mathbf L_J)
       \bigl(w_JX^-_{w_Jw}\times X_w\bigr),\qquad
 C_J=\bigcup_{w\in{}^JW^I}C_{J,w}.
\end{equation}
There is no omitted closure in this notation. Indeed, the product in
parentheses is invariant under $\mathbf L_J\cap wPw^{-1}$, and the
associated bundle
\[
 \mathbf L_J\mathbin{\times^{\mathbf L_J\cap wPw^{-1}}}
       \bigl(w_JX^-_{w_Jw}\times X_w\bigr)
 \longrightarrow X\times X
\]
is projective. Its image is $C_{J,w}$. Brion's component description
also says that this morphism is birational onto its image.

\begin{theorem}\label{r2gr:schubert}
The Hilbert limit of the graphs of $g_m$ is the reduced subscheme with
support
\begin{equation}\label{r2gr:cartan-graph}
 R_{(g_m)}=(k'^{-1},kb)C_J.
\end{equation}
Consequently
\begin{equation}\label{r2gr:exact-schubert-fibre}
 \D^X_{(g_m)}(x)
 =kb\,\{y:(k'x,y)\in C_J\}.
\end{equation}
The $C_{J,w}$ are its distinct irreducible components before
translation. Each has dimension $\dim X$ and coefficient one in the
fundamental cycle. No comparison of two divergent root gaps is needed.
\end{theorem}

\begin{proof}
For clarity, the sign in the torus chart can be checked in Brion's
parametrisation. The group element $h$ in that compactification
parametrises $\Gamma_{h^{-1}}$: the orbit map is
$(a,c)\mapsto ac^{-1}$, whereas its graph is $\Gamma_{ca^{-1}}$.
Thus $\Gamma_{z_m}$ is parametrised by $z_m^{-1}$.
The torus chart has coordinates
\[
 q_\alpha=\alpha(z_m^{-1})=e^{-\alpha(\nu_m)}.
\]
They tend to $1$ for $\alpha\in J$ and to $0$ otherwise. This is
Brion's canonical point of type $J$, which lies on the boundary when $J\ne\Delta$.
Proposition~10 of \cite{r2grBrion2003} identifies the support of its
fibre with \eqref{r2gr:CJ}; Theorem~\ref{r2gr:brion} identifies the
fibre with the reduced scheme on that support.

Since $z_m$ commutes with $b_m$, equation~\eqref{r2gr:action} gives
\[
 \Gamma_{g_m}=(k_m'^{-1},k_mb_m)\Gamma_{z_m}.
\]
The two factors converge, so the algebraic group action on the Hilbert
scheme gives the claimed Hilbert limit. Theorem~\ref{r2gr:destinations}
then proves both the Hausdorff limit and
\eqref{r2gr:exact-schubert-fibre}.

The component indexing and the birational associated bundles are
Brion's description. Their dimensions may also be read at $wP$.
Choose $H_0\in\ak$ with $\alpha(H_0)=0$ for $\alpha\in J$ and
$\alpha(H_0)>0$ for $\alpha\notin J$. The tangent space $T_{wP}X$
splits into positive, zero and negative weight spaces according to
its weights evaluated on $H_0$; zero weights give the Levi directions.
Thus the signs are relative to this positive central direction,
not to an ordering of all characters of the centre of $\mathbf L_J$.
The open part of the associated bundle has as base the Levi orbit
(the zero directions), and as its two factors the negative and
positive unipotent orbits.
Its dimension is their sum, $\dim X$. Theorem~\ref{r2gr:brion} gives
coefficient one by reducedness. The chart limit used only the finite
or infinite alternatives for the simple roots.
\end{proof}

The fibre of \eqref{r2gr:cartan-graph} is the unstable orbit closure
of Theorem~\ref{thm:dest}. We now also know its algebraic family and
its graph components. These have dimension $\dim X$, but their fibres
over a source may have different dimensions.

\begin{corollary}\label{r2gr:regular}
If every simple root gap diverges, then
\[
 C_\varnothing=\bigcup_{w\in W^I}X_w^-\times X_w.
\]
If $k'x\in B^-vP/P$, with $v\in W^I$, then
\[
 \D^X_{(g_m)}(x)=kX_v,
 \qquad \dim_\C\D^X_{(g_m)}(x)=\ell(v).
\]
\end{corollary}

\begin{proof}
Here $\mathbf L_\varnothing=T$, $w_\varnothing=1$, and $b=1$ in
the splitting used above. The products of Schubert varieties are
$T$-stable. A point of $B^-vP/P$ belongs to $X_w^-$ exactly when
$w\le v$ in Bruhat order. Its fibre is therefore
$k\bigcup_{w\le v}X_w=kX_v$.
\end{proof}

For projective space the $X_w$ form a nested family of linear spaces,
and this is the regular case of Proposition~\ref{prop:proj}. For a
general flag variety it is the Bruhat order, rather than a single
projective dimension, that describes the incidences of the fibres.

\begin{example}[Two different surfaces of destinations]\label{r2gr:grassmannian}
Take $X=\operatorname{Gr}_2(\C^4)$ and
\[
 a_m=\operatorname{diag}(e^{3m},e^m,e^{-m},e^{-3m}).
\]
Every simple root gap diverges. Put
$E_j=\langle e_1,\ldots,e_j\rangle$. The Schubert indices are the
two-element subsets $\{i_1<i_2\}$ of $\{1,2,3,4\}$, with
\[
 X_{i_1i_2}=\{U:\dim(U\cap E_{i_a})\ge a\ (a=1,2)\},
 \qquad \dim X_{i_1i_2}=i_1+i_2-3.
\]
At the fixed source $U_{ij}=\langle e_i,e_j\rangle$, the opposite
cell has the same index, so Corollary~\ref{r2gr:regular} gives
\[
 \begin{aligned}
 \D_{(a_m)}(U_{14})&=\{U:E_1\subset U\}\cong\PP^2,\\
 \D_{(a_m)}(U_{23})&=\{U:U\subset E_3\}\cong\PP^2.
 \end{aligned}
\]
Neither surface contains the other. Their intersection is the line
$\{U:E_1\subset U\subset E_3\}$, which is itself the destination
fibre of $U_{13}$. At $U_{24}$ the fibre is the Schubert divisor
$\{U:U\cap E_2\ne0\}$. Thus equal destination dimensions need not
mean equal carriers; the incidence conditions matter.

The Pl\"ucker eigenvalues at $e_1\wedge e_4$ and
$e_2\wedge e_3$ both equal $1$. This equality of two weight values
does not merge the two intrinsic fibres. The graph calculation takes
place on the Grassmannian, where the Pl\"ucker incidence equation is
retained throughout the degeneration.
\end{example}

\section{Residual maps and graph resolutions}
\label{r2res:section}

Let $V$ be a complex vector space, with the full $\PGL(V)$ acting on
$\mathbb P(V)$. The complete-collineation compactification of
Section~\ref{sec:SPvsWond} is the Hilbert graph compactification;
this projective-space identification is due to Thaddeus
\cite[Main Theorem~2.3 and \S8]{rwThaddeus1999}.
Theorem~\ref{r2gr:brion} gives the corresponding model for complex flag
varieties. We first express the classical residual map in a blow-up
chart, then use the block components to calculate multidegrees and
diagonal return degrees.

\subsection{The residual map in the projective-plane blow-up}
\label{r2res:plane-chart}

Let $\dim V=3$, let $Y=\mathbb P(\operatorname{End}V)$, and let
$\Sigma\subset Y$ be the rank-one Segre variety.  By Theorem~\ref{thm:blowup}, the complete-collineation model is
\[
 \pi:\widehat Y=\operatorname{Bl}_{\Sigma}Y\longrightarrow Y.
\]
For a rank-one $A$, put $K=\ker A$ and $L=\operatorname{im}A$.
The point $[A]$ determines the kernel and image. What happens when
a moving source approaches the kernel? The exceptional divisor records the missing
projective residual map from $K$ to $V/L$.

\begin{proposition}
\label{r2res:schur}
There is a canonical identification
\[
 \pi^{-1}([A])=\mathbb P\operatorname{Hom}(K,V/L).
\]
In independent source and target bases in which $A=\operatorname{diag}(1,0,0)$,
write a nearby matrix as
\[
 M=\begin{pmatrix}1&u\\v&W\end{pmatrix},
 \qquad C=W-vu.
\]
The rank-one centre has equations $C=0$, and its blow-up retains $[C]$.
If $E$ is the exceptional divisor and $D$ the strict transform of the
determinantal hypersurface, then
\[
 E\cap\pi^{-1}([A])\cong\mathbb P^3,
 \qquad D\cap\pi^{-1}([A])=\{[B]:\det B=0\}\cong\mathbb P^1\times\mathbb P^1.
\]
The reduced boundary is $E\cup D$, with normal crossings.  The total
transform of the determinant divisor has multiplicities
\begin{equation}\label{r2res:determinant-multiplicity}
 \pi^*\operatorname{div}(\det)=D+2E.
\end{equation}
\end{proposition}

\begin{proof}
Gaussian elimination gives
\begin{equation}\label{r2res:schur-factorization}
 M=
 \begin{pmatrix}1&0\\v&I\end{pmatrix}
 \begin{pmatrix}1&0\\0&C\end{pmatrix}
 \begin{pmatrix}1&u\\0&I\end{pmatrix}.
\end{equation}
Thus the rank-one ideal on this chart is the ideal of the four entries
of $C$.  Blowing it up replaces $C=0$ by their projective direction
$[B]$, so a blow-up chart has $C=tB$.  Intrinsically, the normal space at $[A]$ is
$\operatorname{Hom}(K,V/L)\otimes(\mathbb C A)^*$: the normal
component of a tangent vector is its restriction to $K$, followed by
the quotient to $V/L$.  Tensoring by the one-dimensional space
$(\mathbb C A)^*$ disappears on projectivising, which proves the
canonical fibre identification.

The determinant is $\det M=\det C=t^2\det B$.  The exceptional divisor
is $t=0$ and the strict transform is $\det B=0$, proving
\eqref{r2res:determinant-multiplicity}.  A nonzero singular two-by-two
matrix has rank one; its projective class is determined by its kernel
and image lines.  This is the smooth Segre quadric in $\mathbb P^3$.
The equation $\det B=0$ is smooth on its rank-one locus and independent
of the coordinate $t$, so the two boundary divisors meet transversely.
Outside the exceptional divisor the determinant hypersurface is smooth
along its rank-two locus.  These charts cover the boundary.
\end{proof}

The multiplicity two in~\eqref{r2res:determinant-multiplicity} is a divisor
multiplicity.  It must not be confused with the multiplicity one of the
irreducible graph components in Theorem~\ref{r2gr:schubert}.

Let $p=\mathbb P(L)$, let $q:V\to V/L$ be the quotient, and, for
$0\ne B\in\operatorname{Hom}(K,V/L)$, define
\begin{equation}\label{r2res:residual-relation}
 R_B=
 \bigl(\mathbb P(V)\times\{p\}\bigr)
 \ \cup\ 
 \{([x],[y])\in\mathbb P(K)\times\mathbb P(V):
                     Bx\wedge q(y)=0\}.
\end{equation}
The vanishing equation allows either factor to be zero; this convention
is essential at a residual kernel point.

\begin{corollary}[The residual relation]
\label{r2res:rankone-fibre}
The support of the graph boundary point over $([A],[B])$ is exactly
$R_B$.  If $B$ is invertible, the fibre over a point $[x]\in\mathbb P(K)$
is the line $\mathbb P(q^{-1}(B\mathbb Cx))$ through $p$.  If $B$ has
rank one, write $a=\mathbb P(\ker B)$ and
$H=\mathbb P(q^{-1}(\operatorname{im}B))$.  The fibres over
$\mathbb P(K)\setminus\{a\}$ are $H$, and the fibre over $a$ is the
whole plane.  In both cases the fibre off $\mathbb P(K)$ is $\{p\}$.
\end{corollary}

\begin{proof}
Choose normalized representatives $M_m\to A$ in the Schur chart for
an invertible sequence converging to $([A],[B])$ in the blow-up. Write
$C_m=t_mB_m$, where $t_m\to0$ and $B_m\to B\ne0$. The two outer factors in
\eqref{r2res:schur-factorization} tend to the identity. For the middle
factors $S_m=\operatorname{diag}(1,t_mB_m)$, we have $S_m\to A$ and
$t_m^{-1}\wedge^2 S_m\to B\oplus 0$ in the corresponding splitting
of the second exterior power. This gives convergence of the exterior
data along the whole sequence. If $B$ is
invertible, the middle factors have blocks of dimensions $(1,2)$,
with residual map $B$.  If $B$ has rank one, the ratio of the two
singular values of $B_m$ diverges; the middle factors consequently have
three one-dimensional blocks, with source flag
$\ker B\subset K$ and target flag
$L\subset q^{-1}(\operatorname{im}B)$.  The destination formula of
Theorem~\ref{thm:ext-boundary} gives the stated fibres along the whole
sequence.  Their union is precisely the closed relation
\eqref{r2res:residual-relation}.  Theorem~\ref{r2gr:destinations}
identifies that relation with the graph support.
\end{proof}

This is the coordinate-free form of Example~\ref{ex:d3fibre}.
Equation~\eqref{r2res:residual-relation} is used here as a support
identity; no assertion about an arbitrary incidence scheme follows
from the wedge equation alone.  For example,
$\operatorname{diag}(1,t,t)$ and $\operatorname{diag}(1,2t,t)$ have
the same leading rank-one map but attach different destination lines
to $[e_2+e_3]$.  The residual map, rather than merely the exceptional
divisor containing it, distinguishes the two.

\subsection{Block resolutions and their defects}
\label{r2res:block-resolution-section}

Let $\dim V=d=n+1$, and let $\xi$ be an exterior boundary point
with ordered block sizes $d=b_1+\cdots+b_k$ and endpoints
$r_j=b_1+\cdots+b_j$, $r_0=0$.  Retain its decreasing source flag
$F_\bullet$, increasing target flag $G_\bullet$, and graded maps
\[
 A_j:F_j/F_{j+1}\xrightarrow{\sim}G_j/G_{j-1}.
\]
Let $Z_\xi$ be its Hilbert graph degeneration, which is reduced by
Theorem~\ref{r2gr:brion}, and let $Z_j$ be the reduced
closure of its relation over
$\mathbb P(F_j)\setminus\mathbb P(F_{j+1})$.
Here $F_1=G_k=V$, $F_{k+1}=G_0=0$,
$\dim F_j=d-r_{j-1}$ and $\dim G_j=r_j$.  Write
$Q_j=F_j/F_{j+1}$, and denote the quotient maps by
$q_j:F_j\to Q_j$ and $\pi_j:G_j\to G_j/G_{j-1}$.

Recording the active quotient line gives the block resolution of
Brion \cite[Example~3]{r2grBrion2003}. We write it in the source and
target flags of the destination formula. Over the locus where both
active quotient coordinates vanish, its fibre is the projective space
of active lines. The last assertion below computes the additional
cohomology of this resolution.

\begin{proposition}
\label{r2res:block-resolution}
The $Z_j$ are precisely the irreducible components of $Z_\xi$.
The incidence variety
\[
 \widetilde Z_j=
 \{([v],[w],\ell)\in\mathbb P(F_j)\times\mathbb P(G_j)\times\mathbb P(Q_j):
          q_j(v)\in\ell,\ \pi_j(w)\in A_j\ell\}
\]
is smooth, projective, and irreducible of dimension $n$, and forgetting
$\ell$ gives a proper birational map $\rho_j:\widetilde Z_j\to Z_j$.
It is an isomorphism away from
\[
 E_j=\mathbb P(F_{j+1})\times\mathbb P(G_{j-1}),
\]
and $\rho_j^{-1}(E_j)=E_j\times\mathbb P(Q_j)$.
The set $E_j$ is the locus where both active quotient coordinates
vanish. If $b_j=1$, then $\rho_j$ is an isomorphism everywhere,
even when $E_j\ne\varnothing$. If $b_j>1$, its non-isomorphism
locus in $Z_j$ is exactly $E_j$.
On the complex analytic spaces, for a field of coefficients
$\mathbf k$, put
$\mathscr D_{\rho_j}=\operatorname{Cone}(\mathbf k_{Z_j}\to
R\rho_{j*}\mathbf k_{\widetilde Z_j})$, and let $i_j:E_j\hookrightarrow Z_j$.
Then
\begin{equation}\label{r2res:block-defect}
 \mathcal H^q(\mathscr D_{\rho_j})=
 \begin{cases}
 (i_j)_*\mathbf k_{E_j},&q=2,4,\ldots,2b_j-2,\\
 0,&\text{otherwise}.
 \end{cases}
\end{equation}
Here $\mathbb P(0)=\varnothing$; in particular $E_j=\varnothing$
for each extreme block, $j=1$ or $j=k$.
\end{proposition}

\begin{proof}
Over $\ell\in\mathbb P(Q_j)$ the two projective fibres are
$\mathbb P(q_j^{-1}\ell)$ and $\mathbb P(\pi_j^{-1}A_j\ell)$.
They are the projectivisations of vector bundles of ranks
$d-r_j+1$ and $r_{j-1}+1$.  Their fibre product is therefore smooth,
projective, irreducible, and has dimension
\[
 (b_j-1)+(d-r_j)+r_{j-1}=d-1.
\]
Its locus $q_j(v)\ne0$ maps to exactly the source-stratum piece from
the destination formula.  That locus is dense, so its proper image is
$Z_j$.  Off $E_j$ at least one of $q_j(v)$ and $\pi_j(w)$ is nonzero;
it determines $\ell$ uniquely, with a regular inverse.
On either such open set the incidence equations identify the reduced
image scheme with the quotient-graph equations.  Thus this is an
isomorphism of schemes, not only a bijection of points.  On $E_j$ both
vanish, and every $\ell\in\mathbb P(Q_j)$ is allowed.
If $b_j=1$, then $\mathbb P(Q_j)$ is a point, the incidence conditions
impose no further restriction, and forgetting $\ell$ is an isomorphism
on all of $Z_j$. If $b_j>1$, every fibre over $E_j$ has positive
dimension, so the non-isomorphism locus in the target is exactly $E_j$.

The full relation is the union of these finitely many closures.  Their
source and target projections are $\mathbb P(F_j)$ and
$\mathbb P(G_j)$; the target dimensions are strictly increasing.
Thus the $Z_j$ are distinct irreducible varieties of the same dimension,
none contains another, and they are exactly the graph components.
This recovers the classical block resolution
\cite[Example~3]{r2grBrion2003} directly from the destination formula.

Proper base change for sheaves on the analytic topology
\cite[Tag~09V6]{r2resStacks} identifies the stalks of
$R^q\rho_{j*}\mathbf k$ with the cohomology of its fibres.  These are
points off $E_j$ and $\mathbb P^{b_j-1}$ on $E_j$.  They are connected,
so the unit in degree zero is an isomorphism.  In positive degrees the
sheaves are supported on $E_j$; the displayed product over $E_j$ and
its hyperplane class make them constant there in precisely the even
degrees $2,4,\ldots,2b_j-2$.  Taking the cone removes degree zero and
gives~\eqref{r2res:block-defect}.
\end{proof}

Formula~\eqref{r2res:block-defect} determines the cohomology sheaves
of the derived object.  For type $(1,2,1)$ in vector dimension four, the middle
locus $E_j$ is a point, its inverse image is $\mathbb P^1$, and the
only nonzero defect sheaf is $\mathbf k$ in degree two, supported
at that point. For type $(1,1,1,1)$ every such defect vanishes.

\section{Multidegrees and diagonal returns}
\label{r2res:returns-section}

Intersecting the components with general linear spaces recovers their
block sizes. Intersecting with the diagonal brings us back to returns.

\subsection{The class of an individual component}
\label{r2res:cycles-section}

The whole graph has the diagonal class. For the block components just
described, the following calculation separates that class into
consecutive multidegree intervals and recovers the ordered block
sizes. Intersecting each component with the diagonal gives the degree
$b_j$; this is the intersection number which will be used for returns.
Write $h_1,h_2$ for the source and target hyperplane classes; the
component classes lie in $A^n(\mathbb P^n\times\mathbb P^n)$.

\begin{theorem}
\label{r2res:multidegree}
For each block component,
\begin{equation}\label{r2res:multidegree-formula}
 [Z_j]=\sum_{a=r_{j-1}}^{r_j-1}h_1^a h_2^{n-a}.
\end{equation}
Consequently
\[
 [Z_\xi]=\sum_{a=0}^{n}h_1^a h_2^{n-a}=[\Delta_{\mathbb P^n}],
 \qquad \deg\Delta^![Z_j]=b_j.
\]
Here $\Delta^![Z_j]\in A_0(Z_j\cap\Delta)$ is the refined
intersection class, whether or not the intersection is proper.
In particular each block component contains a dynamically returning
point.  The ordered block sizes are recovered from the multidegree
intervals of the individual components.  The total diagonal cycle
class alone does not recover them.
\end{theorem}

\begin{proof}
Write $[Z_j]=\sum_a c_a h_1^a h_2^{n-a}$.  Its coefficient is
\[
 c_a=\int_{Z_j}h_1^{n-a}h_2^a.
\]
Choose general linear spaces $L_1,L_2\subset\mathbb P^n$ of
codimensions $n-a$ and $a$.  The source projection has dimension
$n-r_{j-1}$ and the target projection dimension $r_j-1$.
Thus the first intersection is empty if $a<r_{j-1}$, and the second
is empty if $a\ge r_j$.

For $r_{j-1}\le a\le r_j-1$, their intersections with the source and
target projection spaces have dimensions $a-r_{j-1}$ and $r_j-a-1$.
Generality makes them avoid $\mathbb P(F_{j+1})$ and
$\mathbb P(G_{j-1})$.  The quotient maps embed them as general linear
spaces in the two active projective quotients.  Identifying these by
$[A_j]$, their dimensions add to $b_j-1$.  They therefore meet
transversely in a single point.  The quotient-graph description of
$Z_j$ gives $c_a=1$ on this interval, proving
\eqref{r2res:multidegree-formula}.

The graph fibre is reduced, and the intervals partition
$\{0,\ldots,n\}$, giving the total class.  Compatibility of refined Gysin maps with proper push-forward
\cite[\S6.2]{rwFulton1998} identifies the degree of the refined
class with $\int_{\mathbb P^n\times\mathbb P^n}[Z_j][\Delta]$.
Equivalently, both $h_1$ and $h_2$ restrict on the diagonal to its
hyperplane class $h$.  Each summand in
\eqref{r2res:multidegree-formula} contributes
$\int_{\mathbb P^n}h^n=1$, and there are $b_j$ summands.
The refined class is supported on $Z_j\cap\Delta$; its nonzero
degree rules out an empty intersection even when this locus has
positive dimension.  By Theorem~\ref{r2gr:destinations}, every pair in the graph support
is a genuine moving-source limit, so these diagonal points return dynamically.
Finally, the nonzero coefficient intervals give their consecutive
endpoints $r_{j-1},r_j$, hence the ordered $b_j$.
\end{proof}

\subsection{Existence and the residual return operator}

Each component also meets the diagonal by a direct linear argument.
Put $H=F_j\cap G_j$ and $b=b_j$. The dimensions of the source and
target projection spaces give
\[
 \dim H\ge(d-r_{j-1})+r_j-d=b.
\]
Write $f=q_j|_H$ and $g=\pi_j|_H$. If $0\ne v\in\ker f$, choose
the active line $A_j^{-1}(\C g(v))$ when $g(v)\ne0$, and any active
line otherwise. The incidence description gives $([v],[v])\in Z_j$.
If $\ker g\ne0$, choose the active line $\C f(v)$ for a nonzero
$v\in\ker g$, or any line if $f(v)=0$. If both maps are injective,
their $b$-dimensional targets force $\dim H=b$ and make both maps
isomorphisms. An eigenline of the complex operator $g^{-1}A_jf$
again gives a diagonal point. This proves the existence of a return in
each component. The field matters here: a real operator need not have
a real eigenline.

The degree $b_j$ is an intersection degree; identifying it with the
length of a finite return scheme requires a proper-intersection
argument. Different components may return to the same point. Under
the quotient identifications in the next proposition, the residual
map determines the entire return scheme as a projective fixed scheme.

\begin{proposition}
\label{r2res:return-operator}
Put $H_j=F_j\cap G_j$.  Suppose explicitly that
\[
 q_j|_{H_j}:H_j\xrightarrow{\sim}Q_j,
 \qquad \pi_j|_{H_j}:H_j\xrightarrow{\sim}G_j/G_{j-1}.
\]
Then the scheme-theoretic intersection $Z_j\cap\Delta$, regarded as
a subscheme of $\mathbb P(H_j)$, is the fixed scheme of the projective
operator
\[
 \mathcal R_j=\bigl[(\pi_j|_{H_j})^{-1}A_j(q_j|_{H_j})\bigr]\in\operatorname{PGL}(H_j).
\]
If this scheme is zero-dimensional, its length is $b_j$.  If
$\mathcal R_j$ has $b_j$ distinct eigenvalues, there are $b_j$ reduced
return points.
\end{proposition}

\begin{proof}
A diagonal point belongs to $\mathbb P(H_j)$.  The two isomorphism
assumptions give $\dim H_j=b_j$ and prevent either quotient from
vanishing at a nonzero vector of $H_j$.  The equations of the quotient
graph on the diagonal are therefore
\[
 A_jq_j(v)\wedge\pi_j(v)=0.
\]
Off $E_j$, Proposition~\ref{r2res:block-resolution} identifies the
component scheme with its quotient-graph incidence equations.  Thus the
displayed equations define the intersection scheme, and are exactly
the eigenline equations of $\mathcal R_j$, including their scheme
structure.  This locus avoids $E_j$, so $Z_j$ is smooth
there by Proposition~\ref{r2res:block-resolution}.  If it is finite, it has the expected codimension $n$ in the smooth
$n$-fold $Z_j$.  The $n$ local generators of the diagonal ideal
therefore form a regular sequence on $Z_j$.  Consequently the refined
intersection is its fundamental zero-cycle
\cite[\S\S6.1,7.1]{rwFulton1998}, and its scheme length is the
degree $b_j$.  If the eigenvalues are distinct, the tangent eigenvalue
ratios at each eigenline are different from one, so each fixed point is
reduced.
\end{proof}

\begin{example}[A double return at one point]\label{r2res:double-return}
In vector dimension three, take the two-block datum realised by
\[
 g_t=\begin{pmatrix}1&1&0\\0&1&0\\0&0&t\end{pmatrix},\qquad
 t\longrightarrow0,\quad t\ne0.
\]
For the first block, $H_1=\langle e_1,e_2\rangle$ and both quotient
restrictions are isomorphisms.  Its return operator is the projective
class of $\left(\begin{smallmatrix}1&1\\0&1\end{smallmatrix}\right)$.
It has a single eigenline $[e_1]$.  In the affine coordinate
$[1:z]$ around that line, the fixed equation is
$z/(1+z)=z$, or $z^2=0$.  Thus the return scheme has length two and
one point of support, in agreement with $b_1=2$.
\end{example}

For $g_m=\operatorname{diag}(e^m,e^m,e^{-2m})$, both quotient
restrictions for the first block are isomorphisms, but $\mathcal R_1$
is the identity.  The entire projective line returns, while the refined
degree is two.  For the translations $z\mapsto z+m$ on $\mathbb P^1$,
both limiting graph components meet the diagonal at $\infty$.
These componentwise returns alone do not imply the existence of a
regular critical seed, recurrence of a whole critical carrier, or the
existence of a largest domain.

\begin{remark}\label{r2res:scope}
The integer $b_j$ is preserved as a refined intersection degree, not as
the number of distinct points of a fibre.  The complete-collineation
model resolves action graphs of the full projective group; an
arrangement model obtained by blowing up a collection of diagonal or
incidence subspaces has different centres and a different universal
problem.  No identification of these models, or of their exceptional
fibres, enters the argument.
\end{remark}

\section{Support, scheme and cycle}\label{r2gr:cycles}

The support $|Z_\xi|$ records destination pairs, the scheme $Z_\xi$
their equations, and the fundamental cycle the components with their
generic lengths. The Hilbert--Chow isomorphism of
Theorem~\ref{r2gr:brion} concerns this graph family, not the scheme
structures of arbitrary intersections with the same support.

For the regular boundary the corresponding equality in the Chow group
is
\begin{equation}\label{r2gr:diagonal-class}
 [\Delta_X]=\sum_{w\in W^I}[X_w^-\times X_w]
       \quad\text{in }A_{\dim X}(X\times X).
\end{equation}
For rational equivalence, choose a regular dominant integral
cocharacter $\lambda$. For $t\ne0$, the graphs of $\lambda(t)^{-1}$
extend to a flat Hilbert family over $\mathbb A^1$ with reduced special fibre
$C_\varnothing$. The fibres at $1$ and $0$ are rationally equivalent,
proving \eqref{r2gr:diagonal-class}. Translates have the same class:
$\mathbf G$ is generated by its torus and root subgroups, whose
one-parameter actions give the required rational equivalences.

\begin{example}\label{r2gr:jump}
In $\PP^1\times\PP^1$, use source coordinates $[x_0:x_1]$ and target
coordinates $[y_0:y_1]$. The graphs of
$\operatorname{diag}(t^{-1},t)$ have equation
\[
 x_0y_1-t^2x_1y_0=0.
\]
Their flat limit is the reduced divisor
\[
 (\{[0:1]\}\times\PP^1)\ \cup\
 (\PP^1\times\{[1:0]\}).
\]
Both components have coefficient one, and this divisor is
Cohen--Macaulay. Its first projection has a point fibre at every
source except $[0:1]$, where the fibre is all of $\PP^1$. In
particular, the projection is not flat.
\end{example}

\begin{example}\label{r2gr:double-return}
Let $g_t$ be the translation $z\mapsto z+t^{-1}$ on $\PP^1$. Its
graph has equation
\[
 t(x_1y_0-x_0y_1)-x_1y_1=0.
\]
At $t=0$ it is the union
$\{\infty\}\times\PP^1\ \cup\ \PP^1\times\{\infty\}$, with both
components reduced and of multiplicity one. Restriction of its
equation to the diagonal gives $x_1^2=0$. The diagonal intersection
is thus a double point supported at $\infty$, of length two. There is
one returning point, and both graph components contribute there.
\end{example}

A reduced Cohen--Macaulay graph can thus have a nonreduced diagonal
intersection. At the expected codimension, the local equations of a
smooth intersecting subvariety form a regular sequence on the graph;
for a finite proper intersection, the resulting lengths compute its
degree. Excess intersection requires a refined class. The set of returning
points does not record these lengths, and the component coefficients
alone do not give the number of distinct returns.

\section{A flag variety inside projective space}\label{r2gr:restriction}

The restriction statement of Chapter~\ref{ch:parabolic} is an inclusion for destination
relations. It may be strict even for a homogeneous projective embedding.
The incidence equations constrain the moving sources before their
images are taken to the limit.

\begin{example}[The Segre surface]\label{r2gr:segre}
Embed $X=\PP^1\times\PP^1$ in $\PP^3$ by the Segre map, and use the
basis
\[
 e_1=u_+\otimes v_+,\quad e_2=u_+\otimes v_-,\quad
 e_3=u_-\otimes v_+,\quad e_4=u_-\otimes v_-.
\]
Consider
\[
 g_m=\left(
 \begin{pmatrix}e^{2m}&0\\0&e^{-2m}\end{pmatrix},
 \begin{pmatrix}e^m&0\\0&e^{-m}\end{pmatrix}\right).
\]
On the tensor product its matrix is
$\operatorname{Diag}(e^{3m},e^m,e^{-m},e^{-3m})$. At $x=[e_3]$,
Proposition~\ref{prop:proj} gives the ambient destination fibre
\[
 \D^{\PP^3}_{(g_m)}(x)=\PP\langle e_1,e_2,e_3\rangle.
\]
The Segre equation is $z_1z_4=z_2z_3$, so the intersection of this
plane with $X$ has support
\[
 \PP\langle e_1,e_2\rangle\ \cup\ \PP\langle e_1,e_3\rangle.
\]
The intrinsic fibre is only
\begin{equation}\label{r2gr:segre-intrinsic}
 \D^X_{(g_m)}(x)=\PP^1\times\{[v_+]\}
                  =\PP\langle e_1,e_3\rangle.
\end{equation}

To verify this along the whole sequence, write a source approaching
$([u_-],[v_+])$ as
$([a_m u_++u_-],[v_++c_m v_-])$, where $a_m,c_m\to0$.
The two images are
\[
 [e^{4m}a_m u_++u_-],\qquad[v_++e^{-2m}c_m v_-].
\]
The second necessarily tends to $[v_+]$. For any finite $s\in\C$,
the choice $a_m=e^{-4m}s$ realises $[su_++u_-]$ in the first
factor. The choice $a_m=e^{-2m}$ realises $[u_+]$. Taking $c_m=0$
proves every destination in \eqref{r2gr:segre-intrinsic}. Hence
$[e_2]$ belongs to the ambient intersection and is not an intrinsic
destination.
\end{example}

The ambient intersection already has too many points. The intrinsic
Hilbert family keeps the Segre constraint on both sources and targets
before taking the flat limit, as does the Schubert description.

\section{Asymptotic graphs of a group}\label{r2gr:group}

We can now return from a single sequence to a group. The preceding
compactification parametrises the graph limits whose union gives the
relation of Chapter~\ref{ch:limit}. Let $\Gamma$ be a countably
infinite group acting through $\rho:\Gamma\to\mathbf G$. Define its
asymptotic graph set by
\[
 \partial_X^\infty\Gamma=
 \bigcap_{F\subset\Gamma\ {\rm finite}}
 \overline{\{[\Gamma_{\rho(\gamma)}]:\gamma\in\Gamma\setminus F\}}^{\,\mathrm{an}}
 \subset\mathcal H_{X,\mathbf G}.
\]
The superscript denotes closure in the complex analytic topology on
$\mathcal H_{X,\mathbf G}$. The definition keeps the original group
elements as labels, even when their images coincide.
An exhaustion of $\Gamma$ by finite sets and a decreasing neighbourhood
basis show that every point of $\partial_X^\infty\Gamma$ is the limit
of graph points along pairwise distinct elements of $\Gamma$.

\begin{corollary}\label{r2gr:group-relation}
The destination relation of the original $\Gamma$-action is
\[
 \D^X_\Gamma=\bigcup_{\xi\in\partial_X^\infty\Gamma}|Z_\xi|.
\]
Hence a $\Gamma$-invariant open set $U\subset X$ is a domain of proper
discontinuity precisely when
\[
 |Z_\xi|\cap(U\times U)=\varnothing
 \quad\text{for every }\xi\in\partial_X^\infty\Gamma.
\]
\end{corollary}

\begin{proof}
For a destination pair choose its pairwise distinct witnesses in
$\Gamma$ and extract a convergent subsequence of Hilbert graph points.
Its limit belongs to every tail closure, and
Theorem~\ref{r2gr:destinations} puts the pair in its support. Conversely,
choose pairwise distinct elements whose graph points tend to a given
$\xi\in\partial_X^\infty\Gamma$. The same theorem realises every
pair of $|Z_\xi|$ along those elements. The properness statement follows
from the relational criterion of Chapter~\ref{ch:limit}: a forbidden
pair in $U\times U$ gives convergent witnesses in compact subsets of
$U$, and failure of properness supplies such witnesses after extraction.
\end{proof}

If $\rho$ has finite kernel and discrete image, this asymptotic set is
the analytic closure of the graph points of $\rho(\Gamma)$ with those
graph points removed.
For a general action it can contain ordinary graphs. In particular,
an infinite kernel puts $\Delta_X$ in $\partial_X^\infty\Gamma$;
then every point is recurrent and no nonempty open set is a domain of
proper discontinuity for the original group. The effective algebraic
quotient simplifies the graphs but cannot erase group labels. Component
classes and intersection lengths have no automatic counterpart on a real
flag manifold.

The supports of the asymptotic graphs recover the destination relation
and hence the criterion for proper discontinuity of
Chapter~\ref{ch:limit}. Over $\C$, the graph family also records
component classes and the scheme structure of diagonal intersections.
Finite diagonal intersections have lengths, and the refined
intersection class also applies in positive dimension. This connects the
limits of moving sources with the algebraic geometry of degenerating graphs.

\backmatter
\chapter*{Declarations}
\addcontentsline{toc}{chapter}{Declarations}

\section*{Funding}

\mbox{}
The research of A.~Cano was partially supported by
SECIHTI--SNII 104023. Marlon M.~López-Flores received partial funding from the
EVAAT-GCN project under FINEP agreement 01.24.0661.00,
ref.~3311/24. José Seade was supported by SECIHTI--SNII 1170 and
DGAPA--PAPIIT--UNAM IN101424.
Andrew Shaw was supported by SECIHTI, Mexico, through
graduate fellowship No.~4076100.

\section*{Acknowledgements}
We are indebted to William Goldman, Alberto~Verjovsky, John Parker, Dennis Sullivan, 
Richard Canary, and Tengren Zhang for several truly inspiring conversations.
We thank A.~Arroyo, C.~Cabrera, L.~A.~Cano Cordero, A.~Guillot,
P.~Makienko, M.~M\'endez, H.~Z\'u\~niga, W.~Barrera, J.~P.~Navarrete,
and I.~Garc\'ia for helpful discussions on complex projective dynamics
and related questions.

\section*{AI Declaration}
\label{rev59:ai-declaration}
This monograph uses assistance from OpenAI ChatGPT, Claude and Codex for language
and structural editing, literature searches, and compuational  experiments.

\end{document}